\documentclass{article}
\usepackage[latin1]{inputenc}
\usepackage[T1]{fontenc} 

\usepackage{graphicx,xcolor}
\usepackage{subfigure}          

\usepackage[english]{babel}
\usepackage[autostyle=true,german=quotes]{csquotes}

\usepackage{epigraph}

\usepackage{placeins}

\usepackage{tabularx}
\usepackage{adjustbox} 	

\usepackage{amsmath}
\usepackage{amsfonts}
\usepackage{amssymb}
\usepackage{mathabx}

\usepackage{lmodern,textcomp}

\usepackage{algorithm}
\usepackage{algpseudocode}

\usepackage{enumerate}

\usepackage{xspace}
\newcommand\largeparbreak{\par\bigskip}

\def\mathrlap{\mathpalette\mathrlapinternal}

\def\mathrlapinternal#1#2{%
	\rlap{$\mathsurround=0pt#1{#2}$}}

\newcommand{\mdtheorem}[2]{\newtheorem{#1}{#2}}
\mdtheorem{Definition}{Definition}
\mdtheorem{Requirement}{Requirement}
\mdtheorem{Remark}{Remark}
\mdtheorem{Theorem}{Theorem}
\mdtheorem{Corollar}{Corollar}
\mdtheorem{Conjecture}{Conjecture}
\mdtheorem{Lemma}{Lemma}
\mdtheorem{Example}{Example}

\newcommand*\tageq{\refstepcounter{equation}\tag{\theequation}}

\newcommand{\No}{\mathbb{N}_0\xspace}

\newcommand{\rg}{\operatorname{rg}}
\newcommand{\rk}{\operatorname{rank}}

\newcommand*{\LargerCdot}{\raisebox{-0.25ex}{\scalebox{1.7}{$\cdot$}}}

\newcommand{\sign}{\operatorname{sign}\xspace}

\newcommand{\R}{\mathbb{R}\xspace}
\newcommand{\N}{\mathbb{N}\xspace}
\newcommand{\C}{\mathbb{C}\xspace}
\newcommand{\Cno}{\mathbb{C}\setminus \lbrace 0 \rbrace \xspace}
\newcommand{\Rno}{\mathbb{R}\setminus \lbrace 0 \rbrace \xspace}

\newcommand{\cM}{\mathcal{M}\xspace}

\newcommand{\cB}{\mathcal{B}\xspace}
\newcommand{\cP}{\mathcal{P}\xspace}
\newcommand{\cQ}{\mathcal{Q}\xspace}
\newcommand{\cG}{\mathcal{G}\xspace}

\newcommand{\cC}{\mathcal{C}\xspace}
\newcommand{\cK}{\mathcal{K}\xspace}
\newcommand{\cN}{\mathcal{N}\xspace}

\newcommand{\cO}{\mathcal{O}\xspace}

\renewcommand{\t}{^{\textsf{T}}\xspace}

\newcommand{\CG}{\textsf{CG}\xspace}

\newcommand{\IDR}{\textsf{IDR($s$)}\xspace}
\newcommand{\IDReins}{\textsf{IDR($1$)}\xspace}
\newcommand{\IDRstab}{\textsf{IDR($s$)stab($\ell$)}\xspace}
\newcommand{\GMstabI}{\textsf{GM($s$)stab1}\xspace}
\newcommand{\GMstabII}{\textsf{GM($s$)stab2}\xspace}
\newcommand{\GMstab}{\textsf{GM($s$)stab}\xspace}
\newcommand{\IDRstabp}[2]{\textsf{IDR(#1)stab(#2)}\xspace}
\newcommand{\IDRp}[1]{\textsf{IDR(#1)}\xspace}
\newcommand{\GMRESkp}[1]{\textsf{GMRES(#1)}\xspace}

\newcommand{\Mstabno}{\textsf{$\mathcal{M}$stab}\xspace}

\newcommand{\IDRstabbiortho}{\textsf{IDR($s$)stab($\ell$)biortho}\xspace}
\newcommand{\IDRstabno}{\textsf{IDRstab}\xspace}

\newcommand{\IDRO}{\textsf{IDR}\xspace}

\newcommand{\SRIDR}{\textsf{SRIDR($s$)}\xspace}

\newcommand{\Mstab}{\textsf{$\mathcal{M}$($s$)stab($\ell$)}\xspace}

\newcommand{\BiCG}{\textsf{BiCG}\xspace}

\newcommand{\BiCGstab}{\textsf{BiCGstab}\xspace}
\newcommand{\BiCGstabL}{\textsf{BiCGstab($\ell$)}\xspace}

\newcommand{\GMRES}{\textsf{GMRES}\xspace}

\newcommand{\GCR}{\textsf{GCR}\xspace}

\newcommand{\GCRODR}{\textsf{GCRO-DR}\xspace}

\newcommand{\CGS}{\textsf{CGS}\xspace}
\newcommand{\IDRbiortho}{\textsf{IDR($s$)biortho}\xspace}
\newcommand{\IDRobio}{\textsf{IDR($s$)obio}\xspace}

\newcommand{\hJ}{{{J}}}

\newcommand{\bA}{\textbf{A}\xspace}

\newcommand{\bB}{\textbf{B}\xspace}
\newcommand{\bG}{\textbf{G}\xspace}
\newcommand{\bU}{\textbf{U}\xspace}
\newcommand{\bV}{\textbf{V}\xspace}
\newcommand{\bZ}{\textbf{Z}\xspace}
\newcommand{\bD}{\textbf{D}\xspace}
\newcommand{\bS}{\textbf{S}\xspace}
\newcommand{\bP}{\textbf{P}\xspace}
\newcommand{\bPi}{\boldsymbol{\Pi}\xspace}
\newcommand{\bQ}{\textbf{Q}\xspace}

\newcommand{\bF}{\textbf{F}\xspace}
\newcommand{\bH}{\textbf{H}\xspace}
\newcommand{\bI}{\textbf{I}\xspace}
\newcommand{\bR}{\textbf{R}\xspace}
\newcommand{\bL}{\textbf{L}\xspace}
\newcommand{\bT}{\textbf{T}\xspace}

\newcommand{\be}{\textbf{e}\xspace}
\newcommand{\bw}{\textbf{w}\xspace}
\newcommand{\bM}{\textbf{M}\xspace}

\newcommand{\bC}{\textbf{C}\xspace}

\newcommand{\bW}{\textbf{W}\xspace}
\newcommand{\bY}{\textbf{Y}\xspace}
\newcommand{\bx}{\textbf{x}\xspace}
\newcommand{\bdx}{\textbf{dx}\xspace}

\newcommand{\by}{\textbf{y}\xspace}
\newcommand{\bb}{\textbf{b}\xspace}

\newcommand{\br}{\textbf{r}\xspace}
\newcommand{\bc}{\textbf{c}\xspace}

\newcommand{\bgam}{\boldsymbol{\gamma}\xspace}
\newcommand{\bchi}{\boldsymbol{\chi}\xspace}
\newcommand{\bzeta}{\boldsymbol{\zeta}\xspace}
\newcommand{\bxi}{\boldsymbol{\xi}\xspace}

\newcommand{\btau}{\boldsymbol{\tau}\xspace}
\newcommand{\hbtau}{\hat{\boldsymbol{\tau}}\xspace}
\newcommand{\bEta}{\boldsymbol{\eta}\xspace}

\newcommand{\bq}{\textbf{\textit{q}}\xspace}

\newcommand{\bu}{\textbf{u}\xspace}
\newcommand{\bv}{\textbf{v}\xspace}

\newcommand{\hbU}{\hat{\textbf{U}}\xspace}
\newcommand{\hbV}{\hat{\textbf{V}}\xspace}

\newcommand{\hbY}{\hat{\textbf{Y}}\xspace}
\newcommand{\hbH}{\hat{\textbf{H}}\xspace}
\newcommand{\hbZ}{\hat{\textbf{Z}}\xspace}
\newcommand{\hhbH}{{\hat{\hbH}}\xspace}

\newcommand{\hbx}{\hat{\textbf{x}}\xspace}

\newcommand{\hbr}{\hat{\textbf{r}}\xspace}

\newcommand{\bp}{\textbf{p}\xspace}

\newcommand{\tbA}{\tilde{\textbf{A}}\xspace}

\newcommand{\tbV}{\tilde{\textbf{V}}\xspace}
\newcommand{\tbx}{\tilde{\textbf{x}}\xspace}
\newcommand{\tbb}{\tilde{\textbf{b}}\xspace}
\newcommand{\tbr}{\tilde{\textbf{r}}\xspace}
\newcommand{\tbv}{\tilde{\textbf{v}}\xspace}

\renewcommand{\d}{^{\dagger}\xspace}
\newcommand{\h}{^{\textsf{H}}\xspace}

\newcommand{\bO}{\textbf{0}\xspace}

\newcommand{\nEqns}{{n_{\text{Systems}}}\xspace}

\newcommand{\ubI}{\underline{\bI}\xspace}
\newcommand{\uubI}{\underline{\ubI}\xspace}

\newcommand{\ubH}{\underline{\bH}\xspace}

\newcommand{\cond}{\operatorname{cond}\xspace}
\newcommand{\opspan}{\operatorname{span}\xspace}

\newcommand{\hp}[1]{^{(#1)}\xspace}
\newcommand{\fp}[1]{_{\lbrace#1\rbrace}\xspace}
\newcommand{\ha}{^{(1)}\xspace}
\newcommand{\hb}{^{(2)}\xspace}
\newcommand{\hc}{^{(3)}\xspace}
\newcommand{\hia}{^{(\iota)}\xspace}

\newcommand{\tol}{\mathrm{tol}\xspace}
\newcommand{\tolabs}{{\mathrm{tol}_\text{abs}}\xspace}
\newcommand{\tolrel}{{\mathrm{tol}_\text{rel}}\xspace}
\newcommand{\epsmachine}{\varepsilon_{\textrm{Machine}}\xspace}

	{\end{cases}\right.}

\title{A restarted \GMRES-based implementation of \IDRstab to yield higher robustness\footnote{This text is a former version of my master thesis that has been submitted on July 20, 2017 to the University.}}
\author{Master thesis of Martin P. Neuenhofen\and First supervisor: Professor Georg May (AICES, RWTH)\and Second supervisor: Professor Chen Greif (CS, UBC)
	}

\date{August 6, 2017}
\begin{document}

\maketitle

\begin{abstract}
In this thesis we propose a novel implementation of \IDRstab that avoids several unlucky breakdowns of current \IDRstab implementations and is further capable of benefiting from a particular lucky breakdown scenario.	
\IDRstab is a very efficient short-recurrence Krylov subspace method for the numerical solution of linear systems.

Current \IDRstab implementations suffer from slowdowns in the rate of convergence when the basis vectors of their oblique projectors become linearly dependent.

We propose a novel implementation of \IDRstab that is based on a successively restarted \GMRES method. Whereas the collinearity of basis vectors in current \IDRstab implementations would lead to an unlucky breakdown, our novel \IDRstab implementation can strike a benefit from it in that it terminates with the exact solution whenever a new basis vector lives in the span of the formerly computed basis vectors.

Numerical experiments demonstrate the superior robustness of our novel implementation with regards to convergence maintenance and the achievable accuracy of the numerical solution.
\end{abstract}

\tableofcontents

\newpage

\FloatBarrier
\newpage

\section{Introduction}

In this master thesis we present a novel highly robust implementation of the numerical method \IDRstab \cite{IDRstab-Paper,GBiCGstab}. \IDRO stands for \emph{induced dimension reduction} and denotes in specific a short-recurrence principle for Krylov subspace methods \cite{IDReins}.

\IDRstab is a parametric short-recurrence Krylov subspace method for the solution of symmetric and non-symmetric large sparse (preconditioned) linear systems. The method constitutes a generalisation of a multitude of commonly used iterative methods, such as \CGS \cite[p.\,215]{Saad1}, \BiCGstab \cite[p.\,217]{Saad1}, \BiCGstabL \cite{BiCGstabL}, \IDR \cite{IDR-report} and -- speaking for our novel implementation -- restarted \GMRES \cite[p.\,167]{Saad1}. Originally, \IDRstab was developed to solve one system of linear equations, i.e.
\begin{align}
	\bA \cdot \bx = \bb\,, \label{eqn:Axb}
\end{align}
where $\bA\in\R^{N \times N},\,\bb\in\R^N,\,N \in \N$ are given with $\bA$ regular and numerical values for $\bx$ are sought with a certain accuracy. However, from generalisations such as \SRIDR \cite{Report15} and \Mstab \cite{Mstab-report,Mstab-paper} it is well-known that \IDRstab variants can also be efficiently applied in order to solve sequences of linear systems, such as they occur e.g. in Quasi-Newton methods. This is why in this thesis we also test of our novel implementation of \IDRstabno as a Krylov subspace recycling method for solving sequences of linear systems of the form
\begin{align}
\bA \cdot \bx\hia = \bb\hia\,, \quad\iota =1 ,...,\nEqns\,.\label{eqn:AxbSequence}
\end{align}

\subsection{Outline}
So far, all current implementations of \IDRstab are based on a \GCR approach, where the \enquote{general conjugation} is a biorthogonalisation of the residual. It is known that \GCR methods have an unlucky breakdown when the residual does not change during one iteration \cite[p.\,5, l.\,10]{GCR-breakdown}, whereas \GMRES does never break down before the solution is found \cite[Prop.\,6.10]{Saad1}. This is why in this thesis we propose the first implementation of \IDRstabno that is not based on a \GCR approach but on an interior restarted \GMRES approach.

Further, we utilise the well-conditioned Krylov subspace basis from the robust Arnoldi iteration of the interior restarted \GMRES method in order to construct well-conditioned basis matrices for the oblique projectors that are required for subsequent iterations of \IDRstabno. By this we ensure that several further breakdown scenarios that are related to (nearly) linearly dependent columns in some oblique projection bases cannot occur.

\subsection{Structure}
In Section\,\ref{sec:Motivation:SR-KSS(R)methods} we review and motivate Krylov methods. In Section\,\ref{sec:IntroIDR} we provide an introduction into \IDRO, which is the mathematical theory for short-recurrence methods that the Krylov method \IDRstab is based on.
Afterwards, in Section\,\ref{Sec:DerivationIDRstab} we derive and analyse the algorithm of \IDRstab by using the formerly introduced mathematical theory.

In Section\,\ref{sec:ImplementationsIDRstab} we lay out the implementation of selected variants of \IDRstab and discuss their numerical properties with some remarks on possible breakdown scenarios. At the end of Section\,\ref{sec:ImplementationsIDRstab} we discuss numerical weaknesses of all these implementations and describe which of these weaknesses could be prevented by implementing \IDRstabno in a more sophisticated way.

In Section\,\ref{Sec:NovelRestartedGMRESbasedImpelementationOfIDRstab} we propose our main contribution, namely the restarted \GMRES-type implementation of \IDRstab. This implementation uses different code blocks for two possible choices of the method parameter $\ell$. Both variants and the geometric properties of their computed vector quantities are derived in detail. Afterwards it is discussed how both code blocks can be merged into one combined practical implementation of \IDRstabno with an adaptive choice of the method parameter $\ell$. Eventually we discuss how this new implementation avoids several weaknesses that the aforementioned implementations from Section\,\ref{sec:ImplementationsIDRstab} suffer from.

Before we evaluate the practical performance of our new implementation in Section\,\ref{sec:NumericalExperiments} by numerical experiments we first give a brief introduction into Krylov subspace recycling for \IDRstabno in Section\,\ref{sec:Mstab}. This is mainly done because we want to show numerical experiments where \IDRstabno is utilised in a Krylov subspace recycling approach. 

The numerical experiments in Section\,\ref{sec:NumericalExperiments} consist of non-preconditioned and preconditioned single linear systems and sequences of linear systems.

The conclusion in Section\,\ref{sec:Conclusion} summarises the main results of our thesis and reviews the practical advantages of our new method as found from the numerical experiments.

\subsection{Notation}
Throughout this subsection let $a,b,c \in \N$. $a \geq b$ and $c$ can be of any relation to $a,b$.

$\|\cdot\|$ denotes the Euclidean norm. $\|\bv\|_\bS:= \sqrt{\bv\t\cdot\bS\cdot\bv}$ is the $\bS$-norm for a symmetric positive definite matrix $\bS$. For a matrix $\bB \in \R^{a \times c}$ the condition $\cond(\bB)$ is defined as the extremal ratio of the singular values of $\bB$.

$\bB\d$ is the Penrose-Moore pseudo-inverse. When it is unclear whether the inverse of $\bB$ exists then we do always use pseudo-inverses. We will never use pseudo-inverses of matrices that have more than $10$ columns (so any concern relating to computational complexity is obsolete).
\largeparbreak

This thesis deals with linear algebra in the $N$-dimensional Euclidean space. The theory, especially that of Induced Dimension Reduction, applies to $\C^N$.\footnote{This is inevitable because of the possibly complex relaxation values for $\omega$.} Besides, this theory for $\C^N$ is used in the derivation of numerical algorithms that are presented in this thesis. However, for real systems $\bA \in \R^{N \times N}, \bb\in\R^N$ these algorithms will never compute any data in $\C^N\setminus\R^N$ but only in $\R^N$. This behaviour is desired since computations on numbers in $\R$ can be performed at a lower computational cost than on numbers in $\C$. This is why in our algorithms we use notation for numbers in $\R$ to emphasize on this fact. However, all the methods may easily be applied to complex systems by replacing transposes, denoted for a matrix $\bB$ by $\bB\t$, by Hermitian transposes, denoted by $\bB\h$.
\largeparbreak

In this thesis there occur many linear vector (sub-)spaces. These are denoted by calligraphic capital letters. For spaces we use orthogonal complements. Considering a space $\cB \subset \C^N$ with basis matrix $\bB$, the orthogonal complement of $\cB$ is denoted by $\cB^\perp$ and defined as $\cB^\perp := \cN(\bB)$, where $\cN(\bB)$ is the null-space of $\bB$. The sum of two spaces $\cB = \rg(\bB)$ and $\cC=\rg(\bC)$ is defined as $\cB+\cC := \rg([\bB,\,\bC])$.

We define the \textit{(block) Krylov subspace} for a general matrix $\bB \in \R^{N \times a}$ as $\cK_d(\bA;\bB) = \rg([\bB,\bA\cdot\bB,...,\bA^{k-1}\cdot\bB])$. The index $d$ is called \textit{degree}. For the case $a=1$ this matches the conventional definition. There is always a degree $d\leq N$ for which the (block) Krylov subspace does not grow any more in dimension, called \textit{grade}. This space is called full Krylov subspace and it is denoted by $\cK_\infty(\bA;\bB)$ since in this thesis we do not care about the particular value of the grade.

For the orthogonal complement of (block) Krylov subspaces $\cK_d(\bA;\bB)$ we use the notation $\cK_d^\perp(\bA;\bB):= \big(\cK_d(\bA;\bB)\big)^\perp$. Since the expression $\cK_\ell^\perp(\bA\h;\bP)$ occurs frequently we write it out for better accessibility:
\begin{align*}
					& & &\bu \in \cK_\ell^\perp(\bA\h;\bP) \\
	\Leftrightarrow & & &\bu\h \cdot [\bP,(\bA\h)\cdot\bP,...,(\bA\h)^{(\ell-1)}\cdot\bP] = \bO\\
	\Leftrightarrow & &\bA^k \cdot &\bu \in \cN(\bP)\quad \forall k=0,...,\ell-1\,.
\end{align*}
\largeparbreak

We describe several algorithms in this thesis. In this, we use some standard linear algebra subroutines, for instance the QR-decomposition. We write $[\bQ,\bR] = \texttt{qr}(\bB)$ for $\bB \in \R^{a \times b}$ to denote the reduced QR-decomposition $\bQ \cdot \bR = \bB$, where $\cond(\bQ)=1$ and $\bR \in \R^{b \times b}$ is upper triangular. This decomposition is to be computed by the modified Gram-Schmidt procedure.

For a matrix $\bB \in \R^{b \times a}$ we define the LQ-decomposition $[\bL,\bQ]=\texttt{lq}(\bB)$ as $\bL = \bB \cdot \bQ$, where $\cond(\bQ)=1$ and $\bL \in \R^{b \times b}$ is lower triangular. The LQ-decomposition shall be computed by a (reduced) QR-decomposition of the transpose of $\bB$.

Further, we define some orthonormalisation routines. $\bQ = \texttt{orth}(\bB)$ for $\bB \in \R^{a \times b}$ constructs $\bQ \in \R^{a \times b}$ such that $\rg(\bQ)=\rg(\bB)$ and $\cond(\bQ)=1$. $\bQ$ can be obtained from the QR-decomposition above. We also define the orthogonalisation for rows of matrices. $\bQ = \texttt{roworth}(\bB)$ for $\bB \in \R^{b \times a}$ generates $\bQ$ with orthonormalised columns, i.e. $\cond(\bQ) = 1$ and $\cN(\bQ)=\cN(\bB)$.
\largeparbreak

Let in the following $m,n,\ell,s\in\N$. In the algorithms there will be basis matrices $\bW \in \R^{N \times n}$, $\bV^{(g)} \in \R^{N \times s}$, $g=-1,0,1,...,\ell$ with columns $\bw_q,\bv_q^{(g)}$, $q=1,...,s$\,. $\bW_{:,d:f}$ is the sub-matrix of $\bW$ consisting of $[\bw_d,...,\bw_f]$, the columns from $d$ to $f \in \N$. We use the same sub-index notation for each $\bV\hp{g}$. In general, for any matrix $\bM \in \R^{m \times n}$ the notation $\bM_{i:j,d:f}$ refers to the sub-matrix of rows $i$ to $j$ and columns $d$ to $f$. If one of the letters $i,d$ respectively $j$ or $f$ is missing then this is a short writing for the number $1$ respectively $m$ or $n$. For example, $\bM_{i:,:f}$ is the submatrix of $\bM$ of rows $i$ to $m$ and columns $1$ to $f$.

\section{Motivation of short-recurrence Krylov subspace (recycling) methods}\label{sec:Motivation:SR-KSS(R)methods}
\IDRstab is a short-recurrence Krylov subspace method \cite{IDRstab-Paper}. It is a very particular method for solving large sparse symmetric and non-symmetric systems of linear equations. In this section we provide an overview of some important classes of methods for solving large sparse systems and explain which role \IDRstab takes among these available methods.

This section has the purpose to motivate short-recurrence Krylov subspace methods and short-recurrence Krylov subspace recycling methods. It is not an introduction into iterative methods and does not review the mathematical foundations of iterative methods and Krylov methods. These can be found e.g. in \cite{Saad1}.

\subsection{Motivation of matrix-free methods}
All the methods presented in this thesis are matrix-free methods. Matrix-free methods are one class of methods for solving linear systems like \eqref{eqn:Axb}. 

As a characteristic property, in {matrix-free} methods the system matrix $\bA$ is only accessed by applying matrix-vector products $\bv:=\bA\cdot\bu$ with it for given vectors $\bu \in \C^N$ \cite{Matrixfree-method}. The benefit of matrix-free methods for large sparse systems is that matrix-vector products can be computed in $\cO(N_\text{nnz})$, where $N_\text{nnz}$ is the number of non-zeros of $\bA$, whereas computing a complete factorization might still require a time complexity in $\cO(N^3)$ and necessity of storage for $\cO(N^2)$ numbers. A matrix-free method in contrast requires at most storage for $\cO(N_\text{nnz})$ numbers (plus storage for the method itself).

\subsection{Motivation of iterative methods}
Iterative methods are the antipode to direct methods. A direct method is for instance Gaussian elimination: It decomposes $\bA$ in one expensive computation into triangular factors and computes forward and backward substitutions. In iterative methods instead, an approximate solution is successively improved in terms of numerical accuracy by repeatedly applying a cheap (compared to Gaussian elimination) computational step.

A (non-generic) framework of an iterative method is:
\begin{algorithmic}[1]
\State Given is the problem $\bA,\bb$ and an initial guess $\bx_0$ for $\bx$.
\For{$j=0,1,2,3,...$}
	\State $\bx_{j+1} := \Phi(\bA,\bb;\bx_j)$
	\If{$\bx_{j+1}$ is accurate enough}
		\State \Return $\bx_{j+1}$
	\EndIf
\EndFor
\end{algorithmic}
In this, $\Phi$ is a computationally cheap function that computes from a given approximation $\bx_j$ to the exact solution $\bx^\star := \bA^{-1} \cdot \bb$ a more accurate approximation $\bx_{j+1}$. A convergence theory for this framework is provided by Banach's fixed point theorem \cite[p.\,414]{Saad1}.

In order to determine whether the currently found numerical approximate solution $\bx$ for the exact solution $\bx^\star$ is accurate enough, a commonly used upper error estimate is
\begin{align*}
\underbrace{\frac{\|\bx^\star-\bx\|}{\|\bx^\star\|}}_\text{relative error} \leq \underbrace{\|\bA\|\cdot\|\bA^{-1}\|}_{\cond(\bA)} \cdot \underbrace{\frac{\|\br\|}{\|\bb\|}}_{\leq \tolrel}\,,\tageq\label{eqn:RelErrEst}
\end{align*}
in which the \textit{residual} $\br = \bb-\bA \cdot \bx$ is the offset in the right-hand side that is caused by the gap between the inaccurate solution $\bx$ and the exact solution $\bx^\star$.

In \eqref{eqn:RelErrEst}, a guess for the right-hand side can be found by estimating the condition number of $\bA$ from, e.g., the rate of convergence of the method itself (since it usually converges faster for a better condition and since further the condition can be estimated from that of the Hessenberg matrix) or problem characteristics of the system to be solved. The residual on the other hand can be directly computed from the numerical solution $\bx$.

The second factor in the right-hand side is called \textit{relative residual norm}. A small positive tolerance value $\tolrel \in \R^{+}$ for the relative residual norm can be used as a stopping criterion for the above for-loop. $\tolrel$ must be chosen with respect to the desired solution accuracy and $\cond(\bA)$. From $\tolrel$ an absolute tolerance $\tolabs$ for $\|\br\|$ can be computed.

Iterative methods can be cheaper than direct methods when only relatively few iterations are required to drive $\|\br\|\leq\tolabs$. This is the case if at least one of the following properties holds: 
\begin{itemize}
\item The initial guess is already accurate.
\item Only a coarse accuracy is required.
\item The average reduction in $\|\br\|$ per iteration is huge.
\end{itemize}
Potential ways to achieve the above properties are discussed in the next subsection.

\subsection{Preconditioners and initial guess}
Iterative methods often converge in fewer iterations to accurate solutions when the system matrix is better conditioned. Further to that, it is beneficial when $\tolabs$ is not much smaller than the initial residual norm.

In order to achieve both of that for a given system \eqref{eqn:Axb}, one can use an initial guess $\bx_0 \in \R^N$ and left and right regular preconditioners $\bM_L,\bM_R \in \R^{N\times N}$ to obtain
\begin{align}
\underbrace{\bM_L^{-1} \cdot \bA \cdot \bM_R^{-1}}_{=:\tbA} \cdot \by = \bc\,, \label{eqn:precAxb}
\end{align}
where $\bc = \bM_L^{-1}\cdot(\bb - \bA \cdot \bx_0)$ and the final numerical solution for $\bx$ can be reconstructed from $\bx = \bM_R^{-1} \cdot \by + \bx_0$. In this situation, the iterative method is applied to solve for $\by \in \R^N$.

When $\bM_L \cdot \bM_R \approx \bA$ then $\cond(\tbA)\approx 1$ and there are good chances that the iterative method can yield huge accuracy improvements with each iteration. Further, when $\|\bb-\bA\cdot\bx_0\|$ is small (i.e. $\bx_0$ is a good guess for the solution) then $\|\bc\|$ is already close to the required absolute tolerance. Both potentially leads to a huge reduction of the required number of iterations.

In order to use preconditioning one applies the matrix-free method on \eqref{eqn:precAxb} instead of \eqref{eqn:Axb} and evaluates the matrix-vector-product with $\tbA$ by a subsequent application of its factors.

For the ease of presentation in the remainder of this thesis we will often use the initial guess $\bx_0=\bO$ without loss of generality since this refers to the initial guess of $\by$ then.
\largeparbreak
In the following three subsections we motivate Krylov subspace methods and further to that Krylov subspace recycling methods \cite{RGCRO,GCRO-DR,PGCRO-DR,RBiCG,Compressing} from basic iterative methods.

\subsection{Basic iterative methods}
Basic iterative methods are matrix-free iterative methods $\Phi$ that use only $\bA,\bb$ and the current numerical solution $\bx$ (and the residual $\br=\bb-\bA\cdot\bx$).

Without any exception, basic iterative methods can be expressed by the update formula
\begin{align*}
	\bx_{j+1} := \bx_j + \omega_j \cdot \br_j\,,\tageq\label{eqn:PrimitiveUpdate}
\end{align*}
which is called \textit{Richardson's method} \cite{RichardsonIteration}. In this, $\bx_j$ is the current solution guess, $\bx_{j+1}$ is the improved solution guess and $\br_j$ is the residual of $\bx_j$. $\omega_j$ is a scalar that is called relaxation parameter, stabilisation coefficient or step size, cf. below. I.e. $\Phi(\bA,\bb;\bx) := \bx + \omega \cdot (\bb - \bA \cdot \bx)$.
\largeparbreak
Other common basic iterative methods are Jacobi's method and the Gauss-Seidel method. Splitting the system matrix $\bA$ into its strictly lower triangular matrix $\bL$, the diagonal matrix $\bD$ and the strictly upper triangular matrix $\bU$ such that
\begin{align*}
	\bA = \bL + \bD + \bU
\end{align*}
the methods fit into the above update scheme of Richardson's method by choosing the following preconditioners and step sizes: For Jacobi's method choose $\omega=1$ and $\bM_L = \bD$. For the Gauss-Seidel method choose instead $\bM_L = \bL + \bD$. For $\omega>1$ the methods are called over-relaxated, for $\omega<1$ under-relaxated \cite{SORmethod}. Also a multi-grid method fits into Richardson's update scheme by defining the preconditioner in such a way that $\br$ is the multi-grid error correction step from the residual of the non-preconditioned system.

\subsection{Motivation of Krylov subspace methods}\label{sec:Motivation_KrylovMethods}
One could try to motivate Krylov subspace methods from the perspective that they compute a somehow optimal solution in the Krylov subspace $\cK_j(\bA;\bb)$, whereas every basic iterative method only computes some solution in $\cK_j(\bA;\bb)$. However, this motivation would not work since most Krylov methods do not compute solutions in $\cK_j(\bA;\bb)$ that are by any means optimal. Consider for instance \BiCG \cite[p.\,211]{Saad1}, that computes biorthogonal residuals. Such residuals do not form the solution to the minimisation of any error norm.

In the following we motivate Krylov subspace methods as projection methods because this causes their superlinear convergence property, which is their actual practical advantage. To this end, the following passages and images used therein are quoted from our paper \cite{Mstab-paper}.
\largeparbreak
{Krylov subspace methods} do still use updates that are either identical or very similar to that of Richardson's method.
However, the advantage of Krylov subspace methods such as, e.g., Conjugate Gradients (\CG) over primitive iterative methods is that with each iteration they also iteratively project the linear system onto a smaller system\footnote{The original solution is then $\bx=\tbx + (\bA^{-1} \cdot \bC) \cdot \bC^\textsf{H}\cdot (\bb-\bA\cdot\tbx)$.} \cite{GCRO-DR}:
\begin{align}
\underbrace{(\bI-\bC\cdot\bC^\textsf{H})\cdot\bA}_{=:\tbA}\cdot\tbx=(\bI-\bC\cdot\bC^\textsf{H})\cdot\bb\,. \label{eqn:ProjSys}
\end{align}
In this, $\bC \in \C^{N \times c}$ is a matrix of orthonormal basis vectors of a \emph{Petrov-space} $\cC \subset \C^N$ that grows iteratively during the iterations. For example, \CG uses the spaces $\cC_0=\lbrace\bO\rbrace$, $\cC_{k+1}=\cC_k+\opspan\lbrace\br_k\rbrace$, $k=0,\,1,\,2,\,...$\ \cite{Saad1}.

Using the above projection approach is equivalent to computing $\bx_j \in \C^N$ such that
\begin{align}
\bx_j &\in \cK_j(\bA;\bb)\quad \wedge \quad \br \perp \cC_j\,. \label{eqn:KrylPrinc}
\end{align}

\begin{figure}
	\centering
	\includegraphics[width=0.55\linewidth]{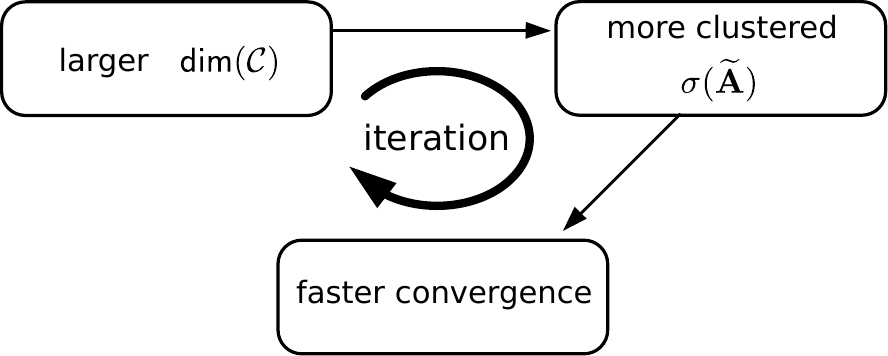}
	\caption{Working principle of superlinear convergence.}
	\label{fig:SuperlinConv}
\end{figure}

In consequence of the projection \eqref{eqn:ProjSys} the singular values of the system matrix $\tbA$ cluster. This in turn improves the rate of convergence of the Richardson iteration that is inherited in the Krylov method. Fig.~\ref{fig:SuperlinConv} depicts this: With each iteration $k$ the dimension $\dim(\cC)$ grows. This leads to a stronger clustering of the singular values $\sigma(\tbA)$, which in turn in the long term improves the rate of convergence for the inherited update scheme \eqref{eqn:PrimitiveUpdate}.
Thus, at a certain stage the rate of convergence probably becomes only faster and faster, which is referred to as superlinear convergence \cite{superlinGMRES}.

Certainly, Krylov subspace methods have their issues. For example, superlinear convergence is not proven in the general case. Besides, in the non-symmetric case there is no method such as \CG that both minimises an error norm and uses only short recurrences \cite{Saad1}. Nevertheless, Krylov methods are considered as one of the most important classes of algorithms \cite{SIAM-top10}.

\subsection{Motivation of Krylov subspace recycling methods}\label{sec:MotivationKSSR}
Krylov subspace recycling methods try to extend the projection approach of Krylov subspace methods for a sequence of more than only one linear system.

Consider solving one after the other the two linear systems
\begin{align}
\bA \cdot \bx\ha=\bb\ha,\quad\quad\bA\cdot\bx\hb=\bb\hb\,. \label{eqn:b1b2}
\end{align}
When solving the first system with a Krylov method in $\star$ iterations, it iteratively builds a Petrov space $\cC\ha_\star$. Now, for solving the second system, there are two options:
\begin{enumerate}[(i)]
	\item Solve for $\bx\hb$ with a Krylov method. However, using a conventional Krylov method for the second system means that we throw $\cC\ha_\star$ and iteratively build a new Petrov space, starting from $\cC_0\hb = \lbrace \bO \rbrace$.
	\item It seems more desirable to have a method that starts the solution of the second system with a Petrov space $\cC_0\hb = \cC\ha_\star$. This is what Krylov subspace recycling methods would ideally do, of course using only short recurrences.
\end{enumerate}
The hope of using a Krylov subspace recycling method is the following: Since the dimension of $\cC$ and thus the projection in \eqref{eqn:ProjSys} is larger right from the beginning, there is hopefully an earlier occurrence of superlinear convergence. The earlier occurrence of a fast rate of convergence in turn leads to a reduction of the number of iterations that is required to achieve the desired solution accuracy.

\begin{figure}
\centering
\includegraphics[width=0.6\linewidth]{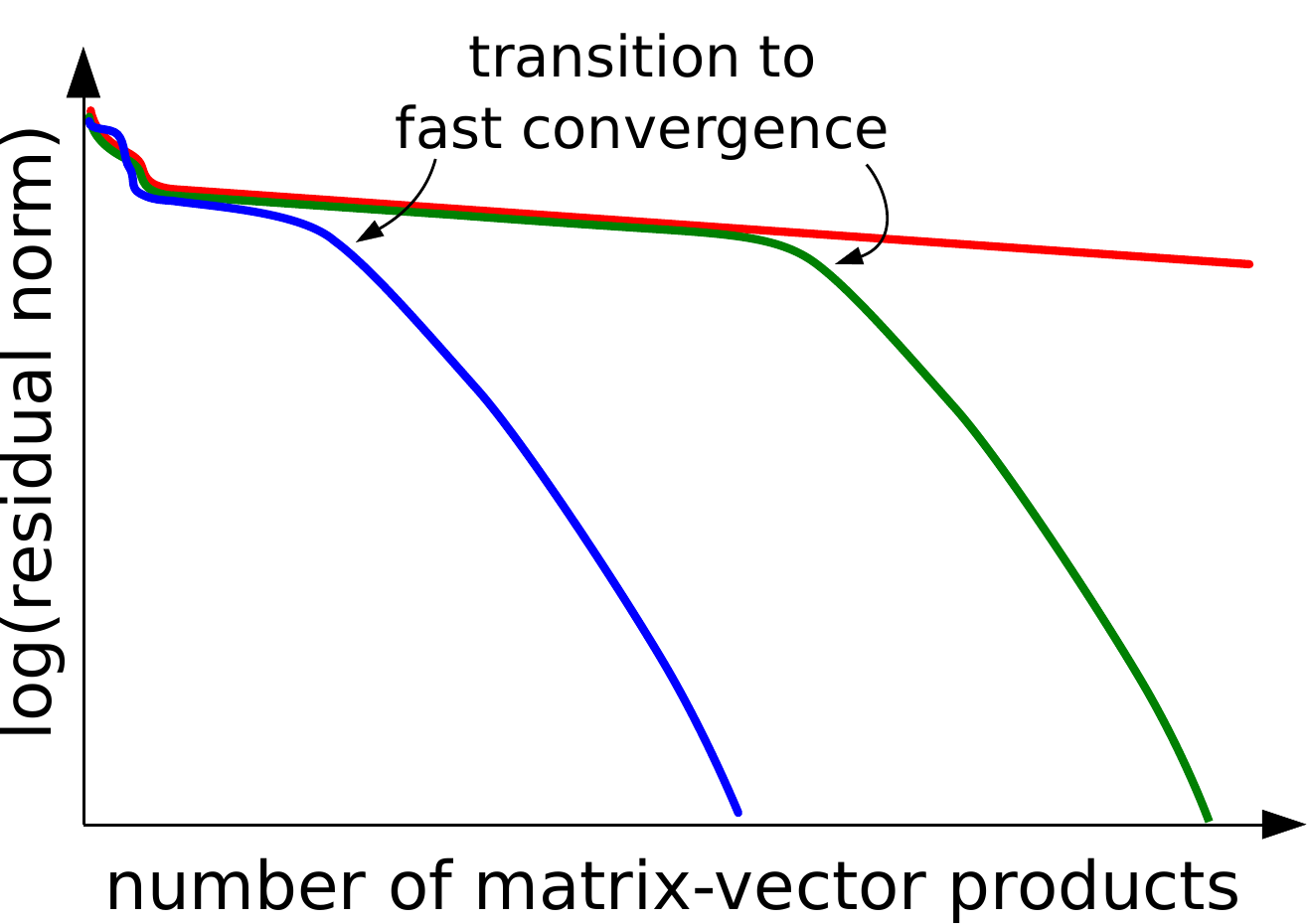}
\caption{Principlal convergence behaviour of a basic iterative method (red), a Krylov subspace method (green) and a Krylov subspace recycling method (blue).}
\label{fig:superlinearConvergenceCurves}
\end{figure}

\largeparbreak
Here ends the citation from our paper \cite{Mstab-paper}. In the following we discuss and sketch in Fig.\,\ref{fig:superlinearConvergenceCurves} the intended convergence behaviour of a Krylov method and a Krylov subspace recycling method in comparison to a basic iterative method:

Considered the case that the basic iterative scheme converges at all it achieves in the limit for a large number of iterations only a linear rate of convergence. Thus, plotting the logarithm of the residual norm over the number of computed matrix-vector products, one obtains -- apart from some erratic behaviour in the beginning -- a straight line. This is sketched in the figure by the red curve.

For the Krylov method during a long fraction of the iterations the convergence behaviour is identical to that of the basic iterative scheme. However, eventually there is a point at which the projection of the linear system provides a sufficiently strong clustering of $\sigma(\tbA)$. This has then a similar effect on the rate of convergence as preconditioning. In the figure this convergence behaviour is sketched by the green curve and the point where the convergence improvement happens is marked by an arrow.

For a Krylov subspace recycling method the convergence behaviour is similar to that of a Krylov subspace method, however the effectively solved system is projected onto a much smaller spatial dimension. This is because the dimension of the space $\cC$ is much larger right from the start. As a consequence of this, the initial erratic behaviour differs because the geometry of the problem is now massively changed. As a further consequence, the transition to fast convergence is encountered earlier since there is already a stronger clustering of $\sigma(\tbA)$ right from the start.

\subsection{Short-recurrence vs. long-recurrence Krylov subspace methods}
In order to realise the iterative projection of the original system \ref{eqn:Axb} onto a smaller system \eqref{eqn:ProjSys} the step direction for the update of $\bx$ is not the residual $\br$ itself but a projected residual $\br - \bC \cdot (\bC\t \cdot \br)$. $\bC$ is a matrix depending on $\cC_j$. For the choice of $\cC_j$, there are in general two classes of Krylov methods, cf. in \cite{Saad1} the introduction of chapter 6:
\begin{enumerate}[(A)]
	\item Long-recurrence methods, that require $\cO(j)$ vector operations (AXPYs and DOTs) to compute the projected residual. These methods are called \textit{long-recurrence methods}.
	\item Short-recurrence methods, that require only $\cO(1)$ of such vector operations to compute the projected residual. These are called \textit{short-recurrence methods}.
\end{enumerate}
One could think of two extremes of sophisticated Krylov subspace methods. The first kind would use long recurrences but ensure that the images of the update directions, i.e. $\bA\cdot (\bI-\bC\cdot\bC\t)\cdot\br_j$, are orthogonal for distinct values of $j$ because then the Krylov principle yields that $\omega_j$ is chosen such that $\|\br_{j+1}\|$ is minimised. So after each iteration the residual is as small as possible. Methods that realise this principle are \GMRES \cite[p.\,158]{Saad1} and \GCR \cite[p.\,183]{Saad1}.

The other extremal would choose $\cC_j$ such that in every step $(\bI - \bC \cdot \bC\t)\cdot\br = \br$ holds because then the projected residual requires no vector operations at all. However, in this case the pairwise angles between the vectors $\bA\cdot (\bI-\bC\cdot\bC\t)\cdot\br_j$ for $j=1,...,N$ are far from orthogonal. As a consequence of this, the Krylov principle would yield that the linear coefficients $\omega_j$, that must be used to combine the residuals $\br_1,...,\br_N$ to the exact solution $\bx^\star$, would not turn out to satisfy that the intermediately found iterates $\bx_j:=\sum_j \omega_j \cdot \br_j$ $\forall j=1,...,N$ are by any means optimal approximations in $\cK_j(\bA;\bb)$ to $\bx^\star$. This is formally proven \cite{FaberManteuffel}.
\largeparbreak
Neither of the above approaches is practical: Long recurrences require storage for $\cO(j)$ column vectors, which is impractical when many iterations $j$ are required or if the system is so large in terms of $N$ that only a few dozen column vectors can be stored. 

On the other hand, searching merely for short recurrences (or as above, recurrences of length zero) without ensuring that the angles between subsequent residuals do not become too small would result in methods with ridiculously large intermediate residuals.\footnote{This statement does only hold for the general, i.e. non-symmetric, case.}

An accepted compromise in practice seems to be to use a stabilised short-recurrence Krylov subspace method such as \BiCGstab \cite[p.\,217]{Saad1}. In such a method, those matrix-vector products that must be computed not for projections but only to make the method transpose-free \cite[p.\,214--221]{Saad1} are utilised to compensate for the potential growth of intermediate residual norms.

In fact, \BiCGstab is the most commonly used short-recurrence Krylov subspace method for general systems \cite[abstract]{IDRstab-Paper}. In the remainder of this thesis we will deal with a parametric generalisation of \BiCGstab that is called \IDRstab \cite{IDRstab-Paper}. \IDRO stands for Induced Dimension Reduction and is a projection theory for short-recurrence Krylov subspace methods \cite{IDR-report}. In \cite{SRPCR,Mstab-report} this theory has been generalised in two unrelated ways to Krylov subspace recycling methods.
\largeparbreak
In the next section we provide the reader with a rigorous introduction into the ideas and the mathematical theory of \IDRO. Eventually we review the method \BiCGstab as a particular \IDRO method.

\section{Introduction into Induced Dimension Reduction}\label{sec:IntroIDR}
In \eqref{eqn:KrylPrinc} a conditional principle for defining a Krylov subspace method is given. We have discussed above why this principle leads to superior convergence over a basic iterative scheme. In this last subsection we have discussed that the choice of the space $\cC$ has some non-trivial effects on the length of the recurrences and especially on the stability of the iterative methods since it affects the step-sizes $\omega$ and the magnitudes of the intermediate residual norms.

The theory of \IDRO uses an equivalent principle to \eqref{eqn:KrylPrinc} by not defining $\cC_j$ but its orthogonal complement $\cG_j:= \cC_j^\perp$, which is called \textit{Sonneveld space} \cite{IDR-report,IDR-Gutkn,IDRbreakdown}. In this setting, one searches $\bx_j \in \C^N$ such that
\begin{align*}
	\bx_j \in \cK_j(\bA;\bb)\quad\land\quad\br_j \in \cG_j\,.
\end{align*}
For increasing values of $j$, the residual is no longer orthogonalised with respect to a Petrov space $\cC_j$ of growing dimension. Instead, the residual is restricted into a Sonneveld space $\cG_j$ of shrinking dimension. In \cite{IDR-as-PetrovGalerkinMethod} the authors show that both approaches, i.e. $\br_j \in \cG_j$ and $\br_j \perp \cC_j$, are equivalent, by deriving a formula for $\cC_j$. In \cite{Mstab-report} we went further and derived explicitly a recursive update formula for $\cC_{j+1}$ from the recursive update formula of $\cG_{j+1}$ that is introduced in a subsequent subsection.
\largeparbreak
\IDRO methods are stabilised short-recurrence Krylov subspace methods that iteratively build a sequence of iterates $\bx_j$, $j=0,1,2,...$ with residuals $\br_j$ that live in Sonneveld spaces $\cG_0\supset\cG_1\supset\cG_2\supset...$ of shrinking dimension. In particular it holds $\br_j \in \cG_j$ for $j=0,1,2,...$\,. The spaces are designed such that for a particular value of $j\in\N$ it holds $\cG_j = \lbrace\bO\rbrace$. At this iterate the numerical solution $\bx_j$ is accurate since $\br_j=\bO$ follows from $\br_j \in \cG_j$.

\subsection{Motivation of \IDRO methods}
In this subsection we motivate \IDRO methods by describing their expected superlinear convergence behaviour and discussing some properties of Sonneveld spaces.
\largeparbreak

We have discussed above that \IDRO methods fit into the framework of Krylov subspace methods. By this we mean that orthogonalising the residual against a growing Petrov space is equivalent to restricting it into a shrinking Sonneveld space. For methods based on the Petrov principle we have discussed why we can expect them to converge in a superlinear fashion such as sketched in Fig.\,\ref{fig:superlinearConvergenceCurves}.

\begin{figure}
\centering
\includegraphics[width=0.9\linewidth]{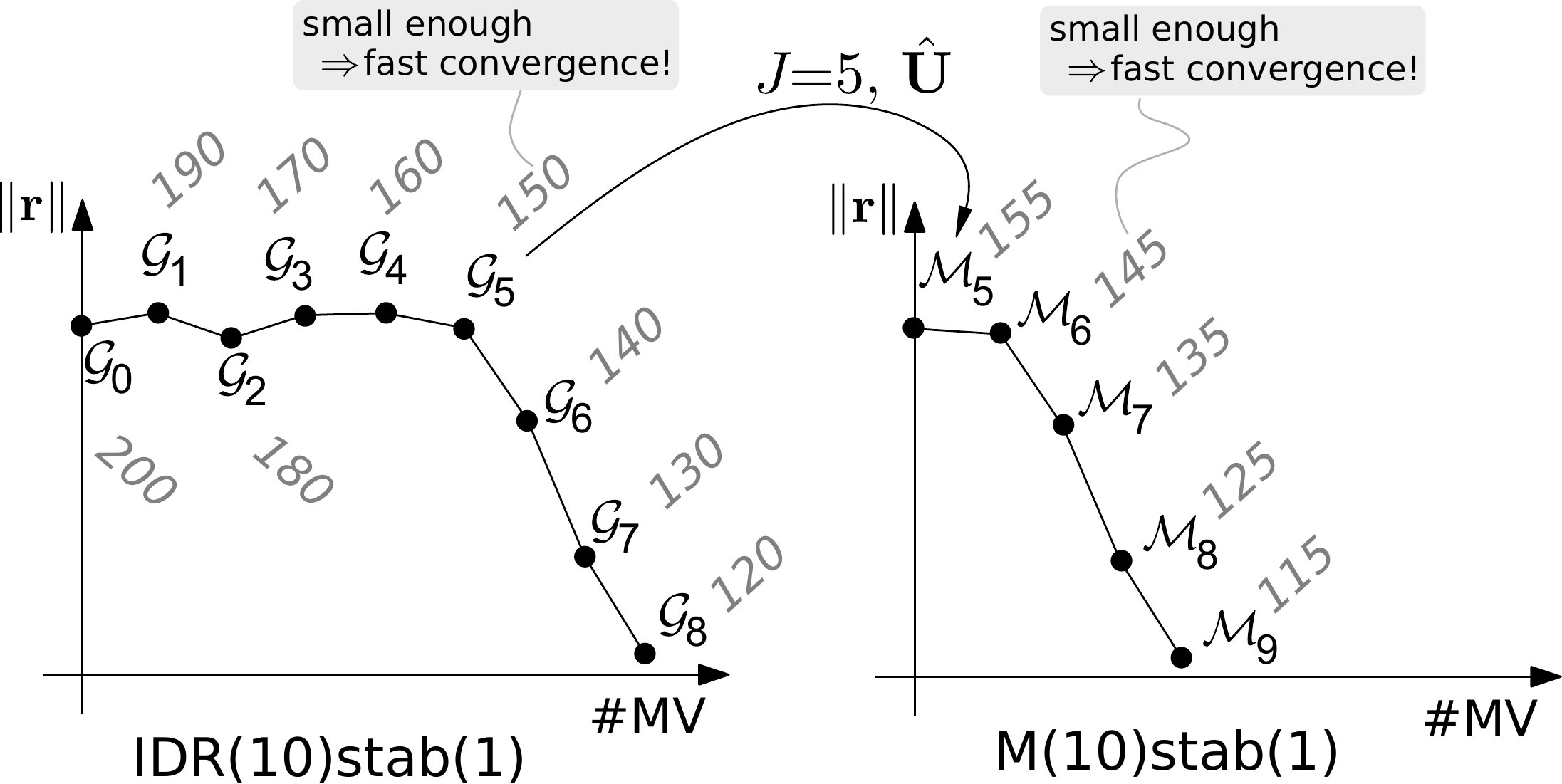}
\caption{Sketch of the principal convergence behaviour of \IDRstab (left) and \Mstab (right).}
\label{fig:MstabMotivation}
\end{figure}

Consequently, we can also expect that \IDRO methods achieve superlinear convergence. Fig.\,\ref{fig:MstabMotivation} shall illustrate this: For an \IDRO method we can assign a residual of its convergence graph to a Sonneveld space of a particular dimension. The dimension of the Sonneveld space is shown in the figure by grey italic numbers. Once the dimension of the Sonneveld space is sufficiently small, the spectrum $\sigma(\tbA)$, where $\tbA$ is again the projected operator onto $\cC^\perp \equiv \cG$, becomes clustered and causes a fast rate of convergence.

In the right part of the figure we sketched the expected convergence behaviour of \Mstab \cite{Mstab-report,Mstab-paper}, which is a Krylov subspace recycling method and an \IDRO method. It is based on the principle to start the iterative scheme of \IDRO from an initial Sonneveld-like space that has a smaller dimension. This leads to a stronger clustering of $\sigma(\tbA)$ right from the beginning and thus to an earlier occurrence of superlinear convergence. The figure illustrates this: A fictional scenario is considered where fast convergence is achieved when the size of the Sonneveld(-like) space is $\leq 150$ dimensions. \Mstabno starts with a residual that lives already in a much smaller Sonneveld-like space, thus it will obtain a fast rate of convergence at an earlier iteration.

\subsection{Mathematical theory of \IDRO}
In this subsection we introduce the spaces and theorems that \IDRO methods are based on. In the subsequent subsection we will derive \BiCGstab as an \IDRO method by using this theory.
\largeparbreak

Originally, \IDRO methods are based on the following spaces \cite{IDReins,IDR-report}.
\begin{Definition}[Sonneveld spaces]\label{def:Sonneveldspace}
	Given $\bA \in \C^{N \times N}$, $\bb \in \C^N$, $\bP \in \C^{N \times s}$, $\rk(\bP)=s$, $\lbrace \omega_j \rbrace_{j \in \N} \subset \Cno$. We define the following sequence of vector spaces.
	\begin{align*}
	\cG_0 &:= \cK_\infty(\bA;\bb)\\
	\cG_{j+1} &:= (\bI - \omega_{j+1} \cdot \bA) \cdot \big(\cG_{j} \cap \cN(\bP)\big)\quad \forall j \in \No\,.\tageq\label{eqn:Sonneveldspace}
	\end{align*}
	The space $\cG_j$ is called Sonneveld space of degree $j$.
\end{Definition}

Sonneveld spaces are nested, i.e. $\cG_0\supset\cG_1\supset\cG_2\supset...$\,. Further, as mentioned above, for an increasing degree of $j\in\No$ the dimension of the Sonneveld space becomes smaller. The following theorem states this.

\begin{Theorem}[\IDRO Theorem]\label{theo:IDR-theo}
	Given the Sonneveld spaces $\lbrace\cG_j\rbrace_{j\in\No}$ from $\bA \in \C^{N \times N}$, $\bb \in \C^N$, $\bP \in \C^{N \times s}$, $\lbrace \omega_j \rbrace_{j \in \N} \subset \Cno$. If $\cG_0 \cap \rg(\bP)$ do not share a non-trivial invariant subspace of $\bA$ then it holds for all $j \in \N$:
	\begin{align*}
	\text{(a)}& & \cG_{j} 		& \subset \cG_{j-1}\\
	\text{(b)}& & \dim(\cG_{j}) & \leq \max\lbrace\, 0,\,\dim(\cG_{j-1}) - s\, \rbrace\,.
	\end{align*}
	\underline{Proof}:\\
	Proposition (a). Proof by induction.
	\begin{enumerate}
		\item Basis: Since $\bA \cdot \cG_0 \subset \cG_0$ it follows $\cG_1 \subset \cG_0$.
		\item Hypothesis: $\exists j\in\N \ : \ \cG_j \subset \cG_{j-1}$.
		\item Induction step: $\cG_j \subset \cG_{j-1} \Rightarrow \cG_{j+1} \subset \cG_j$ is shown.\\
		Choose an arbitrary $\bx \in \cG_{j+1}$. It follows:
		\begin{align*}
		\exists \by \in \cG_j \cap \cN(\bP)\ : \ \bx = \by - \omega_{j+1} \cdot \bA \cdot \by
		\end{align*}
		From the induction hypothesis follows $\by \in \cG_{j-1} \cap \cN(\bP)$, thus $\tbx \in \cG_j$, where
		\begin{align*}
		\tbx := \by - \omega_{j} \cdot \bA \cdot \by\,.
		\end{align*}
		Since $\bx \in \opspan\lbrace\by,\tbx\rbrace \subset \cG_j$ it follows $\bx \in \cG_j$.
	\end{enumerate}
	Proposition (b) follows from Lemma\,\ref{lem:GenRekSonneveldspace}, which is presented later. Using Lemma 2 and assuming further for simplicity that $\rk(q_{0,j}(\bA))=N$ $\forall\,j \in \No$ holds (where $q_{0,j}$ is a polynomial of degree $j$), it follows:
	\begin{align*}
	\dim(\cG_j) &= \dim\Big(q_{0,j}(\bA) \cdot \big(\cG_0 \cap \cK^\perp_j(\bA\h;\bP)\big) \Big)\\
	&=\dim\Big(\cG_0 \cap \cK^\perp_j(\bA\h;\bP)\Big) = \max\lbrace 0\,,\ \dim(\cG_0) - j \cdot s\rbrace \quad\forall\,j\in\No\,.
	\end{align*}
	In this, the third equality follows from the mild condition that $\cG_0$ and $\rg(\bP)$ do not share a non-trivial invariant subspace of $\bA$. Consequently, the dimension is reduced by $s$ for each degree.
	
	When the simplifying assumption $\rk(q_{0,j}(\bA))=N$ $\forall\,j \in \No$ does not hold then the second equality symbol becomes \enquote{$\leq$}. This means that the above presented result for the dimension reduction is only a sharp lower bound, i.e. the worst case. $\boxtimes$
\end{Theorem}

To design an algorithm that produces iterates $\bx_j$ with residuals $\br_j \in \cG_j$ for $j=0,1,2,...$\,, the following template of two computational steps can be used. It provides a scheme to restrict one vector from $\cG_j$ into $\cG_{j+1}$:
\begin{enumerate}[(1)]
	\item Biorthogonalisation: Obliquely project one vector from $\cG_j$ with other vectors from $\cG_j$ into the null-space of $\bP$.
	\item Polynomial step: Multiply the projected vector by $(\bI - \omega_{j+1}\cdot\bA)$ from the left to move it from $\cG_{j} \cap \cN(\bP)$ into $\cG_{j+1}$.
\end{enumerate}

In the following subsection we describe in all detail on the example \BiCGstab how an \IDRO method can be constructed in concrete terms.

\subsection{Derivation of \BiCGstab as an \IDRO method}

\BiCGstab is the simplest \IDRO method \cite{BiCGstab-is-IDReins}\footnote{To avoid any confusion: The reference only says that \BiCGstab is an \IDRO method.}. In this subsection we derive it from the above theory. To this end we first present the algorithm and then lay out with some illustrations the geometric properties of its computed quantities. Just before we do this, we introduce a simple lemma that helps the presentation.
\begin{Lemma}[Remaining of vectors in $\cG_j$]\label{Lem:RemainingVectorsInGj}
	Given the Sonneveld spaces $\lbrace\cG_j\rbrace_{j\in\No}$ from $\bA \in \C^{N \times N}$, $\bb \in \C^N$, $\bP \in \C^{N \times s}$, $\lbrace \omega_j \rbrace_{j \in \N} \subset \Cno$. It holds for all $j \in \No$:
	\begin{align*}
		\bA \cdot \big(\cG_j \cap \cN(\bP)\big) \subset \cG_j\,.
	\end{align*}
	\underline{Proof:}\\
	Choose an arbitrary $\bx \in \cG_j\cap\cN(\bP)$. Then
	\begin{align*}
		(\bI - \omega_{j+1} \cdot \bA) \cdot \bx = \underbrace{\bx}_{\in \cG_j} - \omega_{j+1} \cdot \bA \cdot \bx \in \cG_{j+1} \subset \cG_j
	\end{align*}
	Consequently, $\bA \cdot \bx \in \cG_j$. $\boxtimes$
\end{Lemma}
\largeparbreak
Now we explain the algorithm of \BiCGstab. For this purpose, Alg.\,\ref{Algo:BiCGstab_naive_indexed} shows one implementational variant of \BiCGstab. This variant is formulated with a for-loop from line 6 to 30 that contains the actual iterative scheme.

Before the start of that scheme, from lines 2 to 4 the following quantities are computed: An intial approximate solution $\bx\hp{0}\fp{0}$ and its according residual $\br\hp{0}\fp{0}$, a so-called \textit{auxiliary vector} $\bv\hp{0}\fp{0}$ and its pre-image $\bv\hp{-1}\fp{0} = \bA^{-1} \cdot \bv\hp{0}\fp{0}$ and a matrix $\bP \in \R^{N \times 1}$. $\bP$ plays a role in the definition of the Sonneveld spaces that the method uses. From the way how $\br\hp{0}\fp{0},\bv\hp{0}\fp{0}$ are initialised, it follows the property in line 5.

\begin{algorithm}
	\begin{algorithmic}[1]
		\Procedure{BiCGstab}{$\bA,\bb,\tolrel$}
		\State $\br\fp{0}\hp{0} := \bb$,\quad $\bx\fp{0}\hp{0} := \bO$,\quad$\bv\fp{0}\hp{-1} := \bb$,\quad$\tolabs:=\tolrel\cdot\|\bb\|$
		\State $\bP := $\texttt{orth}$(\texttt{randn}(N,1))$,\quad $\bv\hp{0}\fp{0} := \bA \cdot \bv\hp{-1}\fp{0}$
		\State $\bZ := \bP\t\cdot\bv\hp{0}_0 \in \R^{1 \times 1}$
		\State \textit{// $\br\hp{0}\fp{0},\bv\hp{0}\fp{0} \in \cK_\infty(\bA;\bb) \equiv \cG_0$}
		\For{$j=0,\,1,\,2,\,...$}
		\State \textit{// $\br\hp{0}\fp{j},\bv\hp{0}\fp{j} \in \cG_j$}
		\State \textit{// - - - Biorthogonalisation - - - }
		\State \textit{// Residual}
		\State $\bxi := \bZ\d \cdot (\bP\t\cdot\br\fp{j}\hp{0})$
		\State $\tbr\fp{j}\hp{0} := \br\fp{j}\hp{0} - \bv\fp{j}\hp{0} \cdot \bxi$,\quad$\tbx\fp{j}\hp{0} := \bx\fp{j}\hp{0} + \bv\fp{j}\hp{-1} \cdot \bxi$\quad\textit{// $\tbr\fp{j}\hp{0} \in \cG_j\cap\cN(\bP)$}
		\State $\tbr\fp{j}\hp{1} := \bA \cdot \tbr\fp{j}\hp{0}$\quad\textit{// $\tbr\fp{j}\hp{1} \in \cG_j$}
		\State \textit{// Auxiliary vector}
		\State $\bxi := \bZ\d \cdot (\bP\t\cdot\tbr\fp{j}\hp{1})$
		\State $\tbv\fp{j}\hp{0} := \tbr\fp{j}\hp{1} - \bv\fp{j}\hp{0} \cdot \bxi$,\quad$\tbv\fp{j}\hp{-1} := \tbr\fp{j}\hp{0} - \bv\fp{j}\hp{-1} \cdot \bxi$\quad\textit{// $\tbv\fp{j}\hp{0} \in \cG_j\cap\cN(\bP)$}
		\State $\tbv\fp{j}\hp{1} := \bA \cdot \tbv\fp{j}\hp{0}$
		\State $\bZ := \bP\t\cdot\tbv\hp{1}\fp{j}$
		\State \textit{// - - - Polynomial step - - - }
		\State $\omega_{j+1} := [\tbr\fp{j}\hp{1}]\d \cdot \tbr\fp{j}\hp{0}$\quad\textit{// $\omega_{j+1} = \operatornamewithlimits{argmin}_{\omega \in \C}\lbrace\|\tbr\fp{j}\hp{0}-\omega\cdot\tbr\fp{j}\hp{1}\|\rbrace$}
		\State \textit{// Residual}
		\State $\bx\fp{j+1}\hp{0} := \tbx\fp{j}\hp{0} + \omega_{j+1} \cdot \tbr\fp{j}\hp{0}$
		\State $\br\fp{j+1}\hp{0} := \tbr\fp{j}\hp{0} - \omega_{j+1} \cdot \tbr\fp{j}\hp{1}$\quad\textit{// $\equiv (\bI - \omega_{j+1} \cdot \bA) \cdot \tbr\hp{0}\fp{j} \in \cG_{j+1}$}
		\State \textit{// Auxiliary vector}
		\State $\bv\fp{j+1}\hp{-1} := \tbv\fp{j}\hp{-1} - \omega_{j+1} \cdot \tbv\fp{j}\hp{0}$
		\State $\bv\fp{j+1}\hp{0} := \tbv\fp{j}\hp{0} - \omega_{j+1} \cdot \tbv\fp{j}\hp{1}$\quad\textit{// $\equiv (\bI - \omega_{j+1} \cdot \bA) \cdot \tbv\hp{0}\fp{j} \in \cG_{j+1}$}
		\State $\bZ := -\omega_{j+1} \cdot \bZ$
		\If{$\|\br\hp{0}\fp{j+1}\|\leq\tolabs$}
		\State \Return $\bx\hp{0}\fp{j+1}$
		\EndIf
		\EndFor
		\EndProcedure
	\end{algorithmic}
	\caption{\BiCGstab with lower indices}\label{Algo:BiCGstab_naive_indexed}
\end{algorithm}

\largeparbreak
In the following we discuss the iterative scheme in the for-loop.
At the beginning of the for-loop $\br\hp{0}\fp{j},\bv\hp{0}\fp{j}$ live in $\cG_j$, as is written as a comment in line 7. One iteration of the for-loop performs computations such that afterwards two vectors $\br\hp{0}\fp{j+1},\bv\hp{0}\fp{j+1}$ exist that live in $\cG_{j+1}$.  To this end the above template of the two steps \textit{biorthognalization} and \textit{polynomial step} is applied for both the residual and the auxiliary vector. 

\begin{figure}
\centering
\includegraphics[width=1\linewidth]{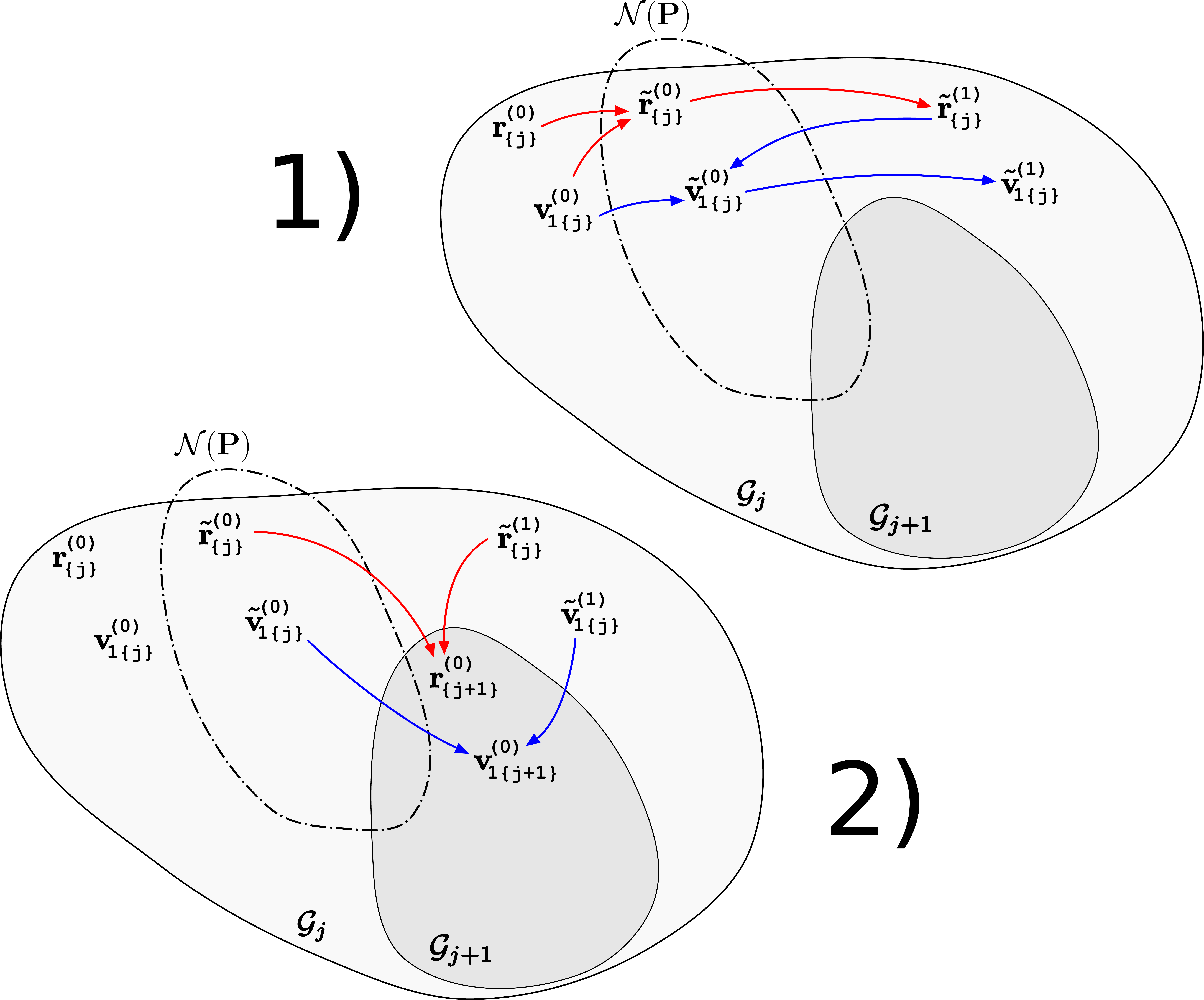}
\caption{Geometric principle of \IDRO in \BiCGstab.}
\label{fig:IDRcycle1}
\end{figure}

The two steps are illustrated in Fig.\,\ref{fig:IDRcycle1}. The figure has two parts. Both parts show spaces and vectors from Alg.\,\ref{Algo:BiCGstab_naive_indexed} that are contained in these spaces. The spaces are $\cG_j$ (light grey), $\cG_{j+1}$ (dark grey, contained in $\cG_j$) and $\cN(\bP)$ (dot-dashed bordered region, intersects in a chaotic way with both $\cG_j,\cG_{j+1}$). Part 1 of the figure illustrates the biorthogonalisation step whereas part 2) illustrates the polynomial step.
\largeparbreak
The iterative scheme starts with the biorthogonalisation. In lines 10 and 11 a new residual $\tbr\hp{0}\fp{j}$ is constructed from $\br\hp{0}\fp{j},\bv\hp{0}\fp{j} \in \cG_j$ such that $\tbr\hp{0}\fp{j} \perp \bp$ holds (where $\bp$ is the column vector of the $N \times 1$ matrix $\bP$). $\tbx\hp{0}\fp{j}$ is the according numerical solution to this residual. Then, in line 11, the image $\tbr\hp{1}\fp{j} = \bA \cdot \tbr\hp{0}\fp{j}$ is computed. Using Lem.\ref{Lem:RemainingVectorsInGj} we find that $\tbr\hp{1}\fp{j} \in \cG_j$ holds. The lines 10--12 are indicated by red arrows in part 1 of Fig.\,\ref{fig:IDRcycle1}.

Afterwards, in lines 14 to 16 the vectors $\tbr\hp{1}\fp{j},\bv\hp{0}\fp{j}$ are used to construct another vector $\tbv\hp{0}\fp{j}$ in $\cG_j \cap \cN(\bP)$ and its image $\tbv\hp{1}\fp{j} \in \cG_j$ (blue arrows in part 1 of the figure). The vector $\tbv\hp{-1}\fp{j}$ is constructed such that it is the pre-image of $\tbv\hp{0}\fp{j}$. With these computations the biorthogonalisation step is completed. 
\largeparbreak
In line 18 the polynomial step begins. It is shown in part 2 of Fig.\,\ref{fig:IDRcycle1}. The polynomial step starts in line 19 with the computation of an appropriate value for the relaxation coefficient $\omega_{j+1}$. This value is chosen such that the norm of $\br\hp{0}\fp{j+1}$ in line 22 will be minimised.

In line 22 a new residual $\br\hp{0}\fp{j+1}$ is computed. Its construction is illustrated by red arrows in Fig.\,\ref{fig:IDRcycle1} part 2. We explain why it holds $\tbr\hp{0}\fp{j+1} \in \cG_{j+1}$:

The figure shows that $\br\hp{0}\fp{j+1}$ is computed from a residual $\tbr\hp{0}\fp{j}$ that lives in $\cG_j \cap \cN(\bP)$ and its image $\tbr\hp{1}\fp{j}=\bA\cdot\tbr\hp{0}\fp{j}$. Consequently, $\br\hp{0}\fp{j+1}$ lives in $(\bI - \omega_{j+1} \cdot \bA) \cdot (\cG_j \cap \cN(\bP))$, which is just $\cG_{j+1}$. The vector $\bx\hp{0}\fp{j+1}$ is constructed such that it is the according numerical solution to $\br\hp{0}\fp{j+1}$.

In the same way as the new residual in $\cG_{j+1}$ has been computed, a new auxiliary vector $\bv\hp{0}\fp{j+1}$ with its according pre-image $\bv\hp{-1}\fp{j+1}$ is constructed in lines 24 to 25.

\subsubsection{Properties of \BiCGstab}
In this subsection we discuss potential geometric breakdowns, finite termination properties and the computational costs of \BiCGstab.

\paragraph{Breakdowns}
First, we notice that \BiCGstab in the implementation provided in Alg.\,\ref{Algo:BiCGstab_naive_indexed} has no formal breakdowns (only since a pseudo-inverse exists for every matrix). However, considering the geometric approach of \IDRO, it can happen that the biorthogonalisation fails so that potentially $\tbr\hp{0}\fp{j}$ and/or $\tbv\hp{0}\fp{j}$ do only live in $\cG_j$ but not in $\cG_j \cap \cN(\bP)$. Besides, it may happen that $\tbr\hp{1}\fp{j} \perp \tbr\hp{0}\fp{0}$ holds, e.g. when $\bA$ is a skew-symmetric matrix \cite{MaintainingConvergenceBiCGstab,EnhancedBiCGstabL}. In this case $\omega_{j+1}=0$, which destroys the property (b) in Theorem\,\ref{theo:IDR-theo}.

On top of these breakdown scenarios, the numerical round-off destroys ab initio that the properties in line 5 hold. In practice we observe that \BiCGstab (and its generalisations) suffer from loss of (superlinear) convergence when a (near) breakdown scenario is encountered.

\paragraph{Finite termination}
Assumed that there was no numerical round-off and no breakdown, the residual $\br\hp{0}\fp{N}$ is zero. This is because $\cG_N = \lbrace \bO \rbrace$ follows from Theorem\,\ref{theo:IDR-theo}. Thus, Alg.\,\ref{Algo:BiCGstab_naive_indexed} terminates after at most $N$ repetitions of the for-loop. Since during each repetition $2$ matrix-vector products are computed the method requires in total at most $2 \cdot N$ matrix-vector products to terminate.

\paragraph{Computational cost}
As a preliminary remark, there are more efficient implementations available of \BiCGstab than Alg.\,\ref{Algo:BiCGstab_naive_indexed}, cf. \cite{BiCGstab-Paper}.
The presented method in Alg.\,\ref{Algo:BiCGstab_naive_indexed} requires storage for $7$ column-vectors in $\R^N$, namely $\bx\hp{0},\br\hp{0},\br\hp{1},\bv\hp{-1},\bv\hp{0},\bv\hp{1},\bp$. The foot-indices ${}\fp{j}$ and the tildes $\tilde{\phantom{\bv}}$ were only introduced to help the derivation but one could as well always overwrite the original seven vectors. Alg.\,\ref{Algo:BiCGstab_naive} provides an according implementation that comes along with storage for only these seven vectors.

During each for-loop $2$ matrix-vector products, $7$ DOTs and $8$ AXPYs must be computed.

\begin{algorithm}
	\begin{algorithmic}[1]
		\Procedure{BiCGstab}{$\bA,\bb,\tolabs$}
		\State $\br\hp{0} := \bb$,\quad $\bx\hp{0} := \bO$,\quad$\bv\hp{-1} := \bb$
		\State $\bP := $\texttt{orth}$(\texttt{randn}(N,s))$,\quad$\bv\hp{0} := \bA \cdot \bv\hp{-1}$
		\State $\bZ := \bP\t \cdot \bv\hp{0} \in \R^{1 \times 1}$,\quad$j:=0$
		\State \textit{// $\br\hp{0},\bv\hp{0}\in \cK_\infty(\bA;\bb) \equiv \cG_0$}
		\While{$\|\br\hp{0}\|>\tolabs$}
		\State \textit{// $\br\hp{0},\bv\hp{0}\in\cG_j$}
		\State \textit{// - - - Biorthogonalisation - - - }
		\State \textit{// Residual}
		\State $\bxi := \bZ\d \cdot (\bP\t\cdot\br\hp{0})$
		\State $\br\hp{0} := \br\hp{0} - \bv\hp{0} \cdot \bxi$\quad\textit{// $\br\hp{0} \in \cG_j\cap\cN(\bP)$}
		\State $\bx\hp{0} := \bx\hp{0} + \bv\hp{-1} \cdot \bxi$
		\State $\br\hp{1} := \bA \cdot \br\hp{0}$\quad\textit{// $\br\hp{1} \in \cG_j$}
		\State \textit{// Auxiliary vector}
		\State $\bxi := \bZ\d \cdot (\bP\t\cdot\br\hp{1})$
		\State $\bv\hp{0} := \br\hp{1} - \bv\hp{0} \cdot \bxi$\quad\textit{// $\bv\hp{0} \in \cG_j\cap\cN(\bP)$}
		\State $\bv\hp{-1} := \br\hp{0} - \bv\hp{-1} \cdot \bxi$
		\State $\bv\hp{1} := \bA \cdot \bv\hp{0}$
		\State $\bZ := \bP\t\cdot\bv\hp{1}$
		\State \textit{// $\br\hp{0},\bv\hp{0} \in \cG_j \cap \cN(\bP)$}
		\State \textit{// - - - Polynomial step - - - }
		\State $\tau := [\br\hp{1}]\d \cdot \br\hp{0}$
		\State \textit{// Residual}
		\State $\bx\hp{0} := \bx\hp{0} + \br\hp{0} \cdot \tau$
		\State $\br\hp{0} := \br\hp{0} - \br\hp{1} \cdot \tau$
		\State \textit{// $\br\hp{0} \in (\bI-\tau\cdot\bA)\cdot\big(\cG_j \cap \cN(\bP)\big) \equiv \cG_{j+1}$}
		\State \textit{// Auxiliary vector}
		\State $\bv\hp{-1} := \bv\hp{-1} - \bv\hp{0} \cdot \tau$
		\State $\bv\hp{0}  := \bv\hp{0}  - \bv\hp{1} \cdot \tau$
		\State \textit{// $\bv\hp{0} \in (\bI-\tau\cdot\bA)\cdot\big(\cG_j \cap \cN(\bP)\big) \equiv \cG_{j+1}$}
		\State $\bZ:=-\tau \cdot \bZ$,\quad$j:=j+1$
		\EndWhile
		\State \Return $\bx\hp{0}$
		\EndProcedure
	\end{algorithmic}
	\caption{\BiCGstab}\label{Algo:BiCGstab_naive}
\end{algorithm}

\section{Derivation of \IDRstab from \BiCGstab}\label{Sec:DerivationIDRstab}
\begin{figure}
\centering
\includegraphics[width=0.6\linewidth]{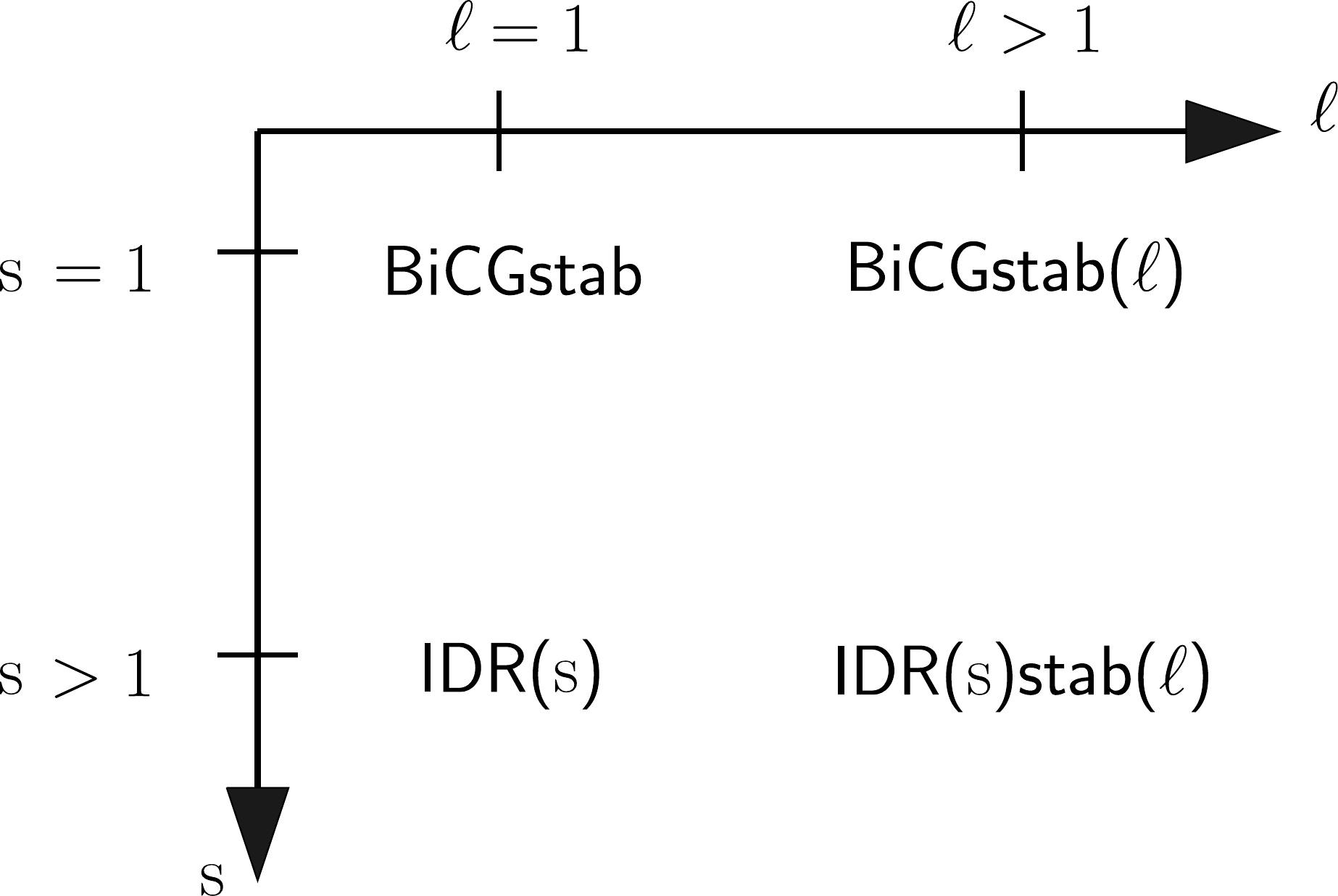}
\caption{Classification of \BiCGstabL, \IDR and \IDRstab as generalisations of \BiCGstab.}
\label{fig:IDRstab_categories}
\end{figure}

In the former section we have rederived \BiCGstab as an \IDRO method. \IDRstab is a twofold generalisation of \BiCGstab \cite[p.\,2688, l.\,33]{IDRstab-Paper}. \BiCGstab is \IDRstabp{$1$}{$1$}. \IDRstab is a twofold generalisation in that sense that it generalises on the one hand the biorthogonalisation step and on the other hand the polynomial step. As a result of the generalised biorthogonalisation step, \IDRstab for $s>1$ terminates earlier than \BiCGstab (which leads for large scale problems usually to an earlier occurrence of superlinear convergence). As a result of the generalised polynomial step, for $\ell>1$ the method converges more reliable for strongly a-symmetric linear systems \cite{BiCGstabL}.
\largeparbreak
Whereas \IDRstab is a twofold generalisation of \BiCGstab, there are also methods that are only generalisations of \BiCGstab in either direction. This is illustrated in Fig.\,\ref{fig:IDRstab_categories}: \BiCGstabL is a \BiCGstab method with only a generalised polynomial step, whereas \IDR is only a generalisation in the biorthogonalisation step. In order to derive \IDRstab, we first derive \BiCGstabL and then \IDR.

\subsection{The method \BiCGstabL}
In this subsection we introduce \BiCGstabL as a generalisation of \BiCGstab. We first motivate the generalisation and provide the theoretical concepts. Afterwards we present an algorithm of \BiCGstabL on which we explain again illustratively the geometric properties of its computed vectorial quantities.

\subsubsection{Motivation of \BiCGstabL}
The following lemma shows that for the recursive construction of Sonneveld spaces there are several equivalent formulations.

\begin{Lemma}[General recursion of Sonneveld spaces] \label{lem:GenRekSonneveldspace}
	The Sonneveld spaces as defined in Definition \ref{def:Sonneveldspace} satisfy for all $j,\ell \in \No$:
	\begin{align*}
	\cG_{j+\ell} = \underbrace{\prod_{i=j+1}^{j+\ell}(\bI-\omega_i \cdot \bA)}_{=:q_{j,\ell}(\bA)} \cdot \big(\cG_j \cap \cK^\perp_{\ell}(\bA\h;\bP)\big)\,.
	\end{align*}
	\underline{Proof:} Choose an arbitrary $\by \in \C^N$. Define $\bx:= q_{0,j}(\bA)\cdot\by$. In the first part of the proof we show that the statements (A) and (B) for $\by,\bx$ are equivalent.
	\begin{align*}
	\text{(A)} & & \by &\in \cG_0 \cap \cK_{j+\ell}^\perp(\bA\h;\bP)\\[10pt]
	\text{(B)} & & \bx &\in \cG_j \cap \cK_{\ell}^\perp(\bA\h;\bP)
	\end{align*}
	(A) and (B) can be reformulated as follows:
	\begin{align*}
	\text{(A)}& & & & \by &\in \cG_0 & &\\
	\land& & & & \langle \bA^k \cdot \by,\bp_q\rangle&=0\quad& &\forall\,q=1,...,s,\ \forall\,k=0,...,j+\ell-1\\[10pt]
	\text{(B)} & & & & \bx &\in \cG_j & &\\
	\land & & & & \langle\bA^k \cdot q_{0,j}(\bA) \cdot \by,\bp_q\rangle&=0\quad& &\forall\,q=1,...,s,\ \forall\,k=0,...,\ell-1
	\end{align*}
	$\bx \in \cG_j$ can be expressed equivalently by $\by \in \cG_0 \cap \cK_j^\perp(\bA\h;\bP)$, cf. \cite{IDR-Gutkn}. Replacing these conditions in (B) we obtain:
	\begin{align*}
	\text{(A)}& & & & \by &\in \cG_0 & &\\
	\land& & & & \langle \bA^k \cdot \by,\bp_q\rangle&=0\quad& &\forall\,q=1,...,s,\ \forall\,k=0,...,j+\ell-1\\[10pt]
	\text{(B)} & & & & \by &\in \cG_0 & &\\
	\land & & & &\langle \bA^k \cdot \by,\bp_q\rangle&=0\quad& &\forall\,q=1,...,s,\ \forall\,k=0,...,j-1\\
	\land & & & & \langle\bA^k \cdot q_{0,j}(\bA) \cdot \by,\bp_q\rangle&=0\quad& &\forall\,q=1,...,s,\ \forall\,k=0,...,\ell-1
	\end{align*}
	Now obviously (A) and (B) are equivalent. Using this, all in all for an arbitrary $\hbx \in \cG_{j + \ell}$ there exist $\bx,\by$ such that
	\begin{align*}
	\hbx = q_{j,\ell}(\bA) \cdot \underbrace{q_{0,j}(\bA)\cdot \by}_{\equiv\bx} \,.
	\end{align*}
	In the formula for $\hbx$ the restrictions for $\bx$ respectively $\by$ can be expressed equivalently either by (A) or (B). $\boxtimes$
\end{Lemma}
\largeparbreak
So far, in \BiCGstab we have iteratively restricted vectors from $\cG_{j}$ into $\cG_{j+1}$. The algorithmic idea of \BiCGstabL is instead to iteratively restrict vectors from $\cG_j$ into $\cG_{j + \ell}$. This has a particular advantage:

For $\bA \in \R^{N \times N}$, in \BiCGstab we have subsequently chosen $\omega_{j+1},...,\omega_{j+\ell} \in \Rno$ such that each time the subsequent residual is minimised. The overall polynomial $q_{j,\ell}(t) = \prod_{k=j+1}^{\ell}(t^0 - \omega_k\cdot t^1)$ has only real roots. For strongly a-symmetric systems with large imaginary parts in the eigenvalues it is however more suitable to use stabilisation polynomials that can also have complex roots.

In \BiCGstabL the values for $\omega_{j+1},...,\omega_{j+\ell}$ can be chosen all at once by constructing the polynomial $q_{j,\ell}(t) = t^0 - \sum_{k=1}^\ell \tau_k \cdot t^k$ by a $\ell$-degree residual minimisation procedure (such as it is done when performing $\ell$ iterations of \GMRES).

The residual minimisation polynomial of a degree $>1$ makes it possible that the roots are complex although all polynomial coefficients $\tau_k$ are real. Thus, the method is able to approximate complex eigenvalues of the system matrix, which is advantageous for the rate of convergence \cite{BiCGstabL}. The possibility of complex roots makes it further more unlikely that the roots $\omega_{j+1},...,\omega_{j+\ell}$ are (close to) zero, since even if the real part is zero the imaginary part might be still distinct from zero.

\subsubsection{Derivation of \BiCGstabL}
\begin{figure}
	\centering
	\includegraphics[width=1\linewidth]{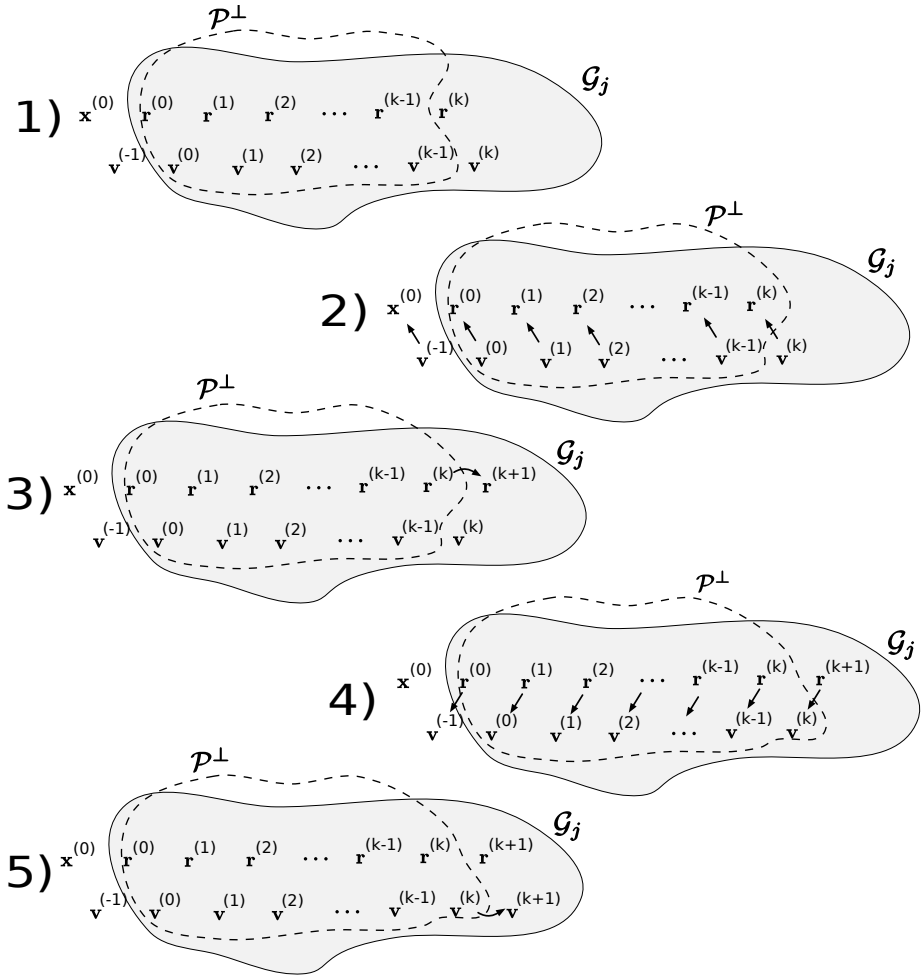}
	\caption{Geometric principle of the biorthogonalisation in \BiCGstabL.}
	\label{fig:Look-ahead2}
\end{figure}

The geometric approach of \BiCGstabL during each repetition of the main loop consists of the following three steps:
\begin{enumerate}[(1)]
	\item Biorthognalisation: Modify $\br\hp{0},\bv\hp{0} \in \cG_j$ such that they live in $\cG_j \cap \cK^\perp_\ell(\bA\h;\bP)$.
	\item Choose $q_{j,\ell}(\cdot)$ such that $\|q_{j,\ell}(\bA)\cdot\br\hp{0}\|$ is minimal.
	\item Polynomial step: Update $\br\hp{0}:=q_{j,\ell}(\bA)\cdot\br\hp{0}$ and $\bv\hp{0}:=q_{j,\ell}(\bA)\cdot\bv\hp{0}$.
\end{enumerate}
The first of the three steps is quite involved and includes a for-loop for $k=0,...,\ell-1$. Fig.\,\ref{fig:Look-ahead2} shows the $k$th repetition of this for-loop of the biorthogonalisation.

\largeparbreak
Our strategy to explain \BiCGstabL is as follows: We first state the algorithm and explain every step except the biorthogonalisation. Finally, we use Fig.\,\ref{fig:Look-ahead2} and the lines of the algorithm to explain the biorthogonalisation procedure.
\largeparbreak

Alg.\,\ref{Algo:BiCGstabL_naive} shows an implementation of \BiCGstabL. It consists of an initialisation phase from lines 2 to 5, which is identical to that of \BiCGstab, and a while-loop from line 6 to 36 that includes the iterative scheme. At the beginning of the while-loop in line 7 the vectors $\br\hp{0},\bv\hp{0}$ live in $\cG_j$. At the end of the loop in line 34 it holds $\br\hp{0},\bv\hp{0} \in \cG_{j+\ell}$.
\begin{algorithm}
	\begin{algorithmic}[1]
		\Procedure{BiCGstabL}{$\bA,\bb,\ell,\tolabs$}
		\State $\br\hp{0} := \bb$,\quad $\bx\hp{0} := \bO$,\quad$\bv\hp{-1} := \bb$
		\State $\bP := $\texttt{orth}$(\texttt{randn}(N,1))$,\quad$\bv\hp{0} := \bA \cdot \bv\hp{-1}$
		\State $\bZ := \bP\t \cdot \bv\hp{0} \in \R^{1 \times 1}$,\quad$j:=0$
		\State \textit{// $\br\hp{0},\bv\hp{0}\in \cK_\infty(\bA;\bb) \equiv \cG_0$}
		\While{$\|\br\hp{0}\|>\tolabs$}
		\State \textit{// $\br\hp{0},\bv\hp{0}\in\cG_j$}
		\State \textit{// - - - Biorthogonalisation - - - }
		\For{$k=0,...,\ell-1$}
		\State \textit{// Residual}
		\State $\bxi := \bZ\d \cdot (\bP\t\cdot\br\hp{k})$
		\State $\br\hp{g} := \br\hp{g} - \bv\hp{g} \cdot \bxi$\quad$\forall\ g=0,...,k$
		\State $\bx\hp{0} := \bx\hp{0} + \bv\hp{-1} \cdot \bxi$
		\State \textit{// $\br\hp{g} \in \cG_j\cap\cN(\bP)$ $\forall\ g=0,...,k$}
		\State $\br\hp{k+1} := \bA \cdot \br\hp{k}$\quad\textit{// $\br\hp{k} \in \cG_j$, cf. Lemma 1}
		\State \textit{// Auxiliary vector}
		\State $\bxi := \bZ\d \cdot (\bP\t\cdot\br\hp{k+1})$
		\State $\bv\hp{g} := \br\hp{g+1} - \bv\hp{g} \cdot \bxi$\quad$\forall\ g=-1,0,...,k$
		\State \textit{// $\bv\hp{g} \in \cG_j\cap\cN(\bP)$ $\forall\ g=0,...,k$}
		\State $\bv\hp{k+1} := \bA \cdot \bv\hp{k}$\quad\textit{// $\bv\hp{k+1} \in \cG_j$, cf. Lemma 1}
		\State $\bZ:=\bP\t\cdot\bv\hp{k+1}$
		\EndFor
		\State \textit{// $\br\hp{0},\bv\hp{0} \in \cG_j \cap \cK_\ell^\perp(\bA\t;\bP)$}
		\State \textit{// - - - Polynomial step - - - }
		\State $\btau := [\br\hp{1},...,\br\hp{\ell}]\d \cdot \br\hp{0}$\quad\textit{// $\btau \equiv (\tau_1,...,\tau_\ell)\t \in \R^s$}
		\State \textit{// $q_{j,\ell}(t) = t^0 - \sum_{k=1}^\ell \tau_k \cdot t^k \equiv \prod_{k=1}^{\ell}(t^0 - \omega_{j+k}\cdot t^1)$}
		\State \textit{// Residual}
		\State $\bx\hp{0} := \bx\hp{0} + \sum_{k=1}^\ell \br\hp{k-1} \cdot \tau_k$
		\State $\br\hp{0} := \br\hp{0} - \sum_{k=1}^\ell \br\hp{k} \cdot \tau_k$
		\State \textit{// $\br\hp{0} \in q_{j,\ell}(\bA)\cdot\big(\cG_j \cap \cK_\ell^\perp(\bA\t;\bP)\big) \equiv \cG_{j+\ell}$, cf. Lemma 2}
		\State \textit{// Auxiliary vector}
		\State $\bv\hp{-1} := \bv\hp{-1} - \sum_{k=1}^\ell \bv\hp{k-1} \cdot \tau_k$
		\State $\bv\hp{0} := \bv\hp{0} - \sum_{k=1}^\ell \bv\hp{k} \cdot \tau_k$
		\State \textit{// $\bv\hp{0} \in q_{j,\ell}(\bA)\cdot\big(\cG_j \cap \cK_\ell^\perp(\bA\t;\bP)\big) \equiv \cG_{j+\ell}$, cf. Lemma 2}
		\State $\bZ:=-\tau_\ell\cdot\bZ$,\quad$j:=j+\ell$
		\EndWhile
		\State \Return $\bx\hp{0}$
		\EndProcedure
	\end{algorithmic}
	\caption{\BiCGstabL}\label{Algo:BiCGstabL_naive}
\end{algorithm}

The lines 8 to 22 realise the biorthogonalisation. After this part, in line 23
the vectors $\br\hp{0},\bv\hp{0}$ live in $\cG_j \cap \cK_\ell^\perp(\bA\t;\bP)$. In addition to that, their \textit{power vectors}
\begin{align*}
	\bv\hp{g} &:= \bA^g \cdot \bv\hp{0} & \forall\,g&=1,...,\ell\\
	\br\hp{g} &:= \bA^g \cdot \br\hp{0} & \forall\,g&=1,...,\ell
\end{align*}
are computed, cf. lines 15 and 20 in the algorithm.

In line 25 the coefficients $\tau_1,...,\tau_\ell \in \R$ for the residual-minimising polynomial $q_{j,\ell}(\cdot)$ are computed. There is no need to compute the actual roots $\omega_{j+1},...,\omega_{j+\ell} \in \Cno$. Thus, even though the roots are likely to be complex, the algorithm still only performs computations with real-valued data.

In lines 27 to 34 the residual and the auxiliary vector are multiplied from the left by $q_{j,\ell}(\bA)$. Since the power vectors of $\br\hp{0},\bv\hp{0}$ have already been computed earlier, the polynomial step consists only of computing $\ell$ AXPYs for each the numerical solution $\bx\hp{0}$, its residual $\br\hp{0}$, the auxiliary vector $\bv\hp{0}$ and its pre-image $\bv\hp{-1}$.
\largeparbreak
In the following we explain the biorthogonalisation in lines 9 to 22:

At line 8 it holds $\br\hp{0},\bv\hp{0} \in \cG_j$. At the beginning of the for-loop of line 9, i.e. at line 10, it holds:
\begin{align*}
	\br\hp{g},\bv\hp{g} &\in \cG_j & &\forall g \in \lbrace n \in \N\ :\ n\geq 0 \ \land \ n\leq k \rbrace\\
	\br\hp{g},\bv\hp{g} &\in \cN(\bP) & &\forall g \in \lbrace n \in \N\ :\ n\geq 0 \ \land \ n\leq k-1 \rbrace
\end{align*}
and the power vectors $\bv\hp{1},...,\bv\hp{k}$ of the auxiliary vector and $\br\hp{1},...\br\hp{k}$ of the residual have been computed already. This scenario is presented in part 1 of Fig.\,\ref{fig:Look-ahead2}: The figure shows all the currently computed vectors and the spaces in which they live. $\cP^\perp$ shall be $\cN(\bP)$. In the following we explain what happens during the $k$th repetition of the for-loop:

In the algorithm, in lines 12 to 13 the residual, its power vectors and its numerical solution are modified by the auxiliary vector, its power vectors and its pre-image, such that $\br\hp{k} \in \cN(\bP)$ holds. In Fig.\,\ref{fig:Look-ahead2} this is illustrated in part 2. Afterwards, in line 15 the next power of the residual is computed. Since $\br\hp{k}\in \cG_j \cap \cN(\bP)$ holds it follows from Lemma\,\ref{Lem:RemainingVectorsInGj} that $\br\hp{k+1}$ lives in $\cG_j$, too. The step in line 15 is shown in the figure in part 3.

Next, the auxiliary vector is biorthogonalised. To this end, in line 18 the auxiliary vector, its power vectors and its pre-image are modified by a linear combination with the residual of the respectively next higher power such that $\bv\hp{k}$ is orthogonal with respect to $\bp$. I.e. afterwards $\bv\hp{k}\in \cG_j \cap \cN(\bP)$ holds. Part 4 of the figure illustrates this step. Finally, the subsequent power vector $\bv\hp{k+1}$ is computed in line 20 of the algorithm respectively part 5 of the figure.

\subsubsection{Properties of \BiCGstabL}
In this subsection we discuss potential geometric breakdowns, finite termination properties and computational costs of \BiCGstabL.

\paragraph{Breakdowns}
First, we notice that \BiCGstabL in the implementation provided in Alg.\,\ref{Algo:BiCGstabL_naive} has no formal breakdowns (only since pseudo-inverses exist for every matrix). However, considering the geometric approach of \IDRO, it can happen that the biorthogonalisation fails so that potentially for some $k\in\lbrace 0,...,\ell-1\rbrace$ the residual $\br\hp{k+1}$ and/or the auxiliary vector $\bv\hp{k+1}$ do only live in $\cG_j$ but not in $\cG_j \cap \cN(\bP)$. Besides, it may happen that the coefficient $\tau_\ell$ of the polynomial step is zero, e.g. when $\bA$ is a particular permutation matrix. In this case $\omega_{j+\ell}=0$ follows, which will destroy the property (b) in Theorem\,\ref{theo:IDR-theo}.

On top of these breakdown scenarios, the numerical round-off destroys ab initio that the properties in line 7 hold. In practice we observe that \BiCGstabL suffers from loss of (superlinear) convergence when a (near) breakdown scenario is encountered. This effect occurs even stronger when $\ell$ is large since then the vectors become collinear (due to convergence properties of the power iteration) and numerical round-off becomes dominant in the overall procedure.

\paragraph{Finite termination}
The termination properties of \BiCGstabL are completely identical to those of \BiCGstab. This is because still the same kind of spaces has been used, only the recursion formula itself has been exchanged by an equivalent one, cf. Lem.\,\ref{lem:GenRekSonneveldspace}. 

Under the assumption that there was no numerical round-off and no breakdown, it holds that the residual $\br\hp{0}\fp{N}$ is zero. This is because $\cG_N = \lbrace \bO \rbrace$ follows from Theorem\,\ref{theo:IDR-theo}. Thus, Alg.\,\ref{Algo:BiCGstab_naive_indexed} terminates after $\lceil N/\ell\rceil$ repetitions of the while-loop. Since during each repetition of the while-loop $2 \cdot \ell$ matrix-vector products are computed the method requires overall $\approx 2 \cdot N$ matrix-vector products to terminate.

\paragraph{Computational cost}
The presented method in Alg.\,\ref{Algo:BiCGstabL_naive} requires storage for $2 \cdot (\ell+1) + 1$ column-vectors in $\R^N$, namely $\bx\hp{0},\br\hp{0},...,\br\hp{\ell},\bv\hp{-1},\bv\hp{0},...,\bv\hp{\ell},\bp$. During each repetition of the while-loop $2\cdot \ell$ matrix-vector products, $3 \cdot \ell + \ell\cdot(\ell+1)/2$ DOTs and $2 \cdot \ell\cdot(\ell+1) + 4 \cdot \ell$ AXPYs must be computed.

Since further to stability issues the cost grows quadratically in $\ell$, the user should make a defensive choice of $\ell \in \N$. Since $\ell=1$ leads to zero complex roots of the residual-minimising polynomial, $\ell=3$ leads still only to two complex roots, and $\ell=4$ is stable by no chance, the value $\ell=2$ seems to be the only reasonable choice.\footnote{This is a personal opinion of the author.}

\subsection{The method \IDR}
From Fig.\,\ref{fig:IDRstab_categories} we have seen so far the methods \BiCGstab and \BiCGstabL. In the introduction of Section\,\ref{Sec:DerivationIDRstab} we have motivated \BiCGstabL as a more reliably converging method when the system matrix has eigenvalues with large imaginary parts. At this stage, it is clear that this improvement in the robustness stems from the capability of \BiCGstabL's polynomial step for $\ell>1$ to approximate these imaginary parts.

In this subsection we introduce \IDR. This method has superior termination properties over \BiCGstab. Whereas \BiCGstab and \BiCGstabL require $\approx 2 \cdot N$ matrix-vector products to terminate, \IDR does only require $\approx (1+1/s) \cdot N$ matrix-vector products. In this, $s\in\N$ is a method parameter. \IDReins is algorithmically equivalent to \BiCGstab.

\subsubsection{Idea of \IDR}
So far, in \BiCGstab and \BiCGstabL we have used the recursion of Sonneveld spaces for $\bP \in \R^{N \times s}$ with $s=1$. However, one can also use $s>1$. To this end, a wider matrix for $\bP$ and $s$ auxiliary vectors instead of only one must be utilised.

The algorithmic recipe for restricting the residual $\br\hp{0}\in\cG_j$ and $s$ auxiliary vectors $\bv\hp{0}_1,...,\bv\hp{0}_s\in\cG_j$ from $\cG_j$ into $\cG_{j+1}$ consists of the following steps:
\begin{enumerate}[(1)]
	\item Biorthogonalisation: Project the residual with the $s$ auxiliary vectors into the null-space of $\bP$. Afterwards, construct $s$ auxiliary vectors that live in $\cG_j \cap \cN(\bP)$, too.
	\item Polynomial step: Choose $\omega_{j+1}$ and perform for all vectors the same polynomial update as in \BiCGstab.
\end{enumerate}

In the following sub-subsection we describe in all detail the computational steps of the algorithm \IDR.

\subsubsection{Derivation of \IDR}
\begin{figure}
\centering
\includegraphics[width=1\linewidth]{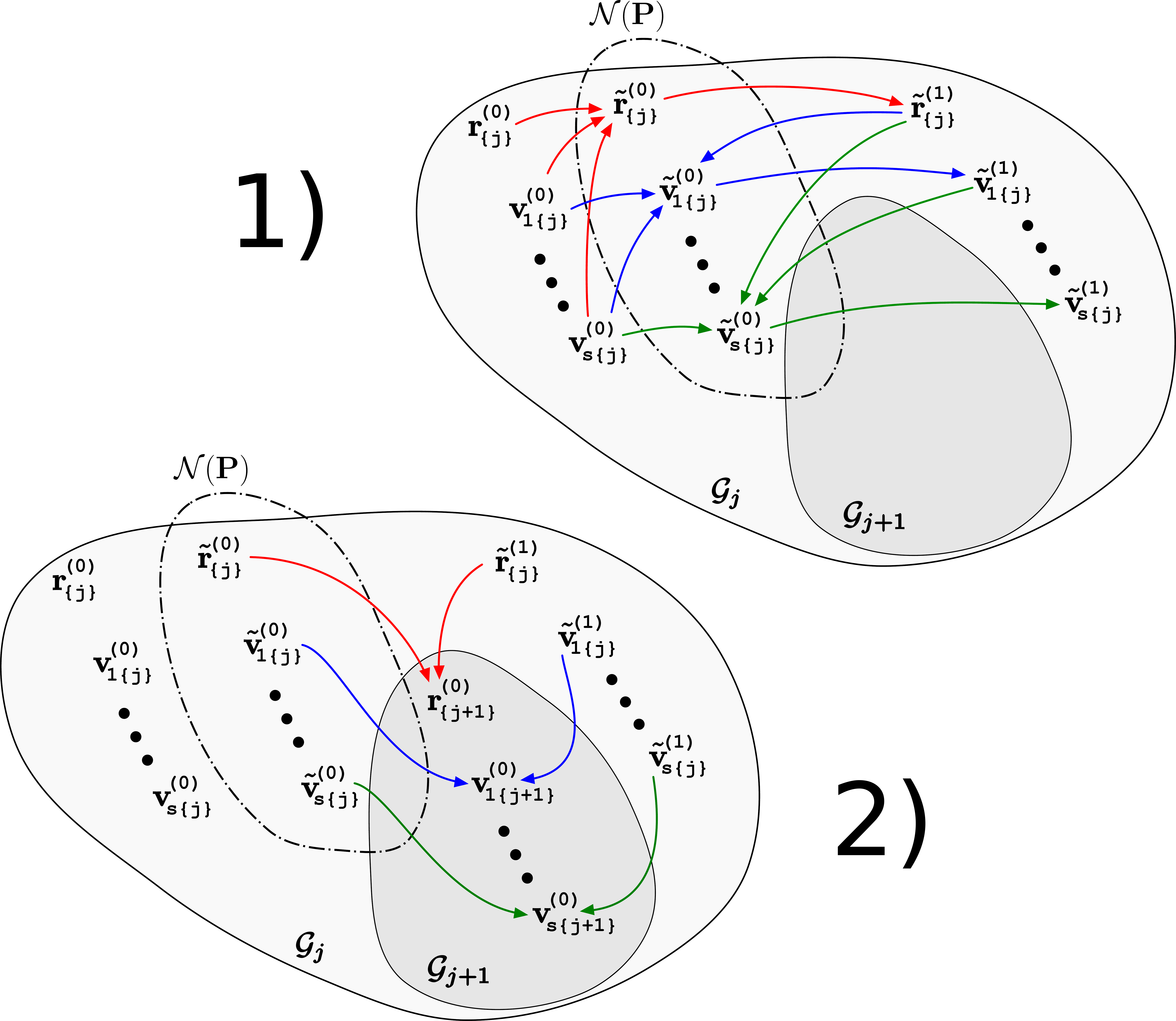}
\caption{Geometric principle of \IDRO in \IDR.}
\label{fig:IDRcycle2}
\end{figure}

Alg.\,\ref{Algo:IDRs_naive} provides an implementation of the method \IDR and Fig.\,\ref{fig:IDRcycle2} shows the computed vectors within \IDR in their respective spaces. In the following we go through the lines of the algorithm and explain with the help of the figure what is done in geometric terms during each step.
\largeparbreak
The first nine lines constitute the initialisation of the method: The initial solution and its residual are computed. Further, a matrix $\bP \in \R^{N \times s}$ for the Sonneveld spaces is initialised. As a difference to what is done in \BiCGstab and \BiCGstabL, this matrix does not only consist of one but $s$ columns, where $s \in \N$ is an arbitrary user-parameter. The columns of $\bP$ are orthonormalised since this leads to a superior conditioning of the matrix $\bZ$ that is used later for the biorthogonalisations.

As a further difference to \BiCGstab and \BiCGstabL, not only one auxiliary vector $\bv\hp{0}$ in $\cK_\infty(\bA;\bb)$ but a list of $s$ auxiliary vectors $\bv_q\hp{0}$, $q=1,...,s$, is built. For each of these vectors, the pre-image $\bv_q\hp{-1} = \bA^{-1}\cdot\bv\hp{0}_q$ is kept.
\largeparbreak
The main loop from line 10 to line 34 consists of the biorthogonalisation step and the polynomial step. During the while-loop, the vectors $\br\hp{0},\bv\hp{0}_1,...,\bv\hp{0}_s$ are moved from $\cG_j$ to $\cG_{j+1}$. Part 1 of Fig.\,\ref{fig:IDRcycle2} shows what is done during the biorthogonalisation: First, using the $s$ auxiliary vectors, the residual is orthogonalised with respect to the columns of $\bP$, cf. lines 14 and 15 of the algorithm. Right afterwards, the image of the biorthogonalised residual is computed in line 16. The biorthogonalisation of the residual and the computation of its image are illustrated by red arrows in part 1 of the figure.

After the residual has been biorthogonalised, the for-loop from line 18 to 24 orthogonalises one after the other each of the auxiliary vectors. For $\bv_1\hp{0}$, the residual's image $\br\hp{1}$ and the auxiliary vectors $\bv_1\hp{0},...,\bv\hp{0}_s$ are used to overwrite the auxiliary vector $\bv\hp{0}_1$ such that it lives in $\cG_j \cap \cN(\bP)$, cf. lines 19 and 20 for $q=1$ in the algorithm. In line 21 and 22 the pre-image and image of this new biorthogonal auxiliary vector are computed, respectively. Part 1 of the figure illustrates the for-loop for $q=1$ with blue arrows.

In order to move the second auxiliary vector $\bv\hp{0}_2$ into $\cG_j \cap \cN(\bP)$ the vectors $\br\hp{1},\bv_1\hp{1},\bv_2\hp{0},...,\bv\hp{0}_s$ are used. And in general, to biorthogonalise the $q$th auxiliary vector, the vectors $\br\hp{1},\bv_1\hp{1},...\bv_{q-1}\hp{1},\bv\hp{0}_q,...,\bv\hp{0}_s \in \cG_j$ are used. This is illustrated in the figure by the green arrows.
\largeparbreak

The biorthogonalisation is followed by the polynomial step from line 25: In line 26 a value $\tau \equiv \omega_{j+1}$ is computed that minimises the norm of the result for $\br\hp{0}$ in the expression in line 29. The lines 27 to 32 are analogous to the polynomial step of \BiCGstab with the only difference that in lines 31 and 32 not only one auxiliary vector but matrices of multiple auxiliary vectors are updated. The polynomial step is illustrated by arrows in part 2 of the figure: Using the respective vector in $\cG_j \cap \cN(\bP)$ and its image, the polynomial step can be computed as a linear combination of both.
\largeparbreak

The matrix $\bZ$ in line 33 is updated in such a way that
\begin{align*}
	\bZ = \bP\t\cdot\bV\hp{0}\tageq\label{eqn:Z=PtV}
\end{align*}
holds for $\bV\hp{0}$ from line 32. This follows from the fact that in the for-loop in line 23 the new columns of $\bZ$ are computed as columns of the matrix $\bP\t\cdot\bV\hp{1}$, whereas at line 30 it holds $\bP\t\cdot\bV\hp{0}=\bO$.

\begin{algorithm}
	\begin{algorithmic}[1]
		\Procedure{IDR}{$\bA,\bb,s,\tolabs$}
		\State $\br\hp{0} := \bb$,\quad $\bx\hp{0} := \bO$,\quad$\bv_1\hp{-1} := \bb$
		\State $\bP := $\texttt{orth}$(\texttt{randn}(N,s))$,\quad$\bv_1\hp{0} := \bA \cdot \bv\hp{-1}_1$
		\For{$q=2,...,s$}
		\State $\bv_q\hp{-1} := \bv\hp{0}_{q-1}$
		\State $\bv_q\hp{0} := \bA \cdot \bv\hp{-1}_q$
		\EndFor
		\State $\bZ := \bP\t \cdot \bV\hp{0} \in \R^{s \times s}$,\quad$j:=0$
		\State \textit{// $\br\hp{0},\bv_1\hp{0},...,\bv_s\hp{0}\in \cK_\infty(\bA;\bb) \equiv \cG_0$}
		\While{$\|\br\hp{0}\|>\tolabs$}
		\State \textit{// $\br\hp{0},\bv_1\hp{0},...,\bv_s\hp{0}\in\cG_j$}
		\State \textit{// - - - Biorthogonalisation - - - }
		\State \textit{// Residual}
		\State $\bxi := \bZ\d \cdot (\bP\t\cdot\br\hp{0})$
		\State $\br\hp{0} := \br\hp{0} - \bV\hp{0} \cdot \bxi$,\quad$\bx\hp{0} := \bx\hp{0} + \bV\hp{-1} \cdot \bxi$\quad\textit{// $\br\hp{0} \in \cG_j\cap\cN(\bP)$}
		\State $\br\hp{1} := \bA \cdot \br\hp{0}$\quad\textit{// $\br\hp{1} \in \cG_j$, cf. Lemma \ref{Lem:RemainingVectorsInGj}}
		\State \textit{// Auxiliary vectors}
		\For{$q=1,...,s$}
		\State $\bxi := \bZ\d \cdot (\bP\t\cdot\br\hp{1})$
		\State $\bv_q^{(0)\phantom{-}} := \br\hp{1} - [\bV\hp{1}_{:,1:(q-1)},\bV^{(0)\phantom{-}}_{:,q:s}] \cdot \bxi$\quad\textit{// $\bv_q\hp{0} \in \cG_j\cap\cN(\bP)$}
		\State $\bv_q\hp{-1} := \br\hp{0} - [\bV\hp{0}_{:,1:(q-1)},\bV\hp{-1}_{:,q:s}] \cdot \bxi$
		\State $\bv_q\hp{1} := \bA \cdot \bv_q\hp{0}$\quad\textit{// $\bv_q\hp{1} \in \cG_j$, cf. Lemma \ref{Lem:RemainingVectorsInGj}}
		\State $\bZ_{:,q} := \bP\t\cdot\bv\hp{1}_q$
		\EndFor
		\State \textit{// - - - Polynomial step - - - }
		\State $\tau := [\br\hp{1}]\d \cdot \br\hp{0}$
		\State \textit{// Residual}
		\State $\bx\hp{0} := \bx\hp{0} + \tau \cdot \br\hp{0}$
		\State $\br\hp{0} := \br\hp{0} - \tau \cdot \br\hp{1}$ \quad\textit{// $\br\hp{0} \in \cG_{j+1}$}
		\State \textit{// Auxiliary vectors}
		\State $\bV\hp{-1} := \bV\hp{-1} - \tau \cdot \bV\hp{0}$
		\State $\bV\hp{0} := \bV\hp{0} - \tau \cdot \bV\hp{1}$ \quad\textit{// $\rg(\bV\hp{0}) \subset \cG_{j+1}$}
		\State $\bZ:=-\tau \cdot \bZ$,\quad$j:=j+1$
		\EndWhile
		\State \Return $\bx\hp{0}$
		\EndProcedure
	\end{algorithmic}
	\caption{\IDR}\label{Algo:IDRs_naive}
\end{algorithm}

\subsubsection{Properties of \IDR}
In this subsection we discuss properties of the above derived method \IDR.

The breakdowns of \IDR are analogous to those of \BiCGstab: The biorthogonalisation can fail when either the residual or one of the auxiliary vectors cannot be projected onto $\cN(\bP)$. Besides, as formerly discussed, the value of $\tau$ can become zero when using the formula in line 26. $\tau=0$ violates the requirements of Theorem\,\ref{theo:IDR-theo}. As a consequence of this violation the Sonneveld spaces would not shrink in dimension any more.

\paragraph{Finite termination} This property is of particular interest in the discussion of \IDR. Reviewing Alg.\,\ref{Algo:IDRs_naive}, during the while-loop the vectors $\br\hp{0},\bv\hp{0}_1,...,\bv\hp{0}_s$ are moved from $\cG_{j}$ to $\cG_{j+1}$.

From Theorem\,\ref{theo:IDR-theo} it follows that $\cG_{j+1}$ is either the null-space or by $s$ dimensions smaller than $\cG_{j}$. Consequently, after at most $\lceil N/s \rceil$ repetitions of the while-loop the method terminates. During each repetition of the while-loop $s+1$ matrix-vector products are computed. Thus, in total the method terminates after at most
\begin{align*}
	\lceil N/s \rceil \cdot (s+1) \approx (1+1/s) \cdot N
\end{align*}
matrix-vector products. Comparing this to \GMRES, the latter must compute $N$ matrix-vector products to terminate with the exact solution in the general case. For mild values of $s>1$, \IDR soon approaches this optimal bound of \GMRES. This is superior over \BiCGstab and \BiCGstabL, which both require $\approx 2 \cdot N$ matrix-vector products.

\paragraph{Computational cost}
The presented method in Alg.\,\ref{Algo:IDRs_naive} requires storage for $3\cdot (s+1) + s$ column vectors of length $N$, namely for the $s$ auxiliary vectors, their images and  pre-images, further for the residual, its image and its according numerical solution, and further for the $s$ column vectors of $\bP$. During each repetition of the while-loop, $s+1$ matrix-vector products, $5+s \cdot (s+1)$ DOTs and $2\cdot s^2 + 4\cdot s + 2$ AXPYs must be computed.

Since further to potential stability issues related to the condition of $\bV\hp{0}$ (i.e. degree of orthogonality of its column vectors) the cost grows quadratically in $s$, the user should make a moderate choice for $s \in \N$. We underline that there are more practical implementations than Alg.\,\ref{Algo:IDRs_naive} that have numerical treatments to improve the conditioning of $\bV\hp{0}$.

Commonly used values for the method parameter $s$ range from $2$ to $8$. From our experience we suggest to use $s=4$ auxiliary vectors for well-conditioned problems and up to $s=8$ for badly conditioned systems or when trying to strike a benefit from Krylov subspace recycling. In the case of Krylov subspace recycling the parameter $s$ has the effect that the method terminates $s-1$ times earlier than \IDR without recycling, cf. \cite[p.\,13 bottom]{Mstab-paper}.

\subsection{The method \IDRstab}
After having presented both kinds of generalisations for \BiCGstab as were discussed along Fig.\,\ref{fig:IDRstab_categories}, in this subsection we now introduce \IDRstab.

\subsubsection{Motivation of \IDRstab}
On the one hand there is \BiCGstabL. It has a treatment for highly a-symmetric systems because it is able to approximate complex eigenvalues of the system matrix. This ability is important to achieve huge reductions in the residual norm during the polynomial step. Unfortunately, \BiCGstabL requires $\approx 2 \cdot N$ matrix-vector products to terminate in the worst case. As a consequence of this one can expect that often it experiences the transition to superlinear convergence only after twice the number of matrix-vector products that \GMRES requires.
\largeparbreak
On the other hand there is \IDR that requires only $\approx (1+1/s) \cdot N$ matrix-vector products to terminate in the worst case. Consequently, it does the transition to superlinear convergence closely after \GMRES. Unfortunately, this method does not have a treatment for system matrices that have large imaginary parts in their eigenvalues. As a consequence of this its convergence behaviour is as unreliable as that of \BiCGstab.
\largeparbreak

\BiCGstabL generalises the polynomial step whereas \IDR generalises the biorthogonalisation. The idea lies near to simply combine both generalisations to combine their advantages. This is done in \IDRstab.

\subsubsection{Derivation of \IDRstab}

The issue with an illustrative derivation of \IDRstab is that one runs out of paper dimensions. Alg.\,\ref{Algo:IDRs_naive} shows an implementation of \IDRstabno. In the following we go through the lines of this implementation and explain the geometric properties of the computed quantities.
\largeparbreak

The first nine lines are identical to the initialisation of \IDR. The while-loop from line 10 to 41 contains the iterative part of the method. During this loop, the vectors $\br\hp{0},\bv\hp{0}_1,...,\bv\hp{0}_s$ are moved from $\cG_j$ to $\cG_{j+\ell}$. 

As always, the iterative scheme consists of two steps: The biorthogonalisation and the polynomial step. In the following we discuss both steps, starting with the biorthogonalisation.
\largeparbreak
The biorthogonalisation works analogous to that of \BiCGstabL: The for-loop from line 13 to 28 in Alg.\,\ref{Algo:IDRstab_naive} has the same functionality as the for-loop from line 9 to 22 in Alg.\,\ref{Algo:BiCGstabL_naive}. In \BiCGstabL, this loop had the purpose that at the end of the $k$th repetition of its loop the following powers of the residual $\br\hp{0}$ and auxiliary vector $\bv\hp{0}$ exist and satisfy 
\begin{align*}
\br\hp{g},\bv\hp{g} &\in \cG_j & &\forall g \in \lbrace n \in \N\ :\ n\geq 0 \ \land \ n\leq k+1 \rbrace\\
\br\hp{g},\bv\hp{g} &\in \cN(\bP) & &\forall g \in \lbrace n \in \N\ :\ n\geq 0 \ \land \ n\leq k \rbrace\,.
\end{align*}
However, in \IDRstab there is not only one auxiliary vector but there are $s$. This is why in the interior of this for-loop with loop-index $k$ there is a nested for-loop in Alg.\,\ref{Algo:IDRstab_naive} from line 21 to 27. For a respective loop-index $q \in \lbrace 1,...,s\rbrace$, this loop biorthogonalises the $k$th power of the $q$th auxiliary vector (cf. line 24) and then computes the subsequent power $\bv_q\hp{k+1}$ in line 25. The geometric procedure of this nested for-loop can be compared to parts 4 and 5 of Fig.\,\ref{fig:Look-ahead2}. However, not only the residual and \textit{one} auxiliary vector are used to do the biorthogonalisation as illustrated in part 4 of the figure but instead the residual and all the \textit{$s$} auxiliary vectors are used.

After the biorthogonalisation the vectors $\br\hp{0},\bv\hp{0}_1,...,\bv\hp{0}_s$ live in $\cG_j\cap\cK^\perp_\ell(\bA\t;\bP)$, cf. line 29 in Alg.\,\ref{Algo:IDRstab_naive}.
\largeparbreak
The polynomial step goes from line 30 to 39. In line 31 the coefficients of the minimised-residual polynomial are computed, cf. for comparison to the comment in Alg.\,\ref{Algo:BiCGstabL_naive} line 26 and Lem.\,\ref{lem:GenRekSonneveldspace}.

Coming back to the \IDRstab implementation in Alg.\,\ref{Algo:IDRstab_naive}, in line 34 the residual is updated to $q_{j,\ell}(\bA)\cdot\br\hp{0}$ in an analogous way to as it is done in the priorly discussed method \BiCGstabL. Line 33 updates the numerical solution in advance such that it is according to the updated residual.

In lines 37 and 38 the auxiliary vectors and their pre-images are updated in an analogous way compared to as it is done for the single auxiliary vector in \BiCGstabL, cf. Alg.\,\ref{Algo:BiCGstabL_naive} lines 32 and 33.

\begin{algorithm}
	\begin{algorithmic}[1]
		\Procedure{IDRstab}{$\bA,\bb,s,\ell,\tolabs$}
		\State $\br\hp{0} := \bb$,\quad $\bx\hp{0} := \bO$,\quad$\bv_1\hp{-1} := \bb$
		\State $\bP := $\texttt{orth}$(\texttt{randn}(N,s))$,\quad$\bv_1\hp{0} := \bA \cdot \bv\hp{-1}_1$
		\For{$q=2,...,s$}
		\State $\bv_q\hp{-1} := \bv\hp{0}_{q-1}$
		\State $\bv_q\hp{0} := \bA \cdot \bv\hp{-1}_q$
		\EndFor
		\State $\bZ := \bP\t \cdot \bV\hp{0} \in \R^{s \times s}$,\quad$j:=0$
		\State \textit{// $\br\hp{0},\bv_1\hp{0},...,\bv_s\hp{0}\in \cK_\infty(\bA;\bb) \equiv \cG_0$}
		\While{$\|\br\hp{0}\|>\tolabs$}
		\State \textit{// $\br\hp{0},\bv_1\hp{0},...,\bv_s\hp{0}\in\cG_j$}
		\State \textit{// - - - Biorthogonalisation - - - }
		\For{$k=0,...,\ell-1$}
		\State \textit{// Residual}
		\State $\bxi := \bZ\d \cdot (\bP\t\cdot\br\hp{k})$
		\State $\br\hp{g} := \br\hp{g} - \bV\hp{g} \cdot \bxi$\quad$\forall\ g=0,...,k$
		\State $\bx\hp{0} := \bx\hp{0} + \bV\hp{-1} \cdot \bxi$
		\State \textit{// $\br\hp{g} \in \cG_j\cap\cN(\bP)$ $\forall\ g=0,...,k$}
		\State $\br\hp{k+1} := \bA \cdot \br\hp{k}$\quad\textit{// $\br\hp{k+1} \in \cG_j$}
		\State \textit{// Auxiliary vectors}
		\For{$q=1,...,s$}
		\State $\bxi := \bZ\d \cdot (\bP\t\cdot\br\hp{k+1})$
		\State $\bv_q\hp{g} := \br\hp{g+1} - [\bV\hp{g+1}_{:,1:(q-1)},\bV\hp{g}_{:,q:s}] \cdot \bxi$\quad$\forall\ g=-1,0,...,k$
		\State \textit{// $\bv_q\hp{g} \in \cG_j\cap\cN(\bP)$ $\forall g=0,...,k$}
		\State $\bv_q\hp{k+1} := \bA \cdot \bv_q\hp{k}$\quad\textit{// $\bv_q\hp{k+1} \in \cG_j$}
		\State $\bZ_{:,q} := \bP\t\cdot\bv_q\hp{k+1}$
		\EndFor
		\EndFor
		\State \textit{// $\br\hp{0},\bv_1\hp{0},...,\bv_s\hp{0} \in \cG_j \cap \cK_\ell^\perp(\bA\t;\bP)$}
		\State \textit{// - - - Polynomial step - - - }
		\State $\btau := [\br\hp{1},...,\br\hp{\ell}]\d \cdot \br\hp{0}$
		\State \textit{// Residual}
		\State $\bx\hp{0} := \bx\hp{0} + \sum_{k=1}^\ell \br\hp{k-1} \cdot \tau_k$
		\State $\br\hp{0} := \br\hp{0} - \sum_{k=1}^\ell \br\hp{k} \cdot \tau_k$
		\State \textit{// $\br\hp{0} \in q_{j,\ell}(\bA)\cdot\big(\cG_j \cap \cK_\ell^\perp(\bA\t;\bP)\big) \equiv \cG_{j+\ell}$}
		\State \textit{// Auxiliary vectors}
		\State $\bV\hp{-1} := \bV\hp{-1} - \sum_{k=1}^\ell \bV\hp{k-1} \cdot \tau_k$
		\State $\bV\hp{0} := \bV\hp{0} - \sum_{k=1}^\ell \bV\hp{k} \cdot \tau_k$
		\State \textit{// $\bv_1\hp{0},...,\bv_s\hp{0} \in q_{j,\ell}(\bA)\cdot\big(\cG_j \cap \cK_\ell^\perp(\bA\t;\bP)\big) \equiv \cG_{j+\ell}$}
		\State $\bZ:=-\tau_\ell \cdot \bZ$,\quad$j:=j+\ell$
		\EndWhile
		\State \Return $\bx\hp{0}$
		\EndProcedure
	\end{algorithmic}
	\caption{\IDRstab}\label{Algo:IDRstab_naive}
\end{algorithm}

\subsubsection{Properties of \IDRstab}\label{sec:PropertiesIDRstab}
\IDRstab combines the properties of \BiCGstabL and \IDR. 

\IDRstabno has a breakdown when $\tau_\ell=0$ or when the biorthgonalisations in line 16 respectively 18 or 23--24 fail.

\paragraph{Finite termination}
We have seen that \BiCGstabL has exactly the same termination properties as \BiCGstab because both use only different recurrences for the same spaces (however, with different values for the $\omega$-values).

In analogy to this, \IDRstab has the same termination properties as \IDR: During the main loop, i.e. the while-loop, the vectors $\br\hp{0},\bv_1\hp{0},...,\bv_s\hp{s}$ are moved from $\cG_j$ into $\cG_{j+\ell}$. Since $\rk(\bP)=s$ it follows that the dimension of $\cG_{j+\ell}$ is either zero, i.e. the method terminates, or $\cG_{j+\ell}$ is by $\ell\cdot s$ dimensions smaller than $\cG_j$. Thus, the method terminates after at most $\lceil N/(\ell\cdot s)\rceil$ repetitions of the main loop.

During each repetition of the main loop $\ell\cdot(s+1)$ matrix-vector products are computed. Consequently, the method terminates after at most $\lceil N/(\ell\cdot s)\rceil \cdot \ell\cdot(s+1)$ matrix-vector products, which is roughly identical to the termination properties of \IDR. Consequently, we can expect that \IDRstab experiences the transition to superlinear convergence roughly at the same matrix-vector product as \IDR.

\paragraph{Computational cost} Here comes the drawback of \IDRstab. We have seen that both \BiCGstabL and \IDR have some overhead compared \BiCGstab. \IDRstab multiplies these overheads as we lay out in the following.

\IDRstab in the above-described implementation from Alg.\,\ref{Algo:IDRstab_naive} requires storage for $(s+1) \cdot (\ell+2) + s$ vectors. This is because for each the residual and the auxiliary vectors $\ell+1$ further vectors must be stored, namely their $\ell$ powers and their pre-images (respectively for the residual the according numerical solution $\bx\hp{0}$). Further to that, the $s$ column vectors of $\bP$ must be stored. Thus, for both directions in which we generalise the method \BiCGstab there occurs a factor in the number of column vectors of length $N$ for which the method requires memory.

The computational cost of the above implementation in DOTs and AXPYs per repetition of the main-loop is determined by counting to
\begin{align*}
\text{\#DOTs} & = \ell\cdot\left(2+s^2 + \frac{\ell+1}{2}\right)\\
\text{\#AXPYs} & = \ell\cdot(s+1)\cdot\left(2+s\cdot \frac{\ell+3}{2}\right)\,.
\end{align*}
Since per repetition of the main-loop the number of matrix-vector products is $\ell\cdot(s+1)$, one obtains per matrix-vector product an overhead that is bilinear in $s$ and $\ell$.

\section{Practical implementations of \IDRstab and \IDR}\label{sec:ImplementationsIDRstab}
In this section we discuss practical implementations of \IDRstab and \IDR. 

\begin{Remark} [Exception for the case $\ell=1$]
On the first glance there is no obvious benefit of discussing an implementation of \IDR separately from \IDRstab since it is covered by an implementation of \IDRstab. However, it turns out, as we show in later subsections, that for the case $\ell=1$ the implementation should massively differ from the case $\ell>1$ to take advantage from algebraic properties of \IDR that \IDRstab for $\ell\neq 1$ does not have.
\end{Remark}

We lay out the structure of this section. In the next subsection we motivate why it is important to discuss more practical implementations of \IDRstab and \IDR. To this end, we review some of the above-mentioned breakdown scenarios and stability issues of \IDRstab and \IDR. Then we give an overview how practical implementations try to overcome these issues.

After a broad motivation and overview, we show some useful algorithmic blocks that help our discussion of the several possible variants in which \IDR and \IDRstab can be implemented. These blocks consist of a standard initialisation scheme, a proper way to compute the coefficients $\btau$ for the polynomial step, and finally some subroutines for the iterative orthogonalisation and biorthogonalisation of the auxiliary vectors. After having prepared some algorithmic blocks we start with the discussion of practical implementations of \IDR and \IDRstab.

\subsection{Motivation}
We have discussed that \IDRstab respectively \IDR both suffer from breakdowns when either the biorthogonalisation fails or when the value $\tau_\ell$ respectively $\tau_1$ becomes zero. This however is only the case when there is no numerical round-off.

When dealing with numerical round-off then in addition we have to make sure that the biorthogonalisation is performed in an appropriate sense to a high relative accuracy and that the respective $\tau$-value is in an appropriate relative sense far enough away from zero.
\largeparbreak
For the biorthogonalisation on the one hand, we can improve the condition number of the matrix $\bZ$ such that the linear-factors $\bxi \in \R^s$ in the biorthogonalisation are computed with a higher accuracy. Although this does not prevent the method from breakdowns when $\bZ$ was singular in exact arithmetic it still makes the method more robust when dealing with numerical round-off.

A desirable goal to achieve more well-conditioned matrices for $\bZ$ is to use not only an orthogonal matrix for $\bP$ but also for $\bV\hp{0}$. This can be explained as follows: Consider the singular value decomposition $\bZ = \bQ_L\t \cdot \bS_Z \cdot \bQ_R$, where $\bQ_L,\bQ_R$ have condition 1 and $\bS_Z$ is symmetric positive semi-definite and diagonal. Then $(\bP\cdot\bQ_L)\t\cdot(\bV\hp{0}\cdot\bQ_R) = \bS_Z$ holds and the diagonal values of $\bS_Z$ are the arc-cosines of the principal angles between the hyperplanes $\rg(\bP)$ and $\rg(\bV\hp{0})$ \cite{principalAngles}. I.e., the condition of $\bZ$ is bounded from below by the geometry of the Sonneveld spaces that $\bP$ determines. Thus, the best achievable condition of $\bZ$ is obtained numerically when both $\bP$ and $\bV\hp{0}$ are unitary.

We will discuss in a subsequent subsection how it can be achieved in an efficient way that the basis matrices remain well-conditioned up to a certain degree.
\largeparbreak

For the avoidance of stability polynomial roots $\omega=0$ on the other hand, there are modified routines for the computation of the polynomial coefficients $\btau$ that can be utilised to enforce a geometrically meaningful off-set from zero. The geometric meaning can be derived by considering biorthogonality properties of the powers of the residual in advance of the execution of the actual polynomial update.

We will discuss the geometric ideas of these sophisticated strategies for the computation of the polynomial coefficients and we will contribute a variant that realises the underlying mathematical idea in a more accurate formula.
\largeparbreak

Besides to the above motivations, there are further good reasons to spend more emphasize on practical implementations of \IDR and \IDRstab. One of these reasons is that one can reduce the computational cost for the above methods: \IDR for instance can be implemented using only storage for $3 \cdot (s+1)$ vectors instead of -- as suggested above -- using storage for $3 \cdot (s+1) + s$ vectors. Further, the authors in \cite{IDR-biorth,IDRstab-biorth} say that both \IDR and \IDRstab can be implemented with fewer AXPYs and DOTs. Finally, \IDR can be implemented in such a way that the relation
\begin{align*}
	\bv\hp{0}_q = \bA \cdot \bv\hp{-1}_q\quad\quad \forall\,q=1,...,s
\end{align*}
holds with maximum possible accuracy, i.e. there is no kind of decoupling for the auxiliary vectors in \IDR.

\largeparbreak
In the rest of this section we discuss all the different implementation variants of \IDR and \IDRstab. However, to keep the discussion efficient and compact, it is beneficial to first prepare a few algorithmic blocks of which all the methods can be composed in a simpler way.

\subsection{Algorithmic blocks}
In this subsection we discuss algorithmic subroutines that are useful for designing implementation variants of \IDR and \IDRstab.

First, we will contribute a new initialisation procedure that generates all the required initial data to start the main-loop of \IDR and \IDRstab. Our novel initialisation has advantages in the condition and structure of some matrices.

Second, we derive and discuss an enhanced scheme due to Fokkema \cite{EnhancedBiCGstabL} and Sleijpen and van der Vorst \cite{MaintainingConvergenceBiCGstab} to compute the polynomial coefficients $\btau$ for an arbitrary value of $\ell \in \N$. We then contribute a refinement of their scheme.


\subsubsection{The Initialisation}
Most \IDRstab respectively \IDR implementations do not spend a lot of effort into the initialisation. For example, in \IDRbiortho \cite{IDR-biorth} the initialisation is not even treated separately.

However, we believe that the initialisation is a very good opportunity to set up the method in a proper way for the following two reasons:
\begin{enumerate}
	\item The initialisation does not add a considerable cost to the method itself since in contrast to the iterative scheme it is only executed once.
	\item When the initialisation is not robust then the damage to the numerical accuracy is already done.
\end{enumerate}

In Algorithm \ref{Algo:Initialisation} we contribute a stable and computationally efficient implementation of the initialisation. We discuss this initialisation by first postulating some mathematical properties of its returned outputs and then describing the geometric constructions that are done in the algorithm in order to obtain these properties.

\paragraph{Properties of the returned outputs}
As inputs, the initialisation receives the system matrix $\bA$, the right-hand side $\bb$, an initial guess $\bx_0$, the \IDRO parameter $s \in \N$ and a stopping tolerance $\tolabs \in \R^+$ for the absolute residual.

It returns the following quantities: 
\begin{itemize}
	\item $\bP \in \R^{N \times s}$. This matrix plays a role in the recursive definition of the Sonneveld spaces that the \IDRO methods use. It has orthonormalised column vectors.
	\item $\bV\hp{-1},\,\bV\hp{0} \in \R^{N \times s}$. These matrices are used for the oblique projections into $\cN(\bP)$. They satisfy:
	\begin{align*}
	\bA \cdot \bV\hp{-1} &= \bV\hp{0}\\
	{\bV\hp{0}}\t\cdot\bV\hp{0} &= \bI
	\end{align*}
	The orthonormality of the columns of $\bV\hp{0}$ is beneficial since it makes sure that $\bZ$ has the best-available condition number subject to the geometric circumstances that the principal angles between the spaces $\rg(\bP)$ and $\rg(\bV\hp{0})$ impose.
	\item $\bZ \in \R^{s \times s}$ is the matrix from \eqref{eqn:Z=PtV}. In our initialisation this matrix is constructed such that it is lower triangular. This is useful in some methods to reduce the number of computed AXPYs and DOTs.
	\item $\bx\hp{0},\br\hp{0} \in \R^N$ are the current numerical solution and its according residual. $\beta$ is the 2-norm of this residual. The residual satisfies $\br\hp{0} \in \cN(\bP)$.
\end{itemize}

Our initialisation is yet the only one that yields that $\bV\hp{0}$ is unitary \emph{and} $\bZ$ is lower triangular. We believe that missing either of these properties is disadvantageous for the following reasons:

A bad condition of $\bV\hp{0}$ leads to a bad condition of $\bZ$. This can be understood by considering the case where one decreases the angle between two column vectors of $\bV\hp{0}$. As a consequence of that, the angle between the same columns in $\bZ$ will decrease. A bad condition of $\bZ$ in turn will probably spoil the superlinear rate of convergence of \IDR respectively \IDRstab.

An argument that speaks for making $\bZ$ lower triangular is that in this case some biorthogonalisations during the iterative scheme can be performed more efficiently in terms of computational cost. In particular, the methods \IDRbiortho \cite{IDR-biorth} and \IDRstabbiortho \cite{IDRstab-biorth} yield a reduction in the number of AXPYs and DOTs per iteration by exploiting that $\bZ$ is lower triangular. Both of these methods are discussed in later sections.

\paragraph{Explanation of the computational steps of the algorithm}
We explain the computational steps of Algorithm \ref{Algo:Initialisation}. In line 4 a \GMRES method of $s$ iterations is performed for the initial guess $\bx_0$ and its according initial residual $\br_0$. If \GMRES already finds a residual that satisfies the tolerance $\tolabs$ then it terminates and returns a sufficiently accurate numerical solution for $\bx$. Otherwise, it returns the matrices $\bW_{:,1:(s+1)}$, $\bQ_H$, $\bR_H$ of an Arnoldi decomposition with the properties from the comment in line 5.

Since after line 5 it is clear that $s$ steps of \GMRES did not converge to a sufficiently accurate solution, all the matrices for \IDRstab must be prepared. Thus, in line 6 the matrix $\bP$ is built. This is done by initialising an $N \times s$ matrix from randomly generated numbers of a normal distribution and then orthonormalising its columns.
\largeparbreak

The lines 7 to 13 realise in an efficient way a biorthogonalisation of the residual. We start the explanation of these lines by using line 11: As the comment says, $\bxi \in \R^s$ is the vector such that
\begin{align}
\br_0 - \bA \cdot \bW_{:,1:s} \cdot \bxi \in \cN(\bP) \label{eqn:Init:BioRes}
\end{align}
holds. In the rest of this paragraph we explain first why the commented equivalence in line 11 holds. Afterwards we explain what is computed in lines 12 to 13. Finally, we explain line 14.

Inserting in this order the decomposition from line 10 and the matrix $\bY$ from line 7 into the expression for $\bxi$ in line 11, it follows:
\begin{align*}
\bxi &= \bR_H^{-1} \cdot \underbrace{\bQ_Z \cdot \bL_Z^{-1}}_{\equiv \bZ^{-1}} \cdot \bEta\\
&= \bR_H^{-1} \cdot \bZ^{-1} \cdot \bEta\\
&= (\underbrace{\bZ}_{\equiv \bY \cdot \bQ_H} \cdot \bR_H)^{-1} \cdot \underbrace{\bEta}_{\equiv \beta \cdot \bY_{:,1}}\\
&= (\underbrace{\bY}_{\equiv \bP\t \cdot \bW_{:,1:(s+1)}} \cdot \bQ_H \cdot \bR_H)^{-1} \cdot (\beta \cdot \underbrace{\bY_{:,1}}_{\equiv\bP\t \cdot \bW_{:,1}})\\
&= (\bP\t \cdot \underbrace{\bW_{:,1:(s+1)} \cdot \bQ_H \cdot \bR_H}_{\equiv\bA\cdot\bW_{:,1:s}})^{-1}  \cdot (\bP\t \cdot \underbrace{\bW_{:,1}	\cdot \beta}_{\equiv\br_0})
\end{align*}
In the last line we inserted the properties of \GMRES{}'s Arnoldi equation.

Now we explain the lines 12 to 13. In line 12 the biorthogonalised residual from \eqref{eqn:Init:BioRes} is computed in the basis $\bW_{:,1:(s+1)}$. From the linear combination $\bxi$ for the update of $\bx_0$ by columns of $\bW_{:,1:s}$, the update of $\br_0$ is $-\bQ_H \cdot \bR_H \cdot \bxi$. The linear combination vector of $\br_0$ in turn is $\beta \cdot \be_1$, since $\bW_{:,1}$ is the normalised initial residual vector. Building the sum of both linear-combination vectors, we obtain $\bc\hp{0}$.

Finally, in line 14 the new projectors are computed such that all the following properties hold:
\begin{align*}
\bV\hp{0} &= \bA \cdot \bV\hp{-1}\\
{\bV\hp{0}}\t \cdot \bV\hp{0} &= \bI\\
\bZ &= \bP\t \cdot \bV\hp{0}\quad\text{is lower triangular}
\end{align*}
We show in the above order for each of these conditions that it holds.
To show the first equation, we insert the expressions from line 14 into the equation and afterwards multiply from the right by $\bQ_Z^{-1} \cdot \bR_H$. As a result we obtain the Arnoldi equation, which obviously holds:
\begin{align*}
& &\bV\hp{0} &= \bA \cdot \bV\hp{-1}\\
\Leftrightarrow & &\bW_{:,1:(s+1)} \cdot \bQ_H \cdot \bQ_Z &= \bA \cdot \bW_{:,1:s} \cdot \bR_H^{-1} \cdot \bQ_Z\\
\Leftrightarrow & &\bW_{:,1:(s+1)} \cdot \bQ_H \cdot \bR_H&= \bA \cdot \bW_{:,1:s} 
\end{align*}

The second equation is trivial to show: $\bV\hp{0}$ is computed as a product of matrices of condition 1, respectively. Each of the matrices from the product has at most as many columns as rows. In consequence, $\bV\hp{0}$ has condition 1.

The third property can be shown by inserting the formula of $\bV\hp{0}$ from line 14 into the equation. However, to avoid confusions we use the names $\bZ,\bQ_Z,\bL_Z$ from line 10, i.e. we show $\bL_Z = \bP\t \cdot \bV\hp{0}$. 
\begin{align*}
& &\bL_Z &= \bP\t \cdot \bV\hp{0}\\
\Leftrightarrow & &\bZ \cdot \bQ_Z &= \bP\t \cdot \bW_{:,1:(s+1)} \cdot \bQ_H \cdot \bQ_Z\\
\Leftrightarrow & &\bZ &= \underbrace{\bP\t \cdot \bW_{:,1:(s+1)}}_{\equiv\bY} \cdot \bQ_H
\end{align*}
From the second to the third line above we have multiplied from the right by $\bQ_Z^{-1}$. We see that the last equation obviously holds since it matches the assignment of $\bZ$ from line 9 in the algorithm.

\largeparbreak
In the remainder of this subsection we explain the advantages of this way of implementing the initialisation:

The advantage of computing the new residual form the columns of $\bW_{:,1:(s+1)}$ is that the computation of $\br\hp{0}$ from $\bc\hp{0}$ is stable. However, $\bc\hp{0}$ in turn cannot be stable because on the one hand its norm can be much smaller than that of $\br_0$ and on the other hand $\bR_H$ can be badly conditioned\footnote{The condition of $\bR_H$ is bounded by that of $\bA$.}. Since $\bc\hp{0}$ is only of dimension $\R^{s+1}$ in a software implementation one could solve a sub-problem where $\bxi$ and $\bc\hp{0}$ are modified in such a way that the error of the equation in line 12 from their digital representations is minimised.

A further advantage of our implementation in comparison to a naive implementation like for instance
\begin{algorithmic}[1]
	\State $\bV\hp{-1} := \bW_{:,1:s} \cdot \bR_H^{-1}$,\quad$\bV\hp{0} := \bW_{:,1:(s+1)} \cdot \bQ_H$
	\State $\bZ := \bP\t \cdot \bV\hp{0}$,\quad$\bEta:=\bP\t\cdot\br_0$
	\State LQ-decomposition: $\bL_Z = \bZ \cdot \bQ_Z$,\quad$\bZ:=\bL_Z$
	\State $\bV\hp{-1}:= \bV\hp{-1}\cdot\bQ_Z$,\quad$\bV\hp{0}:=\bV\hp{0}\cdot\bQ_Z$
	\State $\bxi := \bZ^{-1} \cdot (\bP\t\cdot\br)$
	\State $\bx\hp{0}:=\bx_0 + \bV\hp{-1}\cdot\bxi$,\quad$\br\hp{0}:=\br_0-\bV\hp{0}\cdot\bxi$
	\State $\beta := \|\br\hp{0}\|$
\end{algorithmic}
is that our implementation is simply cheaper in terms of BLAS-1 operations. In particular, the computations in lines 1 and 4 of the naive variant cost $4 \cdot s^2$ AXPYs, whereas in our implementation we combine these two operations in line 14 and thus obtain only $2 \cdot s^2$ AXPYs. 

Further, the additional matrix-products in line 4 introduce an error to the accuracy with which the crucial relation $\rg(\bV\hp{0}) \subset \cK_\infty(\bA;\br_0)$ holds. This is because the round-off from this matrix-product spoils the column vectors of $\bV\hp{0}$ such that the smallest invariant subspace in which $\rg(\bV\hp{0})$ lives is potentially much larger in terms of dimensions than $\cK_\infty(\bA;\br_0)$.

\begin{algorithm}
	\begin{algorithmic}[1]
		\Procedure{Initialisation}{$\bA,\bb,s,\tolabs$}
		\State \textit{// $\bA \in \R^{N \times N},\ \bb\in \R^N,\ s \in \N$}
		\State $\br := \bb$,\quad$\bx:=\bO$,\quad $\beta := \|\bb\|$
		\State [$\bW_{:,1:(s+1)},\,\bQ_H,\,\bR_H$] = \texttt{GMRESm}{($\bA,\bx,\br,\beta,s,\tolabs$)}
		\State \textit{// $\bA \cdot \bW_{:,1:s} = \bW_{:,1:(s+1)} \cdot \bQ_H \cdot \bR_H$,\ $\bQ_H \in \R^{(s+1) \times s},\ \bR_H \in \R^{s \times s}$}
		\State $\bP := $\texttt{orth( randn($N,s$) )}\quad \textit{// $\in \R^{N \times s}$}
		\State $\bY := \bP\t \cdot \bW_{:,1:(s+1)}$
		\State $\bEta := \beta \cdot \bY_{:,1}$\quad\textit{// $\equiv \bP\t \cdot \bb$}
		\State $\bZ := \bY \cdot \bQ_H$\quad\textit{// $\bZ \in \R^{s \times s}$ has best available condition}
		\State [$\bQ_Z,\bL_Z$] = \texttt{lq}($\bZ$)
		\State $\bxi := \bR_H^{-1} \cdot \bQ_Z \cdot \bL_Z^{-1} \cdot \bEta$\quad\textit{// $\equiv (\bP\t \cdot \bA \cdot \bW_{:,1:s})^{-1} \cdot (\bP\t \cdot \br_0)$}
		\State $\bc\hp{0} := \beta \cdot \be_1 - \bQ_H \cdot \bR_H \cdot \bxi$,\quad$\bZ:=\bL_Z$
		\State $\bx\hp{0} := \bx + \bW_{:,1:s}\cdot \bxi$,\quad$\br\hp{0} := \bW_{:,1:(s+1)} \cdot \bc\hp{0}$,\quad$\beta:= \|\bc\hp{0}\|$
		\State $\bV\hp{-1} := \bW_{:,1:s} \cdot (\bR_H^{-1} \cdot \bQ_Z)$,\quad $\bV\hp{0} := \bW_{:,1:(s+1)} \cdot (\bQ_H \cdot \bQ_Z)$
		\State \Return $\bP,\bV\hp{-1},\bV\hp{0},\bZ,\bx\hp{0},\br\hp{0},\beta$
		\EndProcedure
	\end{algorithmic}
	\caption{Initialisation}\label{Algo:Initialisation}
\end{algorithm}

\subsubsection{The stability polynomial coefficients}\label{sec:ThePolynomialCoefficients}

\paragraph{Motivation}
All the \IDRO methods discussed so far have a biorthogonalisation and a polynomial part. The first computational step in the polynomial part is to determine suitable coefficients $\tau_1,...,\tau_\ell$ for the polynomial
\begin{align*}
	q_{j,\ell}(t) = t^0 - \sum_{k=1}^\ell \tau_k \cdot t^k \equiv \prod_{k=1}^{\ell} (t^0 - \omega_{j+k} \cdot t^1)\,.
\end{align*}
Given the residual $\br\hp{0}$, its $\ell$ powers $\br\hp{1},...,\br\hp{\ell}$, and the coefficient vector
\begin{align*}
	\btau = \begin{pmatrix}
	\tau_1\\
	\vdots\\
	\tau_\ell
	\end{pmatrix}\,,
\end{align*}
one uses the following update formula for the residual:
\begin{align*}
	\br\hp{0} := q_{j,\ell}(\bA) \cdot \br\hp{0} = \br\hp{0} - [\br\hp{1},...,\br\hp{\ell}] \cdot \btau
\end{align*}
From the formula it follows that a minimisation of $\|q_{j,\ell}(\bA)\cdot\br\hp{0}\|$ is achieved by solving the following least-squares problem:
\begin{align}
	\btau = \operatornamewithlimits{argmin}_{\bchi \in \C^\ell}\big\|\br\hp{0}-[\br\hp{1},...,\br\hp{\ell}]\cdot\bchi\big\|\label{eqn:tauOpt}
\end{align}
Obviously, if $\rk([\br\hp{0},...,\br\hp{\ell}]) < \ell+1$ then $\exists\, \btau \ : \ q_{j,\ell}(\bA) \cdot \br\hp{0} = \bO$. Since this case is trivial to handle, we can assume in the following that $\rk([\br\hp{0},...,\br\hp{\ell}]) = \ell+1$ and consequently \eqref{eqn:tauOpt} has a unique solution. The issue with this solution is however that the possibility $\tau_\ell = 0$ is not excluded. 

A numerical treatment to this issue is to compute a semi-optimal solution for $\btau$ which yields that $|\tau_\ell|$ is \enquote{sufficiently} larger than zero. In the following we will discuss what a reasonable relative measure for \enquote{sufficiently} is. Then we propose two computational schemes that provide a semi-optimal solution for $\btau$ that satisfies this relative measure.

\paragraph{The geometric approach}

In \cite{MaintainingConvergenceBiCGstab,EnhancedBiCGstabL} it is analysed how the accuracy of the \BiCG coefficients as computed in \BiCGstab and \BiCGstabL is affected by the values of $\omega$ of the polynomial step.

Let us have a look into \cite[p.\,204 eqn.\,(3)]{MaintainingConvergenceBiCGstab} and Alg.\,\ref{Algo:BiCGstab_naive} of our thesis. The scalars $\sigma_k$ of their paper are our matrices $\bZ \in \R^{1 \times 1}$ and their scalars $\rho_k$ are $\bP\t\cdot\br\hp{0}$ in line 10 of Alg.\,\ref{Algo:BiCGstab_naive} of our thesis.

The authors argue that for local maintenance of the convergence one has to make sure that the rounding erros in $\bZ$ and $\bP\t\cdot\br\hp{0}$ are small. In particular, they spend emphasize on the latter quantity. In the following we lay out the idea:
\largeparbreak

Considering the $j$th repetition of the main-loop in \BiCGstab, it holds
\begin{align*}
	\bEta := \bP\t \cdot \br\hp{0}\fp{j} \equiv \underbrace{\bP\t \cdot \tbr\hp{0}\fp{j-1}}_{\equiv \bO} - \omega_j \cdot \bP\t\cdot\tbr\hp{1}\fp{j-1}\,,
\end{align*}
cf. line 10 in Alg.\,\ref{Algo:BiCGstab_naive_indexed} respectively line 10 in Alg.\,\ref{Algo:BiCGstab_naive}. The $:=$ defines what is numerically computed with round-off. The $\equiv$ gives an analytical equivalent expression in exact arithmetic. The approach is to choose $\omega_j$ such that $\bEta$ is computed with a small relative error.

We can decompose the above formula further: We split $\tbr\hp{1}\fp{j-1} = (\tbr\hp{1}\fp{j-1})_t+(\tbr\hp{1}\fp{j-1})_n$, where $(\tbr\hp{1}\fp{j-1})_t$ is the tangent component to $\tbr\hp{0}\fp{j-1}$ and $(\tbr\hp{1}\fp{j-1})_n$ is the normal component. Now, the formula is:
\begin{align*}
\bEta := \bP\t \cdot \br\hp{0}\fp{j} \equiv \underbrace{\bP\t \cdot \tbr\hp{0}\fp{j-1} - \omega_j \cdot (\tbr\hp{1}\fp{j-1})_t}_{\equiv \bO} - \omega_j \cdot \bP\t\cdot(\tbr\hp{1}\fp{j-1})_n
\end{align*}

Only theoretically it holds that the first two terms are zero, however due to numerical round-off there will be a non-zero contribution of these terms to the result of $\bEta$. For very small values of $\omega_j$ this non-zero contribution becomes dominant and destroys the superlinear convergence.

So as an idea, the value of $\omega_j$ must be sufficiently large such that $\|\omega_j\cdot\bP\t\cdot(\tbr\hp{1}\fp{j-1})_n\|$ is much larger than $\|\bP\t\cdot(\tbr\hp{0}\fp{j-1})_t\|$. Since the latter of these norms would be zero in theory and is not intended to be computed once again, one simply says that the vector $\br\hp{0}\fp{j}$ must have a sufficiently small angle with $(\br\hp{1}\fp{j-1})_n$. This suffices to derive a rule for computing the semi-optimal solution for $\btau$. This rule is explained in the following.
\largeparbreak
Let $q^\text{MR}_{j,\ell}(\bA) \cdot\br\hp{0}$ be the polynomial update that leads to a minimisation of the 2-norm (MR = minimal residual). Let further $q^\text{OR}_{j,\ell}(\bA) \cdot\br\hp{0}$ be the polynomial update such that $q^\text{OR}_{j,\ell}(\bA) \cdot\br\hp{0} \perp \opspan\lbrace\br\hp{0},...,\br\hp{\ell-1}\rbrace$ (OR = orthogonal residual). In \BiCGstab the latter vector would be a parallel vector to the component $(\tbr\hp{1}\fp{j-1})_n$. The design approach for the finally used polynomial $q_{j,\ell}(\cdot)$ is to combine $q^\text{MR}_{j,\ell}$ and $q^\text{OR}_{j,\ell}$ to a polynomial $q_{j,\ell}$ such that the angle $\kappa$
\begin{align*}
	\kappa = \angle\big(\, q_{j,\ell}(\bA) \cdot\br\hp{0}\,,\,q^\text{OR}_{j,\ell}(\bA) \cdot\br\hp{0} \,\big)
\end{align*}
is limited from above by a user-defined parameter angle $\alpha$.

\begin{figure}
	\centering
	\includegraphics[width=1\linewidth]{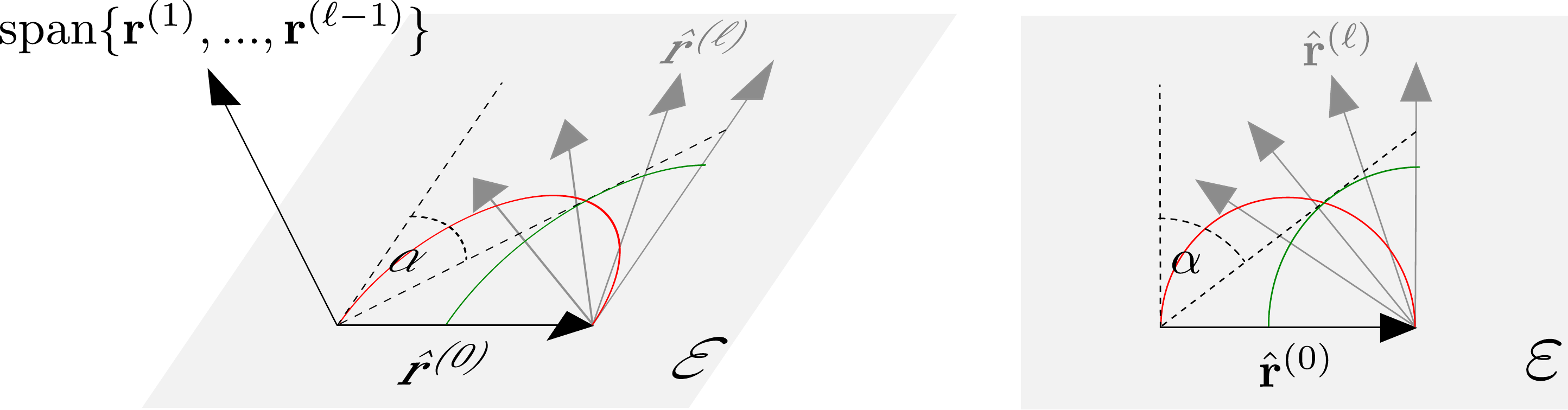}
	\caption{Geometric illustration of the enhanced computation of the stability polynomial coefficients.}
	\label{fig:StabilityStep}
\end{figure}

\paragraph{Computation of the coefficient vector}
Let us consider Fig.\,\ref{fig:StabilityStep}. The left part of the figure shows the two-dimensional plane $\varepsilon$ that is defined as $\varepsilon := \rg([\br\hp{0},...,\br\hp{\ell}])\cap\cN([\br\hp{1},...,\br\hp{\ell-1}])$.

In this two-dimensional plane we have plotted the vectors $\hbr\hp{0},\hbr\hp{\ell}$ that are obtained {after} orthogonalising $\br\hp{0},\br\hp{\ell}$ against $\br\hp{1},...,\br\hp{\ell-1}$. In the right part of the figure the projection plane with the projected vectors $\hbr\hp{0},\hbr\hp{\ell}$ is shown from the top. The vector $\hbr\hp{\ell}$ is shown in grey in various possible angles to $\hbr\hp{0}$.

In the following we draw the bridge from the figure to the above formula for $q^\text{MR}$ and $q^\text{OR}$:
\largeparbreak
It is
\begin{align*}
	\hbr\hp{0} &= \br\hp{0} - [\br\hp{1},...,\br\hp{\ell-1}] \cdot \hbtau\hp{0}\\
	\hbr\hp{\ell} &= \br\hp{\ell} - [\br\hp{1},...,\br\hp{\ell-1}] \cdot \hbtau\hp{\ell}\\
	\hbtau\hp{g} &= [\br\hp{1},...,\br\hp{\ell-1}]\d \cdot \br\hp{g},\quad g\in \lbrace 0,\ell\rbrace\,.
\end{align*}
Considering a linear combination $\hbr\hp{0} + \delta \cdot \hbr\hp{\ell}$, there exists a coefficient vector
\begin{align*}
	\hbtau = \begin{pmatrix}
	1\\
	-\hbtau\hp{0}\\
	0
	\end{pmatrix} + \delta \cdot \begin{pmatrix}
	0\\
	-\hbtau\hp{\ell}\\
	1
	\end{pmatrix} \in \R^{\ell+1}
\end{align*}
such that
\begin{align*}
	\hbr\hp{0} + \delta \cdot \hbr\hp{\ell} \equiv \underbrace{[\br\hp{0},\br\hp{1},...,\br\hp{\ell}]\cdot \hbtau}_{\equiv q_{j,\ell}(\bA)\cdot\br\hp{0}}\,.
\end{align*}
In other words, there exists a polynomial $q_{j,\ell}(\cdot)$ with $q_{j,\ell}(0)=1$ such that the new residual is constructed as $\hbr\hp{0} + \delta \cdot \hbr\hp{\ell}$ and lives in the plane $\varepsilon$. And this is just the way how $q_{j,\ell}(\cdot)$ will be chosen. The only remaining question is how to choose $\delta$. From the figure we will deduce which value for $\delta \in \R$ one should use. Clearly there exists a value for $\delta$ such that  $\hbr\hp{0} + \delta \cdot \hbr\hp{\ell} \equiv q_{j,\ell}^\text{MR}(\bA)\cdot\br\hp{0}$, so the construction that we have made so far does not destroy the method's potential to choose residual-minimal stability-polynomials. Nevertheless, we still have to make sure that $\delta$ is sufficiently far away from zero since it holds $\delta=0 \ \Rightarrow\ \tau_\ell=0 \ \Rightarrow\ \omega_{j+\ell}=0$.
\largeparbreak

In the previous paragraph we have motivated to choose an update for the residual such that the angle between the updated residual and the line $\rg([\br\hp{0},...,\br\hp{\ell}])\cap\cN([\br\hp{0},...,\br\hp{\ell-1}])$ is not larger than a user-defined parameter $\alpha$. Using Fig.\,\ref{fig:StabilityStep}, this means in geometric terms that $\delta$ must be chosen sufficiently large such that the updated residual lives within the two stressed lines.

In the figure we have not only shown the angle $\alpha$ but also a red and a green fraction of a circle. The red curve is a fraction of the Thalis circle. Given a direction of $\hbr\hp{\ell}$, this curve helps finding the optimal contribution $\delta$ of $\hbr\hp{\ell}$ to minimise the length $\|\hbr\hp{0} + \delta \cdot \hbr\hp{\ell}\|$. Unfortunately, some points on this curve do not satisfy the angular requirement because they are below the rightwards stressed line. This means that for these points the value $\delta$ would not be sufficiently large.

Besides, there is a green curve. It is a part of the circle that has the centre in the arrow head of $\hbr\hp{0}$ and the a radius such that it is tangent to the rightwards stressed line.
\largeparbreak

In the following we provide three formulas for $\delta$ w.r.t. $\hbr\hp{1}$,$\hbr\hp{\ell}$ such that the head point of the vector $\hbr\hp{0} + \delta \cdot \hbr\hp{\ell}$ lives on the red, the green or the rightwards stressed line. Given the arc-cosine $\rho$ of the angle $\kappa$ between $\hbr\hp{0}$ and $\hbr\hp{\ell}\phantom{\frac{a}{b}}$
\begin{align*}
	\rho = \frac{(\hbr\hp{\ell})\t \cdot \hbr\hp{0}}{\|\hbr\hp{\ell}\| \cdot \|\hbr\hp{0}\|}\,,
\end{align*}
the formulas for the $\delta$ values are:
\begin{align*}
	\delta^\texttt{red\textcolor{white}{aaaaa}} &= - \frac{\|\hbr\hp{0}\|}{\|\hbr\hp{\ell}\|} \cdot \rho\\
	\delta^\texttt{green\textcolor{white}{aaa}} &= - \frac{\|\hbr\hp{0}\|}{\|\hbr\hp{\ell}\|} \cdot \sign(\rho) \cdot \sin(\alpha)\\
	\delta^\texttt{stressed} &= - \frac{\|\hbr\hp{0}\|}{\|\hbr\hp{\ell}\|} \cdot \sign(\rho) \cdot \frac{\sin(\alpha)}{\sin\Big(\alpha+\frac{\pi}{4}\cdot\big(\sign(\rho)+1\big)-\arccos(\rho)\Big)}
\end{align*}
The value $\delta^\texttt{red}$ yields the value for $\delta$ such that the updated residual is orthogonal to $\hbr\hp{\ell}$. Since $\hbr\hp{0}$ and $\hbr\hp{\ell}$ are orthogonal w.r.t. $\br\hp{1},...,\br\hp{\ell-1}$, it follows that this choice for $\delta$ results in that the updated residual is orthogonal w.r.t. $\br\hp{1},...,\br\hp{\ell}$ Thus, choosing $\btau$ following this rule yields the least-squares solution \eqref{eqn:tauOpt} for $\btau$.
\largeparbreak

When the other values for $\delta$ are close to $\delta^\texttt{red}$ then this means that the updated residual is not so far away from the minimum possible residual. The general strategy is to choose:
\begin{align*}
	\delta := \begin{cases}
		\delta^\texttt{red},& |\rho|\geq \sin(\alpha)\\ \delta^\texttt{green}\text{ or } \delta^\texttt{stressed},& |\rho|< \sin(\alpha)
	\end{cases}
\end{align*}

The green curve hinders the value of $\delta$ from becoming too small when $\rho$ is small. In fact, $\delta^\texttt{green}$ simply replaces $\rho$ by $\sign(\rho) \cdot \sin(\alpha)$.

A formula for an updated residual that lives on the stressed line can be derived by finding an expression for a vector parallel to this line. As a contribution, we propose this formula for $\delta^\texttt{stressed}$ for computing the residual on this line.

\paragraph{Motivation of $\delta^\texttt{stressed}$}
Using the green curve, one must either use a smaller value for $\alpha$ in the formula or one cannot ensure that for $|\rho|\geq \sin(\alpha)$ the new residual is above the stressed line. However, if for the purpose of using the formula of $\delta^\texttt{green}$ one would choose $\alpha$ smaller then in some other regions the residual would be larger than required. This is why we suggest to use $\delta^\texttt{stressed}$ instead of $\delta^\texttt{green}$.

Since the computation of $\delta$ is in $\cO(1)$ whereas the computation of $\hbtau\hp{0},\hbtau\hp{\ell}$ is in $\cO(N \cdot \ell^2)$ it is of no concern that our proposed formula for computing $\delta$ by $\delta^\texttt{stressed}$ is more expensive than $\delta^\texttt{green}$ as proposed in \cite{MaintainingConvergenceBiCGstab}.
\largeparbreak

Alg.\,\ref{Algo:StabCoeffs} shows the overall routine for computing the stability polynomial coefficients $\btau \in \R^\ell$. In the algorithm a matrix $\bS$ is used to compute $\hbtau\hp{0},\hbtau\hp{\ell}$ (as parts of the vectors in line 4), the norms $\kappa_0 \equiv \|\hbr\hp{0}\|$, $\kappa_\ell \equiv \|\hbr\hp{\ell}\|$ (cf. line 5) and the DOT for $\rho$ (cf. line 6). Further, using an $\bS$-norm, the value $\beta$ of the norm of the updated residual can be computed in advance. This is beneficial because no DOT is wasted.

\begin{algorithm}
\begin{algorithmic}[1]
	\Procedure{StabCoeffs}{$\br\hp{0},...,\br\hp{\ell},\beta$}
	\State \textit{// parameter $\alpha \in (0,\pi/2]$ given; typical value: $\alpha = \arctan(0.7)$}
	\State $\bS := [\br\hp{0},...,\br\hp{\ell}]\t \cdot [\br\hp{0},...,\br\hp{\ell}]$\quad\textit{// exploit symmetry and $\bS_{1,1}\equiv \beta^2$}
	\State $\by\hp{0} := \begin{pmatrix} 1\\ -(\bS_{2:\ell,2:\ell})\d \cdot \bS_{2:\ell,1}\\ 0 \end{pmatrix}$, \quad $\by\hp{\ell} := \begin{pmatrix} 0\\ -(\bS_{2:\ell,2:\ell})\d \cdot \bS_{2:\ell,\ell+1}\\ 1 \end{pmatrix}$
	\State $\kappa_0 := \|\by\hp{0}\|_\bS$,\quad $\kappa_\ell := \|\by\hp{\ell}\|_\bS$
	\State $\rho := \frac{(\by\hp{\ell})\t \cdot \bS \cdot \by\hp{0}}{\kappa_0 \cdot \kappa_\ell}$
	\If{$|\rho|\geq\sin(\alpha)$}
		\State $\delta := - \frac{\kappa_0}{\kappa_\ell} \cdot \rho$
	\Else
		\State $\delta := - \frac{\kappa_0}{\kappa_\ell} \cdot \frac{\sin(\alpha)}{\sin\Big(\alpha+\frac{\pi}{4}\cdot\big(\sign(\rho)+1\big)-\arccos(\rho)\Big)}$
	\EndIf	
	\State $\by\hp{0} := \by\hp{0} + \delta \cdot \by\hp{\ell}$
	\State $\btau := \by\hp{0}_{2:(\ell+1)} \in \R^\ell$
	\State $\beta := \|\by\hp{0}\|_\bS$
	\State \Return $\btau,\beta$
	\EndProcedure
\end{algorithmic}
\caption{Computation of the stability polynomial's coefficients}\label{Algo:StabCoeffs}
\end{algorithm}

\subsection{Current practical implementations of \IDR}
In the following we show practical implementations of \IDR. As a motivation of this subsection, we review that there are variants that require less storage, fewer BLAS-1 computations per matrix-vector product and variants that do not decouple the auxiliary vectors' relation
\begin{align*}
	\bV\hp{0} = \bA\cdot \bV\hp{-1}\,.
\end{align*}
The latter variant we call the \textit{decoupling-free} variant of \IDR.

This subsection is organised as follows. First, we introduce a slightly reformulated implementation of \IDR by using our aforementioned code blocks for the initialisation and stability polynomial computation. From this we derive the decoupling-free variant. Afterwards, we re-derive from the decoupling-free \IDR the variant \IDRbiortho, which uses biorthogonal auxiliary vectors, i.e. $\bP\h \cdot \bV\hp{0}$ is lower triangular. Finally, we contribute an extension of \IDRbiortho that ensures in addition that the matrix $\bV\hp{0}$ has condition 1\,.

\largeparbreak
We start by introducting the reformulated version of \IDR. In this version we use the routine \textsc{dirRbio}, which stands for \emph{direct residual biorthogonalisation}. 

\paragraph{Algorithmic sub-block: Direct residual biorthogonalisation}
The direct biorthogonalisation routine is given in Alg.\,\ref{Algo:dirRbio}. Given $\bP$ and $\bV\hp{0}$, the algorithm biorthogonalises the residual with the auxiliary vectors $\bV\hp{0}$ against the columns of $\bP$.

\begin{algorithm}
	\begin{algorithmic}[1]
		\Procedure{dirRbio}{$\bV\hp{-1},\bV\hp{0},\bZ,\bx\hp{0},\br\hp{0},\bEta$}
		\State \textit{$\bP\t\cdot\bV\hp{0}\equiv\bZ$,\quad $\bEta \equiv (\bP\t \cdot \br\hp{0})$}
		\State $\bxi := \bZ\d \cdot \bEta$
		\State $\br\hp{0} := \br\hp{0} - \bV\hp{0} \cdot \bxi$
		\State $\bx\hp{0} := \bx\hp{0} + \bV\hp{-1} \cdot \bxi$
		\State \Return $\bx\hp{0},\br\hp{0},...,\br\hp{k}$
		\EndProcedure
	\end{algorithmic}
	\caption{direct biorthogonalisation of the residual}\label{Algo:dirRbio}
\end{algorithm}

With the help of the direct biorthogonalisation routine we can formulate \IDR in a compacter way in Alg.\,\ref{Algo:IDRs_block}. This algorithm differs in some significant details from the implementation in Alg.\,\ref{Algo:IDRs_naive}. Namely, in line 5, i.e. the beginning of the iterative loop, the residual does not need to be biorthogonalised against the columns of $\bP$. This is because the initialisation has been designed such that the residual is initially already perpendicular to the columns of $\bP$.

As a second difference to the formerly introduced implementation in Alg.\,\ref{Algo:IDRs_naive}, the stability coefficient $\tau \equiv \omega \in \R$ is computed before the for-loop of the biorthogonalisation of the auxiliary vectors is performed, cf. lines 8 and 9. In contrast to that, in Alg.\,\ref{Algo:IDRs_naive} there is first the for-loop for the auxiliary vectors' biorthogonalisation in line 18 and afterwards the computation of $\tau$ in line 26.

Finally, after the polynomial step of the newly introduced implementation in Alg.\,\ref{Algo:IDRs_block}, the updated residual is biorthogonalised using the routine \textsc{dirRbio}. This makes sure that at the beginning of the next repetition of the main-loop the residual already lives in the null-space of $\bP$.

In the following subsection we derive from this implementation of \IDR a variant that uses fewer storage and in which the auxiliary vectors cannot decouple from their pre-images.

\begin{algorithm}
	\begin{algorithmic}[1]
		\Procedure{IDR}{$\bA,\bb,s,\tolabs$}
		\State [$\bP,\bV\hp{-1},\bV\hp{0},\bZ,\bx\hp{0},\br\hp{0},\beta$] = \textsc{Initialisation}($\bA,\bb,s,\tolabs$)
		\State $j:=0$
		\While{$\beta>\tolabs$}
		\State \textit{// $\br\hp{0} \in \cN(\bP)$}
		\State $\br\hp{1} := \bA \cdot \br\hp{0}$\quad\textit{// $\br\hp{1} \in \cG_j$}
		\State $\bEta := \bP\t \cdot \br\hp{1}$
		\State [$\tau,\beta$] = \textsc{StabCoeffs}($\br\hp{0},\br\hp{1},\beta$)
		\For{$q=1,...,s$}
		\State $\bxi := \bZ\d \cdot \bEta$
		\State $\bv_q\hp{-1} := \br\hp{0} - [\bV\hp{0}_{:,1:(q-1)},\bV\hp{-1}_{:,q:s}] \cdot \bxi$
		\State $\bv_q^{(0)\phantom{-}} := \br\hp{1} - [\bV\hp{1}_{:,1:(q-1)},\bV^{(0)\phantom{-}}_{:,q:s}] \cdot \bxi$\quad\textit{// $\bv_q\hp{0} \in \cG_j\cap\cN(\bP)$}
		\State $\bv_q\hp{1} := \bA \cdot \bv_q\hp{0}$\quad\textit{// $\bv_q\hp{1} \in \cG_j$}
		\State $\bZ_{:,q} := \bP\t\cdot\bv\hp{1}_q$
		\EndFor
		\State \textit{// - - - Polynomial step - - - }
		\State \textit{// Residual}
		\State $\bx\hp{0} := \bx\hp{0} + \tau \cdot \br\hp{0}$
		\State $\br\hp{0} := \br\hp{0} - \tau \cdot \br\hp{1}$ \quad\textit{// $\br\hp{0} \in \cG_{j+1}$}
		\State \textit{// Auxiliary vectors}
		\State $\bV\hp{-1} := \bV\hp{-1} - \tau \cdot \bV\hp{0}$
		\State $\bV\hp{0} := \bV\hp{0} - \tau \cdot \bV\hp{1}$ \quad\textit{// $\rg(\bV\hp{0}) \subset \cG_{j+1}$}
		\State $\bZ:=-\tau \cdot \bZ$,\quad$\bEta:=-\tau \cdot \bEta$,\quad$j:=j+1$
		\State \textit{// Biorthgonalize residual, i.e. make $\br\hp{0} \in \cN(\bP)$}
		\State [$\bx\hp{0},\br\hp{0}$] = \textsc{dirRbio}($\bV\hp{-1},\bV\hp{0},\bZ,\bx\hp{0},\br\hp{0},\bEta$)
		\State $\beta := \|\br\hp{0}\|$
		\EndWhile
		\State \Return $\bx\hp{0}$
		\EndProcedure
	\end{algorithmic}
	\caption{\IDR reference implementation}\label{Algo:IDRs_block}
\end{algorithm}

\subsubsection{\IDR without auxiliary decoupling}
So far both of the introduced \IDR implementations require storage for $3 \cdot (s+1) + s$ column vectors. We now introduce a variant that requires only $3 \cdot (s+1)$ column vectors. The idea to this more storage-efficient variant lives in the way how the new auxiliary vectors in $\cG_{j+1}$ can be generated. The following lemma helps our derivation.
\begin{Lemma}[Decoupling-free \IDRO auxiliary vector]\label{Lem:DecFreeIDR}
	Given the Sonneveld spaces $\lbrace\cG_j\rbrace_{j\in\No}$ from $\bA \in \C^{N \times N}$, $\bb \in \C^N$, $\bP \in \C^{N \times s}$, $\lbrace \omega_j \rbrace_{j \in \N} \subset \Cno$.\\	
	Given further an arbitrary $j \in \No$, $\bV\hp{-1},\bV\hp{0} \in \C^{N \times s}$ such that $\bV\hp{0}=\bA\cdot\bV\hp{-1}$, $\rg(\bV\hp{0})\subset\cG_j$ and let $\bZ:=\bP\h\cdot\bV\hp{0}$ have full rank.\\
	Then it follows:
	\begin{align*}
	\bA \cdot \Big(\bI + (-\frac{1}{\omega_{j+1}}\cdot\bV\hp{-1}+\bV\hp{0}) \cdot \bZ\d \cdot \bP\h\Big) \cdot \cG_{j+1} \subset \cG_{j+1}
	\end{align*}
	\underline{Proof:}\\
	Choose an arbitrary $\bx \in \cG_{j+1}$ and define $\hbx$
	\begin{align*}
	\hbx := \bx + (-\frac{1}{\omega_{j+1}}\cdot\bV\hp{-1}+\bV\hp{0}) \cdot \underbrace{\bZ\d \cdot \bP\t \cdot \bx}_{=:\bxi}
	\end{align*}
	Now consider $\bx - \omega_{j+1} \cdot \bA \cdot \hbx$:
	\begin{align*}
	\bx - \omega_{j+1} \cdot \bA \cdot \hbx &= \bx - \omega_{j+1}\cdot\bA\cdot\bx - (\bV\hp{0} - \omega_{j+1}\cdot\bA\cdot\bV\hp{0})\cdot\bxi\\
	&= (\bI - \omega_{j+1}\cdot\bA) \cdot (\bx - \bV\hp{0}\cdot \bxi)
	\end{align*}
	Since $\bxi = (\bP\t\cdot\bV\hp{0})\d\cdot(\bP\t\cdot\bx)$ and $\bZ = \bP\t\cdot\bV\hp{0}$ is regular by requirement, it follows
	\begin{align*}
		\bx - \bV\hp{0} \cdot \bxi \in \cN(\bP)\,.
	\end{align*}
	Due to this it must hold
	\begin{align*}
		\bx - \omega_{j+1} \cdot \bA \cdot \hbx \in \cG_{j+1}\,.
	\end{align*}
	Further, since $\bx \in \cG_{j+1}$, it follows $\bA\cdot\hbx \in \cG_{j+1}$. $\boxtimes$
\end{Lemma}

We now utlise this lemma to derive the decoupling-free \IDR implementation. As a starting point, let us review Alg.\,\ref{Algo:IDRs_block}. Using the above lemma, we can omit the lines 21--23 and replace the lines 10 to 14 by the following:
\begin{algorithmic}
\State $\bw := \br\hp{0} - \tau \cdot \br\hp{1}$
\State $\bv_q\hp{-1} := -\omega_{j+1}\cdot\Big(\bI + (-\frac{1}{\omega_{j+1}}\cdot\bV\hp{-1}+\bV\hp{0}) \cdot \bZ\d \cdot \bP\t\Big) \cdot \bw$
\State $\bv\hp{0}_q := \bA \cdot \bv\hp{-1}_q$
\State $\bZ_{:,q} := \bP\t\cdot\bv_q\hp{0}$
\end{algorithmic}
In the following we explain this in more detail:

When we replace the lines 10--14 by the above code fragment then the code fragment will yield that at line 15 for a respective loop index $q$ it holds $\bv\hp{0}_1,...,\bv\hp{0}_q \in \cG_{j+1}$. This follows from Lem.\,\ref{Lem:DecFreeIDR}: Since $\bw$ is equivalent to the updated residual in line 19, it holds $\bw \in \cG_{j+1}$. In consequence, the vector $\bv\hp{-1}_q$ lives in
\begin{align*}
	\Big(\bI + (-\frac{1}{\omega_{j+1}}\cdot\bV\hp{-1}+\bV\hp{0}) \cdot \bZ\d \cdot \bP\t\Big) \cdot \cG_{j+1}\,.
\end{align*}
The lemma just says that multiplying $\bv\hp{-1}_q$ from the left by $\bA$ yields a vector $\bv\hp{0}_q$ that lives in $\cG_{j+1}$.

Since with the above code fragment new auxiliary vectors are obtained that live already in $\cG_{j+1}$, the polynomial step in lines 21--23 is obsolete. Besides, since $\bw$ is identical to the updated residual in line 19, the computations in lines 18--19 can be moved to earlier in front of line 9.

The above-mentioned variant of \IDR offers two benefits. First, the vectors $\bv\hp{0}_q$ are not modified by AXPYs, thus the accuracy of which the equation $\bv\hp{0}_q = \bA \cdot \bv\hp{-1}_q$ holds is never diminished by rounding errors of further numerical computations. Second, there are no vectors $\bv\hp{1}_q$. As a consequence of this, the method does not need the storage for the $s$ vectors $\bv\hp{1}_1,...,\bv\hp{1}_s$.
\largeparbreak
Alg.\,\ref{Algo:IDRs_noDec} lays out in one implementation all the reformulations described above that lead to the decoupling-free \IDR variant. We see that the computation of the updated residual is moved in front of the inner for-loop (cf. lines 8--9). In line 13 the formula from Lem.\,\ref{Lem:RemainingVectorsInGj} is utilised such that $\bv\hp{0}_q$ lives in $\cG_{j+1}$.

\begin{algorithm}
	\begin{algorithmic}[1]
		\Procedure{IDR}{$\bA,\bb,s,\tolabs$}
		\State [$\bP,\bV\hp{-1},\bV\hp{0},\bZ,\bx\hp{0},\br\hp{0},\beta$] = \textsc{Initialisation}($\bA,\bb,s,\tolabs$)
		\State $j:=0$
		\While{$\beta>\tolabs$}
		\State \textit{// $\br\hp{0} \in \cG_j \cap \cN(\bP)$}
		\State $\br\hp{1} := \bA \cdot \br\hp{0}$\quad\textit{// $\br\hp{1} \in \cG_j$}
		\State [$\tau,\beta$] = \textsc{StabCoeffs}($\br\hp{0},\br\hp{1},\beta$)
		\State $\bx\hp{0} := \bx\hp{0} + \tau \cdot \br\hp{0}$
		\State $\br\hp{0} := \br\hp{0} - \tau \cdot \br\hp{1}$
		\quad\textit{// $\br\hp{0} \in \cG_{j+1}$}
		\State $\bEta := \bP\t \cdot \br\hp{0}$
		\For{$q=1,...,s$}
		\State $\bxi := \bZ\d \cdot \bEta$\quad\textit{// $\bxi \in \R^s$}
		\State $\bv_q\hp{-1} := \bV\hp{-1} \cdot \bxi + \tau \cdot (\br\hp{0} - \bV\hp{0} \cdot \bxi )$\quad\textit{// cf. Lem.\,\ref{Lem:DecFreeIDR}}
		\State $\bv_q^{(0)} := \bA \cdot \bv_q\hp{-1}$\quad\textit{// $\bv_q\hp{0} \in \cG_{j+1}$}
		\State $\bZ_{:,q} := \bP\t\cdot\bv\hp{0}_q$
		\EndFor
		\State $j:=j+1$
		\State [$\bx\hp{0},\br\hp{0}$] = \textsc{dirRbio}($\bV\hp{-1},\bV\hp{0},\bZ,\bx\hp{0},\br\hp{0},\bEta$)
		\State $\beta := \|\br\hp{0}\|$
		\EndWhile
		\State \Return $\bx\hp{0}$
		\EndProcedure
	\end{algorithmic}
	\caption{\IDR without auxiliary decoupling}\label{Algo:IDRs_noDec}
\end{algorithm}

\subsubsection{Iterative biorthogonalisation: \IDRbiortho}
In this sub-subsection we introduce a variant that uses again fewer BLAS-1 operations per repetition than the decoupling-free \IDR. To derive this variant, we use two subroutines; one for biorthogonalising the residual in an iterative way, the other for biorthogonalising the columns of $\bV\hp{0}$ w.r.t. $\bP$ such that $\bZ$ remains lower triangular (as was achieved by our initialisation, cf. Alg.\,\ref{Algo:Initialisation}).

First, we describe the biorthogonalisation subroutine for the residual. This subroutine is called \textsc{itRbio}, which stands for \emph{iterative residual biorthogonalsation}. The name \textit{iterative} means that using this scheme the residual $\br\hp{0}$ is orthogonalised against the columns $\bp_1,\bp_2,...,\bp_s$ subsequently, i.e. one after the other, instead of simultaneously against all of them in one step (as it is done in contrast in \textsc{dirRbio}).

Afterwards we describe the iterative biorthogonalisation of the columns of $\bV\hp{0}$ w.r.t. the columns of $\bP$ so that the lower triangular shape of $\bZ$ is maintained.

\paragraph{Algorithmic sub-block: Iterative residual biorthogonalisation}
The subroutine \textsc{itRbio} consists of three steps: First, the linear coefficient $\xi$ is computed from $\bEta \equiv \bP\t\cdot\br\hp{0}$ and $\bZ_{q,q}\equiv \langle \bp_q,\bv\hp{0}_q\rangle$, cf. line 2. Then, the residual and numerical solution are modified by an AXPY such that the residual becomes perpendicular to $\bp_q$, cf. lines 3--4. Finally, the vector $\bEta$ is modified by an analogous AXPY with the $q$th column of $\bZ$ such that for the biorthogonalised residual it holds again $\bEta = \bP\t\cdot\br\hp{0}$.
\largeparbreak

\begin{algorithm}
	\begin{algorithmic}[1]
		\Procedure{itRbio}{$q,\bV\hp{-1},\bV\hp{0},\bZ,\bx\hp{0},\br\hp{0},\bEta$}
		\State $\xi := \bEta_{q,1} / \bZ_{q,q}$\quad\textit{// $\equiv \langle\bp_q,\br\hp{0}\rangle /\langle \bp_q,\bv_q\hp{0}\rangle$}
		\State $\br\hp{0} := \br\hp{0} - \xi \cdot \bv_q\hp{0}$
		\State $\bx\hp{0} := \bx\hp{0} + \xi \cdot \bv_q\hp{-1}$
		\State $\bEta := \bEta - \bZ_{:,q} \cdot \xi$\quad\textit{// $\bEta \equiv \bP\t\cdot\br\hp{0}$}
		\State \Return $\bx\hp{0},\br\hp{0},\bEta$
		\EndProcedure
	\end{algorithmic}
	\caption{iterative biorthogonalisation of the residual}\label{Algo:itRbio}
\end{algorithm}

\paragraph{Algorithmic sub-block: Iterative auxiliary biorthogonalsation}

This algorithmic sub-block is given in Alg.\,\ref{Algo:itVbio}. It uses a linear combination of the vectors $\bv_1\hp{0},...,\bv\hp{0}_{q-1}$ to orthogonalise the $q$th column of $\bV\hp{0}$ against the columns $\bp_1,...,\bp_{q-1}$. The pre-image $\bv\hp{-1}_q$ is updated accordingly by using the same linear combination for the pre-images $\bv_1\hp{-1},...,\bv\hp{-1}_{q-1}$:
\begin{align*}
\bv\hp{0}_q &:= \bv\hp{0}_q - \sum_{d=1}^{q-1} \bv\hp{0}_d \cdot \xi_d\\
\bv\hp{-1}_q &:= \bv\hp{-1}_q - \sum_{d=1}^{q-1} \bv\hp{-1}_d \cdot \xi_d
\end{align*}
When it comes to the computation of the coefficients $\xi_d$ for $d=1,...,q-1$, one can make use of the fact that $\bv\hp{0}_j \perp \bp_i$ $\forall j>i$ holds. Thus, a Gram-Schmidt procedure can be utilised in which $\bv\hp{0}_q$ is orthogonalised subsequently against $\bp_i$ with $\bv\hp{0}_i$ for $i=1,...,q-1$.

\begin{algorithm}
	\begin{algorithmic}[1]
		\Procedure{itVbio}{$q,\bP,\bV\hp{-1},\bV\hp{0},\bZ$}
		\For{$i=1,...,(q-1)$}
		\State $\xi := \langle \bp_i,\bv_q\hp{0}\rangle / \bZ_{i,i} $\quad\textit{// $\equiv \langle \bp_i,\bv_q\hp{0}\rangle / \langle \bp_i,\bv_i\hp{0}\rangle$}
		\State $\bv\hp{g}_q := \bv\hp{g}_q - \xi \cdot \bv\hp{g}_i$\quad $\forall\ g=-1,0$
		\EndFor
		\State $\bZ_{q:s,q} := \bP_{:,q:s}\t \cdot \bv_q\hp{0}$\quad \textit{$\bZ_{1:(q-1),q}=\bO$}
		\State \Return $\bV\hp{-1},\bV\hp{0},\bZ$
		\EndProcedure
	\end{algorithmic}
	\caption{iterative biorthogonalisation of $\bV\hp{0}$}\label{Algo:itVbio}
\end{algorithm}

\largeparbreak
Using \textsc{itRbio} and \textsc{itVbio} we explain the method \IDRbiortho as presented in Alg.\,\ref{Algo:IDRs_bio}. This algorithm is derived from Alg.\,\ref{Algo:IDRs_noDec}, thus we only describe the differences.

The first difference in \IDRbiortho (Alg.\,\ref{Algo:IDRs_bio}) compared to Alg.\,\ref{Algo:IDRs_noDec} is in lines 12--13. Whereas Alg.\,\ref{Algo:IDRs_noDec} uses all the auxiliary vectors to orthogonalise the residual $\br\hp{0}$ with $\bV\hp{0}$ against $\rg(\bP)$, in \IDRbiortho only the last auxiliary vectors $\bv\hp{0}_q,...,\bv\hp{0}_s$ are used. 

In \IDRbiortho this is possible because the method is built up such that in line 12 for a respective $q$ it holds
\begin{align*}
	\br\hp{0} &\perp \opspan\lbrace\bp_1,...,\bp_{q-1}\rbrace\\
	\bv\hp{0}_q,...,\bv\hp{0}_s & \perp \opspan\lbrace\bp_1,...,\bp_{q-1}\rbrace\,.
\end{align*}
Thus, in lines 12--13 it is possible to only use these latter auxiliary vectors to achieve that $\br\hp{0} - \bV\hp{0}_{:,q:s} \cdot \bxi \in \cN(\bP)$ holds.

A further difference in \IDRbiortho (Alg.\,\ref{Algo:IDRs_bio}) to Alg.\,\ref{Algo:IDRs_noDec} is in lines 15--18: In line 15 the newly computed auxiliary vector $\bv\hp{0}_q$ is modified by the auxiliary vectors $\bv\hp{0}_1,...,\bv\hp{0}_{q-1}$ such that it becomes orthogonal w.r.t. $\bp_1,...,\bp_{q-1}$, cf. the comment in line 16. Afterwards in line 17 the subroutine \textsc{itRbio} is used to orthogonalise $\br\hp{0}$ against $\bp_q$ by using $\bv\hp{0}_q$. Thus, in line 20 it holds that $\br\hp{0}\in \cN(\bP)$.
\largeparbreak

\begin{algorithm}
	\begin{algorithmic}[1]
		\Procedure{IDRsbiortho}{$\bA,\bb,s,\tolabs$}
		\State [$\bP,\bV\hp{-1},\bV\hp{0},\bZ,\bx\hp{0},\br\hp{0},\beta$] = \textsc{Initialisation}($\bA,\bb,s,\tolabs$)
		\State $j:=0$
		\While{$\beta>\tolabs$}
		\State \textit{// $\br\hp{0} \in \cN(\bP)$}
		\State $\br\hp{1} := \bA \cdot \br\hp{0}$\quad\textit{// $\br\hp{1} \in \cG_j$}
		\State [$\tau,\beta$] = \textsc{StabCoeffs}($\br\hp{0},\br\hp{1},\beta$)
		\State $\bx\hp{0} := \bx\hp{0} + \tau \cdot \br\hp{0}$
		\State $\br\hp{0} := \br\hp{0} - \tau \cdot \br\hp{1}$
		\quad\textit{// $\br\hp{0} \in \cG_{j+1}$}
		\State $\bEta := \bP\t \cdot \br\hp{0}$
		\For{$q=1,...,s$}
		\State $\bxi := (\bZ_{q:s,q:s})\d \cdot \bEta_{q:s,1}$\quad\textit{// $\bxi \in \R^{s+1-q}$}
		\State $\bv_q\hp{-1} := \bV_{:,q:s}\hp{-1} \cdot \bxi + \tau \cdot (\br\hp{0} - \bV_{:,q:s}\hp{0} \cdot \bxi )$
		\State $\bv_q^{(0)} := \bA \cdot \bv_q\hp{-1}$\quad\textit{// $\bv_q\hp{0} \in \cG_{j+1}$}
		\State [$\bV\hp{-1},\bV\hp{0},\bZ$] = \textsc{itVbio}($q,\bP,\bV\hp{-1},\bV\hp{0},\bZ,$)
		\State \textit{// $\bv_q\hp{0} \perp \lbrace \bp_1,...,\bp_{q-1}\rbrace$,\quad$\bZ\equiv \bP\t \cdot \bV\hp{0}$}
		\State [$\bx\hp{0},\br\hp{0},\bEta$] = \textsc{itRbio}($q,\bV\hp{-1},\bV\hp{0},\bZ,\bx\hp{0},\br\hp{0},\bEta$)
		\State \textit{// $\br\hp{0} \perp \lbrace \bp_1,...,\bp_{q}\rbrace$,\quad$\bEta \equiv \bP\t \cdot \br\hp{0}$}
		\EndFor
		\State $j:=j+1$,\quad$\beta := \|\br\hp{0}\|$
		\EndWhile
		\State \Return $\bx\hp{0}$
		\EndProcedure
	\end{algorithmic}
	\caption{\IDRbiortho}\label{Algo:IDRs_bio}
\end{algorithm}

The authors in \cite{IDR-biorth} mention that the implementation \IDRbiortho has two benefits:
\begin{enumerate}[(1)]
	\item The columns of $\bV\hp{0}$ are biorthogonal to the columns of $\bP$ in the sense that the matrix $\bZ$ always remains lower triangular. The authors show experiments that indicate that this makes the method numerically more robust when using larger values for $s$. However, this is not obvious since $\bV\hp{0}$ can be still badly conditioned. In fact, there is no proof that justifies why \IDRbiortho should be more robust.
	\item In some papers \cite{IDR-biorth,IDRstab-biorth} it is said that \IDRbiortho is cheaper than Alg.\,\ref{Algo:IDRs_noDec}. However, we cannot verify this. We count $2\cdot (s+1)$ AXPYs in every inner for-loop of \IDRbiortho whereas the decoupling-free variant has $2\cdot s$ AXPYs in every inner for-loop and $2 \cdot s$ additional AXPYs after the for-loop for the direct biorthogonalisation. The number of computed DOTs is $s$ per for-loop repetition for both methods. So in total the computational cost in terms of BLAS-1 operations should not have changed.
\end{enumerate}

\paragraph{Excursus: \IDRobio}
As mentioned above, the authors in \cite{IDR-biorth} observe that \IDRbiortho is numerically more robust for large values of $s$ (they try for instance $s=140$). However, the biorthogonalisation does not achieve that $\cond(\bV\hp{0})$ is bounded. In fact, the contrary holds: Columns of $\bV\hp{0}$ could still become collinear but \IDRbiortho would not terminate with an accurate solution but simply $\bZ$ would become singular and spoil the biorthogonalisations in lines 15 and 17.
\largeparbreak

To overcome the stability issues but keep the biorthogonality relation, we contribute a novel \IDR variant, the so-called \textit{ortho-biortho} variant, or briefer \textit{obio}. This variant is obtained by exchanging the subroutine \textsc{itVbio} by the subroutine \textsc{itVobio}, cf. Alg.\,\ref{Algo:itVobio}.

The subroutine \textsc{itVobio} works as follows: First, in Alg.\,\ref{Algo:itVobio} in line 2 a modified Gram-Schmidt process (named \textsc{itVorth}, given in Alg.\,\ref{Algo:itVorth}) is used to orthonormalise $\bv\hp{0}_q$ against $\bv\hp{0}_1,...,\bv\hp{0}_{q-1}$. After the orthonormalisation we need to biorthogonalise the new column $\bv\hp{0}_q$ against $\bp_1,...,\bp_{q-1}$. However, we cannot simply biorthogonalise it by a modified Gram-Schmidt procedure with $\bv\hp{0}_1,...,\bv\hp{0}_{q-1}$ since this would destroy the formerly achieved orthonormality. In the following we explain how Givens rotations from the right can be utilised to perform the biorthogonalisation without changing $\cond(\bV\hp{0})$:

Consider the following equation for $\bZ$, where a LQ-decomposition for $\bZ$ is computed by Givens rotations:
\begin{align}
	\bP\h\cdot\bV\hp{0} \cdot \bQ_Z = \bZ \cdot \bQ_Z = \bL_Z\label{eqn:Zrot}
\end{align}
Each Givens rotation of which $\bQ_Z$ consists can be applied subsequently onto $\bV\hp{0}$. 

\begin{figure}
\centering
\includegraphics[width=0.3\linewidth]{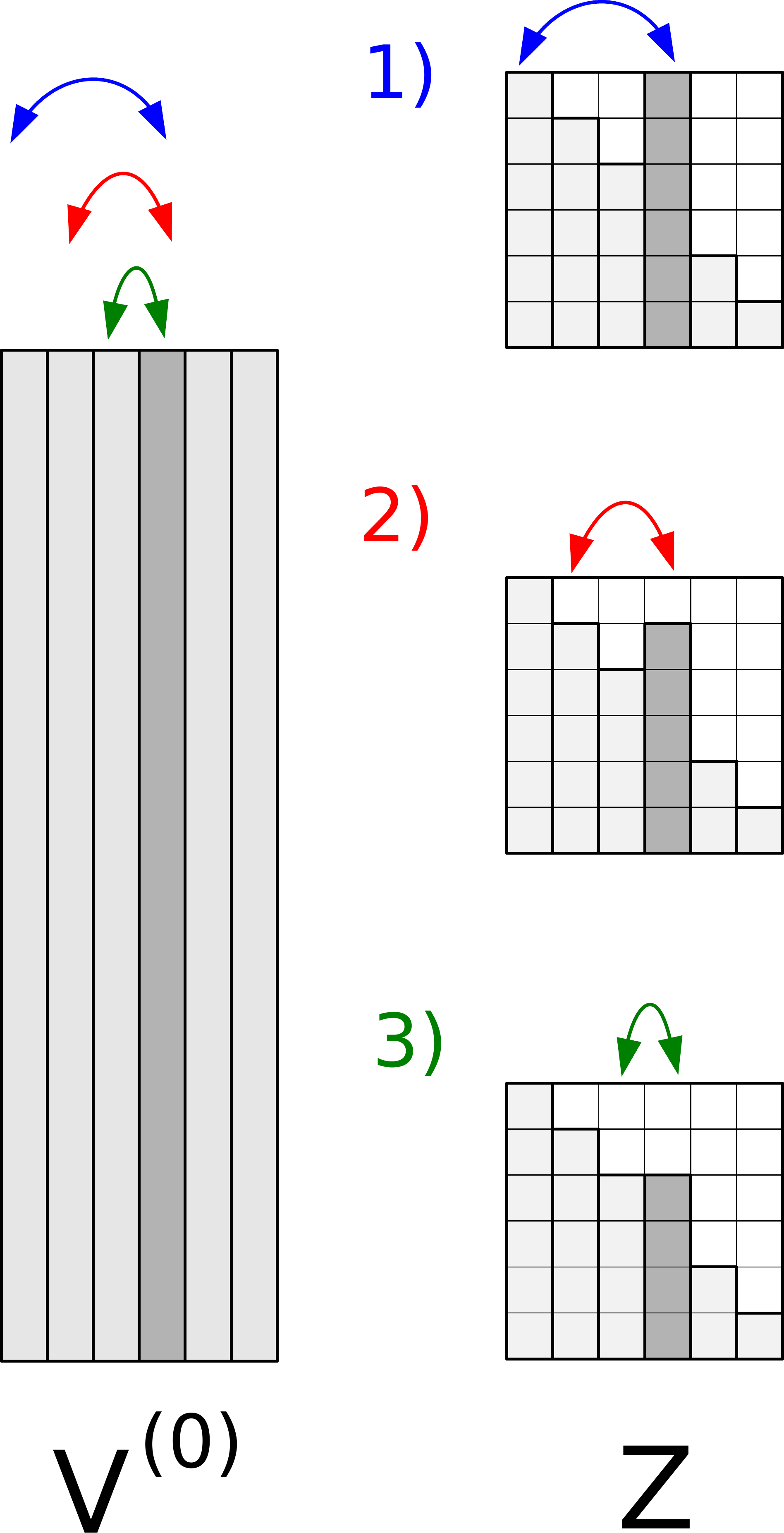}
\caption{Algorithmic principle of the orthogonality-preserving biorthogonalsation.}
\label{fig:itVobio}
\end{figure}

Fig.\,\ref{fig:itVobio} illustrates the rotation approach for the exemplary numbers $q=4$ and $s=6$: Initially, in Alg.\,\ref{Algo:itVobio} after line 3 the fourth column of $\bZ$ can be full (matrix in the upper right part of the figure). By using Givens rotations on the columns of $\bZ$, one can eliminate successively the fist $q-1$ non-zeros of $\bZ_{:,q}$. The figure illustrates with blue, red and green arrows how the rotations act on the matrices $\bZ$ and $\bV\hp{0}$. The figure shows the non-zero pattern of $\bZ$ in advance of each column rotation. The successive biorthogonalisation by Givens rotations is realised by the for-loop from line 4.
\largeparbreak

\begin{algorithm}
	\begin{algorithmic}[1]
		\Procedure{itVorth}{$q,\bV\hp{-1},...,\bV\hp{k}$}
		\For{$i=1,...,(q-1)$}
		\State $\xi := \langle \bv_i\hp{0},\bv_q\hp{0}\rangle$\quad\textit{// $\equiv \langle \bv_i\hp{0},\bv_q\hp{0}\rangle / \langle \bv_i\hp{0},\bv_i\hp{0}\rangle$}
		\State $\bv\hp{g}_q := \bv\hp{g}_q - \xi \cdot \bv\hp{g}_i$\quad $\forall\ g=-1,0,...,k$
		\EndFor
		\State $\xi := 1/\|\bv\hp{0}_q\|$
		\State $\bv_q\hp{g} := \xi \cdot \bv_q\hp{g}$\quad $\forall\ g=-1,0,...,k$
		\State \Return $\bV\hp{-1},...,\bV\hp{k}$
		\EndProcedure
	\end{algorithmic}
	\caption{iterative orthogonalisation of $\bV\hp{0}$}\label{Algo:itVorth}
\end{algorithm}

\begin{algorithm}
	\begin{algorithmic}[1]
		\Procedure{itVobio}{$q,\bP,\bV\hp{-1},\bV\hp{0},\bZ$}
		\State [$\bV\hp{-1},\bV\hp{0}$] = \textsc{itVorth}($q,\bV\hp{-1},\bV\hp{0}$)
		\State $\bZ_{:,q} := \bP\t \cdot \bv_q\hp{0}$
		\For{$i=1,...,(q-1)$}
		\State [$\bQ,\bZ_{[i,q],[i,q]}$] = \texttt{lq}($\bZ_{[i,q],[i,q]}$)
		\State $\bV\hp{g}_{:,[i,q]} := \bV\hp{g}_{:,[i,q]} \cdot \bQ$\quad $\forall\ g=-1,0$
		\EndFor		
		\State \Return $\bV\hp{-1},\bV\hp{0},\bZ$
		\EndProcedure
	\end{algorithmic}
	\caption{iterative orthogonalisation of $\bV\hp{0}$ and biorthogonalisation of $\bV\hp{0}$}\label{Algo:itVobio}
\end{algorithm}

In the remainder we explain how the ortho-biorthogonalsation affects the implementation in Alg.\,\ref{Algo:IDRs_bio}. Replacing \textsc{itVbio} by \textsc{itVobio}, the column vectors of $\bV\hp{0}$ remain orthonormal throughout the whole computation, whereas the matrix $\bZ$ remains again lower triangular. This variant of \IDR we call \IDRobio.

As a final remark, the perfect conditioning of $\bP$ and $\bV\hp{0}$ guarantees that $\bZ$ is not badly conditioned due to some (nearly) collinear columns of $\bV\hp{0}$. However, it does not guarantee that the condition of $\bZ$ is good in general, cf. our above discussion from Sec.\,\ref{sec:ThePolynomialCoefficients} at the end of the paragraph \enquote{The geometric approach}. To our best knowledge, the ortho-biortho variant of \IDR is a novel contribution of this thesis.

\subsection{Current practical implementations of \IDRstab}
In this subsection we review implementation variants of \IDRstab. Except that there is no decoupling-free variant, this subsection describes similar variants to the above subsection, namely a biortho variant and an ortho-biortho variant.
\largeparbreak

This subsection is organised as follows. First, we present a restructured and more comprehensive implementation of \IDRstab in Alg.\,\ref{Algo:IDRstab_block}. Afterwards, we discuss how principles from \IDRbiortho and \IDRobio can be incorporated into Alg.\,\ref{Algo:IDRstab_block}.

\subsubsection{Reference algorithm for \IDRstab}
The below described implementation of \IDRstab in Alg.\,\ref{Algo:IDRstab_block} serves as a starting point for subsequently discussed sophistications. The algorithm uses a small generalisation of the direct residual biorthogonalisation that we discuss in advance of describing Alg.\,\ref{Algo:IDRstab_block}.

\paragraph{Algorithmic sub-block: Direct residual biorthogonalisation}
For a given value of $k$, in this subroutine the $(k+1)$st power of the residual is orthogonalised w.r.t. the columns of $\bP$ by using a linear combination of the columns of $\bV\hp{k+1}$. The residual (powers) $\br\hp{0},...,\br\hp{k}$ and the numerical solution $\bx\hp{0}$ are updated in a consistent way, using the same linear combinations with vectors from $\bV\hp{0},...,\bV\hp{k}$ and $\bV\hp{-1}$.

\begin{algorithm}
	\begin{algorithmic}[1]
		\Procedure{dirRbio}{$\bV\hp{-1},...,\bV\hp{k+1},\bZ,\bx\hp{0},\br\hp{0},...,\br\hp{k+1},\bEta$}
		\State $\bxi := \bZ\d \cdot \bEta$ \quad\textit{// $\bEta\equiv\bP\t \cdot \br\hp{k+1}$}
		\State $\br\hp{g} := \br\hp{g} - \bV\hp{g} \cdot \bxi$\quad $\forall\ g=0,...,k+1$
		\State $\bx\hp{0} := \bx\hp{0} + \bV\hp{-1} \cdot \bxi$
		\State $\bEta := \bO$\quad \textit{// $\in \R^{s \times 1}$}
		\State \Return $\bx\hp{0},\br\hp{0},...,\br\hp{k+1},\bEta$
		\EndProcedure
	\end{algorithmic}
	\caption{direct biorthogonalisation of the residual}\label{Algo:dirRbioL}
\end{algorithm}

\largeparbreak
In the following we briefly point out the main differences between Alg.\,\ref{Algo:IDRstab_block} and Alg.\,\ref{Algo:IDRstab_naive}.

In contrast to Alg.\,\ref{Algo:IDRstab_naive} lines 15--17, in Alg.\,\ref{Algo:IDRstab_block} line 8 the $k$th power of the residual is already orthogonal w.r.t. $\bP$. Thus there is no initial biorthogonalisation of the residual $\br\hp{k}$. The for-loop for the biorthogonalisation of the auxiliary vectors however is identical in both implementations. Afterwards, in contrast to Alg.\,\ref{Algo:IDRstab_naive}, in Alg.\,\ref{Algo:IDRstab_block} in line 19 the $(k+1)$st power of the residual is biorthogonalised.

As a result, after the biorthogonalisation step, i.e. line 21 in Alg.\,\ref{Algo:IDRstab_block}, it holds in addition to the properties from Alg.\,\ref{Algo:IDRstab_naive} line 29 the property:
\begin{align*}
	\br\hp{\ell} \in \cN(\bP)
\end{align*}
This has the effect that the subsequently updated residual in the polynomial step (cf. Alg.\,\ref{Algo:IDRstab_block} line 26) lives not only in $\cG_{j+\ell}$ but also in $\cN(\bP)$.

\begin{algorithm}
	\begin{algorithmic}[1]
		\Procedure{IDRstab}{$\bA,\bb,s,\ell,\tolabs$}
		\State [$\bP,\bV\hp{-1},\bV\hp{0},\bZ,\bx\hp{0},\br\hp{0},\beta$] = \textsc{Initialisation}($\bA,\bb,s,\tolabs$)
		\State $j:=0$
		\While{$\beta>\tolabs$}
		\State \textit{// $\br\hp{0},\bv_1\hp{0},...,\bv_s\hp{0}\in\cG_j$, $\br\hp{0} \in \cN(\bP)$}
		\State \textit{// - - - Biorthogonalisation - - - }
		\For{$k=0,...,\ell-1$}
		\State \textit{// $\br\hp{g} \in \cG_j\cap\cN(\bP)$ $\forall\ g=0,...,k$}
		\State $\br\hp{k+1} := \bA \cdot \br\hp{k}$\quad\textit{// $\br\hp{k+1} \in \cG_j$}
		\State $\bEta := \bP\t \cdot \br\hp{k+1}$
		\State \textit{// Auxiliary vectors}
		\For{$q=1,...,s$}
		\State $\bxi := \bZ\d \cdot \bEta \in \R^s$
		\State $\bv_q\hp{g} := \br\hp{g+1} - [\bV\hp{g+1}_{:,1:(q-1)},\bV\hp{g}_{:,q:s}] \cdot \bxi$\quad$\forall\ g=-1,0,...,k$
		\State $\bv_q\hp{k+1} := \bA \cdot \bv_q\hp{k}$
		\State $\bZ_{:,q} := \bP\t\cdot\bv_q\hp{k+1}$
		\EndFor
		\State \textit{// Residual}
		\State [$\bx\hp{0},\br\hp{0},...,\br\hp{k+1},\bEta$] = ...\quad\quad\quad\quad\quad\quad\quad\quad\quad\quad\quad\quad\quad\quad\quad\quad\quad\quad\quad\quad\quad\linebreak \mbox{\quad\quad\quad\quad\quad\quad\textsc{dirRbio}($\bP,\bV\hp{-1},...,\bV\hp{k+1},\bZ,\bx\hp{0},\br\hp{0},...,\br\hp{k+1},\bEta$)}
		\EndFor
		\State \textit{// $\br\hp{0},\bv_1\hp{0},...,\bv_s\hp{0} \in \cG_j \cap \cK_\ell^\perp(\bA\t;\bP)$, $\br\hp{\ell} \in \cN(\bP)$}
		\State $\beta := \|\br\hp{0}\|$
		\State \textit{// - - - Polynomial step - - - }
		\State [$\btau,\beta$] := \textsc{StabCoeffs}($\br\hp{0},...,\br\hp{\ell},\beta$)
		\State $ {\bx\hp{0}} := \bx\hp{0} + \sum_{k=1}^\ell \br\hp{k-1} \cdot \tau_k$,\quad$\bV\hp{-1} := \bV\hp{-1} - \sum_{k=1}^\ell \bV\hp{k-1} \cdot \tau_k$
		\State $\mathrlap{\br\hp{0}}\phantom{\bx\hp{0}} := \mathrlap{\br\hp{0}}\phantom{\bx\hp{0}} - \sum_{k=1}^\ell \mathrlap{\br\hp{k}}\phantom{\br\hp{k-1}} \cdot \tau_k$,\quad$\mathrlap{\bV\hp{0}}\phantom{\bV\hp{-1}} := \mathrlap{\bV\hp{0}}\phantom{\bV\hp{-1}} - \sum_{k=1}^\ell \mathrlap{\bV\hp{k}}\phantom{\bV\hp{k-1}} \cdot \tau_k$
		\State $\bZ:=-\tau_\ell \cdot \bZ$,\quad$j:=j+\ell$
		\EndWhile
		\State \Return $\bx\hp{0}$
		\EndProcedure
	\end{algorithmic}
	\caption{\IDRstab reference implementation}\label{Algo:IDRstab_block}
\end{algorithm}

In the following sub-subsections we extend the just described \IDRstab variant to incorporate the ideas from \IDRbiortho and \IDRobio.

\paragraph{Attendum: Treatment of $\bZ$} Unfortunately, the computation of $\bZ$ in line 27 leads to an early loss of the rate of convergence. This turns out in the way that the average rate of convergence of the iteratively updated residual becomes a horizontal line at a relative residual of about $10^{-6}$ to $10^{-8}$. $\bZ$ must be recomputed to prevent this, i.e. line 27 must be replaced by $\bZ := \bP\t \cdot \bV\hp{0}$\,. This introduces additional computational costs. These costs could be only prevented when doubling the amount of required storage, cf. the implementations of \IDRstabno in \cite{IDRstab-Paper,GBiCGstab,IDRstab-biorth}.

\subsubsection{The biortho variant of \IDRstab}
To derive \IDRbiortho from the decoupling-free \IDR variant, we simply exchanged the direct residual-biorthogonalisation by an iterative residual-biorthogonalsation. The same is done to derive \IDRstabbiortho in Alg.\,\ref{Algo:IDRstab_biortho} from the above Alg.\,\ref{Algo:IDRstab_block}: Taking Alg.\,\ref{Algo:IDRstab_block} as a reference, we removed line 19 and replaced line 16 by the following code fragment:
\begin{algorithmic}
\State [$\bV\hp{-1},...\bV\hp{k+1},\bZ$] = ...\quad\quad\quad\quad\quad\quad\quad\quad\quad\quad\quad\quad\quad\quad\quad\quad\quad\quad\quad\quad\quad\linebreak
\mbox{\quad\quad\quad\quad\quad\quad\quad\quad\textsc{itVbio}($q,${}$\bP,\bV\hp{-1},...,\bV\hp{k+1},\bZ$)}
\State [$\bx\hp{0},\br\hp{0},...,\br\hp{k+1},\bEta$] = ...\quad\quad\quad\quad\quad\quad\quad\quad\quad\quad\quad\quad\quad\quad\quad\quad\quad\quad\quad\quad\quad\linebreak \mbox{\quad\quad\quad\quad\quad\quad\quad\quad\rlap{\textsc{itRbio}(}\phantom{\textsc{itVbio}($q,\bP,$}$\bV\hp{-1},...,\bV\hp{k+1},\bZ,\bx\hp{0},\br\hp{0},...,\br\hp{k+1},\bEta$)}
\end{algorithmic}

In the following we explain these generalisations for the iterative residual and auxiliary vector biorthogonalsation.

\paragraph{Algorithmic sub-block: Iterative residual biorthogonalsation}
For \IDR we used an iterative biorthogonalsation for the residual $\br\hp{0}$ and the numerical solution $\bx\hp{0}$. During the $q$th interior for-loop the residual was updated using the column $\bv\hp{0}_q$ to orthogonalise $\br\hp{0}$ against $\bp_q$, cf. Alg.\,\ref{Algo:itRbio}.

For \IDRstab the biorthogonalsation part of the whole algorithm does not only have one interior for-loop with index $q$ but an additional for-loop around it with loop-index $k=0,...,\ell-1$\,. As a generalisation of the above biorthogonalisation it turns out in the following paragraphs that in fact a method is needed that biorthogonalises for given loop-indices $q,k$ the residual power vector $\br\hp{k+1}$ against $\bp_q$ by using $\bv_q\hp{k+1}$. Further, the residual (powers) $\br\hp{0},...,\br\hp{k}$ and the numerical solution $\bx\hp{0}$ must be updated consistently, using the vectors $\bv\hp{0}_q,...,\bv\hp{k}_q$ and $\bv\hp{-1}_q$. An appropriate algorithm for doing this is given in Alg.\,\ref{Algo:itRbioL}.

\begin{algorithm}
	\begin{algorithmic}[1]
		\Procedure{itRbio}{$q,\bV\hp{-1},...,\bV\hp{f},\bZ,\bx\hp{0},\br\hp{0},...,\br\hp{f},\bEta$}
		\State $\xi := \bEta_{q,1} / \bZ_{q,q}$\quad\textit{// $\equiv \langle\bp_q,\br\hp{f}\rangle /\langle \bp_q,\bv_q\hp{f}\rangle$}
		\State $\br\hp{g} := \br\hp{g} - \xi \cdot \bv_q\hp{g}$\quad $\forall\ g=0,...,f$
		\State $\bx\hp{0} := \bx\hp{0} + \xi \cdot \bv_q\hp{-1}$
		\State $\bEta := \bEta - \bZ_{:,q} \cdot \xi$\quad\textit{// $\bEta \equiv \bP\t\cdot\br\hp{f}$}
		\State \Return $\bx\hp{0},\br\hp{0},...,\br\hp{f},\bEta$
		\EndProcedure
	\end{algorithmic}
	\caption{iterative biorthogonalisation of the residual}\label{Algo:itRbioL}
\end{algorithm}

\paragraph{Algorithmic sub-block: Iterative auxiliary biorthogonalisation}
This subroutine is given in Alg.\,\ref{Algo:itVbioL}. For given values of $k,q$ it modifies the $q$th column of $\bV\hp{k+1}$ such that it is biorthogonal w.r.t. $\bp_1,...,\bp_{q-1}$. This is done by using a linear combination of the columns $\bv_1\hp{k+1},...,\bv\hp{k+1}_{q-1}$. Since $\bv_j\hp{k+1} \perp \bp_i$ $\forall j>i$ holds one can use again a Gram-Schmidt method to orthogonalise $\bv\hp{k+1}_q$ successively against each column $\bp_i$.

\begin{algorithm}
	\begin{algorithmic}[1]
		\Procedure{itVbio}{$q,\bP,\bV\hp{-1},...,\bV\hp{k+1},\bZ$}
		\For{$i=1,...,(q-1)$}
		\State $\xi := \langle \bp_i,\bv_q\hp{k+1}\rangle / \bZ_{i,i} $\quad\textit{// $\equiv \langle \bp_i,\bv_q\hp{k+1}\rangle / \langle \bp_i,\bv_i\hp{k+1}\rangle$}
		\State $\bv\hp{g}_q := \bv\hp{g}_q - \xi \cdot \bv\hp{g}_i$\quad $\forall\ g=-1,0,...,k+1$
		\EndFor
		\State $\bZ_{q:s,q} := \bP_{:,q:s}\t \cdot \bv_q\hp{k+1}$
		\State \Return $\bV\hp{-1},...,\bV\hp{k},\bZ$
		\EndProcedure
	\end{algorithmic}
	\caption{iterative biorthogonalisation of $\bv_q\hp{k+1}$}\label{Algo:itVbioL}
\end{algorithm}

\begin{algorithm}
	\begin{algorithmic}[1]
		\Procedure{IDRstabLbiortho}{$\bA,\bb,s,\ell,\tolabs$}
		\State [$\bP,\bV\hp{-1},\bV\hp{0},\bZ,\bx\hp{0},\br\hp{0},\beta$] = \textsc{Initialisation}($\bA,\bb,s,\tolabs$)
		\State $j:=0$
		\While{$\beta>\tolabs$}
		\State \textit{// $\br\hp{0},\bv_1\hp{0},...,\bv_s\hp{0}\in\cG_j$}
		\State \textit{// - - - Biorthogonalisation - - - }
		\For{$k=0,...,\ell-1$}
		\State \textit{// $\br\hp{g} \in \cG_j\cap\cN(\bP)$ $\forall\ g=0,...,k$}
		\State $\br\hp{k+1} := \bA \cdot \br\hp{k}$\quad\textit{// $\br\hp{k+1} \in \cG_j$}
		\State $\bEta := \bP\t \cdot \br\hp{k+1}$
		\State \textit{// Auxiliary vectors}
		\For{$q=1,...,s$}
		\State $\bxi := (\bZ_{q:s,q:s})\d \cdot \bEta_{q:s,1} \in \R^{s+1-q}$
		\State $\bv_q\hp{g} := \br\hp{g+1} - \bV\hp{g}_{:,q:s} \cdot \bxi$\quad$\forall\ g=-1,0,...,k$
		\State $\bv_q\hp{k+1} := \bA \cdot \bv_q\hp{k}$
		\State \mbox{[$\bV\hp{-1},\bV\hp{0},...\bV\hp{k+1},\bZ$] = ...}\linebreak
		\mbox{\quad\quad\quad\quad\quad\quad\quad\quad\textsc{itVbio}($q,${}$\bP,\bV\hp{-1},...,\bV\hp{k+1},\bZ$)}
		\State \textit{// $\bv_q\hp{k+1} \perp \lbrace \bp_1,...,\bp_{q-1}\rbrace\,,\quad\bZ\equiv\bP\t\cdot\bV\hp{k+1}$}
		\State [$\bx\hp{0},\br\hp{0},...,\br\hp{k+1},\bEta$] = ...\quad\quad\quad\quad\quad\quad\quad\quad\quad\quad\quad\quad\quad\quad\quad\quad\quad\quad\quad\quad\quad\linebreak \mbox{\quad\quad\quad\quad\quad\quad\quad\quad\rlap{\textsc{itRbio}(}\phantom{\textsc{itVbio}($q,\bP,$}$\bV\hp{-1},...,\bV\hp{k+1},\bZ,\bx\hp{0},\br\hp{0},...,\br\hp{k+1},\bEta$)}
		\State \textit{// $\br\hp{k+1} \perp \lbrace \bp_1,...,\bp_{q}\rbrace\,,\quad\bEta\equiv\bP\t\cdot\br\hp{k+1}$}
		\EndFor
		\EndFor
		\State \textit{// $\br\hp{0},\bv_1\hp{0},...,\bv_s\hp{0} \in \cG_j \cap \cK_\ell^\perp(\bA\t;\bP)$, $\br\hp{\ell}\in \cN(\bP)$}
		\State $\beta := \|\br\hp{0}\|$
		\State \textit{// - - - Polynomial step - - - }
		\State [$\btau,\beta$] := \textsc{StabCoeffs}($\br\hp{0},...,\br\hp{\ell},\beta$)
		\State $ {\bx\hp{0}} := \bx\hp{0} + \sum_{k=1}^\ell \br\hp{k-1} \cdot \tau_k$,\quad$\bV\hp{-1} := \bV\hp{-1} - \sum_{k=1}^\ell \bV\hp{k-1} \cdot \tau_k$
		\State $\mathrlap{\br\hp{0}}\phantom{\bx\hp{0}} := \mathrlap{\br\hp{0}}\phantom{\bx\hp{0}} - \sum_{k=1}^\ell \mathrlap{\br\hp{k}}\phantom{\br\hp{k-1}} \cdot \tau_k$,\quad$\mathrlap{\bV\hp{0}}\phantom{\bV\hp{-1}} := \mathrlap{\bV\hp{0}}\phantom{\bV\hp{-1}} - \sum_{k=1}^\ell \mathrlap{\bV\hp{k}}\phantom{\bV\hp{k-1}} \cdot \tau_k$
		\State $\bZ:=-\tau_\ell \cdot \bZ$,\quad$j:=j+\ell$
		\EndWhile
		\State \Return $\bx\hp{0}$
		\EndProcedure
	\end{algorithmic}
	\caption{\IDRstabbiortho}\label{Algo:IDRstab_biortho}
\end{algorithm}

\paragraph{Attendum: Treatment of $\bZ$} Again, the matrix $\bZ$ must be recomputed to maintain the convergence. However, this time it is not sufficient to only recompute this matrix since it might be no longer lower triangular. Instead, the whole iterative biorthogonalisation must be recomputed. Therefor, line 28 in Alg.\,\ref{Algo:IDRstab_biortho} must be replaced by the following lines:
\begin{algorithmic}[1]
\For {$q=1,...,s$}
	\State [$\bV\hp{-1},\bV\hp{0},\bZ$] = \textsc{itVbio}($q,\bP,\bV\hp{-1},\bV\hp{0},\bZ$)
\EndFor
\end{algorithmic}

\subsubsection{Excursus: The ortho-biortho variant of \IDRstab}
We obtain an ortho-biortho variant from \IDRstabbiortho (Alg.\,\ref{Algo:IDRstab_biortho}) by exchanging the subroutine \textsc{itVbio} by the subroutine \textsc{itVobio}. The latter subroutine does not only biorthogonalise the columns of $\bV\hp{k}$ but also orthonormalises the columns of $\bV\hp{0}$.

The generalisation of \textsc{itVobio} for $\ell>1$ is given in Alg.\,\ref{Algo:itVobioL}. As is done in \textsc{itVobio} for $\ell=1$, the method in Alg.\,\ref{Algo:itVobioL} first uses \textsc{itVorth} to orthonormalise $\bV\hp{0}_{:,1:q}$ by a Gram-Schmidt procedure whilst keeping its powers and pre-image consistent. Afterwards, rotations are applied from the right in an analogous way as discussed along \eqref{eqn:Zrot}, where this time it is $\bZ_{:,1:q} = \bP\t \cdot \bV_{:,1:q}\hp{k+1}$.

\begin{algorithm}
	\begin{algorithmic}[1]
		\Procedure{itVobio}{$q,\bP,\bV\hp{-1},...,\bV\hp{k+1},\bZ$}
		\State [$\bV\hp{-1},...,\bV\hp{k+1}$] = \textsc{itVorth}($q,\bV\hp{-1},...,\bV\hp{k+1}$)
		\State $\bZ_{:,q} := \bP\t \cdot \bv_q\hp{k+1}$
		\For{$i=1,...,(q-1)$}
		\State [$\bQ,\bZ_{[i,q],[i,q]}$] = \texttt{lq}($\bZ_{[i,q],[i,q]}$)
		\State $\bV\hp{g}_{:,[i,q]} := \bV\hp{g}_{:,[i,q]} \cdot \bQ$\quad $\forall\ g=-1,0,...,k+1$
		\EndFor		
		\State \Return $\bV\hp{-1},...,\bV\hp{k+1},\bZ$
		\EndProcedure
	\end{algorithmic}
	\caption{iterative orthogonalisation of $\bv_q\hp{0}$ and biorthogonalisation of $\bv_q\hp{k+1}$}\label{Algo:itVobioL}
\end{algorithm}

\paragraph{Attendum: Treatment of $\bZ$} In analogy to the biortho variant, in the obio variant the line 28 in Alg.\,\ref{Algo:IDRstab_biortho} must be replaced by the following lines:
\begin{algorithmic}[1]
	\For {$q=1,...,s$}
	\State [$\bV\hp{-1},\bV\hp{0},\bZ$] = \textsc{itVobio}($q,\bP,\bV\hp{-1},\bV\hp{0},\bZ$)
	\EndFor
\end{algorithmic}
\largeparbreak
Unfortunately, the ortho-biortho variant is not very useful for the following reason. Even when $\bV\hp{0}$ is constructed such that it is well-conditioned, the condition of $\bV\hp{\ell}$ can still reach the condition of $\bA^\ell$. E.g., for $\ell=2$ one cannot make sure that neither of the matrices $\bV\hp{-1},\bV\hp{2}$ has a condition of $\cond(\bA)^2$, which is usually already a k.o. criterion for the following reason.

For example, it seems reasonable to try to solve systems with a condition of $\cond(\bA)=10^6$ and $\epsmachine=10^{-16}$. The relative residual must then be $<10^{-6}$ according to \eqref{eqn:RelErrEst}. However, if squares of $\cond(\bA)$ occurred then one had only 4 effective digits in some basis matrices, i.e. the required relative residual tolerance is not achievable.

Again, to our best knowledge the obio variant of \IDRstabno is a novel contribution of this thesis.

\subsection{Weaknesses of all the current implementations}

In the former subsections we have frequently pointed out some numerical issues and breakdowns that the methods can encounter. These issues we summarise under the term \textit{weaknesses}. Some weaknesses were related to potential singularities or bad conditioning in $\bZ$. In this section we will motivate why it makes sense to distinguish weaknesses that stem from the underlying geometric principle of the non-symmetric Lanczos process and those that are only related to the particular implementation of the respective \IDRO method.

\subsubsection{Motivation}
Weaknesses in the algorithm can cause a theoretical breakdown. Such breakdowns do not occur in practice. Instead, as a result in practice of these weaknesses, often a loss of superlinear convergence or a delay in the transition to superlinear convergence or an early bound for the achievable final accuracy is observed.

The causes of weaknesses must be distinguished into causes of the underlying mathematical principle, namely \IDRO, and causes of the particular implementation being used. In the following we give reasons why making this distinction is reasonable.
\largeparbreak

On the one hand there are limits of the mathematical concept. That is, the non-symmetric Lanczos procedure can have breakdowns and so can have the recursion of the Sonneveld spaces in the \IDRO framework. No matter which look-ahead strategies as described in \cite{IDRbreakdown,TuneIDR} are used to deal with these breakdowns, they can still occur since there is yet no strategy available that can cover all of them.

On the other hand there are issues that result only from the way how the method is implemented. Arguably, these weaknesses are caused by bad practices in the design of the implementation. As an example, an implementation of \GMRES which uses the power basis has weaknesses that are only related to the implementation and not to the underlying mathematical framework. Issues that are caused by the implementation itself can be removed by reformulating the algorithm in a numerically robuster way. This is exactly what we attempt to do for \IDRstabno in Sec.\,\ref{Sec:NovelRestartedGMRESbasedImpelementationOfIDRstab}.

\subsubsection{Implementation-caused weaknesses}

There is one simple question that lays out that all the so far discussed implementations of \IDR and \IDRstab have weaknesses that are caused only by the way how the respective methods are implemented.
\begin{center}
	What happens when $\bv\hp{0}_{q}$ is (nearly) collinear to $\lbrace\bv\hp{0}_{q+1},...,\bv\hp{0}_{s}\rbrace$?
\end{center}
We know the answer: All the current implementations would (nearly) break down because $\bZ$ becomes (nearly) singular. However, we would actually hope that the method can even strike a benefit from this scenario. In the following we explain why.

Consider the situation that \IDR respectively \IDRstab is in a repetition where $\br\hp{0},\bv_1\hp{0},...,\bv\hp{0}_s \in \cG_j$ and $\dim(\cG_j) < s$. From the theory we know that in this case \IDR respectively \IDRstab should terminate within this repetition of the main-loop. However, in this scenario the above collinearity issue must occur.

Besides, it is a well-known scenario that \GCR has a breakdown when the residual does not change during one iteration \cite{GCR-breakdown}. Consider the case where $\br\hp{0}$ is already orthogonal to the column $\bp_q$ in Alg.\,\ref{Algo:IDRs_bio} line 13 during the $q$th for-loop. In this case the two vectors $\bv\hp{0}_q$ and $\bv\hp{0}_{q+1}$ would be identical, which causes the above-described breakdown of all the so far discussed \IDRstabno variants as an immediate consequence of the unlucky breakdown property of all \GCR-type methods.
\largeparbreak

In the next section we present an implementation of \IDR and \IDRstab (the latter for $\ell=2$ only) that has a proper numerical treatment for the above described scenarios. This is achieved by replacing the nested \GCR by a nested \GMRES.

\section{\GMstab: The novel restarted \GMRESkp{$m$}-based implementation of \IDRstab}\label{Sec:NovelRestartedGMRESbasedImpelementationOfIDRstab}

We have discussed formerly that without any exception all currently available implementations of \IDRstab and \IDR have numerical issues, even when they should actually terminate. In this section we provide a restarted \GMRESkp{$m$}-based implementation of \IDR (i.e. \IDRstab with $\ell=1$ and without auxiliary decoupling) and of \IDRstab (only for $\ell=2$) that has the following two desirable properties:
\begin{itemize}
	\item If during the current repetition of the main-loop there is a sufficiently accurate solution for $\bx$ available with $\br\in \cG_j$ in the basis of the newly computed auxiliary vectors then our novel implementation constructs this solution in a numerically stable way and terminates.
	\item Otherwise, the new auxiliary vectors that will be used in the subsequent repetition of the main-loop are computed in a numerically stable way.
\end{itemize}

This section is organised as follows. In the next subsection we motivate the main idea of our \GMRES-based \IDRO methods. Afterwards we derive the new methods: First \GMstabI, which is the \GMRES-based implementation of \IDR respectively \IDRstab for $\ell=1$; and second \GMstabII, the \GMRES-based implementation of \IDRstab for $\ell=2$.

\subsection{Main idea of \GMstab}
As an introduction to the principal ideas of \GMstab, it helps the presentation to first consider only the case $\ell=1$.

The main idea of \GMstab is to compute the new auxiliary vectors in $\cG_{j+1}$ (cf. Alg.\,\ref{Algo:IDRs_noDec}, line 14) in a different way from simply overwriting the old auxiliary vectors. Instead, the novel approach consists of computing a basis of a suitable Krylov subspace that contains all the required new auxiliary vectors. The basis vectors of this suitable Krylov subspace in turn can be computed by using a reliable numerical scheme.

As a suitable Krylov subspace we choose
\begin{align*}
 	\cK_m(\bA \cdot \bPi ; \br\hp{0})\,,
\end{align*}
where $\bPi := \Big(\bI + (-\frac{1}{\omega_{j+1}}\cdot\bV\hp{-1}+\bV\hp{0}) \cdot \bZ\d \cdot \bP\t\Big)$, i.e. the projector matrix from Lem.\,\ref{Lem:DecFreeIDR}. Taking this matrix has the following advantage: We can choose $\br\hp{0} \in \cG_{j+1}$ and obtain $\cK_m(\bA \cdot \bPi ; \br\hp{0}) \subset \cG_{j+1}$. So new auxiliary vectors $\bv\hp{0}_q \in \cG_{j+1}$ can be obtained from the basis vectors of $\cK_m(\bA \cdot \bPi ; \br\hp{0})$.

The basis of this Krylov subspace is computed by an Arnoldi scheme. For instance, in our novel implementation of \IDR, the value for $m$ is $m=s$. Computing the basis of this particular Krylov subspace by an Arnoldi scheme provides the following trade-off, which can be considered as a positive result:
\begin{enumerate}[(A)]
	\item It is either the case that during the Arnoldi scheme a small value for $|h_{i+1,i}|$, i.e. a sub-diagonal element of the upper Hessenberg matrix, occurs. This means that the principal angles
	\begin{align*}
	\angle\Big(\ \cK_i(\bA \cdot \bPi ; \br\hp{0})\ ,\ \bA\cdot\cK_{i}(\bA \cdot \bPi ; \br\hp{0})\ \Big)
	\end{align*}
	are small, i.e. the Krylov subspace has almost converged in the sense of the principal angles \cite{principalAngles,stableMGS_GMRES}. Thus, an accurate solution for $\br\hp{0}$ can be found in this Krylov subspace by using, e.g., the minimal residual approach of \GMRES.
	\item Or -- instead -- it is the case that during the Arnoldi scheme the values for $|h_{i+1,i}|$ for all $i=1,...,m$ are sufficiently bounded from below. Then the basis vectors for $\cK_m(\bA \cdot \bPi ; \br\hp{0})$, that are obtained by the Arnoldi scheme, are orthonormal and live in $\cK_m(\bA \cdot \Pi ; \br\hp{0})$ with a high numerical accucacy since there is no division by a very small number.
\end{enumerate}
We strongly emphasise that both (A) and (B) are essential because only either of them can hold each time. Both are equally essential because of the following reasons:

When the basis for the Krylov subspace is stable than this means that there is never a division by a small value $|h_{i+1,i}|$, which is correlated to the fact that the method has not converged yet. Thus, using a method that is fully reliant on the case that (B) holds is equivalent in meaning to rely on the assumption that the method will never converge. This is a pointless assumption in itself.

On the other hand, when (B) fails, i.e. the basis of the Krylov subspace is inaccurate, then we cannot proceed. Thus, the only option is to terminate, which is just what (A) offers.

\subsection{The method \GMstabI}
In this subsection we present \GMstabI. This method it totally different from \GMstabII: Whereas \GMstabII is based on an augmented Arnoldi scheme, the method \GMstabI is based on a projected Arnoldi scheme.

This subsection is organised as follows. First, we introduce the projected Arnoldi decomposition to compute a basis of $\cK_m(\bA\cdot\bPi;\br\hp{0})$. Then we present an algorithm for \GMstabI, that utilises this Arnoldi scheme, and explain it line by line.

\subsubsection{The projected restarted \GMRESkp{$m$} method}
Obviously, since we speak of iterative methods, the projected Arnoldi decomposition used by \GMstabI is embedded in a \GMRESkp{$m$} approach. This makes sure that the method terminates immediately whenever a sufficiently accurate solution in the constructed Krylov subspace is available.

The suggested implementation for the projected \GMRES method is given in Alg.\,\ref{Algo:pGMRESm}.

The method computes the following decomposition, cf. lines 6--14:
\begin{align*}
	\bA \cdot (\bW_{:,1:m} - \bV\hp{-1} \cdot \bZ\d \cdot \underbrace{\bY_{:,1:m}}_{:= \bP\h\cdot\bW_{:,1:m}}) = \bW_{:,1:(m+1)} \cdot \bH_{1:(m+1),1:m}
\end{align*}
Further, it computes the reduced QR-decomposition $\bQ \cdot \bR$ of $\bH \in \C^{(m+1)\times m}$.

If there is a sufficiently small residual available in the $i$th iteration of the for-loop in line 5, the method solves the following least-squares problem and terminates:
\begin{align*}
	\bx := \bx\hp{0} + (\bW_{:,1:i} - \bV\hp{-1} \cdot \bZ\d \cdot \bY_{:,1:i}) \cdot \bzeta\,,
\end{align*}
where $\bzeta \in \C^i$ is chosen such that the 2-norm of the residual
\begin{align*}
	\br := \br\hp{0} - \bA \cdot \bx \equiv \bW_{:,1:(i+1)} \cdot (\beta \cdot\be_1 - \bH_{1:(i+1),1:i} \cdot \bzeta)
\end{align*}
is minimised, cf. lines 21--25.

\begin{algorithm}
	\begin{algorithmic}[1]
		\Procedure{pGMRESm}{$\bA,\bP,\bZ,\bV\hp{-1},\bx\hp{0},\br\hp{0},\beta,m,\tolabs$}
		\State \textit{// $\beta \equiv \|\br\hp{0}\|$}
		\State $\bw_1 := 1/\beta \cdot \br\hp{0}$,\quad$\bgam := \beta \cdot \be_1 \in \R^{m+1}$
		\State $\bY := \bO \in \R^{s \times (m+1)}$,\quad$\bQ := \bI \in \R^{(m+1) \times (m+1)}$
		\For{$i=1,...,m$}
		\State $\bY_{:,i} := \bP\t \cdot \bw_i$
		\State $\bxi := \bZ^{-1} \cdot \bY_{:,i} $
		\State $\bw_{i+1} := \bA \cdot (\bw_i - \bV\hp{-1} \cdot \bxi)$
		\For{$j=1,...,i$}
		\State $\bH_{j,i} := \langle \bw_j,\bw_{i+1}\rangle$\quad\textit{// $\equiv \langle \bw_j,\bw_{i+1}\rangle/\langle \bw_j,\bw_{j}\rangle$}
		\State $\bw_{i+1} := \bw_{i+1} - \bH_{j,i} \cdot \bw_j$
		\EndFor
		\State $\bH_{i+1,i} := \|\bw_{i+1}\|$
		\State $\bw_{i+1} := 1/\bH_{i+1,i} \cdot \bw_{i+1}$
		\State $\bR_{:,i} := \bH_{:,i}$
		\For{$j=1,...,(i-1)$}
		\State $\bR_{j:(j+1),i} := \bG\fp{j} \cdot \bR_{j:(j+1),i}$
		\EndFor
		\State [$\bG\fp{i},\bR_{i:(i+1),i}$] = \texttt{givens}($\bR_{i:(i+1),i}$)
		\State $\bQ_{i:(i+1),:} := \bG\fp{i} \cdot \bQ_{i:(i+1),:}$,\quad$\bgam_{i:(i+1),1} := \bG\fp{i} \cdot \bgam_{i:(i+1),1}$
		\If{$|\gamma_{i+1,1}|\leq\tolabs$}
		\State $\bzeta := \bR\d\cdot\bgam_{1:i,1} \in \R^i$
		\State $\bx\hp{0} := \bx\hp{0} + \bW_{:,1:i} \cdot \bzeta - \bV\hp{-1} \cdot \big(\bZ\d \cdot (\bY_{:,1:i} \cdot \bzeta)\big)$,\quad$\beta := |\gamma_{i+1,1}|$
		\State \textit{// $\br\hp{0} := \bW_{:,1:(i+1)} \cdot (\beta \cdot \be_1 - \bH_{1:(i+1),1:i} \cdot \bzeta)$}
		\State \Return $\bx\hp{0},\beta$
		\EndIf
		\EndFor
		\State $\bY_{:,m+1} := \bP\t\cdot\bw_{m+1}$
		\State $\bH := \bH_{1:(m+1),1:m}$,\quad$\bQ:=(\bQ_{1:m,1:(m+1)})\t \in \R^{(m+1) \times m}$
		\State $\bR := \bR_{1:m,1:m}$\quad\textit{// $\bH = \bQ \cdot \bR$}
		\State \Return $\bW,\bY,\bH,\bQ,\bR$
		\EndProcedure
	\end{algorithmic}
	\caption{projected \GMRESkp{$m$}}\label{Algo:pGMRESm}
\end{algorithm}

\subsubsection{The algorithm of \GMstabI}

Alg.\,\ref{Algo:GMstab1} shows the implementation of \GMstabI. The algorithm describes a so-called cycle, which is the interior of the main-loop. I.e. a complete algorithm can be obtained by using these lines of code:
\begin{algorithmic}
	\State \textit{// Given $\bA,\bb,s,\tolabs$}
	\State [$\bP,\bV\hp{-1},\bV\hp{0},\bZ,\bx\hp{0},\br\hp{0},\beta$] = \textsc{Initialisation}($\bA,\bb,s,\tolabs$)
	\While{$\beta>\tolabs$}
		\State \mbox{$[\bZ,\bV\hp{-1},\bV\hp{0},\bx\hp{0},\br\hp{0},\beta]$ = ...}\\
		\quad\quad\quad\quad\textsc{GMstab1\_cycle}($\bA,\bP,\bZ,\bV\hp{-1},\bV\hp{0},\bx\hp{0},\br\hp{0},\beta,\tolabs$)
	\EndWhile
	\State \textit{// return $\bx\hp{0}$}
\end{algorithmic}

Alg.\,\ref{Algo:GMstab1} consists of five steps, as indicated by the comments. In the following we describe each of these steps with reference to the implementation.
\largeparbreak

In the first step, the residual is moved into $\cG_{j+1}$. This is done in an analogous procedure to as is done in Alg.\,\ref{Algo:IDRs_noDec} lines 5--9.

Step 2 consists of modifying the projector. Given the matrices $\bV\hp{-1},\bV\hp{0}$ before line 7, the matrix $\bV\hp{-1}$ that is computed in line 13 can be expressed as
\begin{align*}
	\bV\hp{-1} := (1/\tau \cdot \bV\hp{-1} - \bV\hp{0}) \cdot \bG\,,
\end{align*}
where $\bG \in \C^{s \times s}$ is just some regular matrix that is tuned s.t. $\bV\hp{-1}$ is well-conditioned. In the next paragraph we explain why we overwrite $\bV\hp{-1}$ in this particular way. Anyway, we can already agree at this moment that since we overwrite $\bV\hp{-1}$ it makes sense to do this in such a way that it is well-conditioned. Yielding the well-conditionedness of $\bV\hp{-1}$ after the modification is the only purpose of all the computations of the second step. A construction of $\bG$ such that $\bV\hp{-1}$ obtains a condition number of 1 requires the reduced QR-decomposition in line 8.

Now we explain step 3: The motivation of overwriting the matrix $\bV\hp{-1}$ in the way laid out above can be understood by looking into Lem.\,\ref{Lem:DecFreeIDR} and the projected Arnoldi decomposition: The decomposition is
\begin{align*}
\bA \cdot (\bW_{:,1:m} - \bV\hp{-1} \cdot \bZ\d \cdot \underbrace{\bY_{:,1:m}}_{\equiv \bP\h\cdot\bW_{:,1:m}}) = \bW_{:,1:(m+1)} \cdot \bH_{1:(m+1),1:m}\,.
\end{align*}
Inserting for $\bV\hp{-1}$ the expression by which it was overwritten (namely $1/\tau \cdot \bV\hp{-1}- \bV\hp{0}$), we obtain
\begin{align*}
& &\bA \cdot \Big(\bW_{:,1:m} - (1/\tau \cdot \bV\hp{-1}- \bV\hp{0}) \cdot \bZ\d \cdot \underbrace{\bY_{:,1:m}}_{\equiv \bP\h\cdot\bW_{:,1:m}}\Big) &= \bW_{:,1:(m+1)} \cdot \bH\\
\Leftrightarrow& &\bA\cdot \Big(\bI + (-\frac{1}{\omega_{j+1}}\cdot\bV\hp{-1}+\bV\hp{0}) \cdot \bZ\d \cdot \bP\t\Big) \cdot \bW_{:,1:m} &= \bW_{:,1:(m+1)} \cdot \bH
\end{align*}
The latter equation can be analysed by using Lem.\,\ref{Lem:DecFreeIDR}: When $\bW_{:,1}$, i.e. the first column of $\bW$, lives in $\cG_{j+1}$ then it follows by induction that $\rg(\bW) \subset \cG_{j+1}$ holds.
In consequence, the following relations must hold for the subsequently computed projected Arnoldi decomposition in line 16\,:
\begin{align*}
\bA \cdot \underbrace{(\bW_{:,1:m} - \bV\hp{-1} \cdot \bZ\d \cdot \underbrace{\bY_{:,1:m}}_{\equiv \bP\h\cdot\bW_{:,1:m}})}_{\subset \cG_j \cap \cN(\bP)} &= \underbrace{\bW_{:,1:(m+1)}}_{\subset \cG_{j+1}} \cdot \bH
\end{align*}
Thus, choosing $\bV\hp{-1}$ in this way achieves that the basis matrix $\bW$ from the Arnoldi decomposition can be used to construct new auxiliary vectors $\bv\hp{0}_q \in \cG_{j+1}$ in a numerically robust way. This decomposition is computed for $m=s$.

The steps 4) and 5) in the algorithm consist of computing an updated residual and projectors from the new column vectors of $\bW$.

Step 4) attempts to compute a residual $\br\hp{0} \in \cG_{j+1} \cap \cN(\bP)$. Using a consistent update for $\bx\hp{0}$ and $\br\hp{0}$ of the form
\begin{align*}
\bx\hp{0} &:= \bx\hp{0} + (\bW_{:,1:s} - \bV\hp{-1} \cdot \bZ\d \cdot \bY_{:,1:s}) \cdot \bxi\\
\br\hp{0} &:= \bW_{:,1:(s+1)} \cdot (\beta \cdot\be_1 - \bH_{1:(s+1),1:s} \cdot \bxi)
\end{align*}
one can derive a formula for $\bxi \in \C^s$ such that $\bP\t\cdot\br\hp{0}=\bO$ follows from the following equations:
\begin{align*}
\bxi 
&=\big(\bP\t\cdot\underbrace{\bA\cdot(\bW_{:,1:s}-\bV\hp{-1}\cdot\bZ\d\cdot\bY_{:,1:s})}_{\equiv\bW_{:,1:(s+1)}\cdot\bH}\big)\d \cdot (\bP\t\cdot\underbrace{\br\hp{0}}_{\equiv\bW_{:,1}\cdot\beta}) \\
&=\big(\underbrace{\bP\t\cdot\bW_{:,1:(s+1)}}_{\equiv\bY_{:,1:(s+1)}}\cdot\bH\big)\d \cdot (\underbrace{\bP\t\cdot\bW_{:,1}}_{\equiv\bY_{:,1}}\cdot\beta) \\
&=\big(\underbrace{\bY_{:,1:(s+1)}\cdot\bH}_{\equiv\bL_Z\cdot\bQ_Z\t \cdot \bR_H}\big)\d \cdot (\underbrace{\bY_{:,1}\cdot\beta}_{\equiv\bEta})\\
&=\bR_H\d \cdot \bQ_Z \cdot \bL_Z\d \cdot \bEta
\end{align*}
The expression from the last equation is used in line 20.

Step 5) constructs a matrix for $\bV\hp{-1}$ such that its image $\bV\hp{0}$ satisfies the following two properties:
\begin{align*}
	\cond(\bV\hp{0}) &= 1\\
	\rg(\bV\hp{0}) &\subset \cG_{j+1}
\end{align*}
This can be achieved by computing $\bV\hp{-1}$ and $\bV\hp{0}$ by an expression of the following form, where $\bG$ is again an arbitrary matrix that is used to tune the conditioning:
\begin{align*}
	\bV\hp{-1} &:= (\bW_{:,1:s} - \bV\hp{-1} \cdot \bZ\d \cdot \bY_{:,1:s}) \cdot \bG\\
	\bV\hp{0} &:= {\bW_{:,1:(s+1)} \cdot \bH_{1:(s+1),1:s}} \cdot \bG \equiv\bA \cdot (\bW_{:,1:s} - \bV\hp{-1} \cdot \bZ\d \cdot \bY_{:,1:s}) \cdot \bG
\end{align*}
To achieve $\cond(\bV\hp{0}) = 1$ one must choose $\bG = \bR_H\d$, cf. lines 25--26. Further, to make sure that the subsequent matrix for $\bZ\equiv\bY \cdot \bQ_H$ is lower triangular the LQ-decomposition from line 19 is used and $\bQ_Z$ is multiplied from the right onto $\bG$. I.e., the overall matrix for $\bG$ is $\bG:=\bR_H\d\cdot\bQ_Z$.

\begin{algorithm}
	\begin{algorithmic}[1]
		\Procedure{GMstab1\_cycle}{$\bA,\bP,\bZ,\bV\hp{-1},\bV\hp{0},\bx\hp{0},\br\hp{0},\beta,\tolabs$}
		\State \textit{// - - - 1) Move $\br\hp{0}$ from $\cG_j$ into $\cG_{j+1}$ - - - }
		\State $\br\hp{1}:= \bA\cdot \br\hp{0}$
		\State $[\tau,\beta] := $\textsc{StabCoeffs($\br\hp{0},\,\br\hp{1},\beta$)}
		\State $\bx\hp{0} := \bx\hp{0} + \tau \cdot \br\hp{0}$
		\State $\br\hp{0} := \br\hp{0} - \tau \cdot \br\hp{1}$\quad\textit{// terminate if $\beta\leq\tolabs$}
		\State \textit{// - - - 2) Modify the current projector - - - }
		\State $\left[\ [\bV\hp{0},\bV\hp{-1}]\,,\, \begin{bmatrix}\bI_{s \times s} & \bC \\
		\bO_{s\times s} & \bR
		\end{bmatrix} \vphantom{\begin{matrix}
			a\\
			b\\
			c
			\end{matrix}}\ \right] $ = \texttt{qr}$\Big([\bV\hp{0},\bV\hp{-1}]\Big)$
		\State $\bF :=  
		\begin{bmatrix}
		1/\tau \cdot \bC - \bI_{s \times s} \\
		1/\tau \cdot \bR
		\end{bmatrix} \in \R^{(2\cdot{}s)\times{}s}$
		\State [$\bQ_F,\bR_F$] = \texttt{qr}($\bF$)
		\State $\bZ := -\bZ \cdot \bR_F\d$
		\State [$\bZ,\bQ_Z$] := \texttt{lq}($\bZ$)
		\State $\bV\hp{-1} := [\bV\hp{0},\bV\hp{-1}]\cdot(\bQ_F\cdot\bQ_Z)$
		\State \textit{// - - - 3) Compute new basis vectors - - - }
		\State \textit{// Use the storage of $[\bV\hp{0},\br\hp{0}] \in \R^{N \times (s+1)}$ for $\bW \in \R^{N \times (s+1)}$\,.}
		\State [$\bW,\bY,\bH,\bQ_H,\bR_H$] = \textsc{pGMRESm}($\bA,\bP,\bZ,\bV\hp{-1},\bx\hp{0},\br\hp{0},\beta,s,\tolabs$)
		\State \textit{// - - - 4) Biorthogonalise the residual - - - }
		\State $\bgam := \beta \cdot \be_1 \in \R^{s+1}$,\quad$\bEta := \beta \cdot \bY_{:,1}$
		\State [$\bL_Z,\bQ_Z$] = \texttt{lq}($\bY \cdot \bQ_H$)
		\State $\bxi := \bR_H\d \cdot \big(\bQ_Z \cdot (\bL_Z\d \cdot \bEta)\big)$
		\State $\bx\hp{0} := \bx\hp{0} + \bW_{:,1:s} \cdot \bxi - \bV\hp{-1} \cdot \big(\bZ\d \cdot (\bY_{:,1:s}\cdot\bxi)\big)$
		\State $\bc\hp{0} := \bgam - \bH \cdot \bxi$
		\State $\beta := \|\bc\hp{0}\|$
		\State \textit{// - - - 5) Build the next projector - - - }
		\State $\bV\hp{-1} := \bW_{:,1:s} \cdot (\bR_H\d\cdot\bQ_Z) - \bV\hp{-1} \cdot \Big(\bZ\d \cdot \big(\bY_{:,1:s} \cdot (\bR_H\d\cdot\bQ_Z) \big) \Big)$
		\State $[\bV\hp{0}\,,\,\br\hp{0}]:= \bW_{:,1:s+1}\cdot\big[\bQ_H\cdot\bQ_Z\,,\,\bc\hp{0}\big]$,\quad $\bZ := \bL_Z$
		\State \Return $\bZ,\bV\hp{-1},\bV\hp{0},\bx\hp{0},\br\hp{0},\beta$
		\EndProcedure
	\end{algorithmic}
	\caption{Robust cycle for \IDRstabp{$s$}{$1$}}\label{Algo:GMstab1}
\end{algorithm}

\subsubsection{Discussion of the new method}
Since the new column vectors $\bW$ are computed by an Arnoldi scheme we can be sure there is no scheme that could have computed the new basis vectors for $\bV\hp{0}$ in a more robust and reliable way. Further, we can be sure that in the case where the method terminates with an acceptable solution, this termination is performed during \textsc{pGMRESm} because the biorthogonal residual in lines 21,22,26 is at best only as small as the minimal residual but never smaller.

There is a further benefit of the interiorly used \GMRES method. During the computation of the matrix-vector products, which dominate the cost of the method, the monitored residual norm decreases monotonously. So one can be sure that not a single matrix-vector product is wasted but that instead the method terminates as soon as this is possible.

All the used projectors $\bV\hp{-1},\bV\hp{0}$ and the overwritten projector $\bV\hp{-1}$ in line 13 are well-conditioned. This is enforced by the QR-decomposition in line 8 and the consistent formulas for $\bV\hp{-1},\bV\hp{0}$ in lines 25--26.
\largeparbreak

Clear drawbacks of the new implementation are that the QR-decomposition in line 8 introduces some additional cost in terms of DOTs and AXPYs (although one notices that since $\bV\hp{0}$ has already orthonormalised columns the Gram-Schmidt procedure can be directly started from the columns of $\bV\hp{-1}$). Further to that, the Arnoldi procedure introduces some additional costs for the orthogonalisation of the columns of $\bW$. However, the average number of DOTs and AXPYs per matrix-vector product lives in $\cO(s)$. Thus, for small values of $s$ the cost is still dominated by the matrix-vector products with the system matrix and the preconditioners.
\largeparbreak

We have seen that \IDR, i.e. $\ell=1$, can be implemented such that there is no decoupling of $\bV\hp{-1},\bV\hp{0}$ by means of the accuracy in which the equation
\begin{align*}
\bV\hp{0} = \bA \cdot \bV\hp{-1}
\end{align*}
is satisfied.

\GMstabI has a similar feature: The new matrices $\bV\hp{-1},\bV\hp{0}$ are both computed from $\bV\hp{-1},\bW$. Thus the decoupling of the old projectors $\bV\hp{-1},\bV\hp{0}$ does not influence the decoupling of the new projectors. We emphasize that this is a very useful property as well because it guarantees that the decoupling of the projector matrices cannot amplify itself over several iterations.

\subsection{The method \GMstabII}
In the following we derive \GMstabII. This method uses an interior augmented Arnoldi decomposition. In analogy to the former subsection on \GMstabI, the structure of this subsection is as follows:

First we introduce the augmented Arnoldi procedure that is embedded into a restarted \GMRESkp{$m$} method. Then we present the algorithm of \GMstabII and explain how it works and utilises the augmented Arnoldi decomposition.

\subsubsection{The augmented restarted \GMRESkp{$m$} method}

The augmented Arnoldi scheme computes a decomposition of the following form.
\begin{align*}
	\bA \cdot \bW_{:,1:m} = [\bV\hp{0},\,\bW_{:,1:(m+1)}] \cdot \begin{bmatrix}
		\bZ\d \cdot \bY_{:,1:m}\\
		\bH_{1:(m+1),1:m}
	\end{bmatrix}\,,
\end{align*}
where this time it is $\bY_{:,i} = \bP\h\cdot\bA\cdot\bW_{:,i}$, $\bY \in \C^{s \times m}$. We stress that in the augmented Arnoldi scheme the expression for $\bY_{:,i}$ is different from that of the projected Arnoldi scheme. There it was $\bY_{:,i} = \bP\h\cdot\bW_{:,i}$, $\bY\in\C^{s \times (m+1)}$.

The augmented \GMRESkp{$m$} method is presented in Alg.\,\ref{Algo:augGMRESm}.
We briefly go through the implementation. Lines 6--15 compute the augmented Arnoldi basis. Here, the major difference in comparison to \textsc{pGMRESm} is that in lines 6 and 7 first the matrix-vector product and afterwards the column of $\bY$ is computed. After the computation of the new basis vector there are again Givens rotations to generate a QR-decomposition of $\bH$.

If after the $i$th iteration a sufficiently accurate solution exists in the span of the computed basis vectors of $\bW$ then in lines 22--26 the method computes a residual-optimal solution by using the following formulas:
\begin{align*}
	\bx &:= \bx\hp{0} + (\bW_{:,1:i}-\bV\hp{-1} \cdot \bZ\d\cdot\bY_{:,1:i}) \cdot \bzeta\\
	\br &:= \bW_{:,1:(i+1)} \cdot (\bgam - \bH_{1:(i+1),1:i} \cdot \bzeta)\,,
\end{align*}
where $\bzeta \in \C^i$ minimises $\|\bgam - \bH_{1:(i+1),1:i} \cdot \bzeta\|$. We briefly explain why this leads to a residual-minimal solution: The augmented Arnoldi decomposition is in fact equivalent to a projected decomposition
\begin{align*}
	\bA \cdot (\bW_{1:i} - \bV\hp{-1} \cdot \bZ\d \cdot \bY_{:,1:i}) = \bW_{:,1:(i+1)} \cdot \bH_{1:(i+1),1:i}\,,
\end{align*}
where $\bV\hp{-1} = \bA^{-1} \cdot \bV\hp{0}$. Since $\cond(\bW_{:,1:(i+1)})=1$, one can find $\bzeta$ again by solving the projected least-squares problem in line 23.

\begin{algorithm}
	\begin{algorithmic}[1]
		\Procedure{augGMRESm}{$\bA,\bP,\bZ,\bV\hp{-1},\bV\hp{0},\bx\hp{0},\br\hp{0},\beta,m,\tolabs$}
		\State \texttt{// $\beta \equiv \|\br\hp{0}\|$}
		\State $\bw_1 := 1/\beta \cdot \br\hp{0}$,\quad$\bgam := \beta \cdot \be_1 \in \R^{m+1}$
		\State $\bY := \bO \in \R^{s \times m}$,\quad$\bQ := \bI \in \R^{(m+1) \times (m+1)}$
		\For{$i=1,...,m$}
		\State $\bw_{i+1} = \bA \cdot \bw_i$
		\State $\bY_{:,i} := \bP\t \cdot \bw_{i+1}$
		\State $\bxi := \bZ\d \cdot \bY_{:,i}$
		\State $\bw_{i+1} := \bw_{i+1} - \bV\hp{0} \cdot \bxi$
		\For{$j=1,...,i$}
		\State $\bH_{j,i} := \langle \bw_j,\bw_{i+1}\rangle$\quad\textit{// $\equiv \langle \bw_j,\bw_{i+1}\rangle/\langle \bw_j,\bw_{j}\rangle$}
		\State $\bw_{i+1} := \bw_{i+1} - \bH_{j,i} \cdot \bw_j$
		\EndFor
		\State $\bH_{i+1,i} := \|\bw_{i+1}\|$
		\State $\bw_{i+1} := 1/\bH_{i+1,i} \cdot \bw_{i+1}$
		\State $\bR_{:,i} := \bH_{:,i}$
		\For{$j=1,...,(i-1)$}
		\State $\bR_{j:(j+1),i} := \bG\fp{j} \cdot \bR_{j:(j+1),i}$
		\EndFor
		\State [$\bG\fp{i} \cdot \bR$] = \texttt{givens}($\bR_{i:(i+1),i}$)
		\State $\bQ_{i:(i+1),:} := \bG\fp{i} \cdot \bQ_{i:(i+1),:}$,\quad$\bgam_{i:(i+1),1} := \bG\fp{i} \cdot \bgam_{i:(i+1),1}$
		\If{$|\gamma_{i+1,1}|\leq\tolabs$}
		\State $\bzeta := \bR\d\cdot\bgam_{1:i,1} \in \R^i$
		\State $\bx\hp{0} := \bx\hp{0} + \bW_{:,1:i} \cdot \bzeta - \bV\hp{-1} \cdot \big(\bZ\d\cdot(\bY_{:,1:i}\cdot\bzeta)\big)$,\quad$\beta := |\gamma_{i+1,1}|$
		\State \Return $\bx\hp{0},\beta$
		\EndIf
		\EndFor
		\State $\bH := \bH_{1:(m+1),1:m}$,\quad$\bQ:=(\bQ_{1:m,1:(m+1)})\t \in \R^{(m+1) \times m}$
		\State $\bR := \bR_{1:m,1:m}$\quad\textit{// $\bH = \bQ \cdot \bR$}
		\State \Return $\bW,\bY,\bH,\bQ,\bR$
		\EndProcedure
	\end{algorithmic}
	\caption{augmented \GMRESkp{$m$}}\label{Algo:augGMRESm}
\end{algorithm}

\subsubsection{The algorithm of \GMstabII}

An implementation of \GMstabII is given in Alg.\,\ref{Algo:GMstab2}. The method consists of 4 steps, as indicated by the comments. In the following we explain what is computed during each step.

Step 1 consists of computing the new basis vectors as columns of $\bW \in \C^{N \times (2 \cdot s + 3)}$. The Arnoldi decomposition has the following properties: Since $\bW_{:,1} = 1/\beta \cdot \br\hp{0} \in \cG_{j} \cap \cN(\bP)$, it follows by induction:
\begin{align*}
	\bA \cdot \underbrace{\bW_{:,1:(2 \cdot s+2)}}_{\subset \cG_{j} \cap \cN(\bP)} = \underbrace{ [\overbrace{\bV\hp{0}}^{\subset \cG_j},\,\overbrace{\bW_{:,1:(2\cdot s+3)}}^{\subset \cN(\bP)}] \cdot \begin{bmatrix}
		\bZ\d \cdot \bY_{:,1:(2 \cdot s + 2)}\\
		\hhbH_{1:(2 \cdot s + 3),1:(2 \cdot s + 2)}
	\end{bmatrix} }_{\subset \cG_{j}}
\end{align*}
In particular, the lower brace on the right-hand side results from Lem.\,\ref{Lem:RemainingVectorsInGj}. From the choice of $\bY$ follows $\rg(\bW)\subset\cN(\bP)$, from which in turn with $\rg(\bV\hp{0})\subset\cG_j$ follows the lower brace on the left-hand side.

For the matrix $\hhbH \in \C^{(2 \cdot s + 3) \times (2 \cdot s +2)}$ we use the following notation:
\begin{align*}
	\bH &:= \hhbH_{1:(2\cdot s + 1),1:(2\cdot s)} \in \C^{(2 \cdot s + 1) \times (2 \cdot s)}\\
	 \hbH &:= \hhbH_{1:(2\cdot s + 2),1:(2 \cdot s + 1)} \in \C^{(2 \cdot s + 2) \times (2 \cdot s+1)}\,.
\end{align*}

Step 2 consists of finding a linear combination vector $\bxi$ for $\br\hp{0} = \bW \cdot (\bgam - \bH \cdot \bxi)$ such that $\br\hp{0} \in \cK^\perp_3(\bA\h;\bP)$ holds. We remember from Alg.\,\ref{Algo:IDRstab_block} line 21 that $\br\hp{0} \perp \cK_{\ell+1}(\bA\h;\bP)$ was required in \IDRstab after the biorthogonalisation. This is equivalent to $\br\hp{0},\br\hp{1},\br\hp{2} \in \cN(\bP)$. In the following we derive the formula for $\bxi$ that leads to these biorthogonality properties of the residual:

The idea for the derivation of a formula for $\bxi$ is to find expressions for the updated residual and solution. Imposing conditions on the updated residual, we can then derive conditions for $\bxi$. During the derivation of the formulas for $\bxi$, to distinguish the updated quantities from the current quantities, we use $\hbx\hp{0}$,$\hbr\hp{0}$,$\hbr\hp{1}$,$\hbr\hp{2}$ for the updated solution, residual and its powers, whereas we use $\bx\hp{0}$,$\br\hp{0}$ for the current solution and residual. An approach of the following expression is used for computing $\hbx\hp{0}$.
\begin{align*}
	\hbx\hp{0} := \bx\hp{0} + (\bW_{:,1:(2\cdot s)} - \bV\hp{-1} \cdot \bZ\d \cdot \bY_{:,1:(2 \cdot s)}) \cdot \bxi
\end{align*}
The consistent update for the residual is then given by
\begin{align*}
\hbr\hp{0} := \bW_{:,1:(2 \cdot s + 1)} \cdot \underbrace{(\bgam - \bH \cdot \bxi)}_{=:\bc\hp{0} \in \C^{2 \cdot s + 1}}\,.
\end{align*}
Since $\rg(\bW_{:,1:(2 \cdot s + 1)}) \subset \cG_j \cap \cN(\bP)$ already holds it automatically follows $\br\in\cN(\bP)$. Thus, there follows no condition for $\bxi$. Next, we look into the power vector of $\hbr\hp{0}$:
\begin{align*}
\hbr\hp{1} =& \bA \cdot \bW_{:,1:(2 \cdot s + 1)} \cdot \bc\hp{0}\\
 			=& [\bV\hp{0},\,\bW_{:,1:(2\cdot s + 2)}] \cdot \begin{bmatrix}
 			\bZ\d \cdot \bY_{:,1:(2 \cdot s + 1)}\\
 			\hbH
 			\end{bmatrix} \cdot \bc\hp{0}
\end{align*}
Since $\hbr\hp{1} \in \cN(\bP)$ is required by the geometric principle of \IDRstabno (cf., e.g., Alg.\,\ref{Algo:IDRstab_block} line 21) and it holds $\rg(\bV\hp{0})\nsubset \cN(\bP)$, we have to choose $\bxi$ such that the $s$ conditions $\bY_{:,1:(2 \cdot s+1)} \cdot \bc\hp{0} = \bO$ hold. Satisfying this equation, it holds for $\hbr\hp{1}$:
\begin{align*}
\hbr\hp{1} = \bW_{:,1:(2\cdot s + 2)} \cdot
\underbrace{\hbH \cdot \bc\hp{0}}_{=:\bc\hp{1} \in \C^{2 \cdot s + 2}}
\end{align*}
Consequently, also $\hbr\hp{1} \in \cN(\bP)$ is satisfied because of $\rg(\bW) \subset \cG_j \cap \cN(\bP)$.

Next, we look into the second power vector of $\hbr\hp{0}$:
\begin{align*}
\hbr\hp{2} =& \bA \cdot \hbr\hp{1} \equiv \bA \cdot \bW_{:,1:(2 \cdot s + 2)} \cdot \bc\hp{1}\\
=& [\bV\hp{0},\,\bW_{:,1:(2\cdot s + 3)}] \cdot \begin{bmatrix}
\bZ\d \cdot \bY_{:,1:(2 \cdot s + 2)}\\
\hhbH
\end{bmatrix} \cdot \bc\hp{1}
\end{align*}
Again, in order to achieve $\hbr\hp{2} \in \cN(\bP)$ we can simply choose $\bxi$ such that the $s$ additional equations $\bY_{:,1:(2\cdot s+2)} \cdot \bc\hp{1}=\bO$ hold. The overall equation system for $\bxi \in \C^{2 \cdot s}$ is then given by
\begin{align*}
	\begin{bmatrix}
		\bY_{:,1:(2 \cdot s+1)}\phantom{{} \cdot \hbH {}} \cdot \bH\\
		\bY_{:,1:(2 \cdot s+2)}         {} \cdot \hbH {}  \cdot \bH
	\end{bmatrix} 
	\cdot \bxi = 
	\begin{pmatrix}
		\bY_{:,1:(2 \cdot s+1)}\phantom{{} \cdot \hbH {}} \cdot \bgam_{1:(2\cdot s)\phantom{+1}}\\
		\bY_{:,1:(2 \cdot s+2)} {}\cdot \hbH {}\cdot \bgam_{1:(2 \cdot s+1)}
	\end{pmatrix}\,.
\end{align*}
A variant of this formula with an improved conditioning for the rows of $\bY$ is used in the algorithm in line 8.

Step 3 computes the stability polynomial coefficients and updates $\bx\hp{0}$ consistently to the polynomial update $q_{j,2}(\bA)\cdot\hbr\hp{0}$. As an algorithmic detail, the stability polynomial coefficients can be computed from the linear combination vectors $\bc\hp{0},\bc\hp{1},\bc\hp{2}$. Next, we explain the formula in line 14.
\begin{align*}
	\bx\hp{0} :=& \hbx\hp{0} + [\hbr\hp{0},\hbr\hp{1}] \cdot \btau\\
			  =& \bx\hp{0} + (\bW_{:,1:(2 \cdot s)}-\bV\hp{-1}\cdot \bZ\d\cdot\bY_{:,1:(2 \cdot s)}) \cdot \bxi + \bW_{:,1:(2 \cdot s +2)}\cdot ([\underline{\bc}\hp{0},\bc\hp{1}] \cdot \btau)
\end{align*}
In the latter equation, the linear combination vector for the columns of $\bW$ can be combined. This results in the formula of line 14. The underline-notation, e.g. $\underline{\bc}\hp{0}$, means that the respective quantity is augmented with one row of zeros.
\largeparbreak

Step 4 computes the new projectors and the final residual. In the formerly discussed \IDRstab implementations for $\ell=2$, matrices $\bV\hp{-1},\bV\hp{0},\bV\hp{1},\bV\hp{2}$ are built such that $\bV\hp{g}=\bA^g \cdot \bV\hp{0}$ $\forall g\in\lbrace-1,0,1,2\rbrace$ and $\rg(\bV\hp{g})\subset\cN(\bP)$ $\forall g \in \lbrace 0,1\rbrace$. Then two weighted sums of these matrices are built, cf. e.g. Alg.\,\ref{Algo:IDRstab_block} lines 25--26.

In the following we derive the formulas from Alg.\,\ref{Algo:GMstab2} lines 21--22 for the new projectors from the above principle. For clarity of notation, we use $\tbV\hp{-1},\tbV\hp{0},\tbV\hp{1},\tbV\hp{2}$ for the projectors after the biorthogonalisation and $\hbV\hp{-1},\hbV\hp{0}$ for the updated projectors after the polynomial step.

First, using the augmented Arnoldi decomposition, we can derive formulas for $\tbV\hp{-1},\tbV\hp{0},\tbV\hp{1},\tbV\hp{2}$:
\begin{align*}
	\tbV\hp{-1} &= (\bW_{:,1:(2 \cdot s)} - \bV\hp{-1} \cdot \bZ\d \cdot \bY_{:,1:(2 \cdot s)}) \cdot \bG\\
	\tbV\hp{0} &= \bW_{:,1:(2 \cdot s + 1)} \cdot \bH \cdot \bG\,,
\end{align*}
where $\bG \in \C^{(2 \cdot s) \times s}$ is an arbitrary regular matrix. Choosing $\bG$ such that
\begin{align*}
	\bY_{:,1:(2 \cdot s + 1)} \cdot \bH \cdot \bG = \bO\tageq\label{eqn:ReqbG}
\end{align*}
we can reformulate the equation for $\tbV\hp{0}$ as follows by subtracting a zero:
\begin{align*}
	\tbV\hp{0} &= (\bW_{:,1:(2 \cdot s + 1)} - \bV\hp{-1} \cdot \bZ\d \cdot \bY_{:,1:(2 \cdot s + 1)})\cdot \bH \cdot \bG
\end{align*}
From this equation it is easy to find a representation of $\tbV\hp{1}$ using the basis vectors from the Arnoldi equation.
\begin{align*}
\tbV\hp{1} &= \bW_{:,1:(2 \cdot s + 2)} \cdot \hbH \cdot \bH \cdot \bG
\end{align*}
Multiplying the above equation from the left by $\bA$ and inserting the augmented Arnoldi decomposition again, it follows for $\tbV\hp{2}$ the equation:
\begin{align*}
\tbV\hp{2} &= [\bV\hp{0},\,\bW_{:,1:(2 \cdot s + 3)}] \cdot \begin{bmatrix}
\bZ\d \cdot \bY_{:,1:(2 \cdot s + 2)}\\
\hhbH
\end{bmatrix} \cdot \hbH \cdot \bH \cdot \bG
\end{align*}
The formulas for the updated matrices are:
\begin{align*}
	\hbV\hp{-1} &= \tbV\hp{-1} - \tau_1 \cdot \tbV\hp{0} - \tau_2 \cdot \tbV\hp{1}\\
	\mathrlap{\hbV\hp{0}}\phantom{\hbV\hp{-1}} &= \mathrlap{\tbV\hp{0}}\phantom{\tbV\hp{-1}} - \tau_1 \cdot \tbV\hp{1} - \tau_2 \cdot \tbV\hp{2}
\end{align*}
Inserting all the formulas into one big expression, we obtain
\begin{align*}
	\hbV\hp{-1} =& [\bW_{:,1:(2 \cdot s + 2)},\,\bV\hp{-1}] \\
	&\cdot \Big( 
	\begin{bmatrix}
		\uubI\\
		-\bZ\d \cdot \bY_{:,1:(2 \cdot s)}
	\end{bmatrix} - \tau_1 \cdot
	\begin{bmatrix}
		\ubH\\
		\bO
	\end{bmatrix} - \tau_2 \cdot
	\begin{bmatrix}
		\hbH \cdot \bH\\
		\bO
	\end{bmatrix}
	\Big) \cdot \bG\\
	\mathrlap{\hbV\hp{0}}\phantom{\hbV\hp{-1}} =& [\bW_{:,1:(2 \cdot s + 3)},\,\mathrlap{\bV\hp{0}}\phantom{\bV\hp{-1}}]\\
	 &\cdot \Big( 
	\begin{bmatrix}
		\underline{\underline{\bH}}\\
		\bO
	\end{bmatrix} - \tau_1 \cdot
	\begin{bmatrix}
		\underline{\hbH} \cdot \bH\\
		\bO
	\end{bmatrix} - \tau_2 \cdot
	\begin{bmatrix}
		\hhbH \cdot \hbH \cdot \bH\\
		\bZ\d \cdot \bY_{:,1:(2 \cdot s + 2)} \cdot \hbH \cdot \bH
	\end{bmatrix}
	\Big) \cdot \bG\,.\tageq\label{eqn:hbVhp0}
\end{align*}
The approach is now to find a matrix for $\bG$ that satisfies \eqref{eqn:ReqbG}, leads to $\cond(\hbV\hp{0})=1$ and makes $\bP\h\cdot\hbV\hp{0}$ lower triangular.

In the following we explain how a candidate for $\bG \in \C^{(2 \cdot s) \times s}$ that satisfies all the above requirements is constructed in \GMstabII.

The matrix $\bG$ is expressed as a product of three matrices.
\begin{align*}
	\bG = \bQ_G \cdot \bR_F\d \cdot \bQ_Z
\end{align*}
The first matrix $\bQ_G \in \C^{(2 \cdot s) \times s}$ is unitary and has the property $\bY_{:,1:(2 \cdot s + 1)} \cdot \bH \cdot \bQ_G = \bO$. We construct this matrix from the Householder QR-decomposition of the matrix $(\bQ_\bY \cdot \bH)\h \in \R^{(2 \cdot s) \times s}$. Let the unreduced Householder QR-decomposition be given as:
\begin{align*}
	[\underbrace{\bT}_{\in \C^{(2 \cdot s) \times s}},\underbrace{\bQ_G}_{\in \C^{(2 \cdot s) \times s}}] \cdot \begin{bmatrix}
	\bR_T\\
	\bO
	\end{bmatrix} = (\bQ_\bY \cdot \bH)\h
\end{align*}
We now show that by this construction the matrix $\bG$ already satisfies \eqref{eqn:ReqbG}:
\begin{align*}
	& &\bY_{:,1:(2 \cdot s + 1)} \cdot \bH \cdot \bG &= \bO\\
	\Leftrightarrow& &\bQ_\bY \cdot \bH \cdot \bQ_G \cdot \bR_F\d \cdot \bQ_Z &= \bO\\
	\Leftarrow& &\bQ_\bY \cdot \bH \cdot \bQ_G &= \bO\\
	\Leftrightarrow& &[\bR_T\h,\,\bO] \cdot \begin{bmatrix}
	\bT\h\\
	\bQ_G\h
	\end{bmatrix} \cdot \bQ_G &= \bO\\
	\Leftrightarrow& &[\bR_T\h,\,\bO] \cdot \begin{bmatrix}
	\bO\\
	\bI
	\end{bmatrix} &= \bO
\end{align*}
The latter equality obviously holds.

The second factor in $\bG$ is responsible for the condition number of $\hbV\hp{0}$. Using the matrix $\bC \in \C^{(3 \cdot s + 3) \times (2 \cdot s)}$ as computed in line 18 and the reduced QR-decomposition of $\bQ_W \cdot \bR_W = [\bW_{:,1:(2 \cdot s + 3)},\,\bV\hp{0}]$, the formula \eqref{eqn:hbVhp0} can be reformulated as
\begin{align*}
	\hbV\hp{0} = \bQ_W \cdot \underbrace{\bR_W \cdot \bC \cdot \bQ_G}_{=:\bF \in \C^{(3 \cdot s+3) \times s}} \cdot \bR_F\d \cdot \bQ_Z\,.
\end{align*}
In this formula the matrix $\bQ_Z$ is unitary, too. Thus, the condition number of $\hbV\hp{0}$ is identical to the condition number of
\begin{align*}
	\bF \cdot \bR_F\d\,.
\end{align*}
We compute $\bR_F\d$ from the QR-decomposition of $\bF$, cf. lines 18--20. This guarantees that $\cond(\hbV\hp{0})=1$ will hold.

Given the formula for $\hbV\hp{0}$, we can derive a formula for $\hbZ := \bP\h \cdot \hbV\hp{0}$:
\begin{align*}
	\hbZ &= \bP\h \cdot \hbV\hp{0}\\
			&= [\underbrace{\bP\h \cdot \bV\hp{0}}_{\equiv \bZ},\,\underbrace{\bP\h \cdot \bW_{:,1:(2 \cdot s+3)}}_{\equiv \bO}] \cdot \bC \cdot \bG\\
			&= \bZ \cdot \underbrace{(-\tau_2 \cdot \bZ\d \cdot \bY_{:,1:(2 \cdot s + 2)} \cdot \hbH \cdot \bH)}_{\text{from }\bC} \cdot \bG\\
			&= \underbrace{-\tau_2 \cdot \bY_{:,1:(2 \cdot s + 2)} \cdot \hbH \cdot \bH \cdot \bQ_G \cdot \bR_F\d}_{(\star)} \cdot \bQ_Z
\end{align*}
We use an LQ-decomposition for the under-braced matrix $(\star)$, cf. line 21. In this way the unitary matrix $\bQ_Z \in \C^{s \times s}$ achieves that $\hbZ \equiv \bL_Z$ is lower triangular.

The lines 22--23 overwrite the projector matrices $\bV\hp{-1},\bV\hp{0}$ by using the above formulas and priorly discussed product representation of $\bG$. Finally, the updated residual $\br\hp{0}$ is computed. The computation of the updated residual is performed so late because $\bW$ is written into the array space of $\br\hp{0}$. Thus, both quantities cannot exist at the same time. As is clear from line 3, the algorithm \GMstabII does not require more storage than the \IDRstab implementations that have been discussed in Sec.\,\ref{sec:ImplementationsIDRstab}.

\begin{algorithm}
	\begin{algorithmic}[1]
		\Procedure{GMstab2\_cycle}{$\bA,\bP,\bZ,\bV\hp{-1},\bV\hp{0},\bx\hp{0},\br\hp{0},\beta,\tolabs$}
		\State \textit{// - - - 1) Compute new basis vectors - - - }
		\State \textit{// Use the storage of $[\br\hp{0},\br\hp{1}\bV\hp{1},\br\hp{2},\bV\hp{2}] \in \R^{N \times (2 \cdot s + 3)}$ for $\bW$}
		\State [$\bW,\bY,\hhbH,\bQ_H,\bR_H$] = ...\quad\quad\quad\quad\quad\quad\quad\quad\quad\quad\quad\quad\quad\quad\newline \mbox{\quad\quad\quad\quad\textsc{augGMRESm}($\bA,\bP,\bZ,\bV\hp{-1},\bV\hp{0},\bx\hp{0},\br\hp{0},\beta,2 \cdot s+2,\tolabs$)}
		\State \textit{// - - - 2) Biorthogonalise the subsequent residual}
		\State $\bgam := \be_1 \cdot \beta \in \R^{2 \cdot s + 1}$
		\State $\bQ_{\bY} = \texttt{roworth}(\bY_{:,1:2\,s+1})$,\quad$\bQ_{\hbY} = \texttt{roworth}(\bY_{:,1:2\,s+2})$
		\State $\bxi := \bR_H\d \cdot \left(\vphantom{\begin{matrix}
			a\\
			b\\
			c
			\end{matrix}} \begin{bmatrix}
		\bQ_{\bY}\cdot\bQ_H \\
		\bQ_{\hbY}\cdot\hbH\cdot\bQ_H
		\end{bmatrix}\d \cdot \begin{pmatrix}
		\bQ_{\bY}\cdot\bgam \\
		\bQ_{\hbY}\cdot\hbH\cdot\bgam
		\end{pmatrix} \right)$
		\State \textcolor{blue}{$\bc\hp{0} := \bgam - \bH\cdot\bxi \vphantom{\Big(\Big)}$}\quad \textit{// cf. above for definition of $\bH,\hbH,\hhbH$}
		\State \textcolor{blue}{$\bc\hp{1} := \hbH\cdot\bc\hp{0} \vphantom{\Big(\Big)}$}
		\State \textcolor{blue}{$\bc\hp{2} := \hhbH\cdot\bc\hp{1}\vphantom{\Big(\Big)}$},\quad$\beta := \|\bc\hp{0}\|$
		\State \textit{// - - - 3) Polynomial step - - - }
		\State $[\btau,\beta] := $\textsc{StabCoeffs($\underline{\underline{\bc}}\hp{0},\underline{\bc}\hp{1},\bc\hp{2},\beta$)}
		\State $\bx\hp{0} := \bx\hp{0} + \bW_{:,1:(2\cdot s+2)}\cdot\textcolor{blue}{\Big(
		\underline{\underline{\bxi}} + \tau_1 \cdot \underline{\bc}\hp{0} + \tau_2 \cdot \bc\hp{1} \Big)}$... \newline
		\mbox{\quad\quad\quad\quad\quad\quad\quad\quad$- \bV\hp{-1}\cdot\big(\bZ\d\cdot(\bY_{:,1:(2\cdot s)} \cdot \bxi)\big)$}
		\State \textit{// - - - 4) Build the next projector - - - }
		\State $\left[\ [\bW_{:,1:(2\cdot s+3)},\bV\hp{0}]\,,\, \bR_W \vphantom{\begin{matrix}
			a\\
			c
			\end{matrix}}\ \right] $ = \texttt{qr}$\Big([\bW_{:,1:(2\cdot s+3)},\bV\hp{0}]\Big)$
		\State $\bQ_G := $\,\texttt{nullbasis}($\bH\t\cdot\bQ_\bY\t$)
		\State \textcolor{blue}{$\bC := \Big( 
		\begin{bmatrix}
		\underline{\underline{\bH}}\\[4pt]
		\bO
		\end{bmatrix} - \tau_1 \cdot
		\begin{bmatrix}
		\underline{\hbH} \cdot \bH\\[4pt]
		\bO
		\end{bmatrix} - \tau_2 \cdot
		\begin{bmatrix}
		\hhbH\\
		\bZ\d \cdot \bY_{:,1:(2 \cdot s + 2)} 
		\end{bmatrix} \cdot (\hbH \cdot \bH)
		\Big)$}
		\State \textcolor{blue}{$\bF := \bR_W \cdot \bC \cdot \bQ_G$}
		\State \textcolor{blue}{$[\bQ_F,\bR_F]$ = \texttt{qr}($\bF$)}
		\State \textcolor{blue}{$[\bL_Z,\bQ_Z]$ = \texttt{lq}\textbf{(}$-\tau_2 \cdot \bY_{:,1:(2 \cdot s + 2)} \cdot \left(\hbH \cdot \big(\bH \cdot (\bQ_G \cdot \bR_F\d)\big)\right)$\textbf{)}}
		\State $\bV\hp{-1} := [\bW_{:,1:(2 \cdot s + 2)},\,\bV\hp{-1}] ...$\newline
		\mbox{\textcolor{blue}{$\quad\quad\quad\quad\cdot \left(\Big( 
		\begin{bmatrix}
		\uubI\\
		-\bZ\d \cdot \bY_{:,1:(2 \cdot s)}
		\end{bmatrix} - \tau_1 \cdot
		\begin{bmatrix}
		\ubH\\
		\bO
		\end{bmatrix} - \tau_2 \cdot
		\begin{bmatrix}
		\hbH \cdot \bH\\
		\bO
		\end{bmatrix}
		\Big) \cdot \big(\bQ_G \cdot (\bR_F\d \cdot \bQ_Z)\big)\right)$}}
		\State $\bV\hp{0} := [\bW_{:,1:(2 \cdot s + 3)},\,\bV\hp{0}] \cdot (\bQ_F \cdot \bQ_Z)$,\quad$\bZ := \bL_Z$
		\State $\br\hp{0} := \bW_{:,1:(2 \cdot s + 3)} \cdot \textcolor{blue}{\big(\underline{\underline{\bc}}\hp{0} - \tau_1 \cdot \underline{\bc}\hp{1} - \tau_2 \cdot \bc\hp{2} \big)}$
		\State \Return $\bZ,\bV\hp{-1},\bV\hp{0},\bx\hp{0},\br\hp{0},\beta$
		\EndProcedure
	\end{algorithmic}
	\caption{Robust cycle for \IDRstabp{$s$}{$2$}}\label{Algo:GMstab2}
\end{algorithm}

\subsubsection{Discussion of the new method}
Again, the new basis vectors $\bW$ for the subsequent projectors $\bV\hp{-1},\,\bV\hp{0}$ are computed by an Arnoldi decomposition. Thus, again we can be sure that the following properties hold:
\begin{itemize}
	\item Either the new column vectors are orthonormal to a high numerical accuracy
	\item or the method terminates with a sufficiently accurate solution.
\end{itemize}
Further, when the method terminates then this termination happens during \textsc{augGMRESm} because the biorthogonal residual in lines 14, 24 is at best only as small as the minimal residual from \textsc{augGMRESm} (Alg. \ref{Algo:augGMRESm} lines 23--24) but never smaller.

Again, there is a further benefit from the interiorly used \GMRES method: During the computation of the matrix-vector products, which dominate the cost of the method, the monitored residual norms decrease monotonously. So one can be sure that not a single matrix-vector product is wasted. This is due to the fact that as soon as there is a residual that satisfies the accuracy requirements the method will immediately terminate.

All the used projectors $\bV\hp{-1},\bV\hp{0}$ and their updates are well-conditioned\footnote{i.e. $\cond(\bV\hp{0})=1$ and $\cond(\bV\hp{-1})\leq\cond(\bA)$}. This is enforced by the QR-decomposition in line 16 and the consistent formulas for $\bV\hp{-1},\bV\hp{0}$ in lines 22--23.
\largeparbreak

Clear drawbacks of the new implementation are that the QR-decomposition in line 16 introduces some additional cost in terms of DOTs and AXPYs (although one notices that since $\bW$ and $\bV\hp{0}$ already have orthonormalised columns one can massively reduce cost by using Householder reflectors to orthogonalise the columns of $\bV\hp{0}$ against $\bW$ without destroying their internal orthonormality). Further to that, the Arnoldi procedure introduces some additional cost for the orthonormalisation of the columns of $\bW$. However, the number of DOTs and AXPYs per matrix-vector products lives in $\cO(s)$. Thus, for small values of $s$ the cost is still dominated by the matrix-vector products with the system matrix and the preconditioners.
\largeparbreak

Drawbacks of \GMstabII compared to the method \GMstabI are the potentially bad condition number of the matrix $\bC$ and the fact that a decoupling of $\bV\hp{-1},\bV\hp{0}$ propagates into the decoupling of their replacements in lines 22--23. In the following we discuss whether these two drawbacks are unavoidable.

From the geometric approach of \IDRstabp{$s$}{$2$} it follows that matrices of potentially the condition number of $\bA^2$ must be dealt with. This is because even if $\bV\hp{0}$ has condition 1 the matrix $\bV\hp{2}$ can have the condition of $\bA^2$. However, since we use an Arnoldi decomposition to generate $\bW$, the expression of $\bV\hp{0}$, which is arguably of the kind $\bV\hp{0} \approx \bW \cdot \bH$, introduces a product with $\bH$. This is why in line 18 the matrix $\bC$ has a triple product of $\bH$, which can have the condition of $\bA^3$.

Since we perform in line 23 the construction of the next projector matrix only from matrices of condition 1, the matrix $\bV\hp{0}$ will still be well-conditioned. Nevertheless, the intermediate round-off amplification can lead to a stronger decoupling of the property
\begin{align*}
	\bV\hp{0} = \bA \cdot \bV\hp{-1}\,.
\end{align*}
Notice that a similar problem is apparent for the updated quantities $\bx\hp{0}$ and $\br\hp{0}$. This is because the computation of $\bc\hp{2}$ from $\bc\hp{0}$ in lines 9--11 is potentially badly conditioned.

An option that can reduce the decoupling is to compute all the blue quantities with a larger mantissa length. Since the blue quantities can be all computed in $\cO(s^3)$ this will only negligibly affect the computational expense of the overall method.

\subsection{Global convergence maintenance through flying restarts and adaptive $\ell \in \lbrace 1,2\rbrace$}

In this subsection we present a strategy for a practical implementation of \GMstab. In this we spend emphasize on numerical accuracy issues that are introduced by numerical round-off.

\subsubsection{Motivation}
In principle one could choose either \GMstabI or \GMstabII and plug it into a while-loop in front of which the initialisation routine is called. As a result of this, one would obtain a complete and functioning implementation of \IDRstabno. However, we believe that this is not sufficient to obtain a \textit{practical} implementation of \IDRstabno for several reasons.
\largeparbreak

First, we have seen that both \GMstabI and \GMstabII update the residual $\br\hp{0}$ and the according numerical solution $\bx\hp{0}$ by uncorrelated equations. Thus, over several iterations, there will be a decoupling of the residual in terms of that the equation
\begin{align*}
	\br\hp{0} = \bb - \bA \cdot \bx\hp{0}
\end{align*}
will only hold with a low numerical accuracy. However, to achieve that a solution $\bx\hp{0}$ can be obtained that meets the required accuracy $\|\bb-\bA \cdot \bx\hp{0}\|\leq\tolabs$ the vector $\br\hp{0}$ must somehow be modified during the iterations to improve the accuracy in which the above equation holds.

This issue will be addressed by using \textit{flying restarts}. This is a particular strategy of recomputing the vector $\br\hp{0}$ to achieve both a fast convergence and a high accuracy for the numerical solution.
\largeparbreak

Second, we have discussed during the motivation of \IDRstab and \BiCGstabL that for some problems the value $\ell=1$ is insufficient for fast convergence. Thus, an implementation that is purely based on \GMstabI is insufficient. On the one hand one needs to utilise \GMstabII to yield fast convergence where $\ell>1$ is required. On the other hand \GMstabII is more affected by numerical round-off and suffers from an amplification in the decoupling of $\bV\hp{-1},\bV\hp{0}$, whereas \GMstabI does not.

Using an adaptive value for $\ell \in \lbrace 1,2\rbrace$, the methods \GMstabI and \GMstabII can be used alternatingly to combine their benefits. This is what we propose to do in a later sub-subsection.

\subsubsection{Structure}
This subsection is organised as follows. In the next sub-subsection we describe the flying restart. This provides a framework of how to replace the residual $\br\hp{0}$ in order to achieve high accuracy of the numerical solution whilst keeping the rate of convergence. Afterwards, there is a sub-subsection that describes the overall algorithmic framework of the new proposed implementation of \GMstab with adaptive $\ell$. Finally, we present and motivate a programmatic rule that chooses the value $\ell$ adaptively for each cycle.

\subsubsection{Flying restarts}\label{sec:FlyingRestarts}

Krylov methods with flying restarts \cite{ResReplacement1,ResReplacement2} are an extension of restarted Krylov methods \cite[p.\,153]{Saad1}. In the following, we first discuss restarted Krylov methods and motivate from that the approach of flying restarts.

\paragraph{Restarts}
The idea of restarting a Krylov method consists of the following approach. Given a linear system
\begin{align*}
	\bb =& \bA \cdot \bx^\star\\
	\bb =& \bA \cdot \bx + \br\,,\tageq\label{eqn:ResidualAccuracy}
\end{align*}
where $\bx^\star$ is the accurate solution and $\bx$ a numerical approximation to it as computed by a Krylov method, one can recompute the residual $\br$ in an accurate way by using the expression
\begin{align}
	\br := \bb - \bA \cdot \bx\,.\label{eqn:ResidualReset}
\end{align}
The residual in turn can be used as a right-hand side of subsequent linear system
\begin{align*}
	\br = \bA \cdot \bdx\,.
\end{align*} 
$\bdx$ in turn can be computed by using a Krylov method, too. If $\bx$ and $\bdx$ can be computed each with a relative residual accuracy of $\tolrel<1$ then $\bx+\bdx$ (if computed without numerical round-off) as an approximation to $\bx^\star$ achieves a relative residual accuracy of $\tolrel^2$. Consequently, the restart improves the achievable accuracy limit.
\largeparbreak

The drawback of restarts is a loss of the superlinear convergence: In the above description it is clear that starting the solution for the remaining residual from scratch discards all the dimensions of the Petrov space $\cC$. (The Petrov spaces are responsible for the superlinear convergence as discussed in Sec.\,\ref{sec:Motivation_KrylovMethods} and along the lines of Fig.\,\ref{fig:SuperlinConv},\,\ref{fig:superlinearConvergenceCurves}\,.) Giving away the superlinear rate of convergence seems numerically inefficient.

\paragraph{Recomputation of the residual}
An alternative to the aforementioned restart is to use the expression \eqref{eqn:ResidualReset} instead of, e.g., Alg.\,\ref{Algo:GMstab1} line 6. However, this leads to a loss of superlinear convergence as well, as we explain in the following:

Consider the case where $\bx$ has a residual norm $\|\bb-\bA \cdot \bx\| \leq \sqrt{\epsmachine}$. Lets say further that $\epsmachine=10^{-16}$. In this case, the expression \eqref{eqn:ResidualReset} yields numerical values for the residual in which at least the last $8$ of the 16 digits consist of pure round-off. As an effect of this round-off, in the general case the new computed vector for $\br$ does not live any more in the Krylov subspace $\cK_\infty(\bA;\bb)$. In consequence, all the \BiCG coefficients in subsequent iterations will be spoilt by the round-off in $\br$. Eventually, this is likely to lead to a loss of the superlinear rate of convergence.

\paragraph{Flying restarts as an engineering solution}
We have explained that the residual recomputation can still spoil the superlinear rate of convergence. The reason was that the recomputed residual introduces numerical round-off that is uncorrelated to the Krylov subspace.

An engineering approach to resolve this issue is to reduce the relative magnitude of this round-off during the recomputation of the residual. This can be achieved by replacing in \eqref{eqn:ResidualReset} the vector $\bb$ by a vector that is closer to $\br$.
\largeparbreak

In \cite{ResReplacement2} the author introduces a flying restart procedure. In this procedure the right-hand side $\bb$ is changed with regard to the current value of $\br$. This change is performed in a particular way such that residual recomputations of the form \eqref{eqn:ResidualReset} do only introduce a small relative amount of numerical round-off.

The idea of the flying restart can be described by considering a global and a local system:
\begin{align*}
	\bb &= \bA \cdot \bx_\text{global} + \bb_\text{local}\\
	\bb_\text{local} &= \bA \cdot \bx_\text{local} + \br 
\end{align*}
The first equation is the global system, where $\bx_\text{global}$ is a numerical solution with a residual $\bb_\text{local}$. The second equation is the local system, where $\bx_\text{local}$ is the numerical solution with the residual $\br$. \underline{The Krylov method is applied to solve the local system.} Residual recomputations of the form \eqref{eqn:ResidualReset} are performed w.r.t. the local right-hand side.

Once that $\|\br\| \leq c_\text{restart} \cdot \|\bb_\text{local}\|$ holds, where $0\leq c_\text{restart} \leq 1$ is a constant, the global and local systems are updated as follows:
\begin{align*}
	\bb_\text{local} &:= \bb_\text{local} - \bA \cdot \bx_\text{local}\\
	\bx_\text{global} &:= \bx_\text{global} + \bx_\text{local}
\end{align*}
This update is called \textit{flying restart}.

The benefit of using $\bb_\text{local}$ is that $\br$ can be recomputed by using $\br := \bb_\text{local} - \bA \cdot \bx_\text{local}$ without introducing a large relative error to $\br$. To explain this in more detail, consider the case $c_\text{restart} = 1$. In this case it holds $\bb_\text{local} \equiv \br$. Thus, the update $\br := \bb_\text{local} - \bA \cdot \bx_\text{local}$ is equivalent to the line 6 in Alg.\,\ref{Algo:GMstab1}. Of this line we know that it does not destroy the superlinear convergence since this is the conventional way of how the residual would be updated. On the other hand, when choosing $c_\text{restart} = 0$ then $\bb_\text{local} \equiv \bb$, thus $\br := \bb_\text{local} - \bA \cdot \bx_\text{local}$ is equivalent to \eqref{eqn:ResidualReset}.

The two extremal values of $c_\text{restart}$ lead either to a perfect maintenance of the superlinear convergence ($c_\text{restart}=1$) or to a perfect accuracy of \eqref{eqn:ResidualAccuracy}. For intermediate values of $c_\text{restart}$ we obtain a compromise of both, instead. A recommended value is $c_\text{restart}=0.01$\,.
\largeparbreak

\begin{figure}
\centering
\includegraphics[width=1\linewidth]{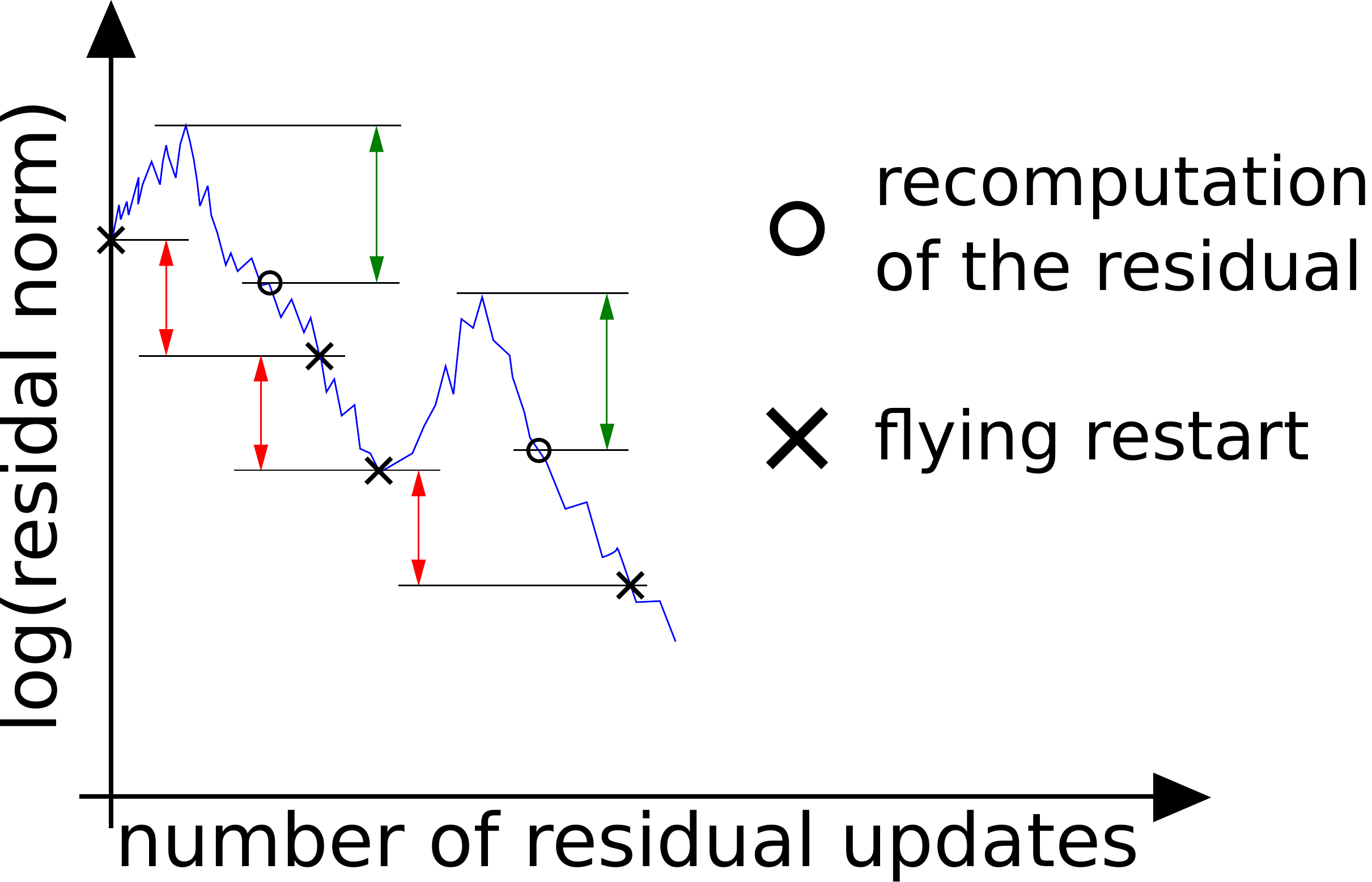}
\caption{Illustration of an exemplary convergence history with residual recomputations and flying restarts.}
\label{fig:FlyingRestart}
\end{figure}

In the following we explain how to implement the flying restart approach by using Fig.\,\ref{fig:FlyingRestart} and the following lines of code.
\begin{algorithmic}[1]
\State $\bx_\text{global}:=\bx_0$,\quad$\bb_\text{local} := \bb-\bA\cdot\bx_0$,\quad$\bx_\text{local} := \bO$
\State \textit{// Initialise the Krylov method, $\beta = \|\br\|$}
\State $\beta_\text{local} := \|\bb_\text{local}\|$,\quad$\beta_\text{max} := \max\lbrace\beta,\,\beta_\text{local}\rbrace$
\While{$\beta > \tolabs$}
	\State \textit{// Perform the iterative scheme of the Krylov method to\newline \mbox{\quad\quad\quad\quad solve for $\bb_\text{local} = \bA \cdot \bx_\text{local}$ with residual $\br$}}
	\If{$\beta > \beta_\text{max}$}
		\State $\beta_\text{max} := \beta$
	\ElsIf{$\beta \leq c_\text{restart} \cdot \beta_\text{local}$}
		\State $\br := \bb_\text{local} - \bA \cdot \bx_\text{local}$,\quad$\beta_\text{max} := \|\br\|$
		\State $\bx_\text{global} := \bx_\text{global} + \bx_\text{local}$,\quad$\bb_\text{local} := \br$,\quad$\bx_\text{local} := \bO$
	\ElsIf{$\beta \leq c_\text{recompute} \cdot \beta_\text{max}$}
		\State $\br := \bb_\text{local} - \bA \cdot \bx_\text{local}$,\quad$\beta_\text{max}:= \|\br\|$
	\EndIf
\EndWhile
\end{algorithmic}
The figure shows an exemplary history of the residual norms as a blue curve. At some points on this curve there are black $\times$ and $\circ$ symbols. A $\times$ is an iterate where a flying restart is performed. A $\circ$ on the other hand is an iterate where a residual recomputation is performed.

As the figure and the code lay out, a flying restart is performed whenever the relative reduction in norm of the residual exceeds a prescribed gap $c_\text{restart}$ w.r.t. the last flying restart (cf. red arrows in the figure and lines 9--10 in the code). Irrespective to that, a residual recomputation is performed whenever $\br$ has experienced a relative reduction in norm that exceeds a prescribed gap $c_\text{recompute}$ (cf. green arrows in the figure and line 12 in the code). Although the figure sketches a scenario where $c_\text{recompute} < c_\text{restart}$ the authors suggest to choose $c_\text{recompute} = c_\text{restart} = 0.01$\,.
\largeparbreak

As the code lays out, the flying restart can be used for any Krylov subspace method that consists of an initialisation phase and an iterative scheme that is repeatedly computed. \IDRstabno and thus \GMstab as a variant fit into this category: The initialsation of \GMstab is given in Alg.\,\ref{Algo:Initialisation} and the iterative scheme is given in Alg.\,\ref{Algo:GMstab1} for \GMstabI respectively Alg.\,\ref{Algo:GMstab2} for \GMstabII. In the next sub-subsection we present an algorithmic framework that uses \GMstab with flying restarts and \textit{both} Alg.\,\ref{Algo:GMstab1} and Alg.\,\ref{Algo:GMstab2}.

\subsubsection{Practical implementation of \GMstab with adaptive $\ell$}

In this sub-subsection we present the practical implementation of our proposed new method \GMstab, which chooses $\ell$ adaptively and uses flying restarts.
\largeparbreak

We have seen that \GMstabII amplifies the decoupling of $\bV\hp{-1},\bV\hp{0}$, whereas \GMstabI resets it. Besides, we have seen that the residual recomputation resets the residual decoupling. Further, we have seen that -- when combined with flying restarts -- the residual recomputation can be used to find a compromise that hopefully yields a small residual gap without spoiling the superlinear convergence behaviour too much.

In the last sub-subsection we have described how flying restarts can be incorporated into a short-recurrence Krylov subspace method. However, what makes \GMstab different from this framework is that it uses distinct iterative schemes, namely one scheme for $\ell=1$ and another one for $\ell=2$.

For our proposed implementation of \GMstab we roll out the flying restarts in a different way. By doing so, we are able to choose the value for $\ell$ depending on whether a restart shall be performed or not.

In the following we first motivate the modified flying restart framework. Then we present the algorithm and explain its computational steps.

\paragraph{Motivation of the modified flying restart}
At some points during the algorithm we want to be able to reset the decoupling of the residual and the auxiliary vectors. For the residual we use a recomputation that is embedded into a flying restart. For the auxiliary vectors we achieve a reset of the decoupling by using a \GMstabI cycle. All in all, the computational strategy is as follows:
\begin{itemize}
	\item Check whether one has to perform a residual recomputation or a flying restart.
	\item If so, then use a \GMstabI cycle and afterwards perform the recomputation or restart.
	\item If not, then use a \GMstabII cycle.
\end{itemize}

However, this strategy is in principle different from the flying restart strategy that was laid out in the former sub-subsection because of the following detail: First, it is analysed whether the residual must be replaced/restarted. Then one still performs an iterative scheme with this residual. Only afterwards the residual is accordingly modified. In the classical flying restart approach instead the residual would be modified immediately after it has been detected that a restart/ recomputation is advised.

What justifies to first compute the cycle and afterwards perform the residual modification is that usually the residual norm converges in a globally smooth way, albeit of some erratic local oscillations.

The reason why we want to compute the cycle of the iterative scheme before the modification (i.e. residual replacement or flying restart) is that in the case of \GMstab we need the auxiliary vectors for this modification of the residual. If now the auxiliary vectors were decoupled then the residual recomputation/ flying restart would become pointless. This is because the decoupling of the auxiliary vectors would lead to a decoupling of the just recomputed residual. Because of that, we first determine whether a recomputation or flying restart must be performed and then choose the cycle accordingly to ensure that $\bV\hp{-1},\bV\hp{0}$ are consistent to a high accuracy whenever a recomputation or restart is actually being computed.

\paragraph{Implementation of the modified flying restart}
The implementation of \GMstab is given in Alg.\,\ref{Algo:GMstab}. In the following we explain this implementation line by line.

In line 2 the flying restart procedure is initialised by setting the right-hand side of the local linear system as the residual of the initial guess. As discussed in Sec.\,\ref{sec:FlyingRestarts} \enquote{Flying restarts as engineering solution} the Krylov method is applied to solve the local system. This is why afterwards in line 3 the initialisation is performed for $\bb_\text{local}$.
\largeparbreak

In line 6 the iterative loop begins. It consists of three phases.

During the first phase only some control parameters are set. These are boolean values and the parameter for $\ell$. $t_\text{restart}$ is true when a flying restart shall be performed and $t_\text{replace}$ indicates that a residual recomputation is required. Lines 8--17 are equivalent to the code fragment that was discussed in Sec.\,\ref{sec:FlyingRestarts} \enquote{Flying restarts as engineering solution}. In lines 18--22 the value for $\ell$ is set to $1$ if the decoupling of the auxiliary vectors shall be reset. This is the case when at least either of the following holds:
\begin{itemize}
	\item A residual recomputation is planned.
	\item A flying restart is planned.
	\item Three consecutive cycles of \GMstabII have been performed. In this case the amplification of the decoupling of the auxiliary vectors might have already become so large that a reset of their decoupling is mandatory.
\end{itemize}

During the second phase the actual iterative scheme is performed. The iterative scheme receives the data $\bA,\bP,\bZ,\bV\hp{-1},\bV\hp{0},\bx\hp{0},\br\hp{0},\beta,\tolabs$ as input and returns updated data for $\bZ,\bV\hp{-1},\bV\hp{0},\bx\hp{0},\br\hp{0},\beta$ as output. For reasons of a compact presentation, in the algorithm the input is abbreviated by \enquote{(...)}.

Finally, during the third phase, the actual flying restart respectively residual recomputation is performed. Since a flying restart implies a recomputation of the residual, the recomputation in line 33 is performed when at least either $t_\text{restart}$ or $t_\text{replace}$ is true. As a speciality, the replacement it followed by a biorthogonalisation of the residual. In theory this was not required since one could achieve by a modified initialisation that $\bb_\text{local}$ always remains in $\cN(\bP)$. However, we believe it is more robust to always re-biorthogonalise the residual in line 35. If a flying restart shall be performed then the local linear system is exchanged in line 38 by defining $\bb_\text{local}$ as the current residual.

It is important that if one of the interior restarted \GMRES subroutines terminates that then its returned solution $\bx\hp{0}$ is the local solution.
In order to obtain the final solution $\bx$ to the original linear system the expression in line 41 must be evaluated.

\begin{algorithm}
	\begin{algorithmic}[1]
		\Procedure{GMstab}{$\bA,\bb,\bx_0,s,\tolabs$}
		\State $\bx_\text{global}:=\bx_0$,\quad$\bb_\text{local} := \bb-\bA\cdot\bx_0$\quad\textit{// $\bx_\text{local} \equiv\bx\hp{0}$}
		\State [$\bP,\bV\hp{-1},\bV\hp{0},\bZ,\bx\hp{0},\br\hp{0},\beta$] = \textsc{Initialisation}($\bA,\bb_\text{local},s,\tolabs$)
		\State $\beta_\text{local} := \|\bb_\text{local}\|$,\quad$\beta_\text{max} := \max\lbrace\beta,\,\beta_\text{local}\rbrace$
		\State $j:=0$,\quad$n_\text{2cycles}:=0$
		\While{$\beta>\tolabs$}
			\State \textit{// - - - 1) Plan flying restart and cycle type - - - }
			\If{$\beta \leq c_\text{restart} \cdot \beta_\text{local}$}
				\State $t_\text{restart} := true$,\quad$t_\text{replace} := false$
			\Else
				\State $t_\text{restart} := false$			
				\If{$\beta < c_\text{recompute} \cdot \beta_\text{max}$}
					\State $t_\text{replace} := true$
				\Else
					\State $t_\text{replace} := false$,\quad $\beta_\text{max} := \max\lbrace\beta,\,\beta_\text{max}\rbrace$
				\EndIf
			\EndIf
			\If{$t_\text{restart}$ \textbf{or} $t_\text{replace}$ \textbf{or} $n_\text{2cycles}>3$}
				\State $\ell:=1$
			\Else
				\State $\ell:=2$
			\EndIf
			\State \textit{// - - - 2) Perform the cycle: Update $\bZ,\bV\hp{-1},\bV\hp{0},\bx\hp{0},\br\hp{0},\beta$ - - - }
			\If{$\ell==1$}
				\State $[\bZ,\bV\hp{-1},\bV\hp{0},\bx\hp{0},\br\hp{0},\beta] = $ \textsc{GMstab1\_cycle}(...)
				\State $n_\text{2cycles}:=0$,\quad$j:=j+1$
			\Else
				\State $[\bZ,\bV\hp{-1},\bV\hp{0},\bx\hp{0},\br\hp{0},\beta] = $ \textsc{GMstab2\_cycle}(...)
				\State $n_\text{2cycles}:=n_\text{2cycles}+1$,\quad$j:=j+2$
			\EndIf
			\State \textit{// - - - 3) Perform the flying restart - - - }
			\If{$t_\text{restart}$ \textbf{or} $t_\text{replace}$}
				\State $\br := \bb_\text{local} - \bA \cdot \bx\hp{0}$
				\State $\bEta := \bP\t \cdot \br\hp{0}$
				\State [$\bx\hp{0},\br\hp{0}$] = \textsc{dirRbio}($\bV\hp{-1},\bV\hp{0},\bZ,\bx\hp{0},\br\hp{0},\bEta$),\quad$\beta:=\|\br\hp{0}\|$
			\EndIf
			\If{$t_\text{restart}$}
				\State $\bb_\text{local} := \br\hp{0}$,\quad$\bx_\text{global}:=\bx_\text{global}+\bx\hp{0}$,\quad$\bx\hp{0}:=\bO$,\quad$\beta_\text{local}:=\beta$
			\EndIf
		\EndWhile
		\State $\bx := \bx_\text{global} + \bx\hp{0}$\quad\textit{// final numerical solution}
		\State \Return $\bx$
		\EndProcedure
	\end{algorithmic}
	\caption{Practical implementation of \GMstab}\label{Algo:GMstab}
\end{algorithm}

\subsection{Discussion of \GMstab}
In this subsection we briefly summarise the ingredients of our new method. Then we discuss the advantages and drawbacks that these ingredients bring to the algorithm.
\largeparbreak

We have developed a new implementation of \IDRstabno that is not based on \GCR but on an interior restarted \GMRES approach. As an effect of this there are no breakdowns in our new method when the columns of $\bW$ become collinear. Instead, the method will converge faster and eventually terminates with a residual-optimal solution w.r.t. the current projectors. In all the other \IDRstabno implementations, since they follow the \GCR approach, there is a breakdown whenever $\bv_q \in \opspan\lbrace\bv_1,...,\bv_{q-1}\rbrace$ occurs: In the non-biortho variants the matrix $\bZ$ becomes singular in this case and in the biortho variant the residual does not change, which leads to collinearity of all the subsequent auxiliary vectors $\bv_{i}$, $i=q+1,...,s$\,.

Besides, our new implementation guarantees that the projector $\bV\hp{0}$ is always unitary and that $\bZ \equiv \bP\h \cdot \bV\hp{0}$ is lower triangular. Consequently, the condition of $\bZ$ is not spoilt in our method by a potentially bad conditioning of $\bV\hp{0}$ as it can be in contrast the case in \IDRbiortho and \IDRstabbiortho. Since $\bZ$ is still lower triangular in our method all linear systems for the biorthogonalsations can be solved in a robust and efficient way.

Further, the new implementation utilises a modified flying restart approach in which the value for $\ell$ is chosen adaptively. The adaptive value for $\ell$ guarantees that the equation $\bV\hp{0}=\bA\cdot\bV\hp{-1}$ always holds to a high accuracy. It further makes sure that numerical round-off in this equation cannot amplify over a long sequence of iterations. This is realised by using occasional cycles of \GMstabI, which act as a reset to the round-off in this equation. Combined with this kind of reset, the flying restarts on the other hand provide an effective handling of the residual decoupling in that they treat two issues: First, it is ensured that the equation $\br\hp{0}=\bb-\bA\cdot\bx\hp{0}$ always holds with a high accuracy. Second, the superlinear convergence is only affected mildly\footnote{This is what we hope. However, in the numerical experiment \textbf{Norris\_torso1} that the flying restarts can still slow down the rate of convergence significantly.} by the recomputations of $\br\hp{0}$.

Finally, our new implementation does not require any additional storage in terms of the number of vectors of length $N$ compared to the so far existing implementations of \IDRstabno.
\largeparbreak

Considering computational cost, our method is clearly more expensive than the cheapest implementation of \IDRstabno, namely the reference \IDRstab. The additional orthogonalisations during the Arnoldi procedures in \textsc{pGMRESm} respectively \textsc{augGMRESm} and the additional QR-decomposition for $[\bV\hp{0},\,\bV\hp{-1}]$ in \GMstabI respectively for $[\bV\hp{0},\,\bW]$ in \GMstabII require $\cO(s)$ additional DOTs and AXPYs per matrix-vector product.

\largeparbreak

At the end of the day the question is always whether the extra cost pays off: Can the number of required matrix-vector products be reduced using our new method because of a better maintenance of superlinear convergence? Or can higher solution accuracies be achieved without a restart? It is possible to solve more difficult problems successfully using our method than with other \IDRstabno implementations? In Sec.\,\ref{sec:NumericalExperiments} we will consider different test problems to find answers to these questions.

\subsection{Comparison of the computational costs}
In the following we provide tables that show the cost in terms of dot-products (DOT) and vector updates (AXPY, acronym for \enquote{$\by := \alpha \cdot \bx + \by$}) for vectors of length $N$, i.e. the system dimension. We show the costs for all the algorithmic sub-blocks and for all \IDR and \IDRstab implementations that we have discussed so far. 

We first present the cost of the algorithmic sub-blocks in Tab.\,\ref{tab:CostSubBlocks}. Presenting the costs of the sub-blocks has mainly two purposes: The user can read the cost that is introduced when replacing one sub-block by another. For example, replacing \textsc{itVbio} by \textsc{itVobio} introduces $q$ additional DOTs in the $q$th interior for-loop of the biorthogonalisation step in \IDRstabno. Second, the formulas from the table help verifying our found costs for the \IDR and \IDRstab implementations and keep the results more transparent. The latter is important since some of our results contradict to propositions from \cite{IDR-biorth} and \cite{IDRstab-biorth}. Namely, the authors propose that the biortho variants are cheaper than the reference implementations, which is not the case as the next table shows:

\begin{table}
	\centering
	\rotatebox{90}{
		\begin{tabular}{||c||c|c|c||}
			\hline
			\hline 		
			       	Name 								& MATVEC 	& DOT 												& AXPY 								\\ 
			\hline
			\hline 
				   	\textsc{Initialisation}  			& $s$  		& $2\,s^2+3\,s+1\vphantom{\Big[\Big]}$ 				& $3\,s^2 +5\,s+1$ 					\\
			       	\textsc{StabCoeffs}   				& $0$  		& $0.5\,(\ell+1)\,(\ell+2)-1\vphantom{\Big[\Big]}$ 	& $0$ 								\\ 
			\hline 
				   	\textsc{QR}($m$) 					& $0$  		& $0.5\,m\,(m+1)\vphantom{\Big[\Big]}$  			& $0.5\,m\,(m+1)$ 					\\
			\hline
				   	\textsc{dirRbio}($k$) 				& $0$  		& $0\vphantom{\Big[\Big]}$ 							& $(k+3)\,s$						\\
				   	\textsc{itRbio}($k$) 				& $0$  		& $0\vphantom{\Big[\Big]}$ 							& $(k+3)$							\\
				   	\textsc{itVbio}($q,k$) 				& $0$  		& $s\vphantom{\Big[\Big]}$ 							& $(k+3)\,(q-1)$					\\
				   	\textsc{itVorth}($q,k$) 			& $0$  		& $q\vphantom{\Big[\Big]}$ 							& $(k+3)\,q$						\\
				   	\textsc{itVobio}($q,k$) 			& $0$  		& $s+q\vphantom{\Big[\Big]}$ 						& $(k+3)\,(3\,q-2)$					\\ 
			\hline
					\textsc{GMRES}($m$) 				& $m$ 		& $0.5\,m^2 + 1.5\,m\vphantom{\Big[\Big]}$ 			& $0.5\,m^2 + 2.5\,m$ 				\\
					\textsc{pGMRES}($m$) 				& $m$ 		& $0.5\,m\,(m+3)+(m+1)\,s\vphantom{\Big[\Big]}$ 	& $0.5\,m\,(m+3)+m\,s$ 				\\
					\textsc{augGMRES}($m$) 				& $m$ 		& $0.5\,m\,(m+3)+m\,s\vphantom{\Big[\Big]}$ 		& $0.5\,m\,(m+3)+m\,s$ 				\\				
			\hline
			\hline
		\end{tabular}
	}
	\caption{Cost of all algorithmic sub-blocks. All these methods work in-place on the array space that they are given by the outer routine that they are called from.}\label{tab:CostSubBlocks}
\end{table}

Tab.\,\ref{tab:CostIDRmethods} shows the cost of the reference, biortho, ortho-biortho (obio) and restarted GMRES (GMstab) implementations of \IDR and \IDRstab. This cost does not include the additional cost for the treatment of $\bZ$ that is required only in \textsf{IDR($s$)stab($\ell$)ref}, \textsf{IDR($s$)stab($\ell$)biorth} and \textsf{IDR($s$)stab($\ell$)obio} to maintain the convergence.

From the table we make the following observations: Whereas the biortho variants are roughly as cheap as the reference implementations, the obio variant has about twice the cost in AXPYs and $1.5$ times as much cost in DOTs. \textsf{GM($s$)stab1} is not only twice but three times as expensive as \textsf{IDR($s$)noDec}, not only in AXPYs but also in DOTs. \textsf{GM($s$)stab2} instead is only about $2.5$ times as expensive as \textsf{IDR($s$)stab($2$)} in DOTs and AXPYs. However, \textsf{GM($s$)stab2} does not need a treatment for $\bZ$ whereas all the other listed \IDRstabno implementations have some additional costs to recompute $\bZ$ that is not considered in the table.

From the numerical experiments in subsequent sections we can compare the costs of the methods indirectly by comparing their runtime: In Fig.\,\ref{fig:BiCGstabL_xpl1_iterres} all the \IDRstabno variants have the same number of matrix-vector products. For $s=4$ and $\ell=2$ the implementations \IDRstabbiortho with $\bZ$-treatment and \GMstab with flying restarts require about twice as much runtime per iteration as the reference implementation of \IDRstabno. The obio variant with $\bZ$-treatment however requires three times as much runtime as \IDRstabno.
\largeparbreak

In Tab.\,\ref{tab:CostIDRmethods} no column space was left to present the number of matrix-vector products (MATVEC) and and memory consumption in the number of stored column vectors of length $N$ that the \IDRO methods require in each cycle. 

For all \IDR variants, i.e. also \GMstabI, the number of matrix-vector products is $(s+1)$ and the number of stored column vectors is $(s+1)\cdot 3$. In contrast to this, for all \IDRstab variants, i.e. also \GMstabII, the number of matrix-vector products is $\ell\cdot(s+1)$ and the number of stored column vectors is $(\ell+3)\cdot(s+1)-1$.

\begin{table}
	\centering
	\rotatebox{90}{
		\begin{tabular}{||c||c|c||}
			\hline
			\hline
					Name 									& DOT 															& AXPY			 														\\
			\hline
			\hline
					\textsc{IDR($s$)noDec} 					& $s^2 + s + 3$ 												& $2\,s^2 + 2\,s + 2$ 													\\
					\textsc{IDR($s$)biortho} 				& $s^2 + s + 3$ 												& $2\,s^2 + 2\,s + 2$ 													\\
					\textsc{IDR($s$)obio} 					& $1.5\,s^2 + 1.5\,s + 3$ 										& $4\,s^2 + 2\,s + 2$ 													\\
			\hline
					\textsc{IDR($s$)stab($\ell$)ref} 		& $\phantom{1.5\,}\ell\,(s^2+s)+0.5\,(\ell+1)\,(\ell+2)$ 		& $0.5\,s^2 \,\ell^2+1.5\,s^2\,\ell+0.5\,s\,\ell^2+2.5\,s\,\ell+4\,s+4$	\\
					\textsc{IDR($s$)stab($\ell$)biortho}	& $\phantom{1.5\,}\ell\,(s^2+s)+0.5\,\ell\,(\ell+3)+1$			& $0.5\,s^2 \,\ell^2+2.0\,s^2\,\ell+0.5\,s\,\ell^2+2.0\,s\,\ell+4\,s+4$	\\
					\textsc{IDR($s$)stab($\ell$)obio} 		& $1.5\,\ell\,(s^2+s)+0.5\,\ell\,(\ell+3)+1$					& $1.0\,s^2 \,\ell^2+4.5\,s^2\,\ell+0.5\,s\,\ell^2+2.5\,s\,\ell+4\,s+4$	\\
			\hline
					\textsc{GM($s$)stab1} 					& $3.0\,s^2 + 3.0\,s + 3$ 										& $6.5\,s^2 + 4.5\,s + 3$												\\
					\textsc{GM($s$)stab2} 					& $5.5\,s^2 + 14.5\,s + 8$ 										& $10.5\,s^2 + 22.5\,s + 10$											\\
			\hline
			\hline
		\end{tabular}
	}
	\caption{Cost per repetition of the main while-loop of various implementations of \IDRO methods.}\label{tab:CostIDRmethods}
\end{table}

\section{\Mstabno: \IDRstabno with Krylov subspace recycling}\label{sec:Mstab}
The original motivation for the novel \IDRstabno implementation that is presented in this thesis was to have a proper implementation to perform the numerical experiments in \cite{Mstab-paper}. It turns out that robustness and convergence maintenance have a more crucial influence on the convergence behaviour when a Krylov subspace recycling variant of \IDRstabno for the solution of multiple linear systems is used. Due to this reason we want to show numerical experiments where \GMstab is applied for solving sequences of linear systems by using it as a Krylov subspace recycling method.

For the sake of a self-contained presentation, before we discuss numerical experiments of this kind in the subsequent section, this section describes in advance the mathematical and algorithmic approach of \IDRstabno with Krylov subspace recycling. However, this presentation gives only a coarse and brief description of the matter since Krylov subspace recycling is not the main topic of this thesis. For further information on the contents presented in this section we refer to our publication \cite{Mstab-paper}.

\subsection{Review of Krylov subspace recycling}
Krylov subspace recycling is a mathematical approach that is used to solve sequences of linear systems
\begin{align*}
	\bA \cdot \bx\hia = \bb\hia,\quad \iota=1,...,\nEqns\,,
\end{align*}
where the systems must be solved in alphanumeric order and one after the other. This is for instance the case when an implicit time-stepping scheme is applied to solve a discretised partial differential equation (PDE) or when a Quasi-Newton method is used.

Whereas in principle each of the $\nEqns \in \N$ systems can be solved individually using a conventional Krylov subspace method, the goal of Krylov subspace recycling methods is to solve subsequent linear systems faster.

In Sec.\,\ref{sec:MotivationKSSR} we have sketched a motivation for Krylov subspace recycling. In the following we review this.
\largeparbreak

Having solved the first linear system $\bA \cdot \bx\hp{1}= \bb\hp{1}$, there is data available in $\cK(\bA;\bb\hp{1})$ that could help solving the second system $\bA\cdot\bx\hp{2}=\bb\hp{2}$ with lower computational efforts.

In general, one tries to reduce the computational cost by improving the rate of convergence of the Krylov method. Fig.\,\ref{fig:superlinearConvergenceCurves} depicts the intended convergence behaviour of Krylov subspace recycling methods in comparison to Krylov subspace methods: Recycling some data from $\cK(\bA;\bb\hp{1})$, the Krylov subspace recycling method shall convergence sooner for the second linear system $\bA\cdot\bx\hp{2}=\bb\hp{2}$.

There are two unrelated ways of how to achieve this faster convergence:
\begin{enumerate}[(A)]
	\item The data from $\cK(\bA;\bb\hp{1})$ can be used to precondition the system matrix $\bA$, for instance by deflation, cf. \cite{RGMRES,RGCRO,GCRO-DR,RBiCGstab}.
	\item The data from $\cK(\bA;\bb\hp{1})$ can be utilised in such a way that the Krylov subspace of the second linear system becomes smaller. This is achieved in the method \Mstab \cite{Mstab-report,Mstab-paper}, which is the Krylov subspace recycling variant of \IDRstab. Another way to achieve this is given by \textit{short representations} \cite{Report15,SRPCR}.
\end{enumerate}

Before we review in the next subsection the motivation of Krylov subspace recycling for \IDRO methods we briefly explain the concept of deflation.
\largeparbreak

Deflation works as follows: One considers the situation that $\bA \cdot \bx\ha = \bb\ha$ has been solved using a Krylov method from which on the fly some Krylov vectors $\bU \in \C^{N \times k}$ and their images $\bC \in \C^{N \times k}$ for some small $k \in \N$, e.g. $k=20$, have been extracted:
\begin{align*}
	\bA \cdot \bU = \bC
\end{align*}
The information from this equation can be useful for solving subsequent linear systems like $\bA \cdot \bx\hb = \bb\hb$. Instead of solving this system conventionally by applying a Krylov method to it, one can solve instead the system
\begin{align}
	\underbrace{(\bI - \bC \cdot \bC\d)\cdot \bA}_{=:\tbA} \cdot \tbx\hp{2} =(\bI - \bC \cdot \bC\d)\cdot \bb\hp{2}\label{eqn:DeflationAxb}
\end{align}
for $\tbx\hp{2}$. Afterwards, the original solution $\bx\hp{2}$ can be reconstructed from $\tbx\hp{2}$ by using the expression
\begin{align*}
	\bx\hp{2} := \tbx\hp{2} + \bU \cdot \bC\d \cdot (\bb - \bA \cdot \tbx\hp{2})\,.
\end{align*}
When for instance $\bU$ approximates an invariant subspace of some outlying eigenvalues (e.g. close to zero) then the matrix $\tbA$ offers a faster rate of convergence than $\bA$. Thus, solving the system \eqref{eqn:DeflationAxb} can be cheaper than solving $\bA \cdot \bx\hp{2}=\bb\hp{2}$.
\largeparbreak

When evaluating Krylov subspace recycling it is important to consider computational overhead. 

We have said that solving \eqref{eqn:DeflationAxb} can lead to a smaller number of iterations. However, the computational cost per iteration is larger for \eqref{eqn:DeflationAxb} than for $\bA \cdot \bx\hp{2}=\bb\hp{2}$. This is because matrix-vector products with $\tbA$ require some overhead compared to those with $\bA$. Thus, one must be sure in advance that the overhead that is introduced by using deflation will massively reduce the iteration count because otherwise it will not pay off.

Unfortunately, deflation is strongly dependent on the spectrum of the linear system. For instance, it can be the case that there are no well-separated eigenvalues. In this and many other situations the deflation brings no improvement to the rate of convergence at all.
\largeparbreak

In the following subsection we describe the second approach (B). As a consequence of its construction, the convergence improvement that can be obtained from this approach does not depend on spectral properties of the system matrix.

\subsection{Idea and motivation of \Mstabno}
During the computation of a solution to $\bA \cdot \bx\hp{1} = \bb\hp{1}$ with \IDRstab it is possible to fetch some vectors $\hbU \in \C^{N \times s}$ from $\cK(\bA;\bb\hp{1})$ on the fly such that the second linear system $\bA \cdot \bx\hp{2}=\bb\hp{2}$ can be solved by searching a solution in a massively reduced (in terms of dimensions) solution space.
\largeparbreak

In Fig.\,\ref{fig:MstabMotivation} we have sketched this for a system size of $N=200$ and a method parameter $s=10$: 

The first linear system is solved conventionally with \IDRstabno, where $\dim(\cK_\infty(\bA;\bb\hp{1}))=200$. The grey numbers give the dimension of the Sonneveld space in which each respective residual on the convergence graph lives. During the computation of the solution some data $\hbU$ is written out. We see that as the spaces become smaller there is eventually a point where the dimension is so small that suddenly a fast rate of convergence occurs.

Using the data $\hbU \in \C^{200 \times 10}$, the solution for the second linear system can be computed from a space that has only dimension $155$, as can be proven \cite[p.\,13]{Mstab-paper}. Starting from that smaller dimensioned space, it is likely that the turning point to superlinear convergence is reached earlier, as is sketched in the figure.

\subsection{Mathematical approach of Krylov subspace recycling for \IDRO methods}

How can it be possible at all that by only using the $10$ vectors $\hbU \in \C^{200 \times 10}$ the problem space of the second linear system can be reduced from $200$ dimensions to $155$? In order to sketch this, we want to lay out in this subsection the geometric idea that underlies the theory of \Mstab.
\largeparbreak

So far we have discussed Sonneveld spaces $\cG_0,\cG_1,\cG_2,...$\,. They followed the recursion
\begin{align*}
	\cG_{\mathrlap{0}\phantom{j+1}} &= \cK_\infty(\bA;\bb\hp{1})\,,\\
	\cG_{j+1} &= (\bI - \omega_{j+1} \cdot \bA) \cdot \big(\cG_j \cap \cN(\bP)\big)\quad \forall j \in \No\,.
\end{align*}
In order to use these spaces in a short-recurrence Krylov method, it was crucial that $\cG_{j+1} \subset \cG_j$ held because otherwise one would run out of vectors for the biorthogonalisation of subsequent vectors once after the first vector has been moved from $\cG_j$ into $\cG_{j+1}$. The second crucial property was the dimension reduction of the Sonneveld spaces since this led to the superlinear rate of convergence.

We have given arguments that with an increasing degree $j$ of the Sonneveld spaces the methods converge faster and faster. We have further explained how this is related to the dimension reduction of the Sonneveld spaces.
\largeparbreak

For Krylov subspace recycling within the \IDRO theory the mathematical approach is that for some large value of $j$ the Sonneveld spaces can be reused in the following way: If the residual $\br$ of the second linear system $\bA \cdot \bx\hp{2} = \bb\hp{2}$ lived in a Sonneveld space $\cG_j$ of large degree than we could simply use the iterative scheme of \IDRstab to restrict it iteratively into further successors of $\cG_j$. Since $\cG_j$ is already of small dimension we could be sure that the residual of the second system converged very quickly.

The question is of course if and how we can achieve that a property similar to $\bb\hp{2} \in \cG_J$ for some large $J \in \N$ can be achieved. In the next subsection we introduce a generalisation of Sonneveld spaces that gives an answer to this question.

\subsection{Mathematical theory of \Mstabno: $\cM$-spaces}

In this subsection we introduce only a particular variant of $\cM$-spaces that is tailored for the application of solving two subsequent linear systems $\bA \cdot \bx\ha = \bb\ha$ and $\bA\cdot\bx\hb=\bb\hb$.

\begin{Definition}[$\cM$-space variant]\label{def:Mspace}
Given $\bA \in \C^{N \times N}$ regular, $\bb\hp{1} \in \C^N$, $\bP \in \C^{N \times s}$ with $\rk(\bP)=s$, $\lbrace \omega_j \rbrace_{j \in \N} \subset \Cno$. Using Def.\,\ref{def:Sonneveldspace} we obtain a sequence of Sonneveld spaces $\lbrace\cG_j\rbrace_{j \in\No}$ for the given data.

Given further $\bb\hp{2}\in\C^N$ and $J \in \N$ such that $\dim(\cG_J) \geq s$. We define the following vector spaces:
\begin{align*}
	\cM_{\mathrlap{0}\phantom{j+1}} &:= \cK_\infty(\bA;\bb\hp{1}) + \cK_\infty(\bA;\bb\hp{2}) & &\\
	\cM_{j+1} &:= (\bI - \omega_{j+1}\cdot\bA) \cdot \big(\cM_j \cap \cN(\bP)\big) + \opspan\lbrace\bb\hp{2}\rbrace,\quad & j&=0,...,J-1\,,\\
	\cM_{j+1} &:= (\bI - \omega_{j+1}\cdot\bA) \cdot \big(\cM_j \cap \cN(\bP)\big),\quad & j&=J,J+1,J+2,...\,.
\end{align*}
The space $\cM_{j}$ is called $\cM$-space of degree $j$.
\end{Definition}

The $\cM$-spaces are nested and reduce their dimension for increasing degree $j$, just like Sonneveld spaces do. This is proven in a subsequent theorem.

$\cM$-spaces can be utilised to solve the two linear systems in a Krylov subspace recycling method that follows the approach described in the former subsection. As we will show, the properties $\cG_J \subset \cM_J$ and $\bb\hp{2} \in \cM_J$ hold. 

The first property is useful because it says that \IDRstab{}'s $s$ auxiliary vectors from $\cG_J$ (from the solution process for $\bA\cdot\bx\ha=\bb\ha$) do also live in $\cM_J$. Further it holds $\bb\hp{2}=:\br \in \cM_J$. So everything is available (namely $s$ auxiliary vectors and a residual in $\cM_J$) to launch \IDRstab to iteratively move $\br$ from $\cM_J$ into its successors $\cM_J\supset\cM_{J+1}\supset\cM_{J+2}...$\,.
\largeparbreak

\begin{figure}
	\centering
	\includegraphics[width=1\linewidth]{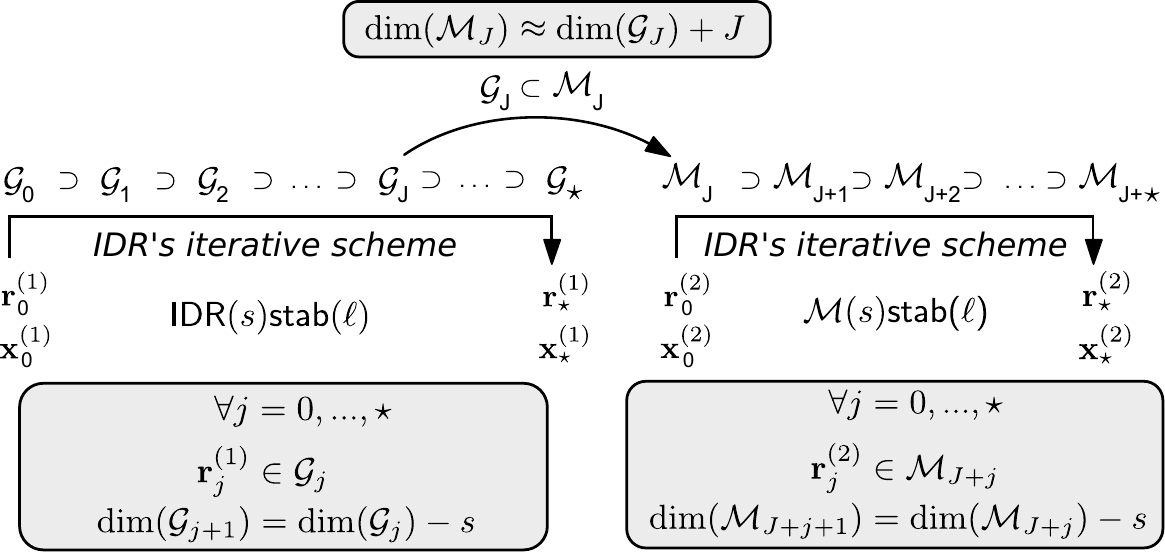}
	\caption{Algorithmic approach of Krylov subspace recycling for \IDRstabno.}
	\label{fig:MstabPrinciple}
\end{figure}

We explain the algorithmic approach in more detail, using Fig.\,\ref{fig:MstabPrinciple}. The figure illustrates the algorithmic approach for Krylov subspace recycling for \IDRstab applied to the two linear systems:

First, on the left-hand side of the figure, the linear system $\bA \cdot \bx\hp{1} = \bb\hp{1}$ is solved using \IDRstab. At some point during the solution procedure there are auxiliary vectors $\bV\hp{-1},\bV\hp{0} \in \C^{N \times s}$ for which $\rg(\bV\hp{0})\subset\cG_J$ holds for some value $J \in \N$ that can be chosen arbitrarily. At this moment we save the matrix $\hbU := \bV\hp{-1}$ on the hard-drive of the computer. Then \IDRstab proceeds until the first system has been solved successfully.

Afterwards, the second system shall be solved. For this we load the data $\hbU$ from memory and initialise the projectors $\bV\hp{-1}:=\hbU$, $\bV\hp{0} := \bA \cdot \bV\hp{-1}$ and the numerical solution $\bx := \bO$ and residual $\br := \bb\hp{2}$. Having done this, it holds
\begin{align*}
	\rg(\bV\hp{0}) &\subset \cM_J\,,\\
	\br &\in \cM_J\,.
\end{align*}
As the figure indicates, the same iterative scheme as is used in \IDRstab can be utilised to iteratively restrict the residual and the auxiliary vectors from $\cM_J$ into subsequent $\cM$-spaces. For a distinction by name, when applying \IDRstab to iterate over $\cM$-spaces we call the method \Mstab.

The benefit of \Mstabno over \IDRstabno can be derived from the grey boxes in the figure: It is likely that the $\cM$-spaces are of much smaller dimension. Thus, it is likely that \Mstabno converges much faster than \IDRstabno. Not only is it likely that the method converges faster. Further, one can prove that finite termination of \Mstab is achieved $s-1$ times earlier than for \IDRstab when choosing $J$ appropriately, cf. \cite[p.\,13]{Mstab-paper}.

\largeparbreak
In the remainder we formulate the above postulated theoretical properties of $\cM$-spaces as a theorem.

\begin{Theorem}[Selected properties of some particular $\cM$-spaces]\label{Theo:Mspace}
	Consider the Sonneveld spaces $\cG_0,\cG_1,...,\cG_J$ that arise from $\bA$, $\br\ha_0$, $\bP$, $\lbrace\omega_j\rbrace_{j\in\N}$, and the $\cM$-spaces that arise from $\bA$, $\br\ha_0$, $\br_0\hb$, $\bP$, $\lbrace\omega_j\rbrace_{j\in\N}$, $J$. If $\rg(\bP)$ and $\cM_0$ do not share a non-trivial invariant subspace of $\bA$ then it holds:
	\begin{enumerate}[(a)]
		\item $\cM_j \subset \cM_{j-1}$\quad $\forall j \in \N$
		\item $\dim(\cM_j) \leq \max\lbrace 0,\,\dim(\cM_{j-1})-(s-1)\rbrace$ \quad $\forall j\in\lbrace0,\,1,\,2,\,...,\,J\rbrace$\\
		$\dim(\cM_j) \leq \max\lbrace 0,\,\dim(\cM_{j-1})-s\rbrace$ \quad $\forall j=J+1,\,J+2,\,...$
		\item $\cG_J \subset \cM_J$\,,\ $\br_0\hb \in \cM_J$
	\end{enumerate}
	\underline{Proofs:}\\
	For ease of notation we use in the following proofs that the constituting recursion for the $\cM$-spaces can be expressed as
	\begin{align}
	\cM_j = (\bI-\omega_j \cdot \bA) \cdot \big(\cM_{j-1} \cap \rg(\bP)^\perp\big) + \cQ_{j} \quad \forall j\in\N \label{eqn:RecM}
	\end{align}
	with the following spaces for $\lbrace\cQ_j\rbrace_{j\in\N}$.
	\begin{align*}
	\cQ_{j} :=
	\begin{cases}
	\opspan\lbrace\br_0\hb\rbrace &\text{for}\ \ j\leq J\,,\\
	\lbrace\bO\rbrace &\text{for}\ \ j > J\,.
	\end{cases}
	\end{align*}
	Further, we define $\chi_j:=1$ for $j\leq J$ and $\chi_j:=0$ for $j> J$, from which follows $\dim(\cQ_j)=\chi_j$ $\forall j \in \N$. We further use $\bQ_j$ as a basis matrix of $\cQ_j$ and $\cP := \rg(\bP)$.\\[7pt]
	Proposition (a). Proof by induction.
	\begin{enumerate}
		\item Basis: Since $\cM_0$ is the sum of two full Krylov spaces, it holds $\bA \cdot \cM_0 \subset \cM_0$. Thus, $\cM_1 \subset \cM_0 + \bA \cdot \cM_0 = \cM_0$.
		\item Hypothesis: $\cM_j \subset \cM_{j-1}$ holds for a specific $j \in \N$.
		\item Induction step: \enquote{$\cM_{j} \subset \cM_{j-1}$ $\Rightarrow$ $\cM_{j+1} \subset \cM_j$} is shown.\\
		Choose an arbitrary $\bx \in \cM_{j+1}$. Then
		\begin{align*}
		\exists \by \in \cM_j \cap \cP^\perp \ \wedge \ \bq \in \cQ_{j+1} \ : \ \bx = (\bI - \omega_{j+1} \cdot \bA)\cdot \by + \bq\,.
		\end{align*}
		From $\cM_j \subset \cM_{j-1}$ follows $\cM_j \cap \cP^\perp \subset \cM_{j-1}\cap\cP^\perp$ and from the latter in turn $\by \in \cM_{j-1}\cap\cP^\perp$, from which we can construct a vector $\tbx$ in $\cM_j$ as follows:
		\begin{align*}
		\tbx := (\bI-\omega_j \cdot \bA) \cdot \by \in \cM_j\,.
		\end{align*}
		Since $\bq \in \cQ_{j+1} \subset \cQ_j \subset \cM_j$, it follows $\by,\tbx,\bq \in \cM_j$. Since $\bx \in \opspan\lbrace\by,\tbx,\bq\rbrace \subset \cM_j$ holds for an arbitrary $\bx \in \cM_{j+1}$ it is shown that $\cM_{j+1} \subset \cM_{j}$.
	\end{enumerate}
	Proposition (b).
	\begin{enumerate}
		\item[] We state the following proposition:
		\begin{align*}
		\dim(\cM_{j+1}) \leq \dim(\cM_j) - \dim(\cP) + \dim(\cQ_j)\quad \forall j \in \No\,.
		\end{align*}
		\item[] When the proposition holds then (b) follows in immediate consequence.
		\item[] We show the above proposition by using Lemma 4, which is given below. In this we consider two cases for the range of the degree of the $\cM$-space due to the definition of $\lbrace\cQ_j\rbrace_{j\in\No}$.
		\item The case $\cM_j$ for $j \in \lbrace0,...,J\rbrace$\,:
		\item[] Using Lemma 4, it holds
		\begin{align*}
		\cM_{j} \subset p_{0,j}(\bA) \cdot \big(\cM_0 \cap \cK_j(\bA^\textsf{H};\bP)\big) + \cK_j(\bA;\br^{(0)})\,.
		\end{align*}
		It follows
		\begin{align*}
		\dim(\cM_{j}) \leq \underbrace{\dim\big(\cM_0 \cap \cK^\perp_j(\bA^\textsf{H};\bP)\big)}_{{}^{(*)}=\max\lbrace 0, \dim(\cM_0) - j \cdot \dim(\cP)\rbrace } + \underbrace{\cK_j(\bA;\br^{(0)})}_{\leq j}\,.
		\end{align*}
		In this the equation (*) holds under the mild requirement that $\cP$ and $\cM_0$ do not share a non-trivial invariant subspace of $\bA$, cf. \cite{IDR-Gutkn}.
		\item[] From this follows the above proposition for $j<J$.
		\item The case $\cM_{J+j}$ for $j\in\N$\,:
		\item[] Using Lemma 4, it holds
		\begin{align*}
		\cM_{J+j} \subset p_{J,j}(\bA) \cdot \big(\cM_J \cap \cK^\perp_j(\bA^\textsf{H};\bP)\big) + \cK_j(\bA;\bO)\,.
		\end{align*}
		It follows
		\begin{align*}
		\dim(\cM_{J+j}) \leq \underbrace{\dim\big(\cM_J \cap \cK_j(\bA^\textsf{H};\bP)\big)}_{{}^{(**)}=\max\lbrace 0, \dim(\cM_J) - j \cdot \dim(\cP)\rbrace}\,.
		\end{align*}
		In this the equation (**) holds analogously when $\rg(\bP)$ and $\cM_J$ do not share a non-trivial invariant subspace of $\bA$. Since $\cM_J \subset \cM_0$ holds the mild conditions on $\cP$ for (*) to hold, which are commonly used for Sonneveld spaces, do imply that equation (**) holds, too.
		\item[] From this follows the above proposition for $j\geq J$. In total, proposition (b) has been shown.
	\end{enumerate}
	Proposition (c) follows by induction over the degree of the two spaces $\cG_j,\cM_j$ for $j=0,...,J$ and is left to the reader $\boxtimes$.
\end{Theorem}

\begin{Lemma}[A superspace of $\cM$-spaces]
	Consider the $\cM$-spaces from Definition\,2. It holds $\forall j,d \in \N$:
	\begin{align*}
	\cM_{j+d} \subset \underbrace{\left(\prod_{k=1}^d (\bI - \omega_{j+d}\cdot\bA)\right)}_{=:p_{j,d}(\bA)}\cdot\left(\cM_j \cap \cK_d^\perp(\bA^\textsf{H};\bP) \right) + \cK_d(\bA;\bQ_j)
	\end{align*}
	\underline{Proof:}\\
	The proposition is shown by complete induction over $d=0,1,2,...$ for a fixed arbitrary $j \in \No$. The proposed equation we identify by the tupel $\lbrace j,d \rbrace$.
	\begin{enumerate}
		\item Basis: The proposition holds obviously for $\lbrace j,0 \rbrace$ since $\cM_j \subset \cM_j$.
		\item Hypothesis: $\lbrace j,d \rbrace$ holds for a specific $j,d \in \No$.
		\item Induction step: "When $\lbrace j,d \rbrace$ holds then $\lbrace j,d+1 \rbrace$ follows." is shown.\\
		We insert the induction hypothesis in the recursion formula of $\cM_{j+d+1}$:
		\begin{align*}
		&\cM_{j+d+1} \\
		=&(\bI-\omega_{j+d+1}\cdot\bA)\cdot\left(\cM_{j+d}\cap\cP^\perp\right) + \cQ_{j+d+1}\\
		\subset&(\bI-\omega_{j+d+1}\cdot\bA)\cdot\bigg( \Big(p_{j,d}(\bA) \cdot \big(\cM_{j}\cap\cK_d^\perp(\bA^\textsf{H};\bP)\big)+\cK_d(\bA;\bQ_j)\Big)\cap\cP^\perp\bigg)\\
		& + \cQ_{j+d+1}
		\end{align*}
		We can bound the last expression from above by neglecting the intersection of $\cK_d(\bA;\bQ_j)$ with $\cP^\perp$:
		\begin{align*}
		&\cM_{j+d+1} \\
		\subset&(\bI-\omega_{j+d+1}\cdot\bA)\cdot\bigg( \Big(p_{j,d}(\bA) \cdot \big(\cM_{j}\cap\cK_d^\perp(\bA^\textsf{H};\bP)\big)\Big)\cap\cP^\perp+\cK_d(\bA;\bQ_j)\bigg)\\
		& + \underbrace{\cQ_{j+d+1}}_{\subset \cQ_j}\\
		\subset&(\bI-\omega_{j+d+1}\cdot\bA)\cdot\bigg( \Big(p_{j,d}(\bA) \cdot \big(\cM_{j}\cap\cK_d^\perp(\bA^\textsf{H};\bP)\big)\Big)\cap\cP^\perp\bigg)\\
		& +\underbrace{(\bI-\omega_{j+d+1}\cdot\bA)\cdot\cK_d(\bA;\bQ_j) + \cQ_{j}}_{\equiv \cK_{d+1}(\bA;\bQ_{j})}
		\end{align*}
		Finally, we can use a formula that Gutknecht \cite{IDR-Gutkn} has used to show a particular structure of Sonneveld spaces. Due to the analogous nested recursion formula to Sonneveld spaces, it holds:
		\begin{align*}
		&(\bI-\omega_{j+d+1}\cdot\bA)\cdot\bigg( \Big(p_{j,d}(\bA) \cdot \big(\cM_{j}\cap\cK_d^\perp(\bA^\textsf{H};\bP)\big)\Big)\cap\cP^\perp\bigg)\\
		=& p_{j,d+1}(\bA) \cdot \big(\cM_{j}\cap\cK_{d+1}^\perp(\bA^\textsf{H};\bP)\big)
		\end{align*}
		Inserting this identity into the above estimate of $\cM_{j+d+1}$ completes the induction step.$\boxtimes$
	\end{enumerate}
\end{Lemma}

\subsection{The implementation of \Mstabno}\label{sec:ImplementationMstab}

We present a framework for the implementation of an \IDRstab variant that follows the approach from Fig.\,\ref{fig:MstabPrinciple}.

An implementation of \IDRstab for solving the first system can be sketched by the following lines of code:
\begin{algorithmic}[1]
\State Given: $\bA,\,\bb\hp{1},\,s,\,\tolabs$
\State $\bx:= \bO$,\quad$\br:=\bb\hp{1}$
\State Perform the initialsation. Use for instance \textsc{Initialisation}.
\While{the termination condition is not satisfied}
	\State Perform the iterative scheme; e.g. \GMstabI or \GMstabII.
	\If{$\bV\hp{-1}$ shall be written out}
		\State $\hbU := \bV\hp{-1}$
	\EndIf
\EndWhile
\State \Return $\bx,\hbU,\bP$
\end{algorithmic}
In line 6 an appropriate condition is to overwrite $\hbU$ by the current columns of $\bV\hp{-1}$ whenever one can be sure that $\dim(\cG_j)$ for the current Sonneveld space $\cG_j$ is larger than $s$. This is likely to be the case when the residual is sufficiently large because then one can be sure that the method has not terminated yet. So an appropriate replacement for lines 6--8 is
\begin{algorithmic}
\If{$\|\br\|>\tol_2 \cdot \|\bb\hp{1}\|$}
	\State $\hbU := \bV\hp{-1}$
\EndIf
\end{algorithmic}
with $\tol_2>\tolrel$. E.g., $\tol_2=10^{-3}$ to $\tol_2=10^{-2}$ seems a reasonable choice.

So far we have discussed how to implement the scheme that is depicted on the left-hand side in Fig.\,\ref{fig:MstabPrinciple}. In the following we sketch how the part from the right-hand side of the figure can be implemented.
\largeparbreak

An implementation of \Mstab for solving the second system can be sketched by the following lines of code:
\begin{algorithmic}[1]
	\State Given: $\bA,\,\bb\hp{2},\,\bP,\,\hbU,\,\tolabs$
	\State $\bx:= \bO$,\quad$\br:=\bb\hp{2}$
	\State \textit{// Initialisation:}
	\State $\bP$ is given. Set $\bV\hp{-1}:=\hbU$, $\bV\hp{0}:=\bA \cdot \bV\hp{-1}$
	\State \textit{// $\br \in \cM_J$ and $\rg(\bV\hp{0}) \subset \cM_J$ for some $J \in \N$}
	\While{the termination condition is not satisfied}
	\State \textit{// move data from $\cM_{J+j}$ into $\cM_{J+j+\ell}$. $j:=j+\ell$}
	\State Perform the iterative scheme; e.g. \GMstabI or \GMstabII.
	\EndWhile
	\State \Return $\bx$
\end{algorithmic}
As the comments lay out, during the while-loop the method iteratively moves the residual and the auxiliary vectors into subsequent $\cM$-spaces.
\largeparbreak

The question remains how to choose $\hbU$ when solving a third system $\bA \cdot \bx\hp{3} = \bb\hp{3}$ with \Mstab. In \cite{Mstab-report,Mstab-paper} we suggest to simply reuse the matrix for $\hbU$ that has been used for solving the second system. For longer sequences (i.e. $\nEqns$ is large) instead it is suggested to compute $\hbU$ periodically from scratch by solving e.g. each $10$th system with \IDRstabno instead of \Mstabno.

\subsection{Examples on the convergence and termination behaviour of \Mstabno}\label{sec:ExamplesConvergenceMstab}
In the subsections so far we have not given too strong of a motivation for \Mstab. This is once because as mentioned this thesis is rather concerned with the robust implementation of the underlying iterative scheme of \IDRstabno. On the other hand we have not motivated \Mstabno too exhaustively because for the numerical experiments in the subsequent section it will not be too much of importance to know exactly at which matrix-vector product the method should terminate in theory. Instead, in the subsequent numerical experiments section it is rather of interest if one particular implementation of \IDRstab or \Mstab can converge faster than all the others.
\largeparbreak
Nevertheless, since we have introduced a Krylov subspace recycling approach for \IDRstab in this section, in the following we briefly want to demonstrate the practical effects that this recycling approach has on the way how the method converges. In order to do so we consider one test case for two kinds of convergence behaviours, respectively: First, we look into a test problem where \IDRstab terminates after a finite number of iterations. Afterwards we consider a problem where \IDRstab experiences a superlinear convergence improvement after a few hundred iterations.

Parts of the problem description and the figures are quoted from \cite[Sec.\,5]{Mstab-paper}. These test cases consider a sequence \eqref{eqn:AxbSequence} where $\bA \in \R^{N \times N}$ arises from a central finite-difference discretisation of a Dirichlet problem
\begin{subequations}\label{eqn:PDE}
	\begin{align}
	-\epsilon \cdot \Delta u + \vec{\alpha}^\textsf{T} \LargerCdot \vec{\nabla}u - \beta \cdot u &=f\quad &\forall\vec{x} &\in \phantom{\partial}\Omega = (0,1)^d\\
	u &=u_\text{D} \quad &\forall\vec{x} &\in \partial\Omega 
	\end{align}
\end{subequations}
on a uniform Cartesian grid of mesh size $h \in 1/\N$, yielding $N=(1/h-1)^d$ equations. The experiments are performed in Matlab with $\epsmachine \approx 10^{-16}$ and $\hJ$ is chosen by the rule proposed in Sec.\,\ref{sec:ImplementationMstab} with $\tol_2$. All problems are solved with the initial guess $\bx_0 = \boldsymbol{0}$.

\paragraph{Example for finite termination}

For this test case $\bA$ is chosen from \eqref{eqn:PDE} with $d=2$, $h=1/351$ ($N=122500$), for $\epsilon=1$, $\vec\alpha = 1000/\sqrt{2}\cdot\vec{1}$, $\beta=1000$, where $f,u_\text{D}$ are chosen such that the discrete solution to $\bb^{(1)}$ is $u(x,y)=x\cdot y \cdot (1-x) \cdot (1-y)$ on the grid points. We choose $\bb\hb = \bA^{-1} \cdot \bb\ha$, i.e. the solutions to $\bb\ha$ and $\bb\hb$ live in the same Krylov subspace.

\begin{figure}
	\centering
	\includegraphics[width=1\linewidth]{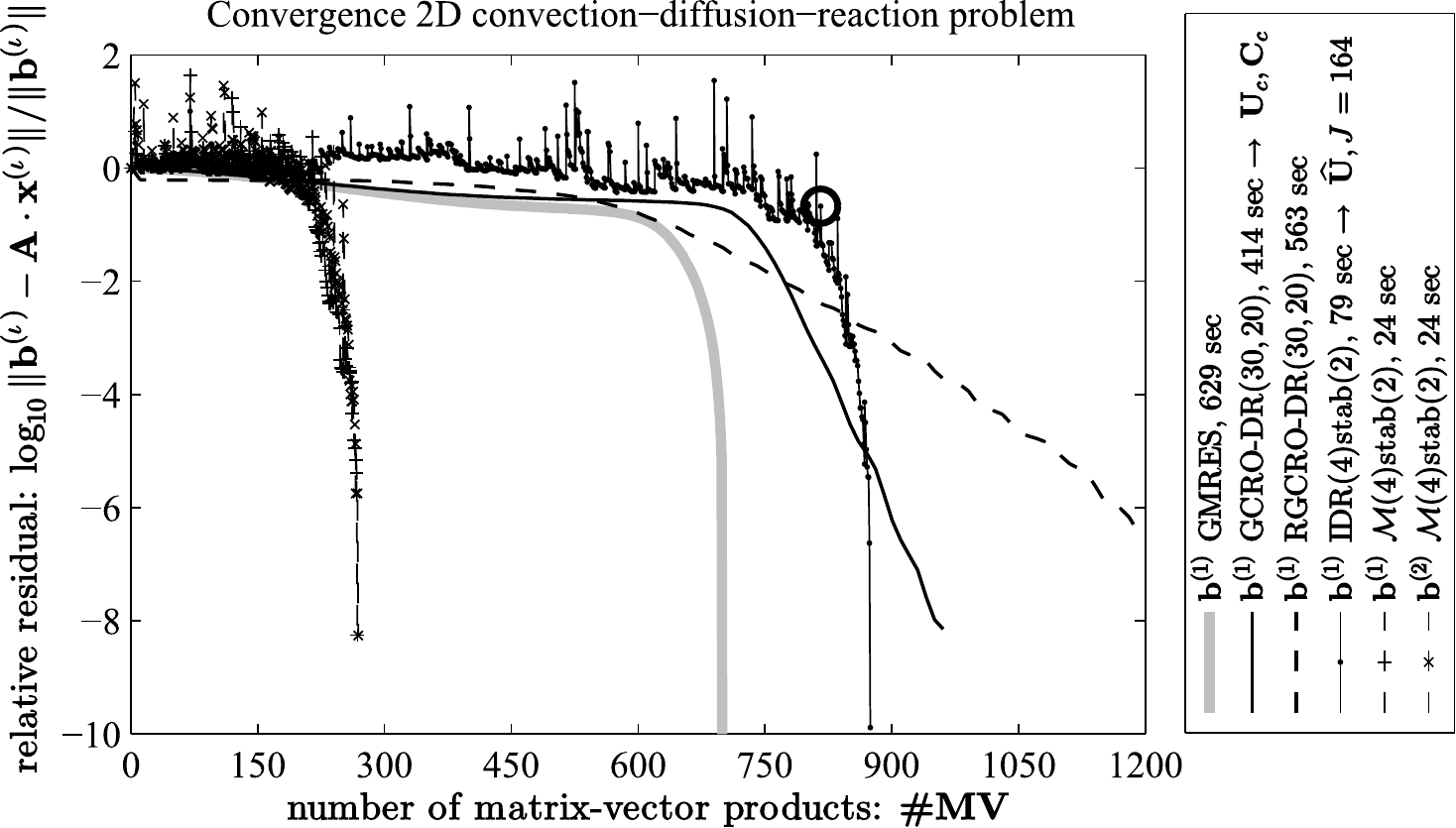}
	\caption{Example of finite termination behaviour of \IDRstabno and \Mstabno.}
	\label{fig:CDR2D}
\end{figure}

Fig.\,\ref{fig:CDR2D} shows the convergence over the number of iterations for \GMRES, \GCRODR \cite{GCRO-DR} (which is a commonly referred Krylov subspace recycling method based on deflation), \IDRstab and \Mstab. The latter two solvers use the \GMstab implementation. The circle on the convergence curve of \IDRstabno is the iterate at which the recycling data $\hbU$ is written out.

We see from \GMRES's convergence graph that the full Krylov subspace of $\bb\ha$ has $\approx 700$ dimensions. Using the theory discussed in Sec.\,\ref{sec:PropertiesIDRstab} one can show that \IDRstabp{$4$}{$\ell$} must terminate in about $700 \cdot (1+1/s) = 875$ matrix-vector products. For \Mstabno in turn one can show that it must terminate in about $875/(s-1) \approx 292 $ matrix-vector products, cf. \cite[Sec.\,3.2]{Mstab-paper}. As the figure shows, despite the numerical round-off the methods \IDRstabno and \Mstabno terminate precisely in the number of matrix-vector products that the theory predicts.

\begin{figure}
\centering
\includegraphics[width=1\linewidth]{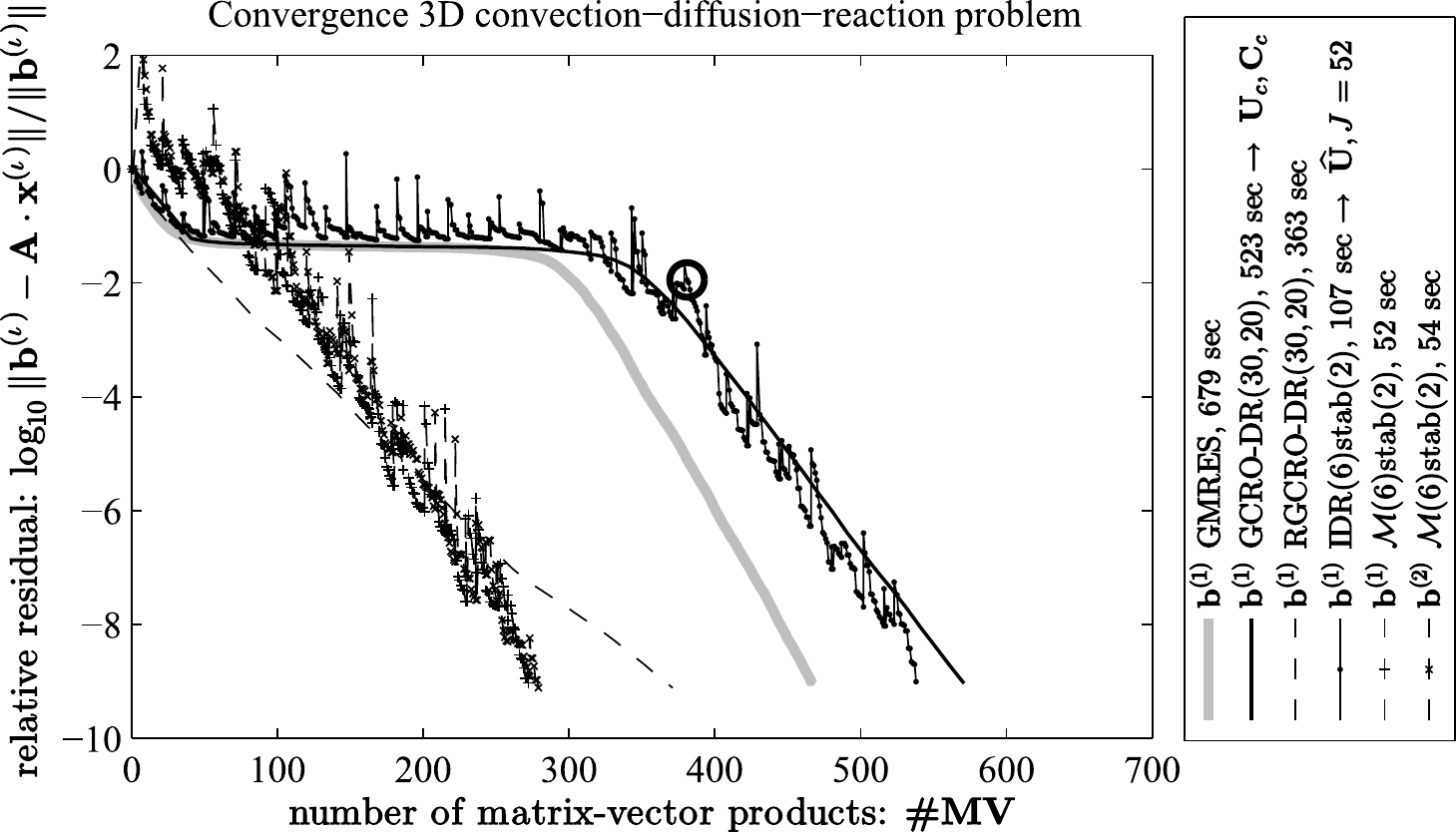}
\caption{Example of superlinear convergence behaviour of \IDRstabno and \Mstabno.}
\label{fig:CDR3D}
\end{figure}

\paragraph{Example for superlinear convergence}
In this case $d=3$, $h=1/61$ ($N=216000$), $\epsilon=1$, $\vec{\alpha}=\vec{x}$, $\beta=-10$ (i.e. in \eqref{eqn:PDE} there is a positive shift) is used. $f$ and $u_\text{D}$ are chosen such that $\bb\ha=\bA\cdot\boldsymbol{1}$. The second right-hand side is chosen as $\bb\hb=\bA^{-1}\cdot\bb\ha$. The results are given in Figure\,\ref{fig:CDR3D}\,.

The same implementations as in the former paragraph have been tested on this problem, except that we have increased the method parameter $s=6$. We can see that \GMRES transitions to superlinear convergence at 300 matrix-vector products. Thus, \IDRstabp{$6$}{$2$} should transition at $300\cdot(1+1/s)=350$ iterations, which it roughly does. Afterwards, its slope of convergence should be only $1+1/s$ times slower than \GMRES's rate, which is roughly the case.

Also the convergence behaviour of \Mstabno meets our expectations: It transitions to superlinear convergence after about $350/(s-1)=50$ matrix-vector products. The convergence curve of \Mstab looks similar to the predicted behaviour that was described in Fig.\,\ref{fig:SuperlinConv} and Fig.\,\ref{fig:MstabMotivation}.

\section{Numerical experiments}\label{sec:NumericalExperiments}

The numerical experiments are organised as follows. We compare our new implementation \GMstab against the three \IDR  and three \IDRstab algorithms from Section \ref{sec:ImplementationsIDRstab}, namely: \textsf{IDR($s$)noDec}, \textsf{IDR($s$)biorth}, \textsf{IDR($s$)obio}, \textsf{IDR($s$)stab($\ell$)ref}, \textsf{IDR($s$)stab($\ell$)biorth}, \textsf{IDR($s$)stab($\ell$)obio}. We first show experiments where the methods are applied to solve a single linear system. Afterwards we evaluate the practical efficiency of the methods with regards to Krylov subspace recycling when solving sequences of multiple right-hand sides.

\subsection{Software implementation for the numerical experiments}
The codes that we have used to run the numerical experiments are written from scratch by the author and published on \texttt{www.MartinNeuenhofen.de} $\rightarrow$ \texttt{GMstab}. They are implemented in Matlab 2011a and were executed on an ASUS Eee Slate 121b (2011) under Windows\,7.

Under the download link the reader finds a zipped file \texttt{Solvers}. This file contains three subfolders:
\begin{itemize}
	\item \texttt{BlockFunctions}: This folder contains all the algorithmic sub-routines such as the \texttt{lq}, \texttt{roworth}, \textsc{Initialisation} and \textsc{StabCoeffs}. This folder does also hold a class \texttt{PerfMeasure}.
	\item \texttt{LinearSystemSolvers}: This folder contains the eight solvers \GMRES, \textsf{IDR($s$)noDec}, \textsf{IDR($s$)biorth}, \textsf{IDR($s$)obio}, \textsf{IDR($s$)stab($\ell$)ref}, \textsf{IDR($s$)stab($\ell$)biorth}, \textsf{IDR($s$)stab($\ell$)obio} and \GMstab. All these routines use structures \texttt{in} and \texttt{out} for their in- and outputs.
	\item \texttt{Testproblems}: This folder contains two subfolders:
	\begin{itemize}
		\item \texttt{\_classes}: This subfolder contains functionality for generating and running a testcase. To this end it holds the class \texttt{Testcase}.
		\item \texttt{Testcases}: This subfolder contains data bases to generate the linear systems for the numerical experiments and holds a bunch of scripts to run the respective experiments.
	\end{itemize}
\end{itemize}

A test problem is generated by instantiating a class object of \texttt{Testcase}. The constructor expects a string that is the name of the test case. After being instantiated, one can push a system to the \texttt{Testcase}-object by passing a Matlab function-handle for the matrix-vector product with $\bA$, the right-hand side $\bb$ and an initial guess $\bx_0$. Afterwards one can push solvers to the \texttt{Testcase}-object that shall be used to solve this system. The solver is pushed as a data structure that specifies the solver name and optional solver parameters. Default solver parameters can be pushed directly to the \texttt{Testcase}-object. Optional parameters can be $s,\ell$ or $\tol_2$ or the information that recycling data $\hbU$ from the output of the $k$th solver from the $\iota$th system shall be used. Pushing systems and solvers can be repeated in arbitrary order. The pushed solvers are applied to the respectively latest pushed system.

Once the instantiation of the testcase is completed, there is a function \texttt{Run} that runs the test case by calling all the solvers on the respective systems.

After the test case has run it returns a structure \texttt{plotData} that can be used by further subroutines to plot the experimental results.

For example, Fig.\,\ref{fig:2D_CDRb_mv} can be generated executing the following script in Matlab:

\begin{verbatim}
%% initialise testcase
% add the file 'Solvers' and all subfolders to the Matlab path
myTestcase   = Testcase('CDR\_2Db\_v4');
% myTestcase is a class object of Testcase.
% The argument is the figure title.
sSaveName    = 'Testcase_CDR_2Db_v4';  % file name to save 
                                       %experimental result data

%% generate problem instance
[A_fun,B,X0] = GenerateLinearSystem(); % subroutine to generate
                                       % the linear systems
[terminate]  = GenerateTermination();  % a data structure that 
                                       % defines tolabs, 
                                       % maxiter, etc.
[param]      = GenerateDefaultParam(); % solver parameters 
                                       % such as s,L

%% 1st system
% A, b, x0 for iota=1
myTestcase.pushSystem( A_fun, B(:,1), X0(:,1) ); 
myTestcase.setTerminate( terminate );
myTestcase.setDefaultParam( param );
% % make solvers for this system
cIDRstabbio = makeSolver(myTestcase,'IDRstabbiortho');
cIDRstabbio.recyclingparam.tolabs2 = 1e-3;
myTestcase.pushSolver(cIDRstabbio);
% systemindex=1, solverindex=1 (*)
cGMstab     = makeSolver(myTestcase,'GMstab');
cGMstab.recyclingparam.tolabs2     = 1e-3;
myTestcase.pushSolver(cGMstab);
% systemindex=1, solverindex=2 (**)

%% 2nd System
myTestcase.pushSystem( A_fun, B(:,2), X0(:,2) );
myTestcase.setTerminate( terminate );
myTestcase.setDefaultParam( param );
% % make solvers for this system
cIDRstabbio = makeSolver(myTestcase,'IDRstabbiortho');
cIDRstabbio.recyclingparam.systemindex  = 1;
cIDRstabbio.recyclingparam.solverindex  = 1;
% M(4)stab(2) with \hat{U} from IDR(4)stab(2)bio
%    from the system iota=1 is used. Cf. (*)
myTestcase.pushSolver(cIDRstabbio);
cGMstab     = makeSolver(myTestcase,'GMstab');
cGMstab.recyclingparam.systemindex      = 1;
cGMstab.recyclingparam.solverindex      = 2;
% M(4)stab(2) with \hat{U} from GM(4)stab(2)
%    from the system iota=1 is used. Cf. (**)
myTestcase.pushSolver(cGMstab);

%% Run testcase
myTestcase.Run();

%% Save plot data
plotData = myTestcase.getPlotData();
save(['TestResults/',sSaveName,'_data.mat'],'plotData');

% % Plot
window = struct('xL',0,'xR',1200,'dx',100,...
   'yL',-12,'yR',4,'dy',2);
% window defines the axis intervals (xL,xR) and ticks (dx)
plot_iterres_over_matvec( plotData,window,...
   ['TestResults/',sSaveName,'_iterres_matvec.pdf']);
\end{verbatim}
The figure is generated completely automatically with all labels and saved as a PDF-file.

\largeparbreak
The runtimes and true residuals are measured by the \texttt{PerfMeasure} class. When the \texttt{Testcase}-object calls the respective Krylov method it passes an object of \texttt{PerfMeasure}. This object can make snapshots during the execution of the Krylov method. One snapshot contains information on the number of computed matrix-vector products, the runtime, the runtime spent in matrix-vector products only, the relative residual norm of the iteratively updated residual and the true residual norm. To this end the \texttt{PerfMeasure} class has four routines (plus several more that are not important for this presentation):
\begin{itemize}
	\item \texttt{[]=resume()}: This (re)starts the stopwatch that evaluates the execution time of the Krylov method.
	\item \texttt{[$\bv$]=matvec($\bu$)}: This is the function handle that evaluated the matrix-vector product $\bv = \bA \cdot \bu$. \texttt{PerfMeasure} counts the number of computed matrix-vector products.
	\item \texttt{[]=stop()}: This pauses the stopwatch. This functionality is required because the time for computing the true residuals should not occur in the measured execution time of the Krylov method.
	\item \texttt{[]=read($\bx$,$\beta$)}: This routine triggers to save a snapshot of the current state of the Krylov method. From $\bx$ the \texttt{PerfMeasure} object computes the true residual norm for this snapshot.
\end{itemize}

\subsection{Experiments with single systems}
In this subsection we apply all the discussed \IDR and \IDRstab implementations to solve single linear systems.

\subsubsection{Experiments without preconditioning}
We have four experiments for single linear systems without preconditioning. The latter of these experiments is parametric and consists of four different systems of which each is solved separately.

\paragraph{Experiment BiCGstabL\_xpl1}
The first test problem is taken from \cite[example 1]{BiCGstabL}. This problem is \eqref{eqn:PDE} for the parameters $d=3$, $\epsilon=-1$, $\vec{\alpha}=(1000,\,0,\,0)\t$, $\beta=0$, $h=1/51$. However, whereas in \cite{BiCGstabL} finite volumes are used, we use central finite differences. $f$ and $u_\text{D}$ are chosen such that the solution on the mesh-points is $u(x,y,z) = \exp(x\cdot y\cdot z)\cdot\sin(\pi \cdot x)\cdot\sin(\pi \cdot y)\cdot\sin(\pi \cdot z)$. This makes a non-symmetric system of size $N=125000$.

In the referred paper the methods \BiCG, \BiCGstab and \BiCGstabL for $\ell=2$ are tested. In the referenced paper it is observed that all methods except \BiCGstabL for $\ell>1$ struggle to converge. We have chosen this test case because it helps demonstrating that in general \IDR is not competitive to \IDRstab.

\begin{figure}
\centering
\includegraphics[width=1\linewidth]{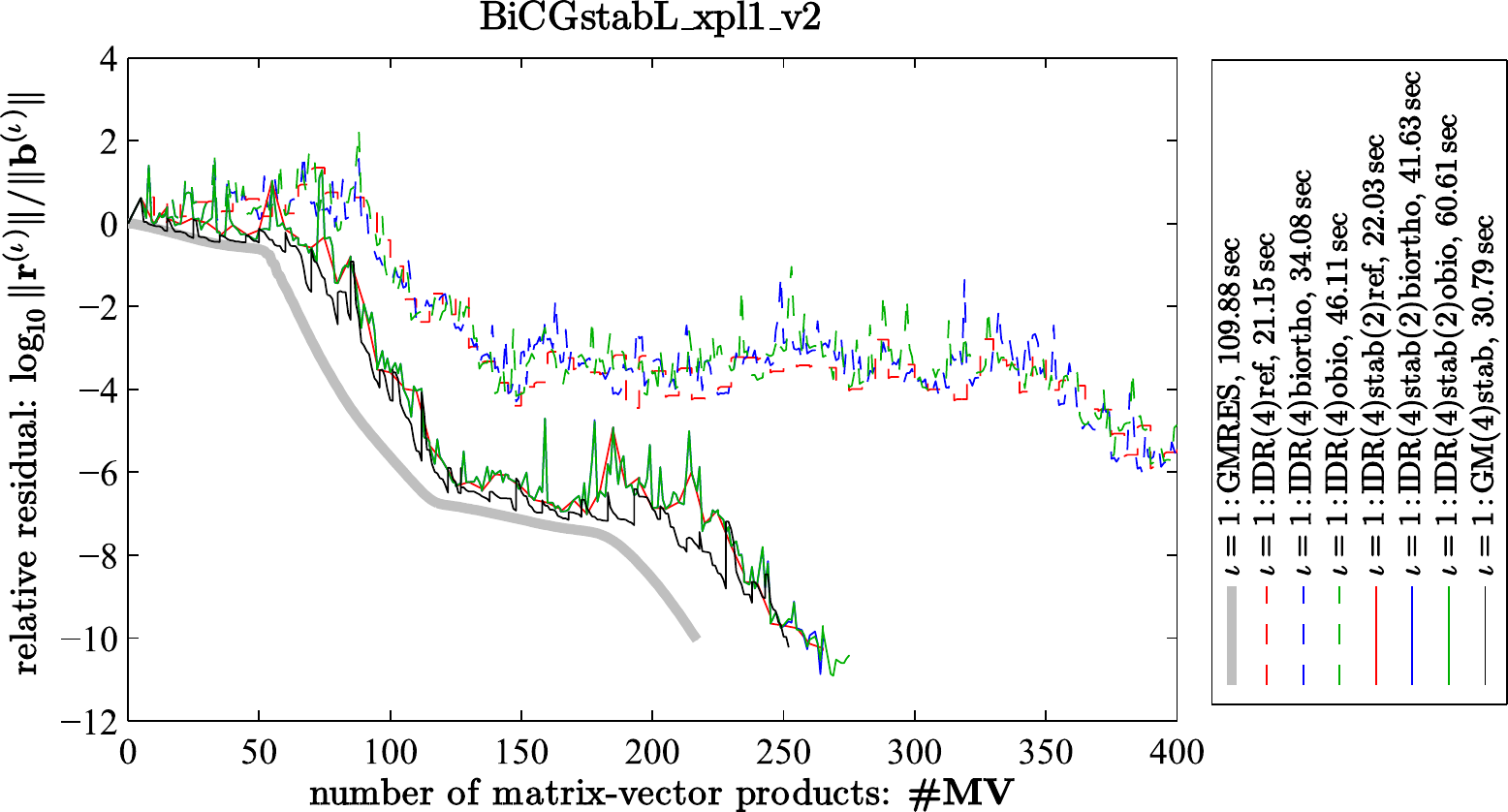}
\caption{\IDRp{$4$} is not competitive to \IDRstabp{$4$}{$2$}.}
\label{fig:BiCGstabL_xpl1_iterres}
\end{figure}

Fig.\,\ref{fig:BiCGstabL_xpl1_iterres} shows the iteratively computed residuals over the number of matrix-vector products for \GMRES (only for a theoretical evaluation of how close the convergence curves of the \IDRO methods are to the optimal rate of convergence of \GMRES) and all the discussed \IDR and \IDRstab variants. For this experiment the plotted iteratively updated residual $\br$ matches accurately in more than $10$ digits with the true residual $\bb-\bA \cdot \bx$ for all the tested methods.

From the figure we make the following observations:
\begin{itemize}
	\item All the \IDRp{$4$} variants are uncompetitive to the \IDRstabp{$4$}{$2$} variants because they do not converge in a reliabe way. The author considers a reliable way of convergence as a graph that is closely behind \GMRES. In fact, this discards all the \IDRp{$4$} variants as black-box linear system solvers.
	\item For this test problem all \IDRstabp{$4$}{$2$} variants converge reliably.
	\item Among the \IDRstab variants, the reference implementation is the cheapest in terms of computation time per matrix-vector product. The second-cheapest variant is \GMstab. This can be easily seen since all the \IDRstab variants require roughly the same number of matrix-vector products.
\end{itemize}

\paragraph{Experiment BiCGstabL\_xpl3}
This test problem is adapted from \cite[example 3]{BiCGstabL}. This problem is \eqref{eqn:PDE} for the parameters $d=2$, $h=1/201$ and
\begin{align*}
 	\epsilon&=0.1\,,\\
 	\vec{\alpha}(x,y)&=\begin{pmatrix}
 	 	4 \cdot x \cdot (x-1) \cdot (1-2\cdot y)\\
 	 	4 \cdot y \cdot (1-y) \cdot (1 - 2 \cdot x)
 	\end{pmatrix}\,,\\
 	\beta&=0\,.
\end{align*}
Again, whereas in \cite{BiCGstabL} finite volumes are used, we use central finite differences. $f$ and $u_\text{D}$ are chosen such that the solution on the mesh-points is $u(x,y) = \sin(\pi \cdot x) + \sin(\pi \cdot y) + \sin(13 \cdot \pi \cdot x) + \sin(13 \cdot \pi \cdot x)$. This makes a non-symmetric system of size $N=40000$.

\begin{figure}
\centering
\includegraphics[width=0.7\linewidth]{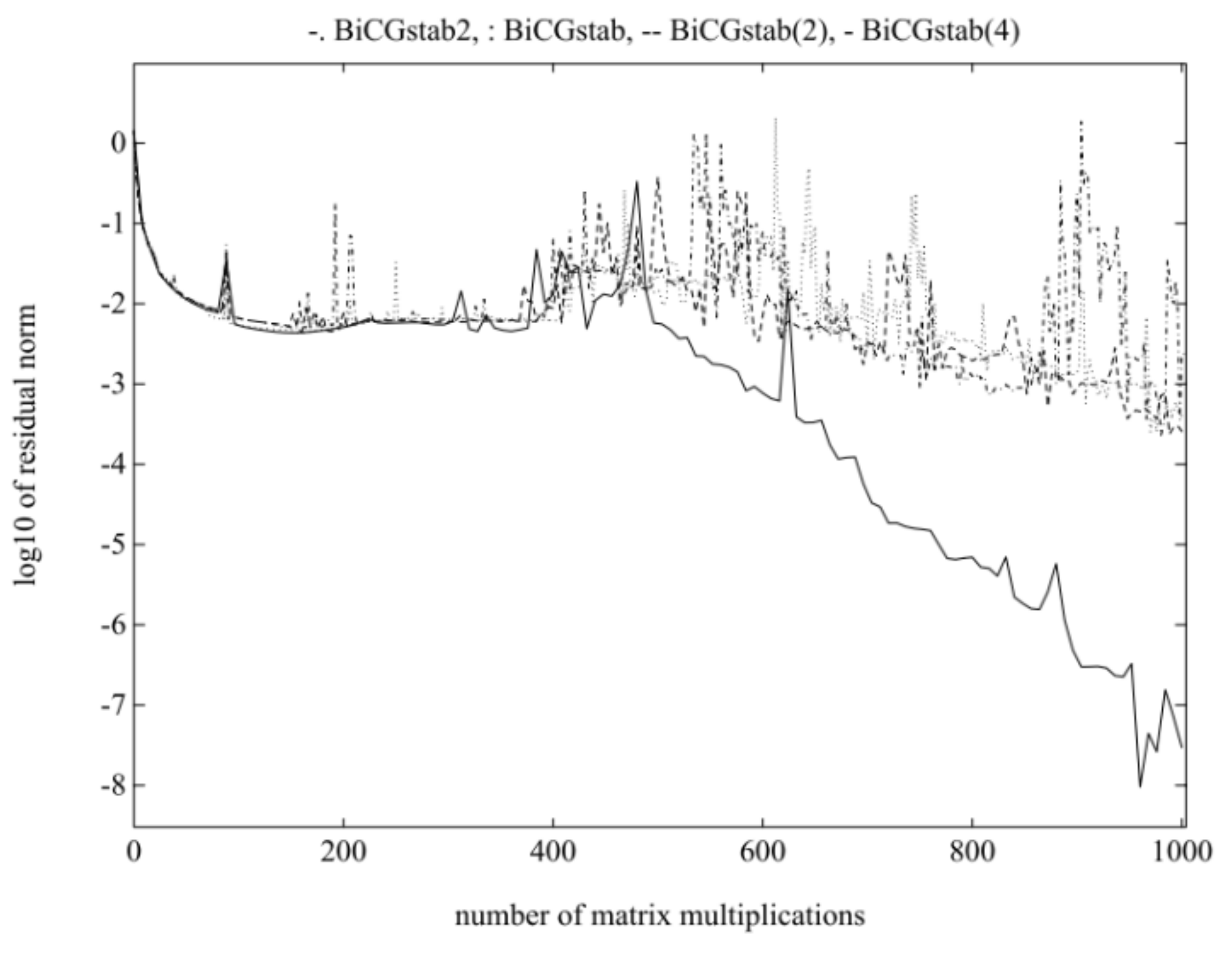}
\caption{Image taken from Fig. 5.3 from \cite{BiCGstabL}. This image is vital to motivate our test problem BiCGstabL\_xpl3.}
\label{fig:BiCGstabL_xpl3_fromReference}
\end{figure}

In the referred paper the methods \BiCGstabL for $\ell \in \lbrace 1,\,2,\,4\rbrace$ are compared. The convergence graph from this paper is quoted in Fig.\,\ref{fig:BiCGstabL_xpl3_fromReference}. The figure shows that after the transition to faster convergence, which is at about 450 matrix-vector products, the variant for $\ell=4$ converges twice as fast as all the other methods.

In the last experiment we have determined whether \IDR is competitive to \IDRstabp{$s$}{$2$}. Now with this experiment instead we want to investigate in the question whether \IDRstab with $\ell>2$ is ever required. Since Fig.\,\ref{fig:BiCGstabL_xpl3_fromReference} indicates that only $\ell=4$ can achieve a fast rate of convergence, we believe that this is a very good test problem for our investigation.

\begin{figure}
\centering
\includegraphics[width=1\linewidth]{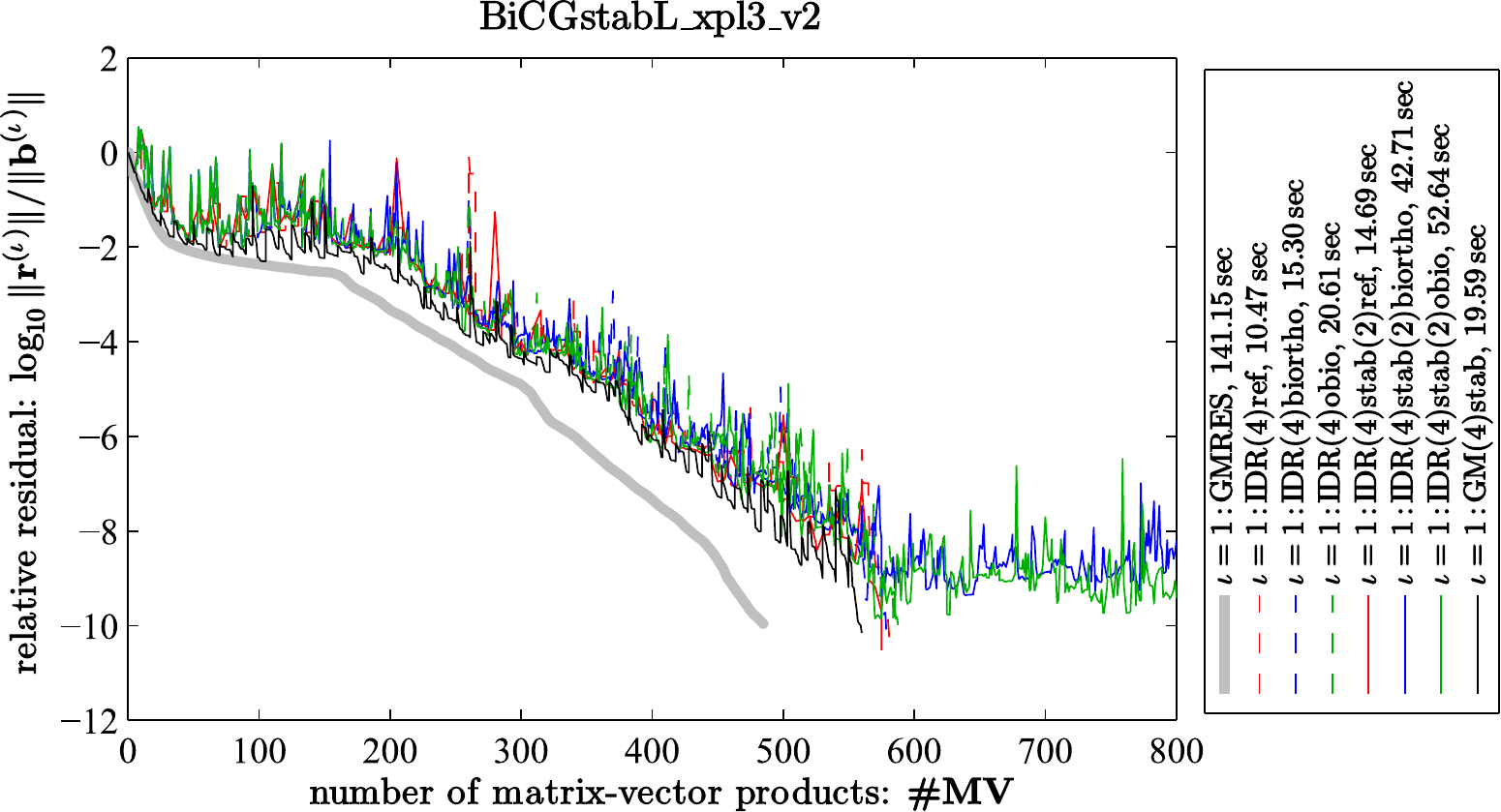}
\caption{\IDRstabp{$4$}{$2$} is not far behind \GMRES, so $\ell>2$ is not required. However, the birotho and obio variant fail to maintain the rate of convergence.}
\label{fig:BiCGstabL_xpl3_iterres}
\end{figure}

In Fig.\,\ref{fig:BiCGstabL_xpl3_iterres} we compare all the formerly discussed implementations of \IDR and \IDRstab for $\ell=2$ by plotting their iteratively updated residual. Also for this test case the residual gap between the iteratively updated vector $\br$ and $\bb-\bA \cdot \bx$ is $<10^{-10}$.

Since all our methods use $\ell<4$ we had expected according to Fig.\,\ref{fig:BiCGstabL_xpl3_fromReference} that they converge at most half as fast as \GMRES. However, in contrast to our expectation, all the compared methods converge closely behind \GMRES. Further, their transition to superlinear convergence is already at $\approx 150$ matrix-vector products and not only after $400$.

The results of our experiment indicate that the observed convergence behaviours from Fig.\,\ref{fig:BiCGstabL_xpl3_fromReference} do not indicate an improved convergence for $\ell=4$. Instead, it rather seems that the quoted figure does only show the erratic behaviour of four methods that do not converge robustly for this problem.

To our best knowledge, there is only one other test problem given in the literature where an \IDRstabno variant with $\ell>2$ converges faster than the according variant for $\ell=2$. This is the test problem that we discuss in the next paragraph.

Further to the above, we make the following obervation.
\begin{itemize}
	\item The biortho and obio variant of \IDRstabno suffer from loss of superlinearity: At $\approx 580$ matrix-vector products the slope of their convergence graphs falls back to horizontal. Since the iteratively updated residual and the true residual match to a far higher accuracy this effect of loss of convergence is not related to the residual gap at all.
	\item Consequently, the biortho and obio variant of \IDRstabno are unable to maintain any rate of convergence up to the desired accuracy, irrespective to the residual gap. In fact, this discards them as black-box linear system solvers.
\end{itemize}

\paragraph{Sherman5}
Sherman5 from the Harwell-Boeing collection\footnote{\texttt{http://math.nist.gov/MatrixMarket/collections/hb.html}} is a highly indefinite system of linear equations, which makes it difficult to solve by iterative methods \cite[Sec.\,6.1]{IDRstab-Paper}. We want to investigate on this test case whether \IDRstabp{$4$}{$2$} is competitive to \IDRstabp{$4$}{$4$}. We wonder about that since our new implementation \GMstab is limited to $\ell=2$.

\begin{figure}
\centering
\includegraphics[width=0.6\linewidth]{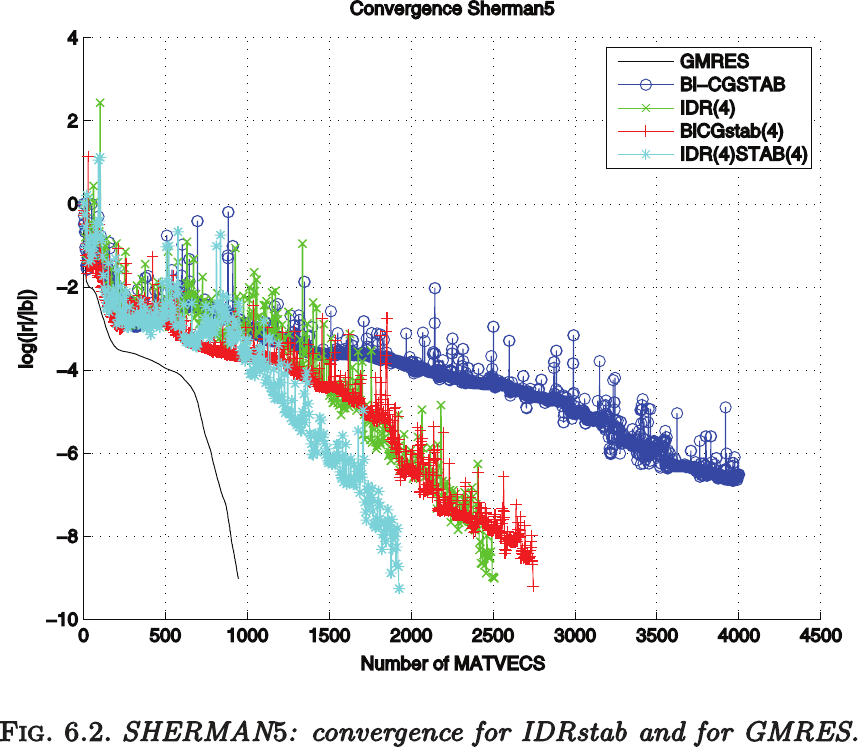}
\caption{Convergence of the iteratively updated residuals of \IDRstab vs. \GMRES, quoted from \cite[Fig.\,6.2]{IDRstab-Paper}. The true residuals of the computations presented in this figure do all stagnate above $10^{-8}$.}
\label{fig:Sherman_reference}
\end{figure}

In \cite[Fig.\,6.2]{IDRstab-Paper}, quoted in Fig.\,\ref{fig:Sherman_reference}, we see that their \IDRstabp{$4$}{$4$} converges in 2000 matrix-vector products to a relative residual of $10^{-9}$, whereas their \IDRstabp{$4$}{$1$} requires 2500 matrix-vector products to do so.

In the following we test the reference and biortho variants of the here discussed \IDRp{$4$} and \IDRstabp{$4$}{$2$} implementations and compare them against \GMstab. We then evaluate whether our implementations with $\ell=2$ converge noticeably slower than \IDRstabp{$4$}{$4$} from \cite[Fig.\,6.2]{IDRstab-Paper}.

\begin{figure}
\centering
\includegraphics[width=1\linewidth]{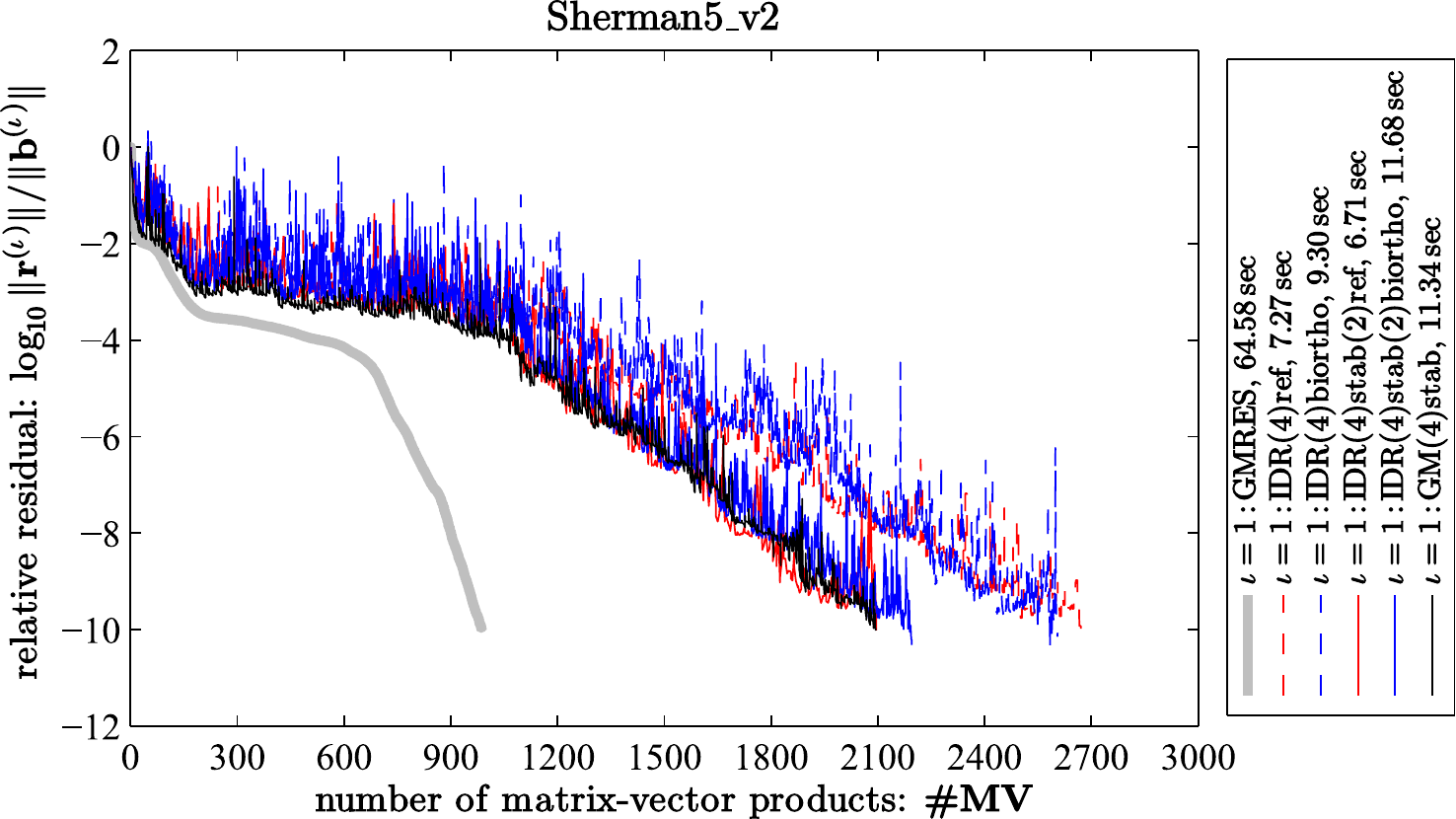}
\caption{\IDRstabno with $\ell=2$ converges just as well as with $\ell=4$.}
\label{fig:Sherman5_iterres}
\end{figure}

Fig.\,\ref{fig:Sherman5_iterres} shows the results. Since we use $\tol=10^{-10}$ whereas the authors in \cite[Fig.\,6.2]{IDRstab-Paper} have chosen $\tol=10^{-9}$, we require about a hundred iterations more to converge. All in all, the convergence graphs of all \IDRstabp{$4$}{$2$} implementations in Fig.\,\ref{fig:Sherman5_iterres} look very similar to the curve of 
\IDRstabp{$4$}{$4$} in \cite[Fig.\,6.2]{IDRstab-Paper}. As a conclusion, also for this difficult test problem a value for $\ell>2$ does not seem to be necessary or useful.

Regarding the methods from this thesis, it occurs that \GMstab converges fastest in terms of matrix-vector products. As a drawback, \GMstab requires more computation time per iteration. Its cost per iteration seems comparable to the cost of \IDRstabbiortho. Clearly, in terms of matrix-vector products the \IDR methods stand behind the \IDRstabno variants that use $\ell=2$.

\paragraph{CDR\_2Dparam($c_1,c_2$)}
This parametric test problem uses a central finite difference discretisation of \eqref{eqn:PDE} for $d=2$, $h=1/350$ ($N=122500$) for the following parameters.
\begin{align*}
	\epsilon &= 1\,,\\
	\vec{\alpha} &= c_1 \cdot \begin{pmatrix}
	1000/\sqrt{2}\\
	1000/\sqrt{2}
	\end{pmatrix}\,,\\
	\beta &= c_2 \cdot 1000\,.
\end{align*}
$f$ and $u_\text{D}$ are chosen for each choice of $c_1,c_2$ such that the solution on the mesh-points is $u(x,y) = x \cdot y \cdot (1-x) \cdot (1-y)$. In the following we consider the four systems where the parameter tuple $(c_1,c_2)$ is $(0,0)$, $(1,0)$, $(0,1)$ and $(1,1)$. For each of the four experiments we compare the reference and biortho implementation of \IDR and \IDRstab against \GMstab. We omit the obio variants since both of them are slower than \GMstab in terms of runtime per matrix-vector product and in order to reduce the number of graphs in the diagrams.

For $(0,0)$ we obtain a Poisson problem (symmetric positive definite), for $(1,0)$ a convection-diffusion problem (close to skew-symmetric), for $(1,0)$ a diffusion-reaction problem (symmetric but strongly indefinite) and for $(1,1)$ a convection-diffusion-reaction problem (strongly a-symmetric and indefinite). The problem is adapted from \cite[Sec.\,6.4]{IDRstab-Paper}.

\begin{figure}
\centering
\includegraphics[width=1\linewidth]{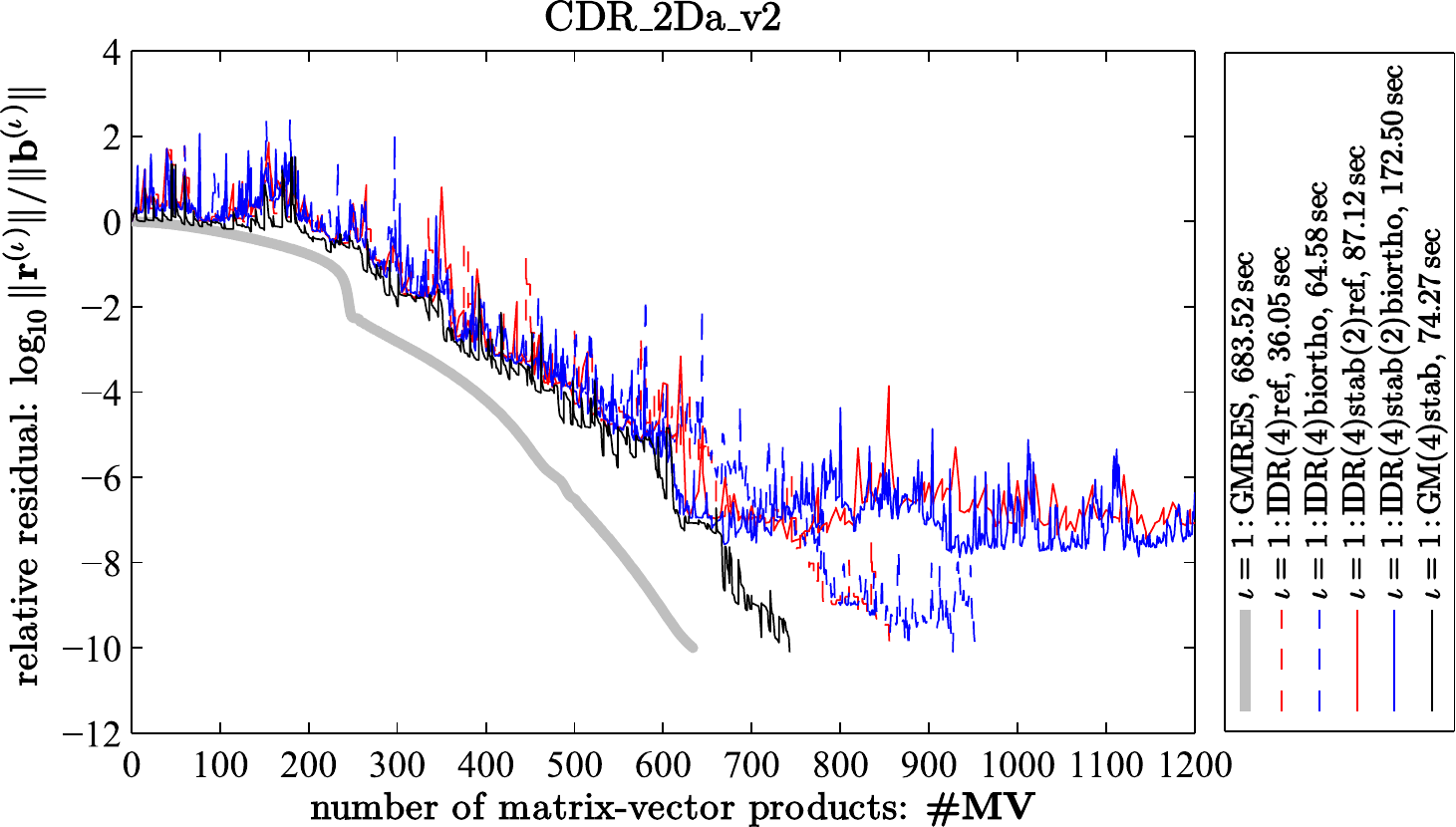}
\caption{CDR\_2Dparam($0,0$). Diffusion.}
\label{fig:CDR_2Da_iterres}
\end{figure}

\begin{figure}
	\centering
	\includegraphics[width=1\linewidth]{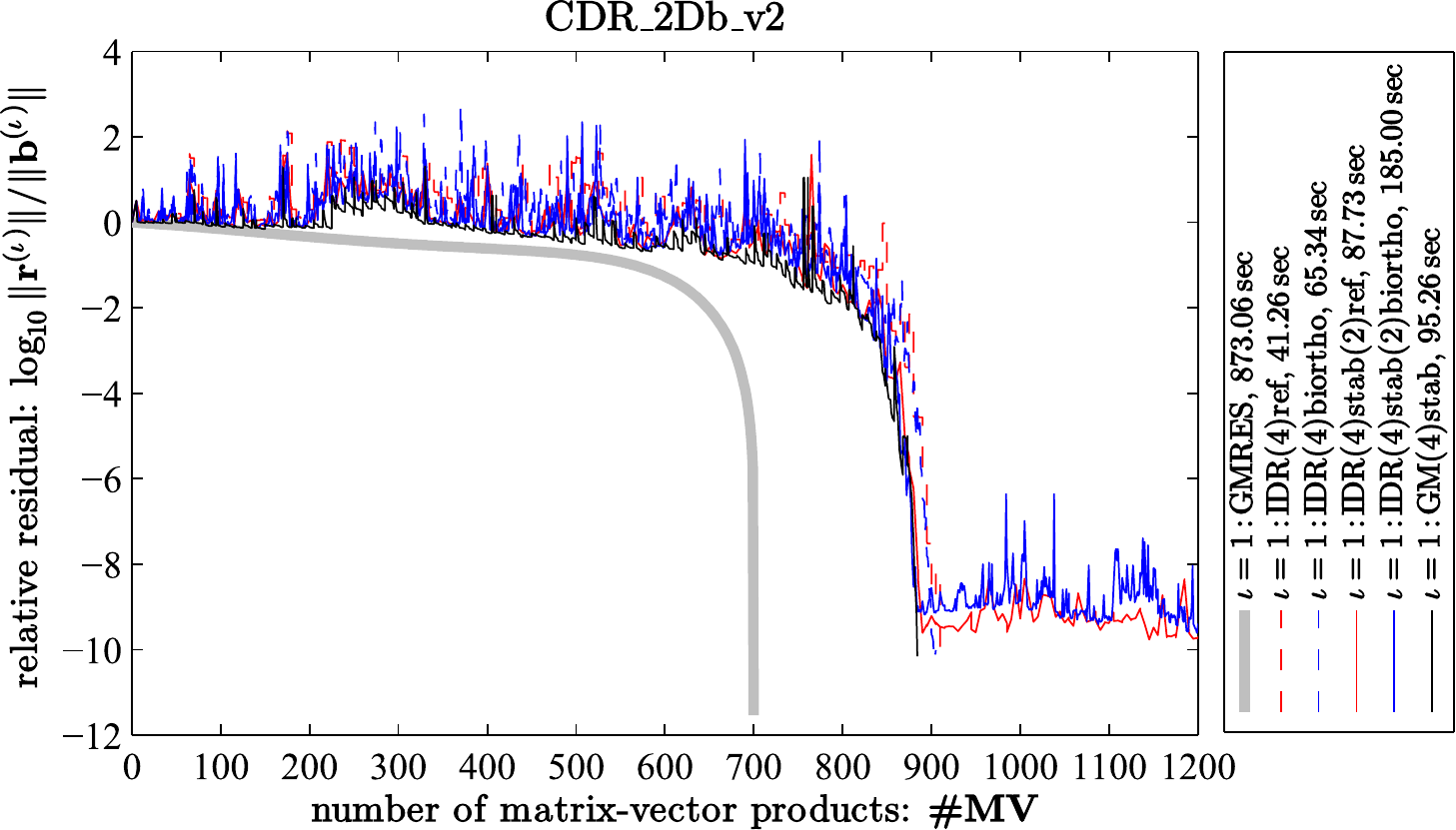}
	\caption{CDR\_2Dparam($1,0$). Convection-Diffusion.}
	\label{fig:CDR_2Db_iterres}
\end{figure}

\begin{figure}
	\centering
	\includegraphics[width=1\linewidth]{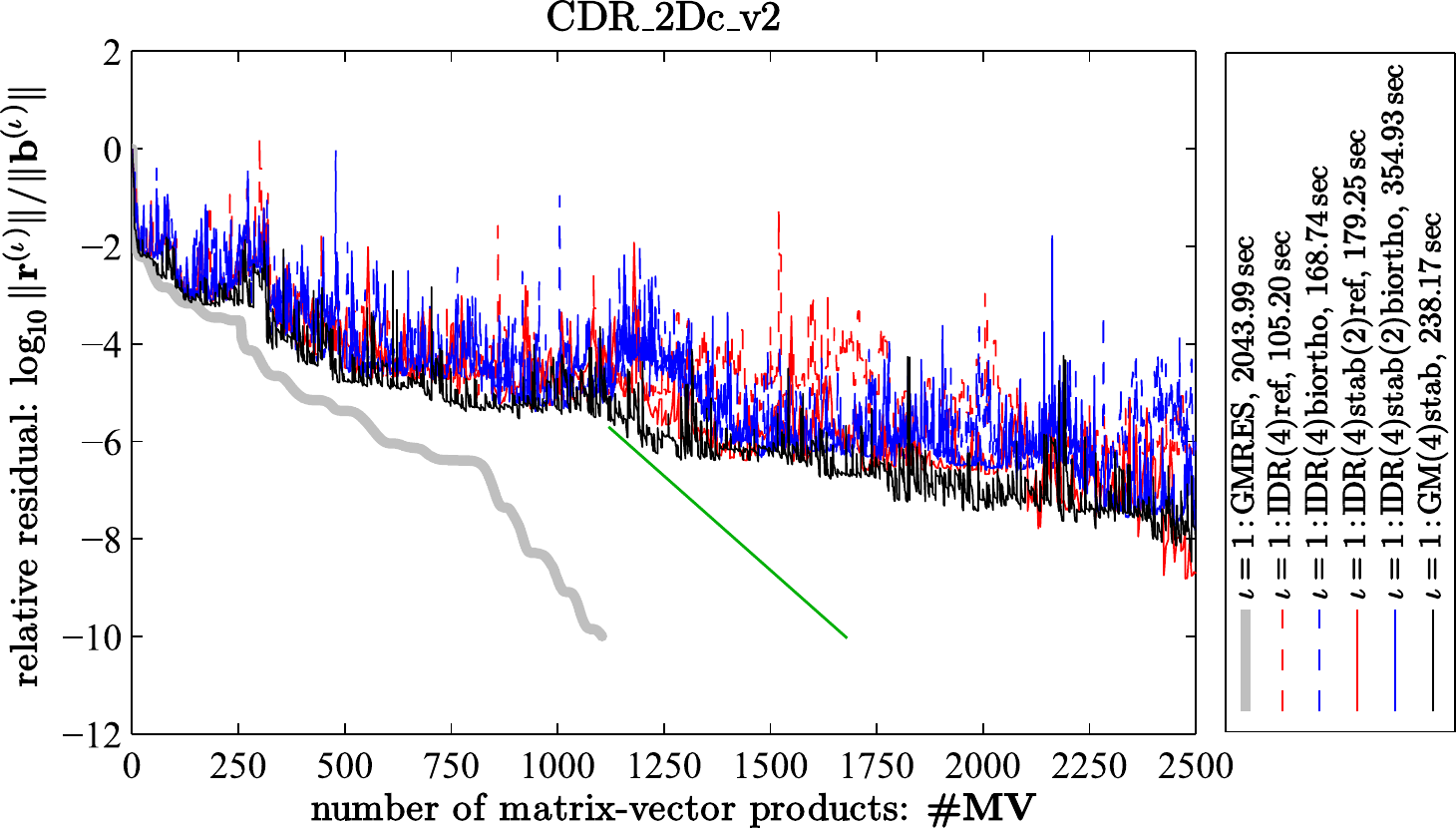}
	\caption{CDR\_2Dparam($0,1$). Diffusion-Reaction.}
	\label{fig:CDR_2Dc_iterres}
\end{figure}

\begin{figure}
	\centering
	\includegraphics[width=1\linewidth]{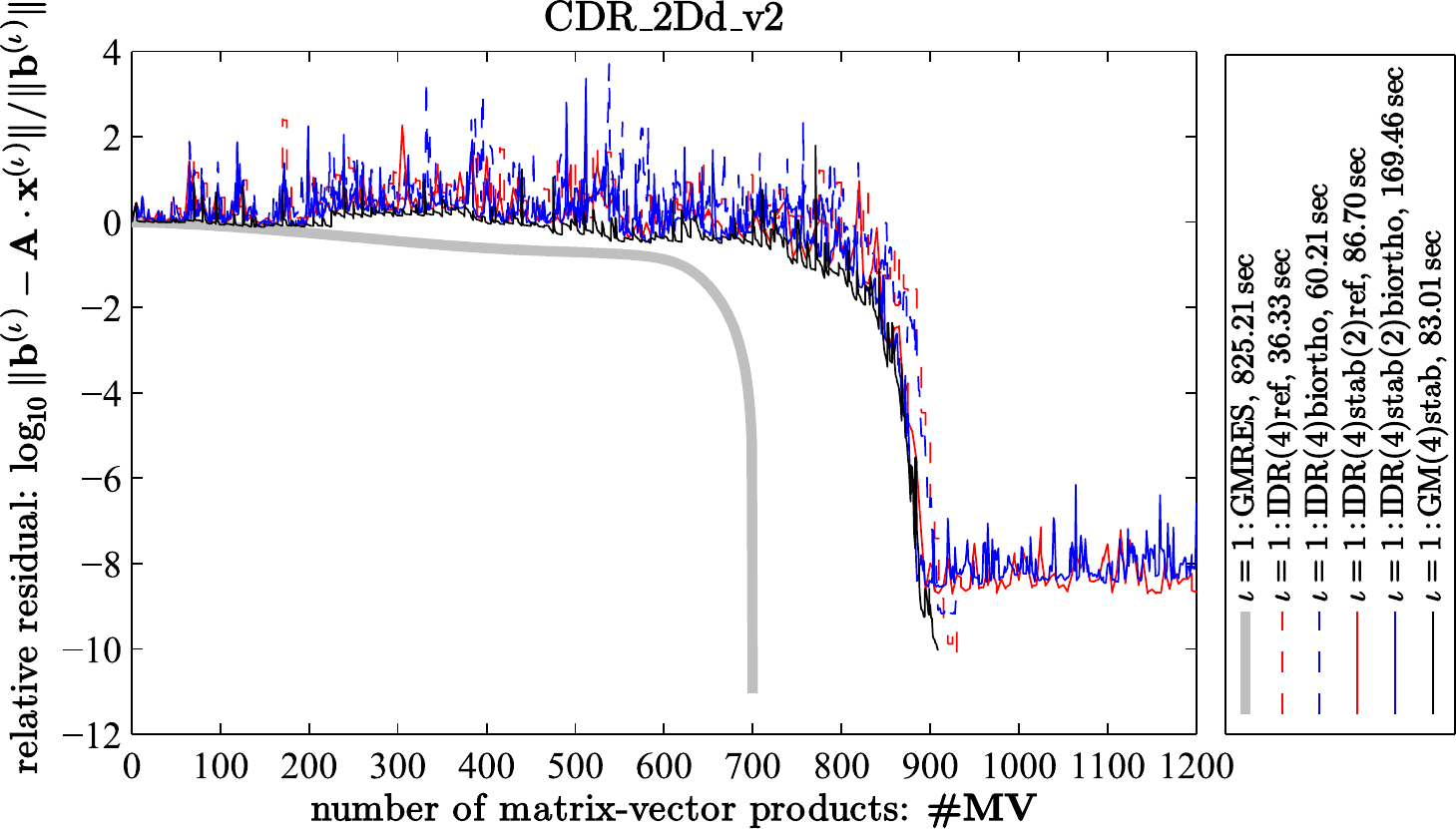}
	\caption{CDR\_2Dparam($1,1$). Convection-Diffusion-Reaction.}
	\label{fig:CDR_2Dd_trueres}
\end{figure}

\begin{figure}
	\centering
	\includegraphics[width=1\linewidth]{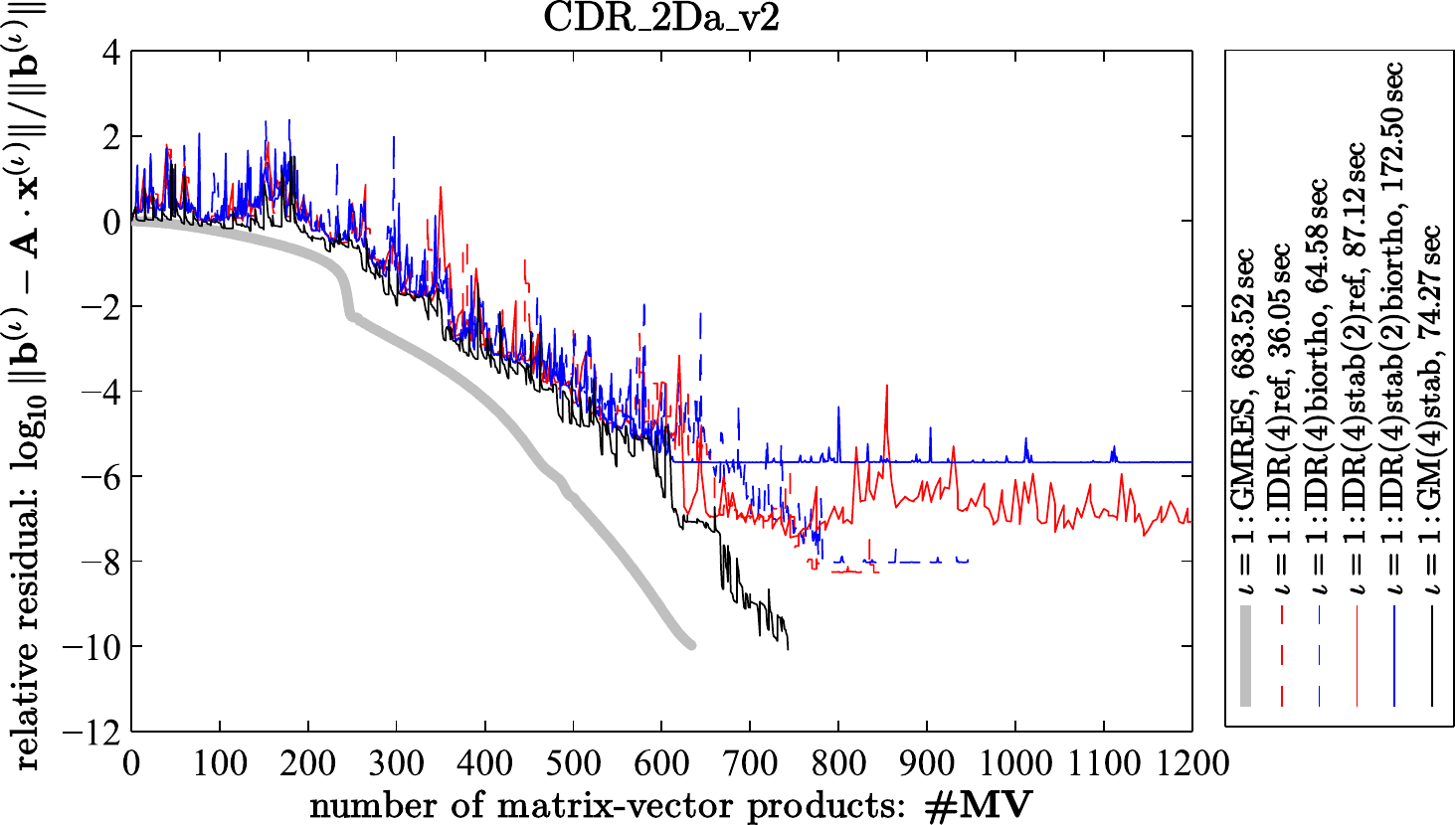}
	\caption{CDR\_2Dparam($0,0$). Diffusion. True residuals.}
	\label{fig:CDR_2Da_trueres}
\end{figure}

The convergence diagrams for the four test problems are given in Fig.\,\ref{fig:CDR_2Da_iterres}--\ref{fig:CDR_2Dd_trueres}.
We emphasize that the latter of these figures shows the true residuals. Additionally, we show in Fig.\,\ref{fig:CDR_2Da_trueres} the convergence of the true residuals for the Poisson problem. For all but the first problem the residual gap is $<10^{-10}$.

We make the following observations:
\begin{itemize}
	\item No method converges for the diffusion-reaction problem ($c_1=0$, $c_2=1$). Convergence in an efficient way would mean to the author that a short-recurrence method achieves at least the rate of convergence that we have drawn in green in Fig.\,\ref{fig:CDR_2Dc_iterres}.
	
	We believe that all the methods loose their superlinearity at $\approx 1250$ matrix-vector products because the large intermediate residual norms amplify rounding errors in the \BiCG-coefficients that underlay the \IDRO recurrence.
	\item Neither the reference nor the biortho implementation of \IDRstabno achieve the desired residual accuracy for any of the four test cases. In all except one case (namely \IDRstabbiortho on the Poisson problem) the residual gap is below $10^{-10}$. Thus, in general this loss of superlinearity is unrelated to the residual gap.
\end{itemize}

The author does not have an explanation why the third problem is much harder to solve than the fourth problem. We have not investigated this further since we do not believe that this is a relevant test problem. This is because of the reason that this system is symmetric, so optimal methods exist anyway.

With regards to all except the third system, the \IDR and \IDRstab variants do all converge reliable and reasonably close behind \GMRES. This is interesting because in \cite{IDRstab-Paper} the experiments indicate that \IDR does not converge well for the second and the fourth problem. This happened because the authors in \cite{IDRstab-Paper} did not use a convergence-maintaining scheme to compute the $\tau$-values, cf. Alg.\,\ref{Algo:StabCoeffs} and \cite{EnhancedBiCGstabL}.

\subsubsection{Experiments with preconditioning}
In the following we show three experiments with preconditioning. To keep things simple, in each case we use a preconditioning that returns two regular triangular matrices $\bL,\bR \in \R^{N \times N}$. Given a system
\begin{align*}
	\tbA \cdot \tbx = \tbb\,,
\end{align*}
we solve the preconditioned system
\begin{align*}
	\bA \cdot \bx = \bb\,,
\end{align*}
where for an arbitrary vector $\bv \in \R^N$ we use the operator $\bA \cdot \bv := \bL^{-1} \cdot \big(\tbA \cdot (\bR^{-1} \cdot \bv)\big)$ and the right-hand-side $\bb := \bL^{-1} \cdot \tbb$. The triangular systems are solved by forward and backward substitution. We pass the preconditioned system to the solvers by passing a function handle that computes the preconditioned matrix-vector product with $\bA$ and by passing the right-hand side $\bb$.

\begin{figure}
\centering
\includegraphics[width=1\linewidth]{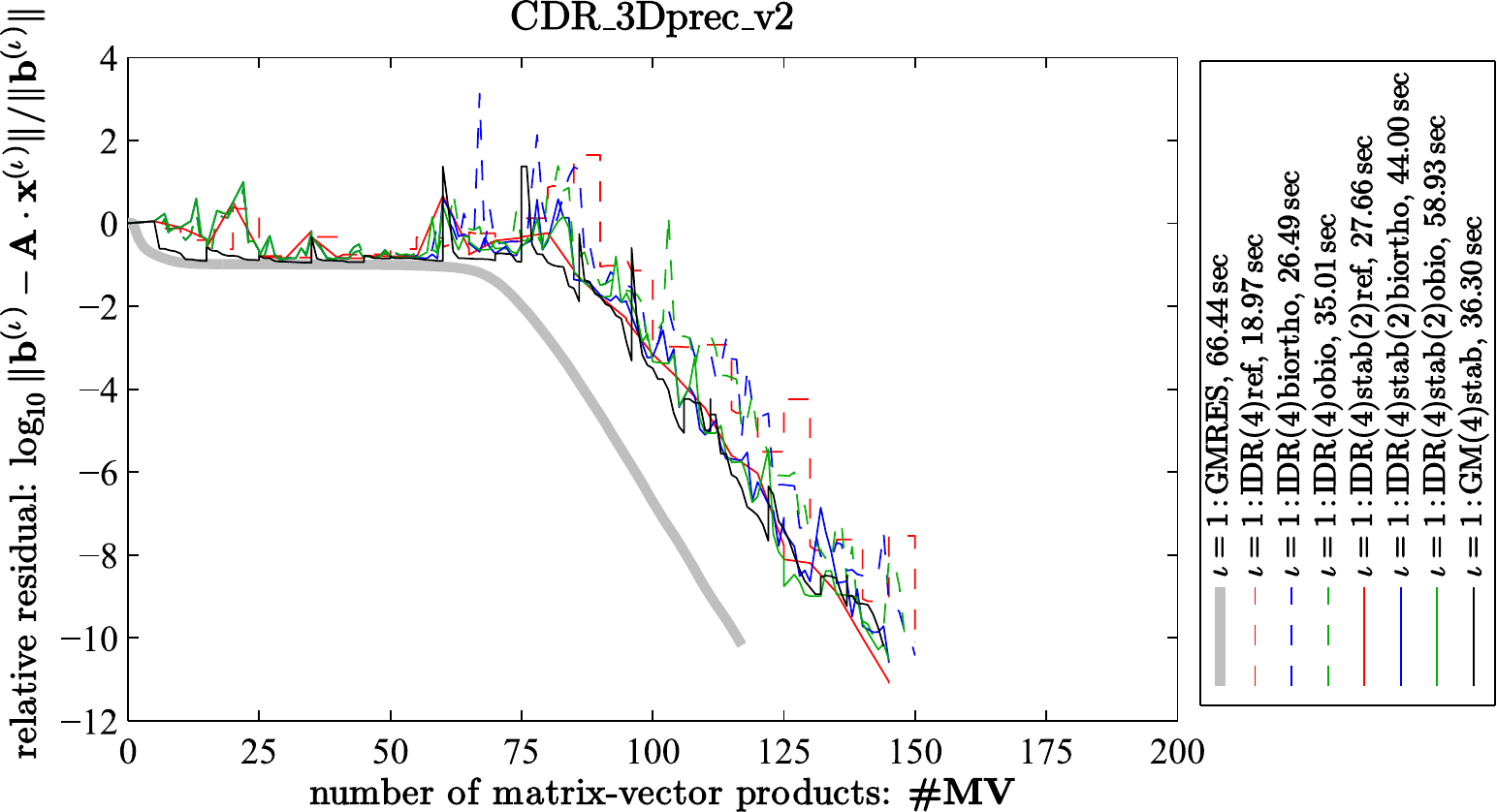}
\caption{If the system is sufficiently well-conditioned then all methods converge well.}
\label{fig:CDR_3Dprec_trueres}
\end{figure}

\paragraph{CDR\_3Dprec}
We consider again the 3D convection-diffusion-reaction problem from Sec.\,\ref{sec:ExamplesConvergenceMstab} and Fig.\,\ref{fig:CDR3D}. We obtain the preconditioners $\bL,\bR$ by applying Matlab's built-in ILU(0) preconditioner to the original system matrix $\tbA$. The results are given in Fig.\,\ref{fig:CDR_3Dprec_trueres}.

The preconditioning makes the system easier to solve in terms of the number of matrix-vector products and less numerical round-off in the \BiCG-coefficients. In consequence, all methods converge in a similar way.

With regards to the runtime, \GMstab is $1/4^\text{th}$ slower than the reference implementation of \IDRstab. The runtimes of the other methods are irrelevant since the \IDR implementations on the one hand are uncompetitive due to $\ell<2$ (cf. our experiment BiCGstabL\_xpl1) and the other \IDRstab variants on the other hand are uncompetitive since they are both slower and have worse-conditioned bases than \GMstab.

\largeparbreak

In this test case the number of non-zeros per row in $\bA$ and both preconditioners is at most 7. Further, all computations are performed in a sequential way. Thus, the time to compute the preconditioned matrix vector products does not yet dominate the runtime over that of the DOTs and AXPYs. However, in a realistic scenario, i.e. high-order finite elements discretisations with computation- and communication-intensive preconditioners, the preconditioned matrix-vector products would completely dominate the runtime. 

This is why for the subsequent test cases we have chosen systems where the amount of non-zeros in $\bA$ is larger so that the relative runtime spent in the matrix-vector products becomes larger. The following two test cases are from \cite{FloridaSparse}.

\paragraph{Norris\_torso1}
This problem is provided by S.\,Norris (\texttt{http://www.esc.auckland.ac.nz/People/Staff/Norris}). The matrix \texttt{torso1} is a coupled finite difference and boundary element discretisation of an electro-physiological model of a human torso in two spatial dimensions. The data provides a real non-symmetric system matrix $\tbA$ only. $\tbA$ has size $N=116158$ and $8516500$ non-zeros. The non-zero pattern is symmetric.

We generate a test problem from this matrix as follows. We reorder $\bA$ by using Matlab's built-in symmetric reverse Cuthill-McKee algorithm. Afterwards we set $\tbb = \tbA \cdot \mathbf{1}$. Since the diagonal of $\tbA$ is non-zero we can compute the preconditioners $\bL,\bR$ by applying ILU(0) directly on $\tbA$.

\begin{figure}
\centering
\includegraphics[width=1\linewidth]{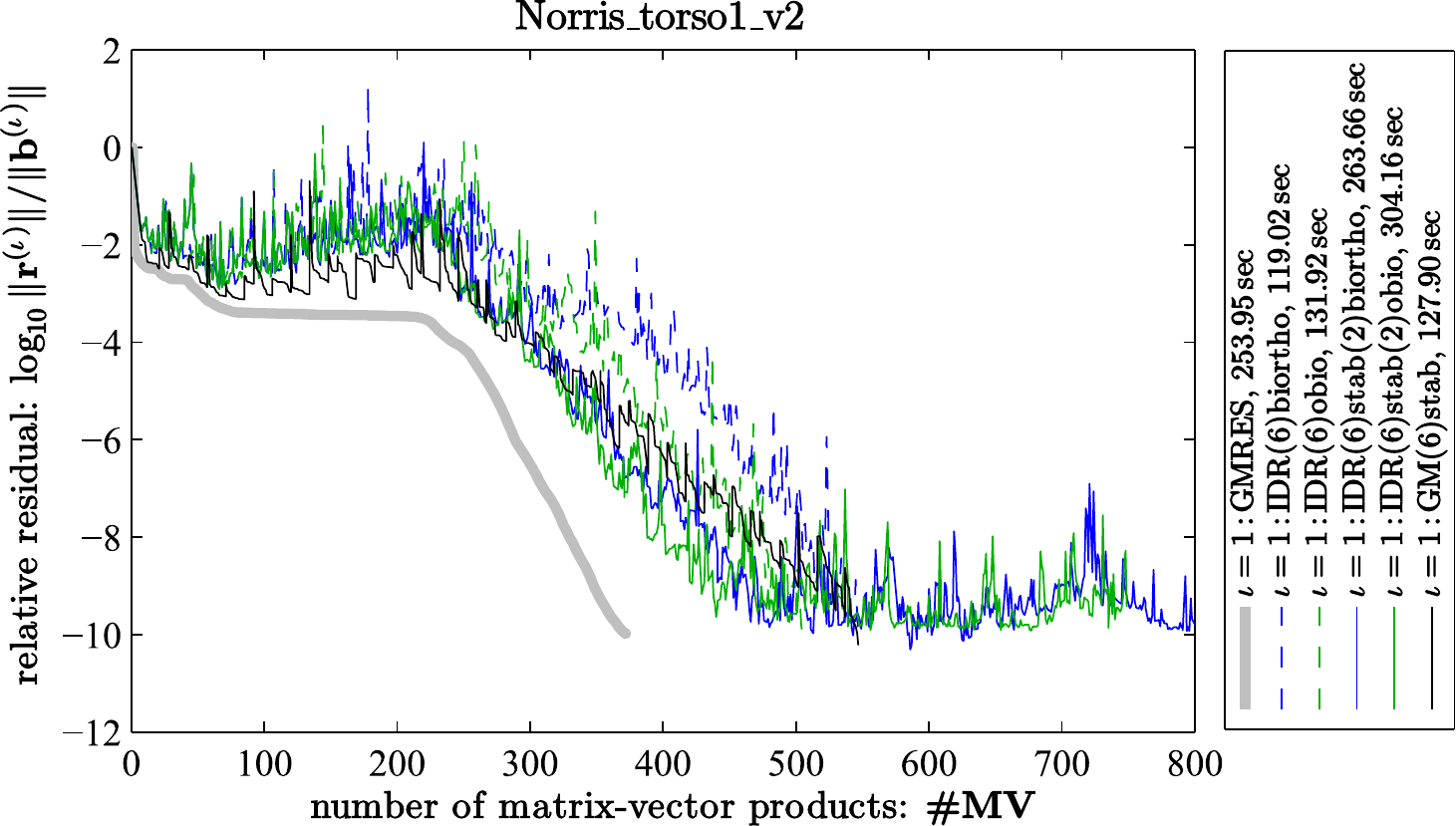}
\caption{Iterative residual vs. number of matrix-vector products for Norris\_torso1.}
\label{fig:Norris_torso1_iterres}
\end{figure}

Fig.\,\ref{fig:Norris_torso1_iterres} shows the iteratively computed residuals. Since we did not want to plot 8 graphs this time we excluded the reference implementations of \IDR and \IDRstab.

\GMstab and \IDRstabbiortho converge both after 550 matrix-vector products. For this test case we have chosen $s=6$, i.e. larger than in the former experiments. In consequence of this increase, the computational cost in \GMstab is now larger than that of \IDRbiortho, as can be seen from the runtime in seconds.

Another expected observation is that again the biortho and obio \IDRstab variants have a bad convergence maintenance: Their iteratively updated residual stagnates above the desired residual accuracy.

\largeparbreak

However, what strikes our eyes much more is the fact that at around the 400th matrix-vector product the rate of convergence of \GMstab in terms of the slope is only half as fast as that of all the other methods. This is unfortunate since good convergence properties with respect to the number of matrix-vector products has been the top 1 goal of our work on \GMstab. This is why in the following we want to look deeper into whether the other methods do truly achieve a better rate of convergence.

\begin{figure}
\centering
\includegraphics[width=1\linewidth]{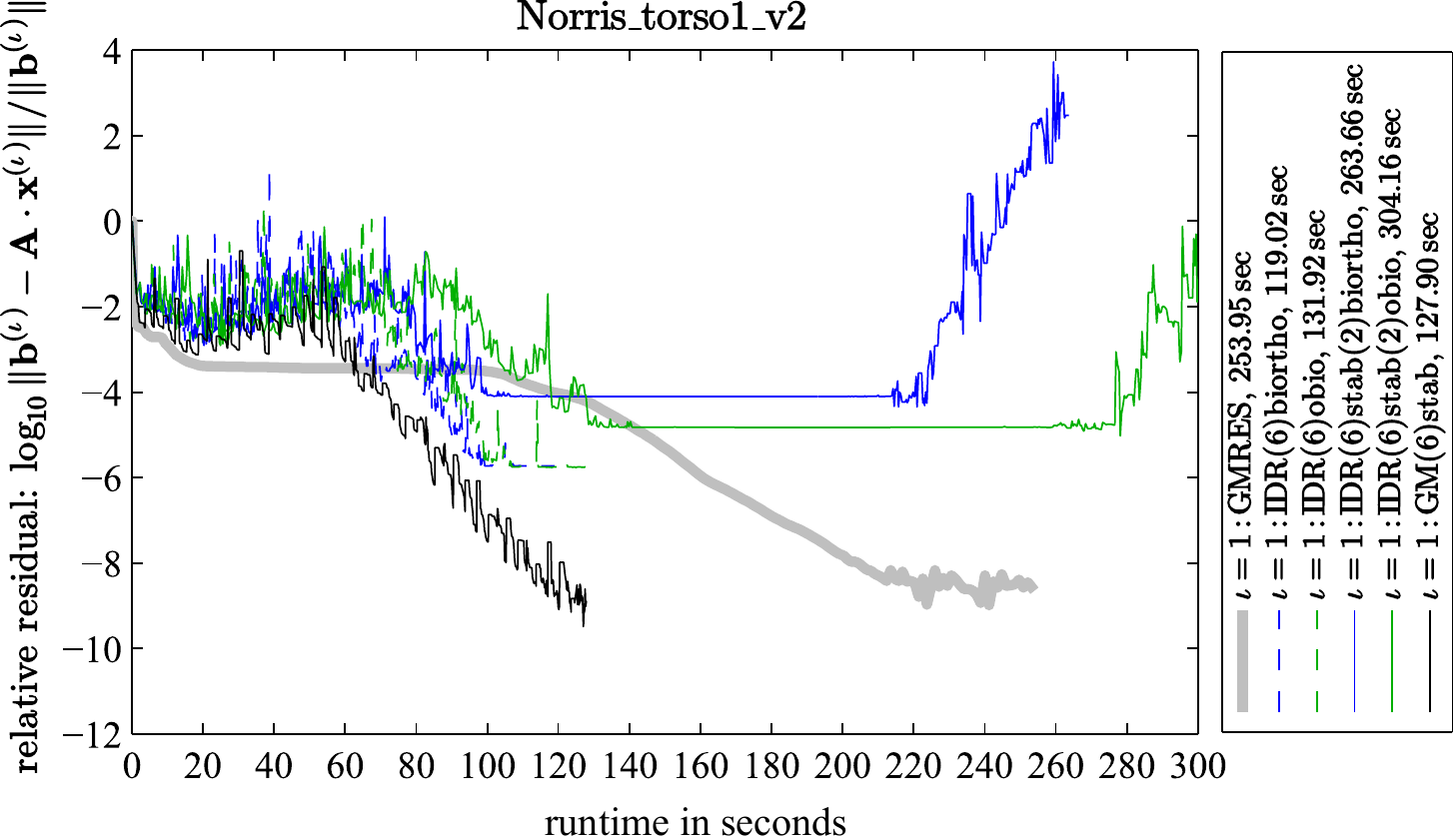}
\caption{True residual vs. runtime for Norris\_torso1.}
\label{fig:Norris_torso1_truerun}
\end{figure}

For this purpose we plot in Fig.\,\ref{fig:Norris_torso1_truerun} the true residual over the runtime in seconds. This shows two things:
\begin{enumerate}
	\item Regarding the \IDRstab variants: From the 330th matrix-vector product, i.e. a relative residual of $10^{-4.5}$, the rate of convergence of the biortho and obio variant supersedes that of \GMstab. However, from Fig.\,\ref{fig:Norris_torso1_truerun} we find that at a relative residual of $10^{-4.5}$ the residuals of these methods are already fully decoupled. Thus, these methods do never really achieve this rate of convergence for the true residuals.
	\item Regarding the \IDR variants: These methods have a huge offset to \GMstab. It seems that at about a relative residual of $10^{-6}$ they can keep up again with \GMstab. However, from Fig.\,\ref{fig:Norris_torso1_truerun} we find that at a relative residual of $10^{-6}$ the residuals of these methods are actually already fully decoupled.
\end{enumerate}
After all, this problem seems to have a strong decoupling effect in the sense that the residual gap grows quickly with respect to the number of iterations. This is why we want to evaluate what rate of convergence our \GMstab would have achieved if we had not used flying restarts. Since the system seems to have a strong decoupling effect, the flying restarts are likely to spoil the rate of convergence of \GMstab.

\begin{figure}
\centering
\includegraphics[width=0.7\linewidth]{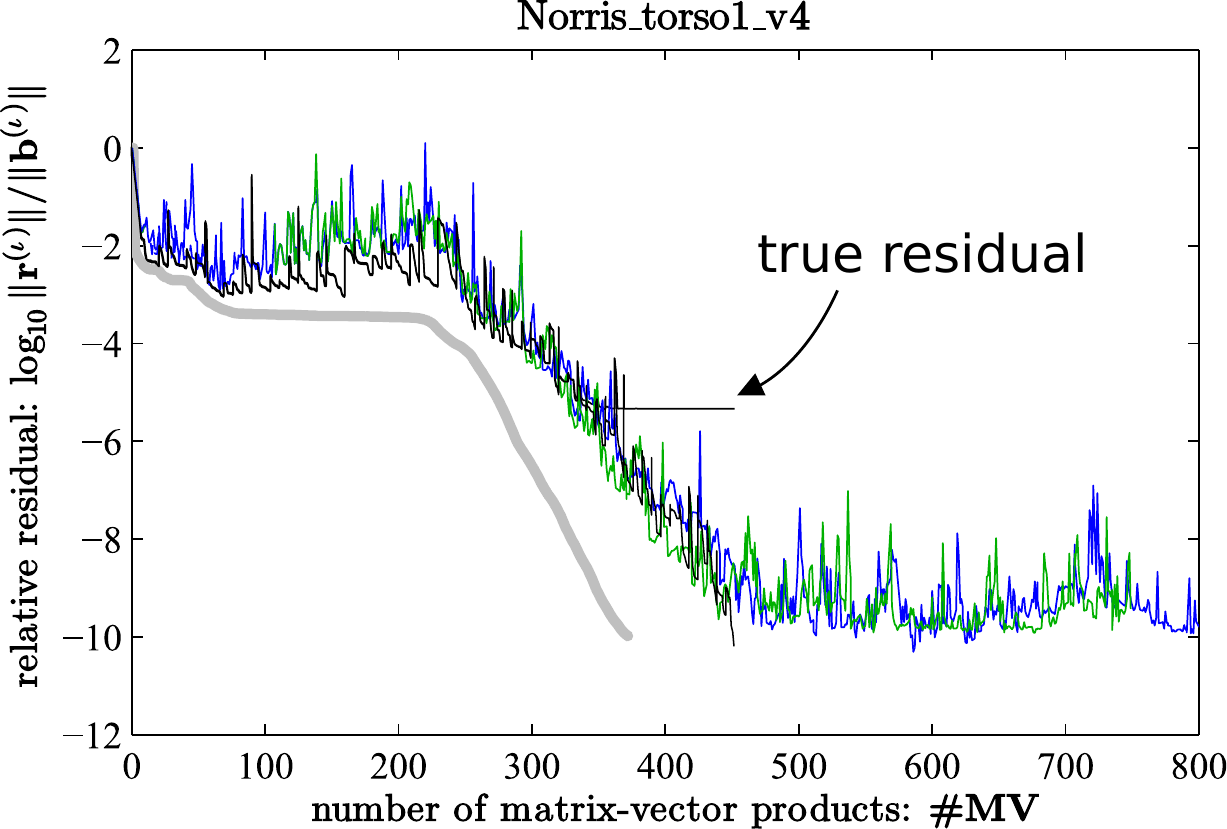}
\caption{Convergence of \GMstab without flying restart vs. \IDRstabno.}
\label{fig:Norris_torso1_decouple}
\end{figure}

Fig.\,\ref{fig:Norris_torso1_decouple} shows the convergence of \GMstab without flying restarts against the graphs of the two \IDRstab variants from Fig.\,\ref{fig:Norris_torso1_iterres}. We observe that under these fairer circumstances the rate of convergence of \GMstab is as good as that of the other two \IDRstab variants\footnote{pride re-established}. However, the true residual now stagnates above $10^{-6}$ also for \GMstab.

\paragraph{Sandia\_ASIC\_320ks}
The real non-symmetric system matrix $\tbA$ of this problem is a circuit simulation matrix from the Sandia National Laboratory. The size of the system is $N=321821$. It has 1316085 non-zeros. The sparsity pattern is non-symmetric. 75 diagonal elements of $\tbA$ are zero.

We generate a test problem from this matrix as follows: We compute $\tbb = \tbA \cdot \mathbf{1}$. We reorder the matrix by using Matlab's built-in approximate minimum degree reordering algorithm \texttt{amd}. The preconditioner matrices $\bL,\bR$ are computed by applying ILU(0) on the reordered and shifted matrix $(\tbA + \lambda \cdot \bI)$ for $\lambda=1$. The shift is necessary since ILU(0) on $\tbA$ leads to a preconditioner that is close to singular. Despite the remarkable shift the preconditioner is still effective.

\begin{figure}
\centering
\includegraphics[width=1\linewidth]{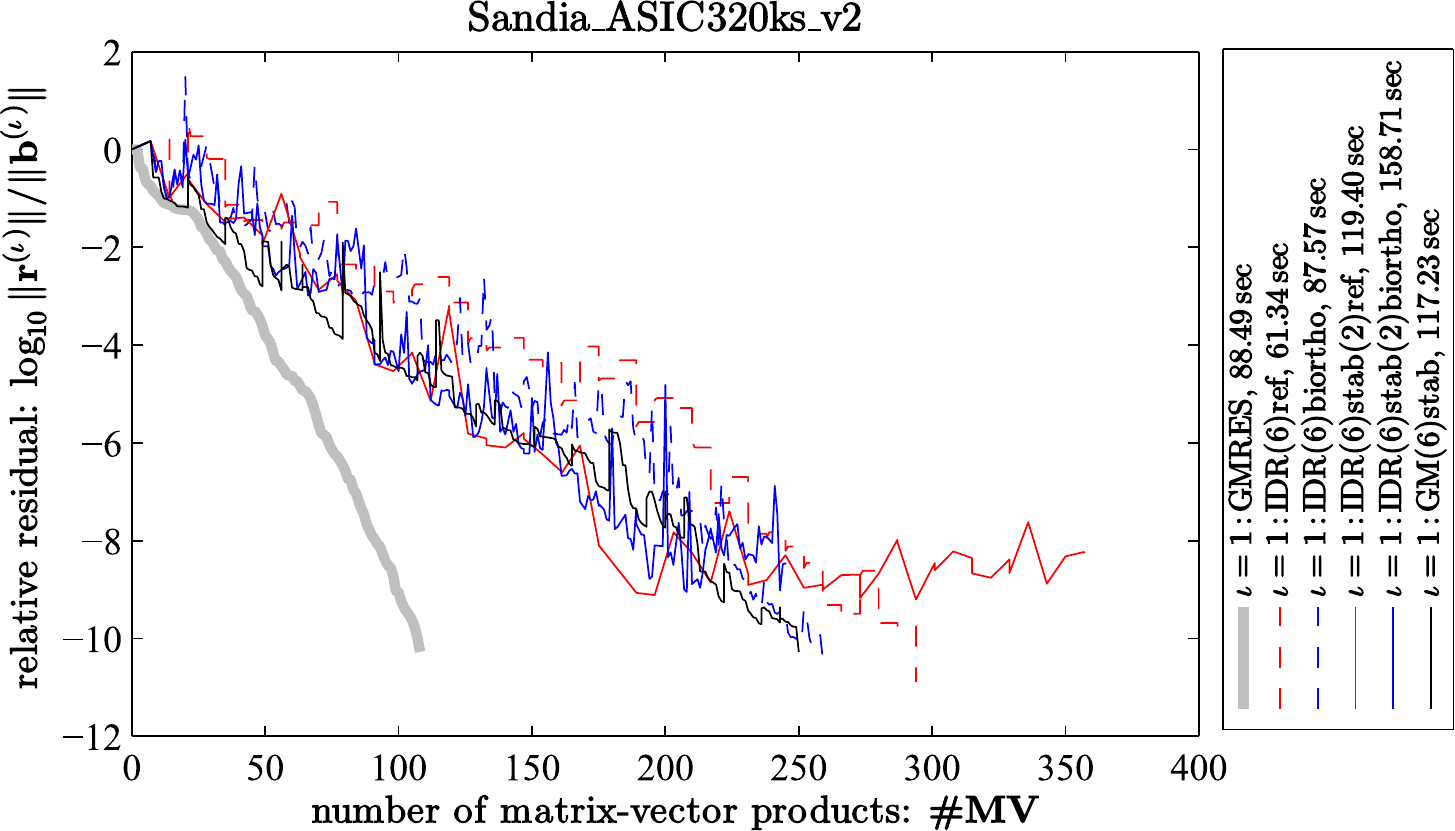}
\caption{All \IDRO methods converge three times slower than \GMRES.}
\label{fig:Sandia_ASIC320ks_iterres}
\end{figure}

Fig.\,\ref{fig:Sandia_ASIC320ks_iterres} shows the convergence of the iteratively updated residual. The true residuals of all methods except \GMstab stagnate above $10^{-8}$. However, it is more of a concern that all the methods do only converge linearly and none of them achieves a faster rate of convergence such as \GMRES. A possible reason could be that columns or rows of $\bA = \bL^{-1} \cdot \tbA \cdot \bR^{-1}$ are badly scaled. In this case the rounding errors in the \BiCG-coefficients become large and no convergence improvement can be made.

\subsection{Experiments with sequences of linear systems}

\largeparbreak





In this subsection we apply \Mstabno for solving sequences of systems of linear equations \eqref{eqn:AxbSequence}. As implementations, we use the codes from the above \IDR and \IDRstab variants and only exchange the initialisation routine.

\subsubsection{Experiments without preconditioning}

\paragraph{CDR\_2Dparam} We solve again the parametric test problem of convection-diffusion-reaction problems for the same parameters as in the latter subsection. However, this time we solve a sequence of systems for each parameter set $(c_1,c_2)$. The first right-hand side $\bb\ha$ is chosen as the right-hand side $\bb$ from the original test problem. The second right-hand side $\bb\hb$ is chosen as $\bA^{-1}\cdot \bb\ha$, i.e. as the solution to the first right-hand side. A similar right-hand side to $\bb\hb$ occurs when solving a transient convection-diffusion-reaction problem using an implicit time-stepping scheme. Thus, this sequence of right-hand sides is representative for a realistic scenario.

\begin{figure}
\centering
\includegraphics[width=1\linewidth]{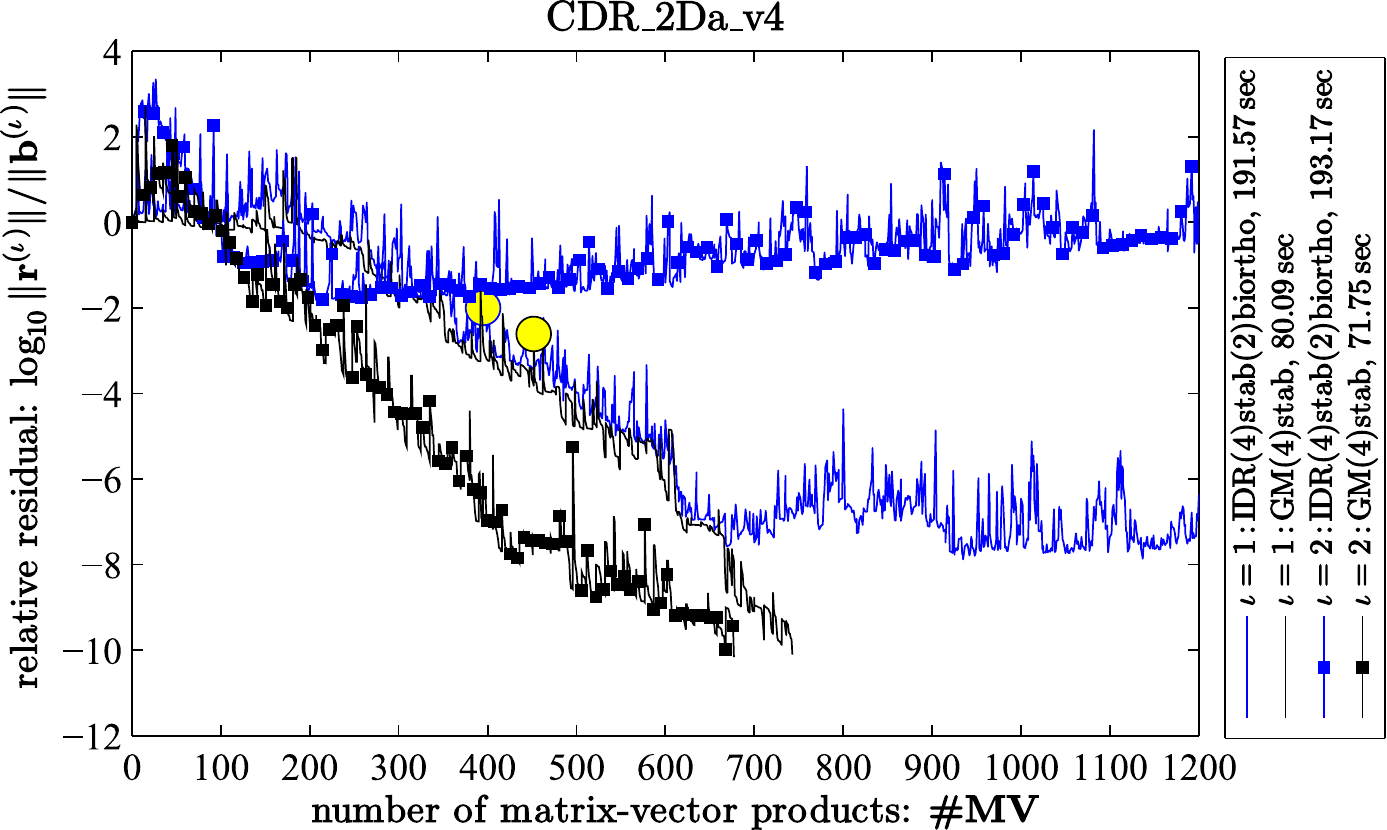}
\caption{CDR\_2Dparam($0,0$). Diffusion.}{}
\label{fig:2D_CDRa_mv}
\end{figure}

\begin{figure}
	\centering
	\includegraphics[width=1\linewidth]{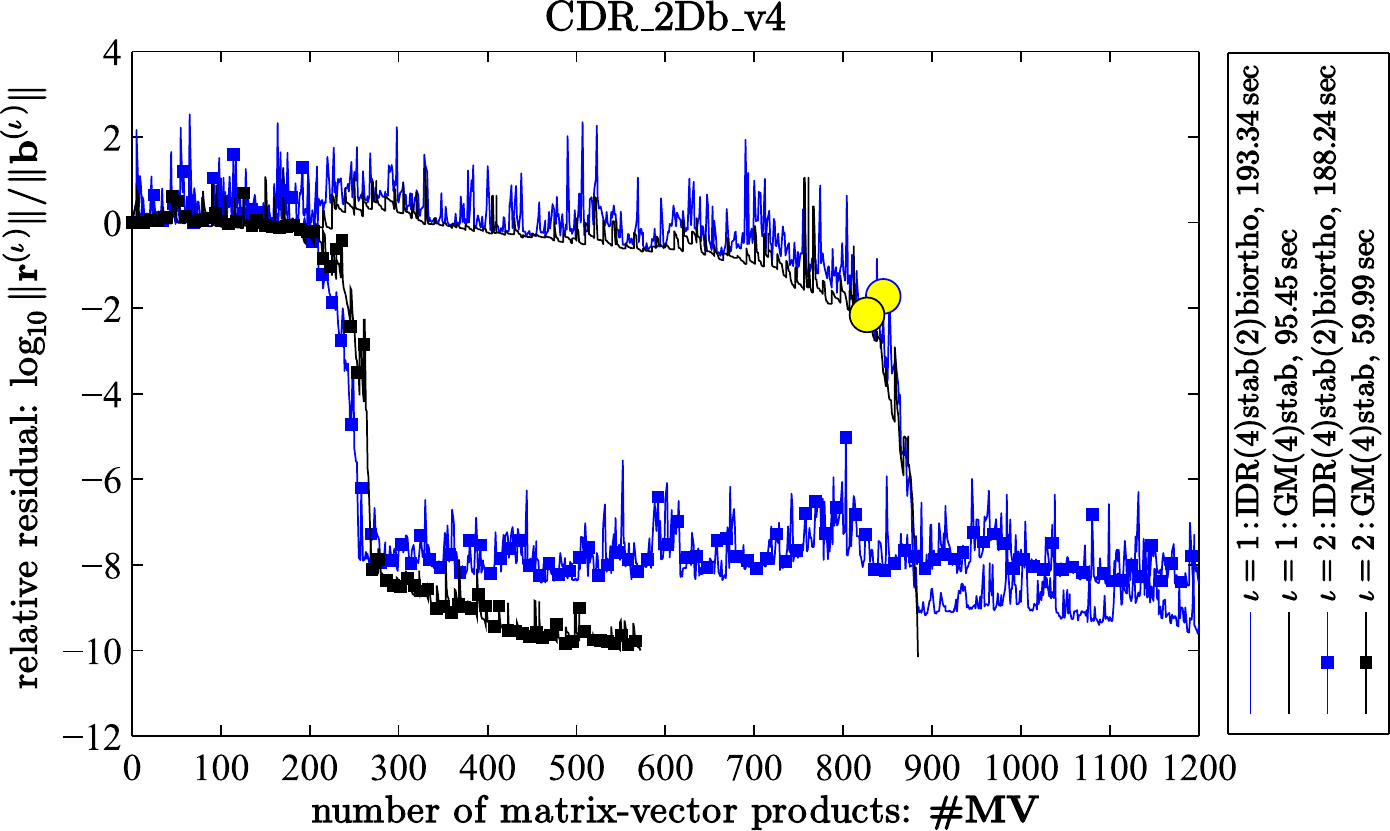}
	\caption{CDR\_2Dparam($1,0$). Convection-Diffusion.}
	\label{fig:2D_CDRb_mv}
\end{figure}

\begin{figure}
	\centering
	\includegraphics[width=1\linewidth]{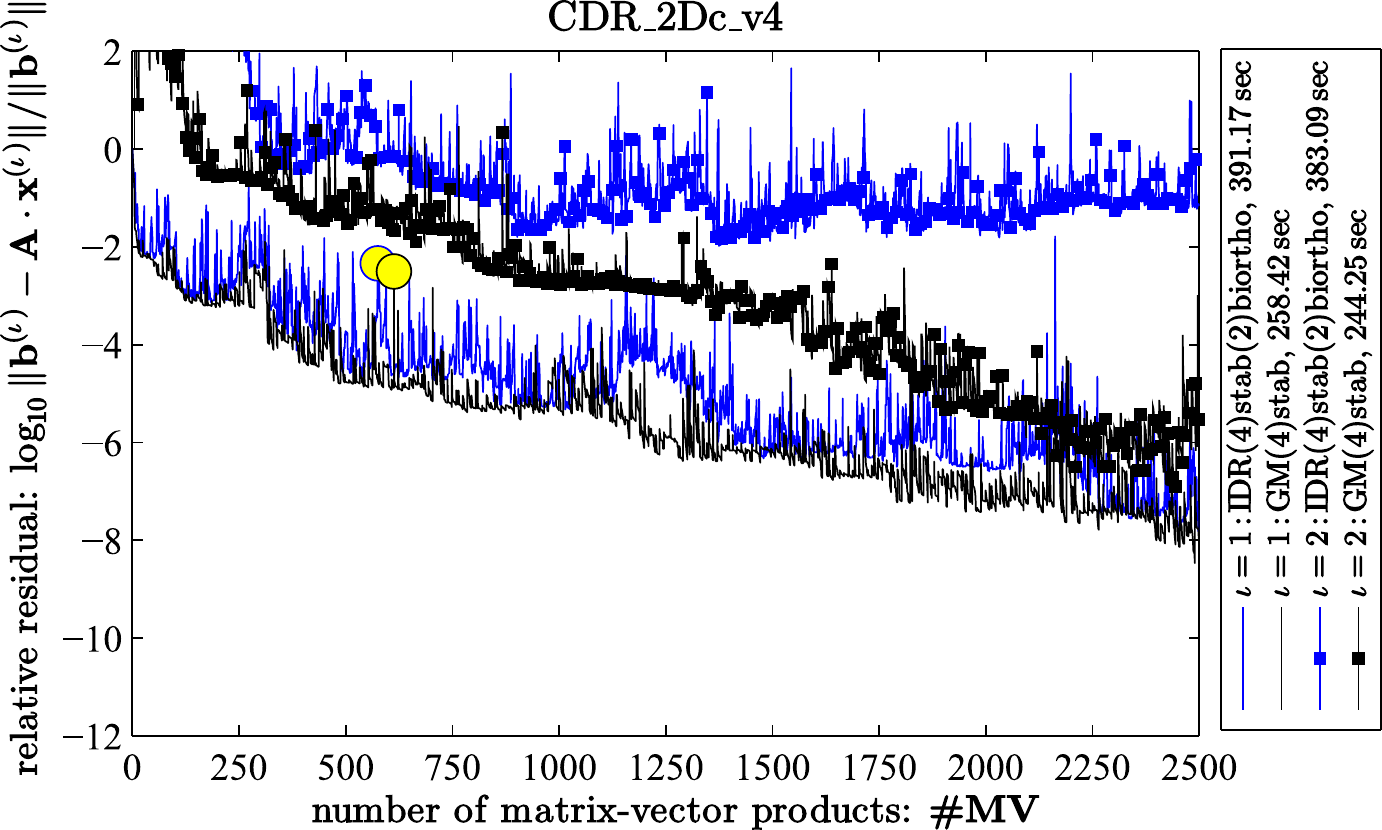}
	\caption{CDR\_2Dparam($0,1$). Diffusion-Reaction.}
	\label{fig:2D_CDRc_mv}
\end{figure}

\begin{figure}
	\centering
	\includegraphics[width=1\linewidth]{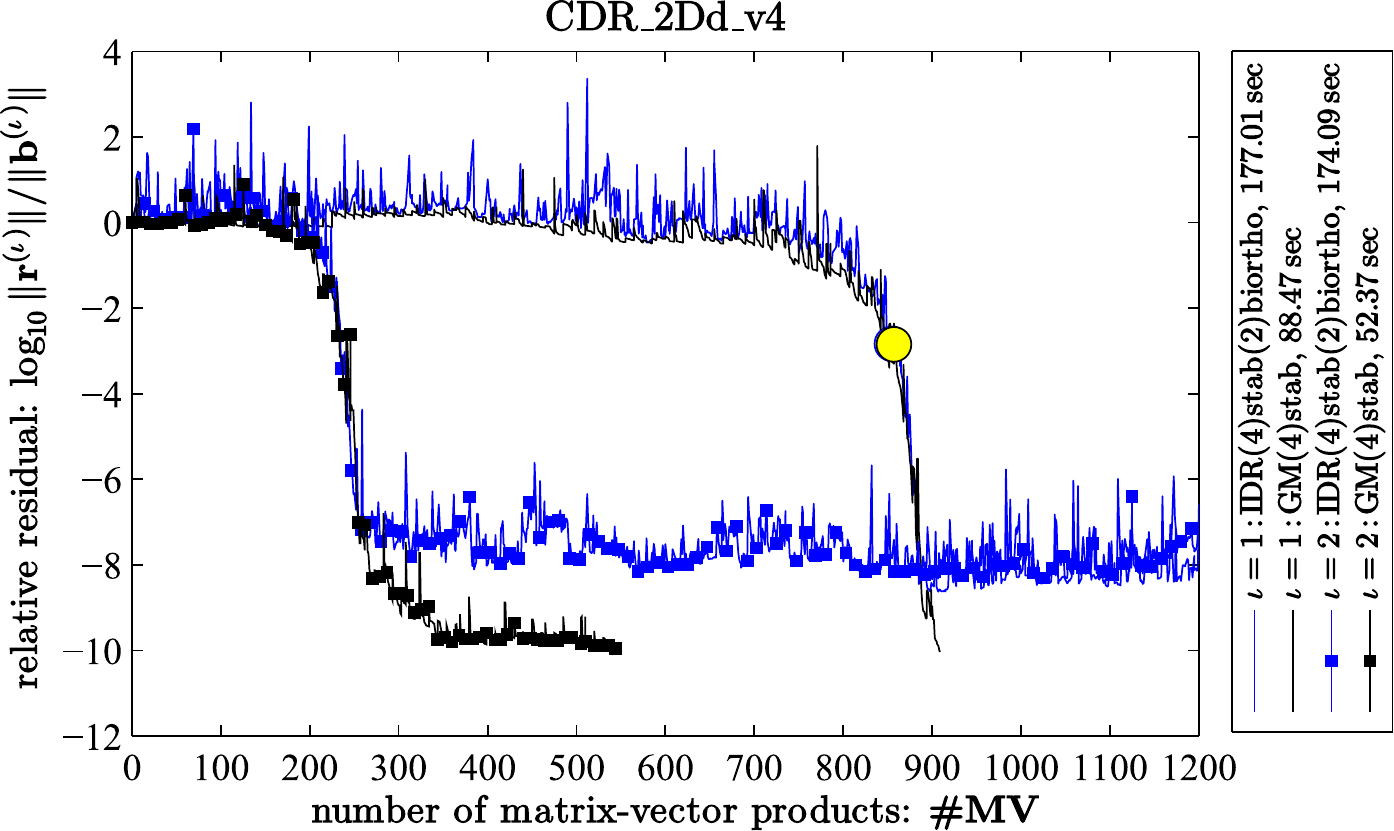}
	\caption{CDR\_2Dparam($1,1$). Convection-Diffusion-Reaction.}
	\label{fig:2D_CDRd_mv}
\end{figure}

In Fig.\,\ref{fig:2D_CDRa_mv}--\ref{fig:2D_CDRd_mv} we have solved each of the four sequences as follows: We use the biortho and restarted GMRES variants of \IDRstabno to solve $\bb\ha$. Then we apply \Mstabno in the respective implementations to solve for $\bb\hb$. The iteration at which the data $\hbU$ is written out in the respective \IDRstabno variant is marked in the figures by a yellow dot that is encircled in the color of the respective implementation. We stress that in \GMstab the relative residual gap is always smaller than $10^{-10}$.

From the figures we make the following observations:
\begin{itemize}
	\item The first and the third test problem, i.e. $(c_1,c_2)=(0,0)$ and $(c_1,c_2)=(0,1)$, do not converge superlinearly. Thus, \Mstabno cannot yield an earlier improvement of the rate of the convergence because there is no convergence improvement attainable at all.
	\item The third test problem $(c_1,c_2)=(0,1)$ is too difficult to solve for both implementation variants: Neither of them convergences in twice the number of matrix-vector products of \GMRES.
	\item The second and the fourth test problem, i.e. $(c_1,c_2)=(1,0)$ and $(c_1,c_2)=(1,1)$, offer a huge convergence improvement. \Mstab can utilise this potential and achieves the transition to superlinear convergence $3$ times earlier than \IDRstabno. Since we have used $s=4$, the speed-up $3$ does perfectly match with the predicted speed-up from our convergence theory in \cite{Mstab-report,Mstab-paper}.
	\item Due to numerical round-off, \Mstab has issues with the convergence maintenance. This can be observed from the figures in that the iteratively updated residuals loose their rate of convergence at relative residual norms of about $10^{-8}$. The biortho implementation suffers significantly more from this issue than the \GMstab implementation.
\end{itemize}

So far, the experiments indicate that the \GMstab implementation can be useful in some cases whereas the biortho implementation is not useful at all for Krylov subspace recycling because it cannot maintain its rate of convergence up to a sufficiently small relative residual.

\largeparbreak
In the following we investigate how the convergence maintenance behaviour of the methods changes when $s$ is increased. To this end, we choose $s=7$ since for this value the theory says that we should precisely obtain twice the speed-up compared to before. Since the first and the third test problem do not offer any potential for Krylov subspace recycling we do only repeat the second and the fourth experiment.

\begin{figure}
	\centering
	\includegraphics[width=1\linewidth]{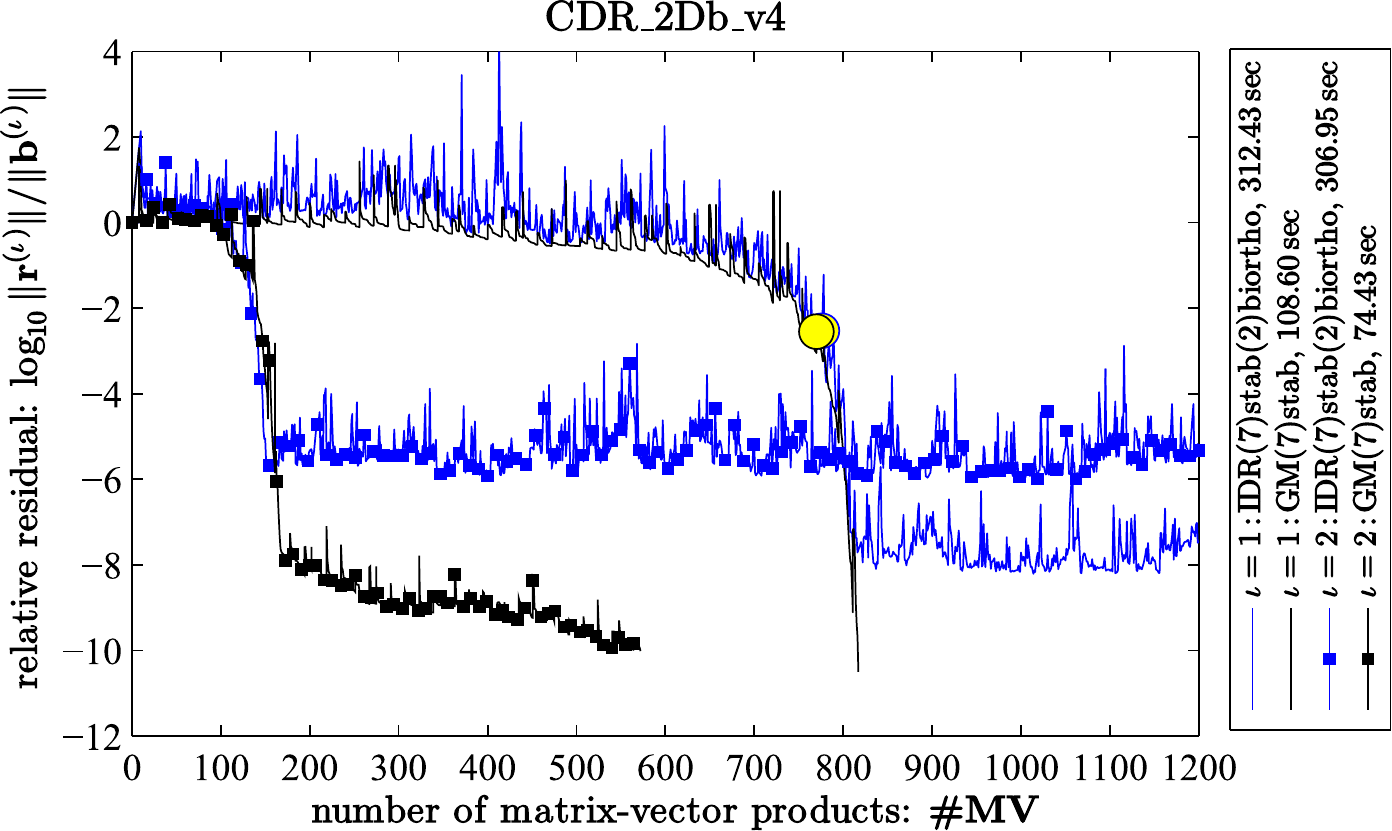}
	\caption{CDR\_2Dparam($1,0$). Convection-Diffusion.}
	\label{fig:2D_CDRb7_mv}
\end{figure}

\begin{figure}
	\centering
	\includegraphics[width=1\linewidth]{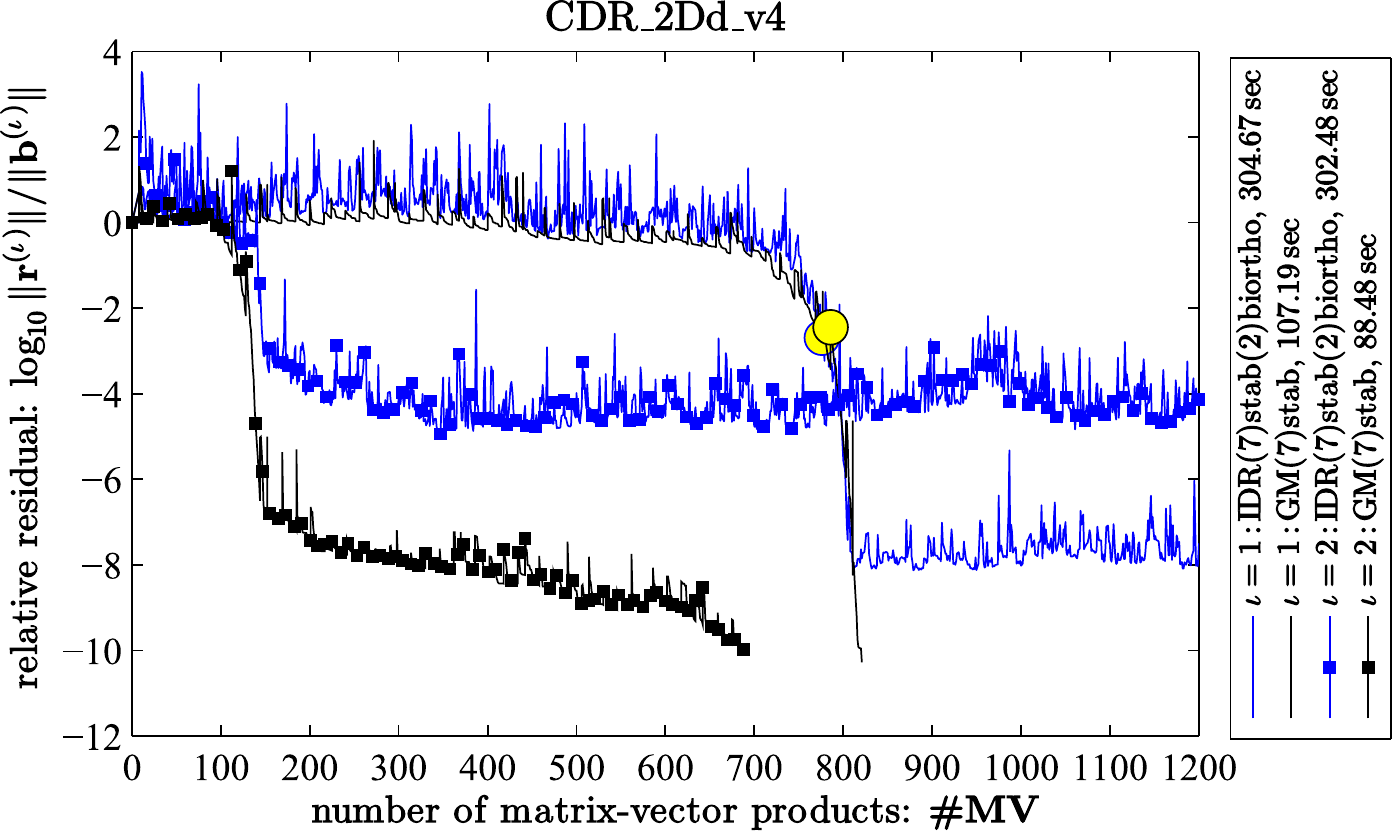}
	\caption{CDR\_2Dparam($1,1$). Convection-Diffusion-Reaction.}
	\label{fig:2D_CDRd7_mv}
\end{figure}

The results of these experiments are given in Fig.\,\ref{fig:2D_CDRb7_mv}--\ref{fig:2D_CDRd7_mv}. We make the following observations:
\begin{itemize}
	\item Indeed, the matrix-vector product at which the respective \Mstab implementation would seemingly terminate in theory has moved from $\approx 250$ to $\approx 125$.
	\item Undesirably, the increase of $s$ has magnified the convergence maintenance problems of both the biortho and the restarted GMRES variant.
\end{itemize}

The conditions of $\bA$ are $\approx 10^3$ and $\approx 10^4$. Thus, the results of the biortho implementation of \Mstab for the test problem in Fig.\,\ref{fig:2D_CDRd7_mv} are not usable. Instead, the results of the restarted GMRES implementation of \Mstab are still useful: Since it maintains the fast convergence until a relative residual of $10^{-7}$, it returns a solution that is at least accurate to the third digit.

\paragraph{CDR\_3D}
So far we have seen how \Mstab can yield dramatic convergence speed-ups of up to $6$ when the problem instance offers very steep superlinear convergence. In the following test case instead we consider a problem where the final rate of convergence is much steeper than at the beginning but not vertical. We consider the second test problem from Sec.\,\ref{sec:ExamplesConvergenceMstab}. On this problem, we first compare the implementations \GMstab and \IDRstab (in the reference implementation) and then \GMstab and \IDRstabbiortho.

\begin{figure}
	\centering
	\includegraphics[width=1\linewidth]{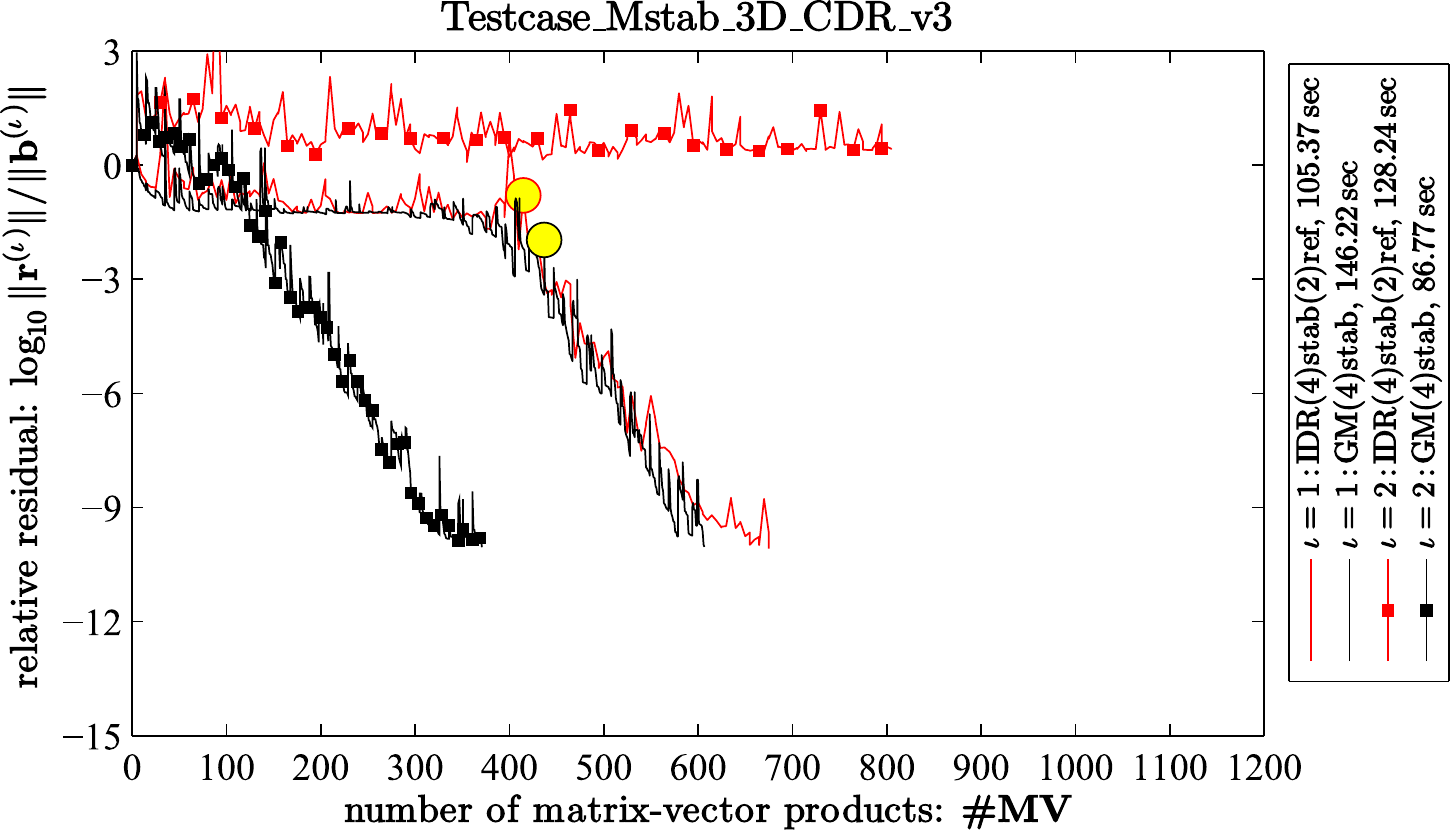}
	\caption{The reference implementation of \IDRstabno is unable to maintain any convergence when being used as \Mstabno.}
	\label{fig:3D_CDR_mv_ref}
\end{figure}

Fig.\,\ref{fig:3D_CDR_mv_ref} compares the reference implementation with \GMstab. We observe that the reference implementation completely looses the convergence when solving for $\bb\hb$. This cannot be a programming error from our side since all implementations use precisely the same initialisation routine. In fact, this absence of any rate of convergence discards the reference implementation as a potential implementation for Krylov subspace recycling.


\largeparbreak

\begin{figure}
	\centering
	\includegraphics[width=1\linewidth]{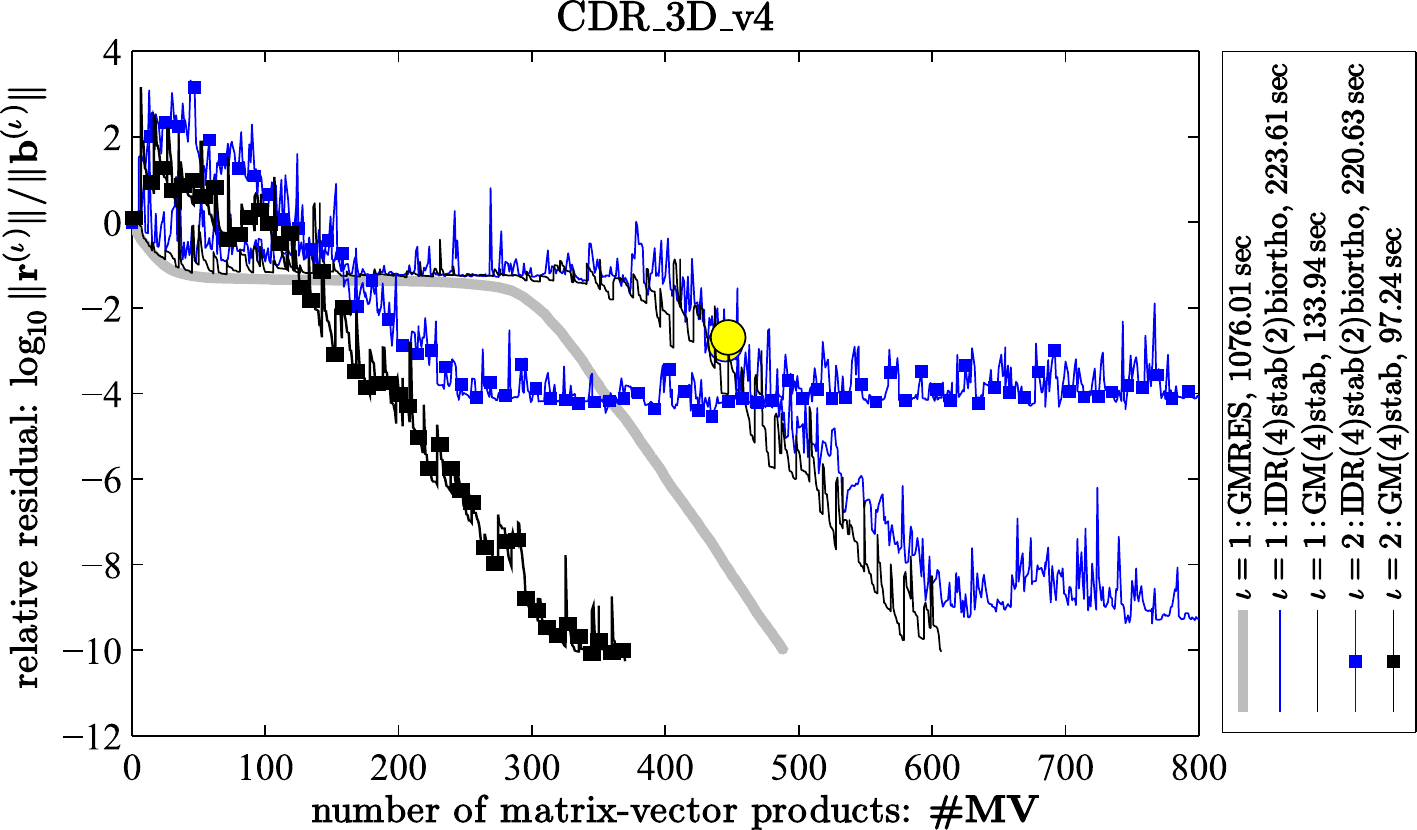}
	\caption{The biortho implementation of \IDRstabno is unable to maintain the convergence long enough when being used as \Mstabno.}
	\label{fig:3D_CDR_mv}
\end{figure}

In Fig.\,\ref{fig:3D_CDR_mv} we compare \GMstab against \IDRstabbiortho. \IDRstabbiortho, which uses in the \Mstabno variant the same modified initialisation routine as the reference \IDRstabno implementation, can yield a similarly improved rate of convergence as \GMstab. However, again it does not maintain this rate of convergence for sufficiently many iterations.

For this test problem the condition is again in the order $10^4$. Thus, again the results that the biortho recycling method generates for $\bb\hb$ are not usable.

\subsubsection{Experiments with preconditioning}
In the following we show two experiments with preconditioning.

\paragraph{CDR\_2Dd\_prec}
This problem is a discretisation of the above parametric two-dimensional convection-diffusion-reaction problem \textbf{CDR\_2Dparam}($1,1$) on a finer mesh, namely $h=1/501$. This makes a system of $N=250\,000$ equations. For the preconditioner $\bL,\bR$ we use ILU(0).

On this system we want to demonstrate the effect that the preconditioner can have on the potential for Krylov subspace recycling. In \cite{SRPCR,Mstab-paper} we discussed that the Krylov subspaces of the right-hand sides $\lbrace \bb\hia\rbrace_{\iota=1,...,\nEqns}$ must overlap in many dimensions because otherwise the computed Krylov information from the first linear system is not useful for the solution of the subsequent linear systems. The consequence in such a scenario is that the recycling brings no benefit to the rate of convergence.

The preconditioner affects the Krylov subspace. Consequently, it has a huge influence on the efficiency of the Krylov subspace recycling. In order to see this from the up-coming experiment, we solve the preconditioned system in the following way:

Given the unpreconditioned system matrix $\tbA \in \R^{N \times N}$ and the two original right-hand sides $\tbb\ha,\tbb\hb \in \R^N$ (as they come from the finer discretisation of the test problem \textbf{CDR\_2Dparam}($1,1$)\, of the former subsection), we consider the following preconditioned linear systems:
\begin{align*}
\underbrace{\bL^{-1} \cdot \tbA \cdot \bR^{-1}}_{=:\bA} \cdot \bx\ha &= \underbrace{\bL^{-1} \cdot \tbb\ha}_{=:\bb\ha}\\
\underbrace{\bL^{-1} \cdot \tbA \cdot \bR^{-1}}_{=:\bA} \cdot \bx\hb &= \underbrace{\bA^{-1} \cdot \bL^{-1} \cdot \tbb\ha}_{=:\bb\hb\equiv \bA^{-1} \cdot \bb\hp{1}}\equiv \bR \cdot \tbb\hb\\
\underbrace{\bL^{-1} \cdot \tbA \cdot \bR^{-1}}_{=:\bA} \cdot \bx\hc &= \underbrace{\bL^{-1} \cdot \tbb\hb}_{=:\bb\hc}\,.
\end{align*}
This has the following effect: The second right-hand side $\bb\hb$ lives in the search space that is built when solving for $\bb\ha$. The third right-hand side $\bb\hc$ however lives in a space that is in general distinct from $\cK_\infty(\bL^{-1}\cdot\tbA\cdot\bR^{-1};\bL^{-1}\cdot\bb\ha)$. Thus, it is not clear whether this right-hand side offers any potential in terms of faster convergence for Krylov subspace recycling methods.
\largeparbreak

\begin{figure}
	\centering
	\includegraphics[width=1\linewidth]{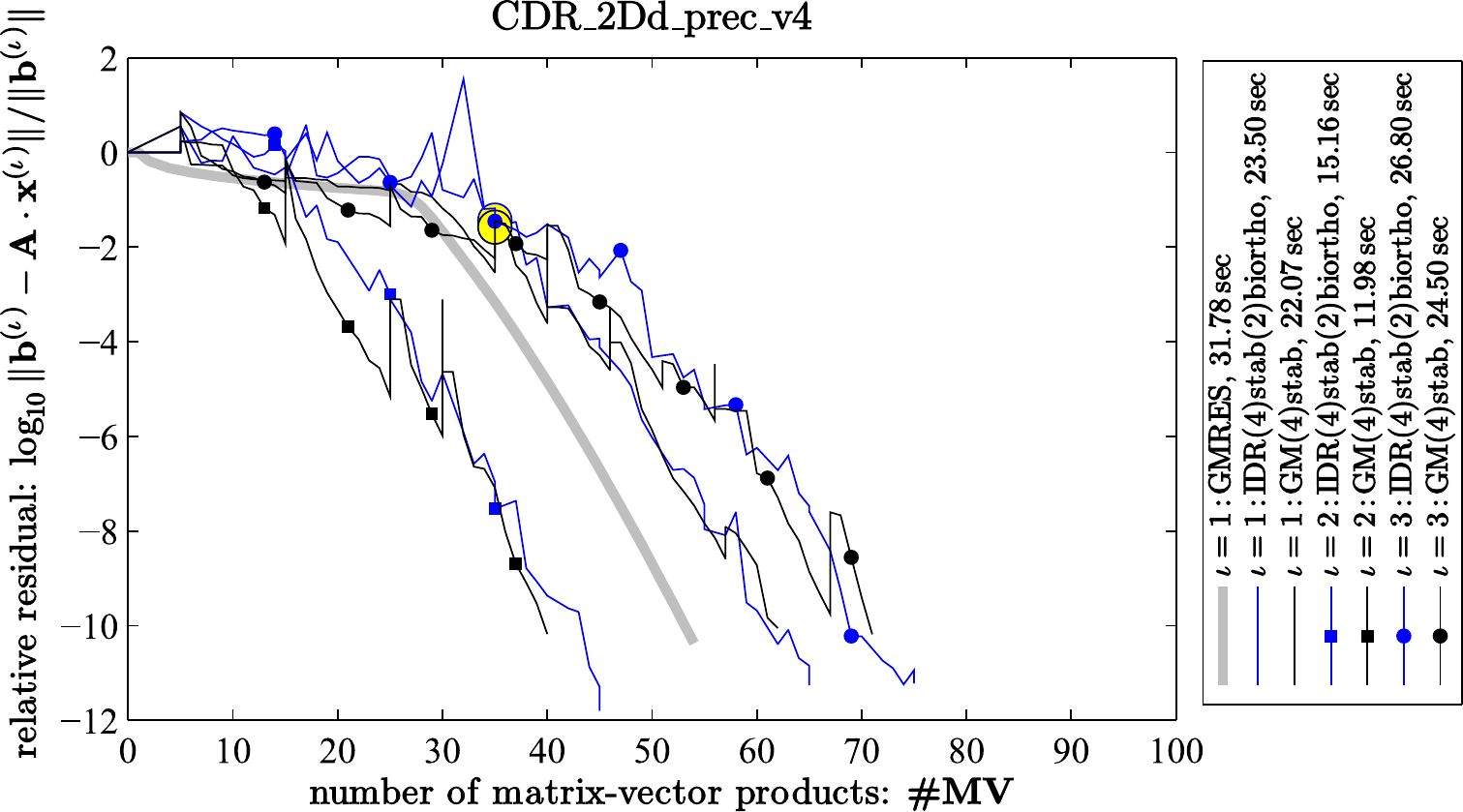}
	\caption{Convergence behaviour of different \Mstabno variants for a preconditioned convection-diffusion-reaction problem.}
	\label{fig:2D_CDRd_mv_prec}
\end{figure}

In Fig.\,\ref{fig:2D_CDRd_mv_prec} we present the convergence graphs of the biortho implementation and the restarted GMRES implementation of \Mstabno. One can see clearly that for $\bb\hb$ the methods converge in $\approx 40$ matrix-vector products, whereas for $\bb\ha$ they require $\approx 60$ matrix-vector-products. Thus, the recycling achieves a speed-up in computation time of $\approx 1.5$\,. In contrast to that, for $\bb\hp{3}$ the methods need a few more iterations than for $\bb\ha$. This is likely because of the fact that the Krylov subspace of $\bb\hc$ is too distinct from that of $\bb\ha$ to yield convergence improvements from the recycling data.

\paragraph{Ocean\_v1}
In the former experiment we have seen that the benefit of recycling is highly sensitive to the preconditioners and the way how the right-hand sides are correlated. The following example motivates that despite all these dependencies Krylov subspace recycling can be still very efficient. The system at hand has a non-symmetric system matrix $\bA \in \R^{N \times N}$ for $N=169850$. The matrix arises from a finite-elements discretisation of the ocean flow of planet earth and is taken from \cite{Ocean}. Twelve right-hand sides are provided that arise from a month-dependent wind-field model.

The system is preconditioned using an ILU(0) preconditioner. Since the finite-elements mesh is based on grid-lines of constant longitude and altitude, some grid areas have a much finer resolution than others (especially the poles). As a result, the columns and rows of the system matrix have strongly distinct scales, resulting in a condition number beyond $10^{6}$. Thus, the preconditioner has columns and rows of strongly distinct scales, too. This could affect the correlation of the right-hand sides in a way that is disadvantageous for Krylov subspace recycling. Nevertheless, it seems that in each right-hand side of the preconditioned system the same modes dominate. This is concluded from the fact that \Mstabno (in an appropriate implementation) can yield convergence speed-ups of $1.6$ for each right-hand side.

The dominating modes (i.e. eigenvectors) characterise the full Krylov subspace. So when the dominating modes of all the right-hand sides are similar (e.g. only low-frequent modes for all right-hand sides) then their full Krylov subspaces are roughly identical -- regardless of the preconditioning that has been used. The author believes that for this test problem this is the reason why the Krylov subspace recycling works so well despite the strong preconditioning and the non-trivial correlation between the right-hand sides.

\begin{figure}
\centering
\includegraphics[width=1\linewidth]{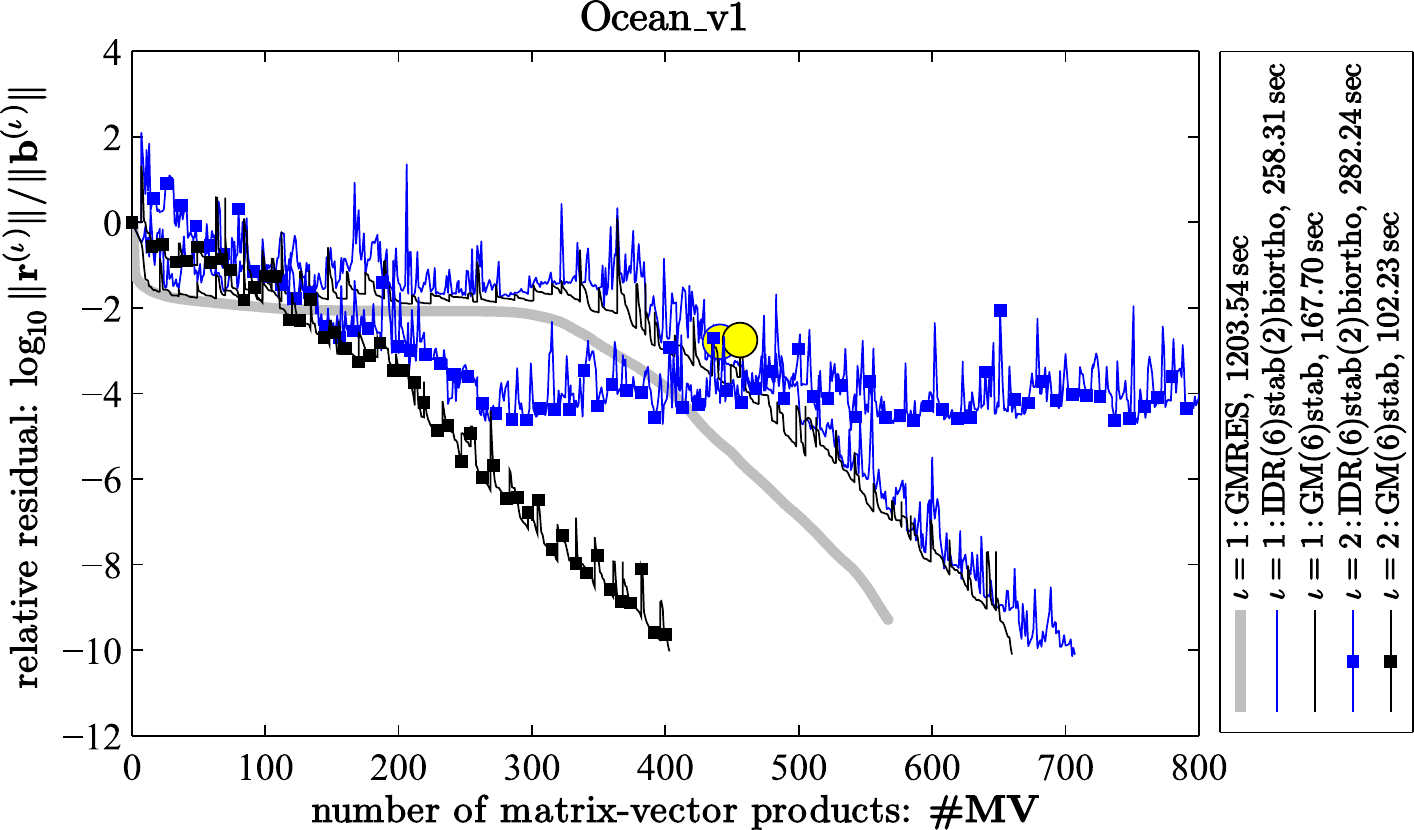}
\caption{Convergence behaviour of different \Mstabno variants for a preconditioned sequence of convection-diffusion problems.}
\label{fig:Ocean_mv}
\end{figure}

Fig.\,\ref{fig:Ocean_mv} shows the convergence of the biortho and restarted GMRES implementation of \Mstabno for the second right-hand side of the sequence. We observe that \Mstabno offers a convergence improvement already after about 100 matrix-vector products. However, only the \GMstab implementation is capable of maintaining this improved rate of converge up to the desired accuracy.

\section{Conclusion}\label{sec:Conclusion}

We have presented a new implementation of \IDRstabp{$s$}{$2$}, called \GMstab. To this end, we first reviewed with many illustrations the motivation of Krylov subspace methods and the theory of Induced Dimension Reduction methods. Then, we re-derived many known \IDR and \IDRstab implementations, which we extended by contributing the ortho-biortho variants.

Afterwards, we pointed out why from the point of theory all these implementations have weaknesses with regards to a robust termination and with regards to the condition of their oblique projection bases. From these insufficiencies we motivated the restarted GMRES implementations of \IDRstabno for $\ell=1$ and $\ell=2$. We derived these implementations from projected and augmented Arnoldi-decompositions. We further discussed the algorithmic steps and the compromises in computational cost that must be made to achieve well-conditioned bases.

We summarised the costs of all the \IDRstabno variants in tables. We also measured the wall-clock time for all methods on problems on which they take a similar iteration count in order to deliver a cost comparison that is close to a realistic scenario. 

In numerous computational experiments we have shown that the new contributed implementation \GMstab is superior with regards to convergence maintenance. Further, in runtime the new implementation occurs computationally cheaper than \IDRstabbiortho, at least for $\ell=2$, $s=4$. All in all, we strongly believe that the reader finds from the numerical experiments, in particular those with Krylov subspace recycling (\Mstabno), that the new implementation \GMstab is absolutely crucial to make the incredible convergence improvements, that the theory of \Mstab offers, accessible for practical computational applications that suffer from numerical round-off.

\paragraph{Outlook}
If we were about to continue the research on \GMstab we would investigate in these three questions.
\largeparbreak
First, we have seen from the the comparison of Fig.\,\ref{fig:Norris_torso1_iterres} and Fig.\,\ref{fig:Norris_torso1_decouple} that the flying restarts still considerably spoil the rate of convergence. So the question arises whether there is a better alternative to flying restarts. Further, the flying restarts do not treat one particular question: How can the round-off with regards to $\pi:=\operatorname{dist}\big(\br,\cK_\infty(\bA;\br_0)\big)$ be reduced? With this is meant that due to round-off a directional component in the residual arises that does not live in the full Krylov subspace of the initial residual. If the grade of the Krylov subspace of this round-off is by far larger than the grade of the original Krylov subspace then not only the residual decoupling but also the magnitude of $\pi$ can spoil the rate of convergence.
\largeparbreak
Second, in Fig.\,\ref{fig:MstabMotivation} we have secretly presumed that only the dimension of the Petrov space respectively Sonneveld(-like) space determines the iteration number at which the transition to superlinear convergence is encountered. This presumption does reflect well the practice as the numerical experiments show. However, there is yet no theoretical proof that justifies it. We believe it would be desirable if someone could formulate this conjecture in precise mathematical terms and show under which assumptions it holds. This implies the following tasks: 1) Find a requirement for system matrices $\bA$ such that superlinear convergence of Krylov methods is guaranteed. 2) Characterise this superlinearity and postulate how it is achieved earlier when recycling (for ease of theoretical accessibility) the full Petrov space. 3) Finally prove that for randomly initialised right-hand sides $\bb\ha\in\C^N$ and $\bb\hb \in \cK_\infty(\bA;\bb\ha)$ this postulation holds with a likelihood of 100\%.
\largeparbreak
Third, we had the goal to provide the robustest implementation of \IDRstabno and \Mstabno that one could think of. However, what is the cheapest possible implementation that is still sufficiently robust to be useful in practice? Having a cheaper compromise that is often as robust as \GMstab but has lower computational overhead might be advantageous in some cases.

\FloatBarrier


\end{document}